\documentclass[12pt]{format}

\usepackage[toc,page]{appendix}
\usepackage{lineno}
\usepackage{times}
\usepackage{bm}
\usepackage[numbers]{natbib}
\usepackage{amsmath}
\usepackage{amssymb}
\usepackage{extarrows}
\usepackage {mathtools}
\usepackage[ruled,vlined]{algorithm2e}
\usepackage{comment}
\usepackage{float}
\usepackage{subfigure}
\usepackage{hyperref}
\usepackage{natbib}
\usepackage{multicol}
\usepackage{color}
\usepackage[svgnames]{xcolor}
\usepackage{appendix}
\usepackage{booktabs}

\usepackage{amsmath}
\usepackage{amssymb}

\usepackage{times}
\usepackage{bm}
\usepackage{graphicx,caption}

\newcommand{\bR}{{\text{R}}}
\newcommand{\bI}{{\text{I}}}
\newcommand{\bu}{{\text{u}}}
\newcommand{\tr}{\text{tr}}
\newcommand{\VEC}{\text{vec}}

\newcommand{\diag}{\mathrm{diag}}

\renewcommand {\thesection}{\arabic{section}.\hspace{-0.25cm}}
\renewcommand {\thesubsection}{\thesection\hspace{+0.25cm}\arabic{subsection}}
\renewcommand{\theequation}{\thesection\hspace{+0.25cm}\arabic{equation}}

\newtheorem{thm}{Theorem}
\newtheorem{lem}{Lemma}

\begin{document}

\title{Directional testing for high-dimensional multivariate normal distributions}
%

\author{
Caizhu Huang \\
\tt{\it\small caizhu.huang@phd.unipd.it }\\
{\it\small Department of Statistical Sciences, University of Padova, 35121 Padova, Italy} \\
Claudia Di Caterina \\
\tt{\it\small claudia.dicaterina@unibz.it}\\
{\it\small Faculty of Economics and Management, Free University of Bozen-Bolzano, 39100 Bolzano, Italy}\\
Nicola Sartori \\
\tt{\it\small nicola.sartori@unipd.it}\\
{\it\small Department of Statistical Sciences, University of Padova, 35121 Padova, Italy}
}

\date{}
\maketitle

\vspace{-0.8cm}




\noindent {\bf Abstract.} 
Thanks to its favorable properties, the multivariate normal distribution is still largely employed for modeling phenomena in various scientific fields. However, when the number of components $p$ is of the same asymptotic order as the sample size $n$, standard inferential techniques are generally inadequate to conduct hypothesis testing on the mean vector and/or the covariance matrix. Within several prominent frameworks, we propose then to draw reliable conclusions via a directional test. We show that under the null hypothesis the directional $p$-value is exactly uniformly distributed even when $p$ is of the same order of $n$, provided that conditions for the existence of the maximum likelihood estimate for the normal model are satisfied. Extensive simulation results confirm the theoretical findings across different values of  $p/n$, and show that under the null hypothesis the directional test outperforms not only the usual first and higher-order finite-$p$ solutions but also alternative methods tailored for high-dimensional settings. 
Simulation results also indicate that power performance of the different tests depends on the specific alternative hypothesis.




\section{Introduction}
Hypothesis testing on the multivariate normal distribution is a subject of great interest in multivariate statistical analysis \citep[see, e.g.,][]{anderson1958,muirhead1982}. It is widely  applied  in various fields, such as  social sciences,  biomedical sciences and finance.  Under  fixed dimension $p$ of the observation vector and large sample size $n$, standard asymptotic results are available for testing hypotheses on the mean vector and/or the  covariance matrix.
For instance, for sufficiently large $n$,  the classical log-likelihood ratio statistic, its Bartlett correction \citep{bartlett:1937}, and the large-deviation modification to the log-likelihood ratio statistic proposed by Skovgaard \cite{skovgaard:2001} have an approximate  $\chi^2_d$ null distribution, with $d$ equal to the number of constraints on the parameters  imposed under the null hypothesis. 
Yet in many modern applications the dimension $p$, while being lower than the sample size, is large and often comparable with $n$. 

Taking inspiration from a classification recently proposed by Battey and Cox \citep{battey2022some}, we distinguish between three asymptotic regimes where $p$ and $n$ diverge, namely low dimensional, high dimensional and ultra-high dimensional.
In the first case  $p/n$ goes to zero, in the second case $p/n \to \kappa \in (0,1]$, while in the third case $p/n$ tends to infinity or to a limit greater than one.
It is well known that inferential problems arise in regression settings where the number of covariates increases with the sample size, both in low and high dimensional regimes; see, for instance, \cite{fan2019nonuniformity}, \cite{sur:2019} and \cite{sur2019likelihood} for results on logistic regression, and \cite{tang2020modified} for results on exponential family models.

Even for the special $p$-variate normal distribution case,
the classical  likelihood-based testing procedures 
 may be already invalid in low dimensional settings. 
Indeed, He, Meng, Zeng and Xu \cite{He2020B} considered the case where $p=o(n^\alpha)$, $0 < \alpha \le 1$, and showed that the log-likelihood ratio statistic's distribution approximates to a $\chi^2_d$ if and only if  $p=o(n^{1/2})$, while the analogous condition for the Bartlett correction is $p=o(n^{2/3})$.
On the other hand, for the high dimensional setting Jiang and Yang \cite{jiang:2013} derived a central limit theorem that allows to construct reliable tests for hypotheses on the  mean vector and/or the  covariance matrix. Specifically,  they proved that the distribution of  the log-likelihood ratio statistic, suitably standardized, converges to a standard normal when both $p$ and $n$ tend to infinity, provided that $p/n \to \kappa \in (0,1]$.

Under the same high dimensional setting of Jiang and Yang \cite{jiang:2013},
we propose the use of a directional approach for testing general hypotheses on multivariate normal distributions. 
Directional inference on a vector parameter of interest was first introduced by Fraser and Massam \cite{fraser:1985} and later developed by Skovgaard \cite{skogaard:1988} and Cheah, Fraser and Reid \cite{fraser1994} using saddlepoint approximations for the distribution and analytical approximations for the required tail probability integrals in regular asymptotic scenarios. More recently, still assuming the classical scenario with fixed parameter dimension and increasing sample size, Davison, Fraser, Reid and Sartori \cite{sartori:2014} and Fraser, Reid and Sartori \cite{sartori:2016} proposed to compute the directional $p$-value by replacing the analytical approximations with one-dimensional numerical integration.
 The numerical evidence in these papers showed the excellent performance  of the proposed method in many examples of practical interest. In fact, later the directional $p$-value was  proven to coincide with that of an exact $F$-test in a few prominent modeling frameworks \citep{sartori:2019Ftest}.

  For multivariate normal distributions, a first example of the directional  approach was given in \cite[][Example 5.3]{sartori:2014}, who considered testing some conditional independence among the components starting from the full dependence structure. Although their simulations  illustrated that the empirical extreme accuracy of directional inference is not limited to simple low-dimensional situations, those results were not formally justified within the asymptotic framework of the high dimensional regime.
  
  We prove that the directional $p$-value is exact when testing a number of hypotheses on the multivariate normal distribution, even in the high dimensional scenario. Precisely, it is only required that $n \ge p+2$, which is the condition for the existence of the maximum likelihood estimate for the covariance matrix. Consequently, we shall not deal here with the ultra-high dimensional regime. The exact uniform null distribution of the directional $p$-value  follows from the exactness of the normalized saddlepoint approximation to the distribution of the canonical sufficient statistic in the multivariate normal model. 
  We focus here on  hypotheses that compare multiple multivariate normal distributions, while results on hypotheses regarding a single distribution are reported in the Supplementary Material.

  Several simulation studies are conducted for comparing  the proposed method with the usual $\chi^2_d$ approximations for the log-likelihood ratio statistic, the Bartlett correction, the modification of the log-likelihood ratio statistic by Skovgaard \cite{skovgaard:2001} and the central limit theorem method by Jiang and Yang \cite{jiang:2013}. The results confirm the theoretical properties of the directional test under the null hypothesis. Indeed,  the method proves uniformly more accurate than the alternative approaches. Only the central limit theorem test  gives a comparable accuracy when the number of components $p$ is large, but it is less reliable for small to moderate values of $p$. The various methods are assessed also in terms of power, after adjusting for Type I error. 
The approach leading to higher power depends on the specific alternative setting.
  In this respect the directional test, which does not need any adjustment, is competitive with the central limit theorem test and overperforms the other candidates across various alternative scenarios.




\section{Background}\label{section2:background}


\subsection{ Notation}\label{section2.1:notation}

Assume that the model for the data $y = (y_1,\dots,y_n)^T$ is an exponential family with canonical  parameter $\varphi=\varphi(\theta)$ and  canonical  sufficient statistic $u=u(y)$. The distribution of $y$ can then be expressed as 
\begin{eqnarray}
  f(y;\theta) &=& \exp\left[ \varphi(\theta)^T u(y) - K\{\varphi(\theta)\}\right] h(y),\nonumber
\end{eqnarray}
with corresponding  log-likelihood function $\ell (\theta;y) = \log f(y;\theta)$. When the dimension $q$ of $\theta$ is equal to the dimension of $\varphi$ and $\varphi(\theta)$ is one-to-one, the statistic $u(y)$ has a full natural exponential family distribution in the canonical  parameterization. Hence,  we can write $f(u;\varphi) = \exp\left\{ \varphi^T u - K(\varphi)\right\} \tilde{h}(u)$ with  associated log-likelihood $\ell(\varphi;u) = \varphi^T u - K(\varphi)$. However, $\tilde{h}(u)$ can rarely be derived explicitly.

For the development of the directional $p$-value in Section \ref{section2.2}, it is notationally convenient to center the sufficient statistic at the observed  data point $y^0$. Hence we let $s=u-u^0$, with $u^0 = u(y^0)$,  and write 
\begin{eqnarray}
  \ell(\varphi;s) = \varphi^T s + \ell^0(\varphi) =  \varphi^T (u-u^0) + \ell(\varphi;u^0), 
  \label{log-likelihood:centering}
\end{eqnarray}
where $\ell^0(\varphi) = \ell(\varphi;s=0_q)=\ell(\varphi;u=u^0)$, with $0_q$ denoting the $q$-dimensional vector of zeroes.
This centering  ensures that the observed value of $s$ is $s^0=0_q$.

Suppose the parameter vector is partitioned as $\varphi=(\psi^T, \lambda^T)^T$, and thus (\ref{log-likelihood:centering}) can be written as
\begin{eqnarray}
  \ell(\varphi;s) = \psi^T s_1 + \lambda^T s_2 + \ell^0(\psi,\lambda), 
  \label{log-likelihood: s1 and s2}
\end{eqnarray}
 where $\psi$ is a $d$-dimensional  parameter of interest, $\lambda$ is a $(q-d)$-dimensional nuisance parameter, and $(s_1^T, s_2^T)^T$ is the corresponding partition of $s$.
Assume we are interested in testing  the hypothesis $H_{\psi}$: $\psi(\varphi)  = \psi$. 
 To a first order of approximation, a parameterization-invariant measure of departure of $s$ from $H_{\psi}$ is given by the log-likelihood ratio statistic
\begin{eqnarray}
  W = 2\{\ell(\hat{\varphi})-\ell(\hat{\varphi}_\psi)\},\nonumber
  \label{log-likelihood ratio statistic}
\end{eqnarray}
where  $\hat{\varphi}$ is the maximum likelihood estimate and  $\hat {\varphi}_{\psi} = ({\psi}^T, \hat{\lambda}_{\psi}^T)^T$  denotes the constrained maximum likelihood estimate of $\varphi$ under $H_{\psi}$.  When $q$ is fixed and  $n \to +\infty$, the statistic $W$ follows asymptotically a  $\chi_d^2$ distribution  with relative error of order $O(n^{-1})$ under $H_\psi$.

Higher-order improvements of likelihood inference for a  vector parameter of interest are available. A first proposal is the Bartlett correction \citep{bartlett:1937,lawley1956}, which rescales the log-likelihood ratio statistic by its expectation under the null hypothesis, i.e.
 \begin{eqnarray}
   W_{BC}=\frac{d}{E(W)}W, \nonumber
 \end{eqnarray}
 and has a $\chi^2_d$ asymptotic null distribution  with relative error of order $O(n^{-2})$ (\cite{davison:2003},  Section 7.3; \cite{mccullage:2018},  Section 7.4). 
 Since the calculation of expectation $E(W)$ is generally not feasible, Lawley \cite{lawley1956} gave an asymptotic expansion for the exact expectation $E(W)$ under $H_{\psi}$  with error of order $O(n^{-1})$. However,  the accuracy of Bartlett correction can be lost when the exact expectation is substituted with such asymptotic expansion \citep{skovgaard:2001,sartori:2014}. Increasing the computational cost, it is possible to replace analytical expansions by parametric bootstrap approximations \citep[][Section 2.7]{cordeiro:2014}. In the present framework, $E(W)$ can be computed exactly and  the condition for validity of the $\chi^2_d$ approximation for the distribution of the Bartlett correction of $W$ is  $ p =o(n^{2/3})$ \citep{He2020B}.   The quantities needed for the Bartlett correction can be found in \cite{jiang:2013}.

Starting from the extremely accurate $r^*$ statistic for inference on a scalar $\psi$ \citep{barndorff1986}, for  a vector parameter of interest \cite{skovgaard:2001} proposed   two   modifications of $W$  designed to maintain high accuracy in the tails of the distribution:
\begin{eqnarray}
W^{*} = W\left(1-\frac{1}{W}\log \gamma \right)^2 \quad \text{and} \quad W^{**} = W-2\log \gamma.
\label{skovgaard's statistics}
\end{eqnarray}
The statistics $W^*$ and $W^{**}$   are  generally easier to calculate than the Bartlett correction and, under standard regularity conditions \citep[see, e.g.,][Section 3.4]{Severini:2000}, they are  also approximately distributed as $\chi_d^2$ when the null hypothesis holds.  Even though the relative error is of order $O(n^{-1})$, as for $W$, they exhibit higher accuracy due to large-deviation properties of the saddlepoint approximation involved in their derivation \citep{skovgaard:2001}.  Among the two forms, $W^{*}$ has the advantages of being always non-negative and  of reducing to the square of Barndorff-Nielsen's $r^*$ statistic when $d=1$ \citep{skovgaard:2001}.  
The general expression of the correction factor $\gamma$ in exponential families \citep[][(13)]{skovgaard:2001} simplifies to
\begin{eqnarray}
\gamma =\frac{\left\{(s-s_\psi)^T J_{\varphi\varphi}(\hat{\varphi}_\psi)^{-1} (s-s_\psi)\right\}^{d/2}}{W^{d/2-1}(\hat{\varphi}-\hat{\varphi}_\psi)^T (s-s_\psi)} \left\{\frac{|{J}_{\varphi \varphi}(\hat{\varphi}_\psi)|}{|{J}_{\varphi \varphi}(\hat{\varphi})|}\right\}^{1/2},
\label{skovgaard's correction factor}
\end{eqnarray}
where $s_\psi$ is the expected value of the sufficient statistic $s$ under $H_\psi$ and  $J_{\varphi \varphi} (\varphi)=- \partial^2 \ell(\varphi;s)/ \partial \varphi $ $\partial \varphi^T $ is the  observed information matrix for $\varphi$, which coincides with the expected Fisher information matrix since $\varphi$ is the canonical parameter. In order to calculate the 
$p$-value, the quantity (\ref{skovgaard's correction factor}) is evaluated at $s=s^0 = 0_q$, corresponding to $y=y^0$.


 \subsection{Directional tests in linear exponential families}\label{section2.2}
Directional tests for a vector parameter of interest in exponential family models were considered  in \cite{skogaard:1988}, \cite{fraser1994} and \cite{sartori:2014}.  In particular, the latter proposed to compute the directional $p$-value via one-dimensional integration. Directional tests in linear exponential families are essentially developed in two dimension-reduction steps, since the sufficient statistic has the same dimension of the canonical parameter  $\varphi$. Specifically, the first step consists of  reducing the  dimension of the sufficient statistic  from $q$ to the dimension of the parameter of interest $d$; indeed,  the conditional distribution of the component relative to $\psi$ of the sufficient statistic in (\ref{log-likelihood: s1 and s2}), $s_1$, given the component relative to $\lambda$, $s_2$, can be accurately approximated using saddlepoint approximations. The second step further reduces the $d$-variate conditional distribution to a one-dimensional conditional distribution given the direction indicated by the observed data point. We review here the key methodological steps, already detailed in \cite{sartori:2014}, to derive the directional $p$-value in linear exponential families.

The simplicity of exponential families makes conditional inference a practicable strategy. In particular, the theory  guarantees that the conditional distribution  of the component of interest $s_1$ of the canonical sufficient statistic,  given $s_2$, depends only on $\psi$ \citep[see, e.g.,][Theorem 5.6]{salvan:1997}. Indeed, we have
 \begin{eqnarray}
  f(s_1 | s_2;\psi) &=& \exp\left\{ \psi^T s_1 - {K}_{s_2}(\psi)\right\} {h}_{s_2}(s_1),\nonumber
\end{eqnarray}
where the cumulant generating function ${K}_{s_2}(\psi)$ and the marginal density ${h}_{s_2}(s_1)$ depend on the conditioning value  $s_2$  and can rarely be  derived explicitly. However, as in \cite{sartori:2014}, a saddlepoint approximation  \citep[see, e.g.,][Section 10.10]{salvan:1997} can be used instead.
Under the null hypothesis $H_\psi$,  the saddlepoint approximation to the density of $s_1$ given $s_2$ is expressed as 
\begin{eqnarray}
  h (s;\psi) = c \exp\left[ \ell(\hat{\varphi}_{\psi};s) - \ell\{\hat{\varphi}(s);s\}\right]  |J_{\varphi \varphi} \{\hat{\varphi}(s);s\}|^{-1/2}, \quad s \in L_{\psi}^0,
  \label{saddle L0}
\end{eqnarray}
where $c$ is a normalizing constant and $L_{\psi}^0$ is a $d$-dimensional plane defined by setting $s_2$ to its observed value, i.e. $s_2={{0}_{q-d}}$.  All values of $s$ in $L_{\psi}^0$ have the same constrained maximum likelihood estimate $\hat{\lambda}_\psi = \hat{\lambda}_\psi^0 $, while they have unconstrained maximum likelihood estimate $\hat{\varphi}(s)$.

We construct a directional test for $H_{\psi}$ by considering the  one-dimensional model based on the magnitude of $s$, $||s||$, conditional on its direction. This is done by  defining a line $L_{\psi}^*$ in $L^0_\psi$ through the observed value of $s$, $s^0 = {{0}_q}$,  and the expected value of $s$ under $H_\psi$, $s_{\psi}$, which depends on the  observed data point $y^0$, i.e. 
\begin{eqnarray}
s_{\psi} = - \ell_{\varphi}^0 \left(\hat{\varphi}_{\psi}^0 \right) = \left\{
-  \ell_{\psi}^0 \left(\hat{\varphi}_{\psi}^0 \right)^T, {{0}_{q-d}^T}\right\}^T.
\label{expected value }
\end{eqnarray}

We parameterize this line by $t \in  \mathbb{R}$, namely $s(t) = s_{\psi} + t (s^0 - s_{\psi})$. In particular, $t=0$ and $t=1$ correspond, respectively, to the expected value $s_\psi$ and to the observed value $s^0$.
The conditional distribution of $||s||$  given the unit vector $s/||s||$ is obtained from (\ref{saddle L0}) by a change of variable from $s$ to $(||s||,s/||s||)$. The Jacobian of the transformation is proportional to $t^{d-1}$. 
The directional  $p$-value to measure the departure from $H_{\psi}$ along the line $L_{\psi}^*$  is defined as the probability that $s(t)$ is as far or farther from $s_\psi$ than is the observed value $s^0$. In mathematical notation,  
\begin{eqnarray}
    p(\psi) = \frac{\int_{1}^{t_{sup}} {t^{d-1}h\{s(t);\psi\}} dt}{\int_{0}^{t_{sup}} {t^{d-1}h\{s(t);\psi\}}dt},
    \label{directed p-value}
\end{eqnarray}
 where the denominator is a normalizing constant. See \citep[][Section 3.2]{sartori:2014} for more details.  The upper limit of the integrals in (\ref{directed p-value}) is the largest value of $t$ for which the maximum likelihood estimate $\hat{\varphi}(t)$  corresponding to $s(t)$ exists; depending on the case, it can be found analytically or approximated numerically.
 The scalar integrals in (\ref{directed p-value}) can be accurately computed via numerical integration.  The error in (\ref{directed p-value}) is therefore essentially given by the error from the saddelpoint approximation used in (\ref{saddle L0}). Some results on such an error when $p$ increases with $n$ are given in \cite{tang2021laplace}. However, in all settings described in Section \ref{section:multivariate normal distribution},  (\ref{saddle L0}) holds exactly, up to the normalizing constant $c$. Thus, since $c$ simplifies in the ratio (\ref{directed p-value}), also the directional $p$-value  is exact.  The results are formally derived in Section \ref{section3:main results}.


 \subsection{Multivariate normal distribution}\label{section:multivariate normal distribution}
Let $y_1,\dots,y_n$ be a sample of independent observations from a multivariate normal distribution $N_p(\mu,\Lambda^{-1})$, where both the mean vector $\mu$ and the concentration matrix $\Lambda$, symmetric and positive definite, are unknown.  Let $y=[y_1 \cdots y_n]^T$ denote the $n \times p$ data matrix and $\tr(M)$ denote the trace of a matrix $M$. Define by $\VEC (M)$  the operator which transforms a matrix $M$ into a vector by stacking its columns one underneath the other. For a symmetric matrix $M$ it is useful to consider also $\text{vech}(M)$, which is obtained from  $\VEC(M)$ by eliminating all supradiagonal elements of $M$.  The two operators satisfy $D_p \text{vech}(M) = \VEC (M)$, where $D_p$ is the duplication matrix \citep[][Section 3.8]{magnus1999}. The log-likelihood for the parameter $\theta =\{\mu^T ,\text{vech}( \Lambda^{-1})^T\}^T$  is
\begin{eqnarray}
\ell(\theta;y) &=& \mu^T \Lambda y^T {{1}}_n-\frac{1}{2} \tr(\Lambda y^Ty)+\frac{n}{2}\log |\Lambda|  -\frac{n}{2}\mu^T\Lambda\mu.\nonumber
\end{eqnarray}

The canonical parameter in this exponential family model is given by $\varphi = \{\xi^T ,$ $ \text{vech}(\Lambda)^T\}^T $ $=\{\mu^T \Lambda, \text{vech} ( \Lambda)^T\}^T$ with canonical sufficient statistic $u=\{n\bar{y}^T, $ $-\frac{1}{2} \text{vech}(y^Ty)^TD_p^T D_p \}^T$, and the corresponding log-likelihood is
\begin{eqnarray}
\ell(\varphi;y)&=&n\xi^T\bar{y} -\frac{1}{2} \tr(\Lambda y^Ty) +\frac{n}{2}\log |\Lambda|  -\frac{n}{2}\xi^T\Lambda^{-1}\xi \nonumber \\
&=& \xi^T n\bar{y} - \text{vech} ( \Lambda)^T \left\{\frac{1}{2}D_p^T D_p \text{vech}(y^Ty)\right\} +\frac{n}{2}\log |\Lambda|  -\frac{n}{2}\xi^T\Lambda^{-1}\xi, \nonumber
\end{eqnarray}
where $\bar{y}=y^T{{1}}_n/n$ with ${1}_n$ a $n$-dimensional vector of ones. 
The score function with respect to the canonical parameter $\varphi$  is 
\begin{eqnarray}
\ell_\varphi (\varphi) &=& \left\{\ell_{\xi}(\varphi)^T,\ell_{\text{vech} (\Lambda)} (\varphi)^T \right\}^T \nonumber \\
&=& \left\{n\bar{y}^T-n\xi^T\Lambda^{-1},\; \frac{n}{2} \text{vech} \left( \Lambda^{-1} - y^T y /n  + \Lambda^{-1} \xi \xi^T \Lambda^{-1}\right)^T \right\}^T.\nonumber
\end{eqnarray}
The maximum likelihood estimates for $\mu$ and $\Lambda^{-1}$ are $\hat{\mu} = \bar{y}$ and $\hat{\Lambda}^{-1} = y^Ty/n - \bar{y}\bar{y}^T$, respectively; thus, $\hat{\xi} = \hat{\Lambda}\hat{\mu}$.  Moreover,  the observed information matrix for components $\xi$ and $\text{vech} (\Lambda)$ of $\varphi$ can be written in block form as 
\begin{eqnarray}
 J_{\varphi \varphi} (\varphi)=  \left[\begin{array}{cc}
   n\Lambda^{-1}   & \quad-n(\xi^T \Lambda^{-1} \otimes \Lambda^{-1} ) D_p \\
   -n D_p^T (\Lambda^{-1} \xi \otimes \Lambda^{-1})    & \quad \frac{n}{2} D_p^T \{ \Lambda^{-1} (\bI_p + 2\xi \xi^T \Lambda^{-1}) \otimes \Lambda^{-1}\} D_p
 \end{array} \right],\nonumber
\end{eqnarray}
 where $ \otimes $ denotes the Kronecker product  \citep[see, e.g.,][Section 5.1]{Lauritzen:1996}. Finally,  the determinant of the observed information matrix, appearing  in  (\ref{saddle L0}), satisfies $ |J_{\varphi \varphi} (\varphi)| \propto |\Lambda^{-1}|^{p+2}$ (see Supplementary Material S1.1).




\section{Directional test for multiple-sample hypotheses}\label{section3:main results}

We consider now testing two hypotheses on the parameters of the multivariate normal model presented in Section \ref{section:multivariate normal distribution}. In particular, we concentrate here on: (I) equality of  covariance matrices in $k$ independent groups; (II) equality of multivariate normal distributions in $k$ independent groups.  We also obtained analogous theoretical results for four one-sample hypothesis about:  (III) sphericity of the covariance matrix; (IV) block-independence; (V)  complete-independence;  (VI) specified values for the mean vector and the covariance matrix. The detailed results for cases (III)-(VI) are available in the Supplementary Material S2.
In all hypotheses, it is shown that the saddlepoint approximation (\ref{saddle L0}) is exact, consequently leading to  an exact directional $p$-value, up to the error from the scalar numerical integrations in (\ref{directed p-value}).




\subsection{Testing the equality of  covariance matrices in $k$ independent groups} \label{case4}

Suppose $y_{i1},\dots,y_{in_i} $,  for $i 
\in \{1,\dots,k\}$, $k\ge2$, are independent realizations of $N_p(\mu_i, \Lambda^{-1}_i)$.   
We focus on testing the null hypothesis
\begin{eqnarray}
 H_{\psi}: \Lambda_1 =\dots = \Lambda_k.
\label{hypothesis case4}
\end{eqnarray}

In the following, with a slight abuse of notation,
 let $y_i$ denote the $n_i \times p$ data matrix of the $i$-th group. We then have $\bar{y}_i=y_{i}^T  \text{1}_{n_i}/n_i$ and   $A_i = y_i ^T y_i -n_i\bar{y}_i\bar{y}_i^T$. The unconstrained maximum likelihood estimates for all $i \in \{1,\dots,k\}$ are  $\hat{\mu}_i = \bar{y}_i$  and $ \hat{\Lambda}_i^{-1} = A_i/n_i $;  the constrained maximum likelihood estimates are instead  $\hat{\mu}_{0i} =\bar{y}_i$ and $ \hat{\Lambda}_0^{-1} = \sum_{i=1}^k A_i/n$, where $n=\sum_{i=1}^k n_i$. Bartlett \cite{bartlett:1937} suggested to use the modified maximum likelihood estimator of $\Lambda^{-1}$,
  that is $\tilde{\Lambda}_i^{-1} =  A_i/(n_i-1)$ and $\tilde{\Lambda}_0^{-1}  = \sum_{i=1}^k A_i/(n-k)$.
 The modified log-likelihood ratio statistic is then equal to
\begin{eqnarray}
\tilde{W}  = \sum_{i=1}^k -(n_i-1) \log |\tilde{\Lambda}_i^{-1} \tilde{\Lambda}_0 |. \nonumber
\label{LRT case4}
\end{eqnarray}
The null distribution of $\tilde{W}$ is approximately  $\chi^2_d$ with $d=kp(p+1)/2-p(p+1)/2= p(p+1)(k-1)/2$ if and only if $p=o(n_i^{1/2})$, and the analogous condition for the Bartlett correction is  $ p =o(n_i^{2/3})$, for all $i \in \{1,\dots,k\}$ with finite $k$ \citep{He2020B}.  The expression for Skovgaard's modifications \cite{skovgaard:2001} can be found in Supplementary Material S1.2.


For the directional $p$-value, under $H_\psi$,  the expectation of $s$ has components
\begin{eqnarray}
 -\left\{{{0}}_p^T,\; \frac{n_i}{2}\text{vech}  \left( \hat{\Lambda}_0^{-1}  - \hat{\Lambda}_i^{-1}\right)^T \right\}^T, \quad i \in \{1,\dots,k\}, \nonumber
\end{eqnarray}
and the tilted log-likelihood, by group independence, can be written as  $\ell(\varphi;t)=\sum_{i=1}^k \ell_i(\varphi_i;t)$ with  the $i$-th group's contribution
\begin{eqnarray}
 \ell_i(\varphi_i;t) 
&=& n_i \xi_i^T \bar{y}_i  - \frac{n_i}{2} \tr \left[\Lambda_i\left\{\frac{y_i^Ty_i}{n_i}+(1-t)\left(\hat{\Lambda}_0^{-1}-  \hat{\Lambda}_i^{-1}\right)\right\}\right]\nonumber\\
 &&  + \frac{n_i}{2} \log |\Lambda_i| -\frac{n_i}{2} \xi_i^T \Lambda_i^{-1}\xi_i. \nonumber
\end{eqnarray}
Maximizing  the tilted log-likelihood leads to the estimates $\hat{\mu}_i(t) = \bar{y}_i$ and $\hat{\Lambda}_i(t)^{-1} = (1-t) \hat{\Lambda}_0^{-1} +t \hat{\Lambda}_i^{-1}$,  $i \in \{1,\dots,k\}$. 
Hence, the saddlepoint approximation (\ref{saddle L0}) takes the form
\begin{eqnarray}
  h\{s(t);\psi\}
  &=& c \exp\left\{ \sum_{i=1}^k \frac{n_i-p-2}{2} \log |\hat{\Lambda}_i^{-1}(t)| \right\}, \nonumber
  \label{saddle case4}
\end{eqnarray}
where $c$ is a normalizing constant. 

The value $t_{sup}$ in (\ref{directed p-value}) is the largest $t$ for which $\hat{\Lambda}_i(t)^{-1}$ is positive definite and is equal to $\{1- \mathop{\min}\limits_{1\le i \le k} \nu_{(1)}^i \}^{-1}$ where $\nu_{(1)}^i$ is the smallest eigenvalue of $\hat{\Lambda}_0 \hat{\Lambda}_i^{-1}$ (see Lemma \ref{lemmacase4} in Section \ref{section:determination of tsup}).

{Since  $\bar{y}_i$ and $\hat{\Lambda}_i^{-1}$ in the multivariate normal distribution are independent, we have that  the  saddlepoint approximation (\ref{saddle L0})  to the density of $s$ is exact, and therefore the directional $p$-value follows exactly a uniform  distribution, even in high dimensional settings with $p$ allowed to grow  with $n_i$. The exact condition for the validity of this result is given in the following theorem. }

\begin{thm}
Assume that $p=p_n$ such that $n_i \ge p+2$ for all $n_i \ge 3, i \in \{1,\dots,k\}$, with $k$ fixed. Then, under the null hypothesis $H_\psi$ (\ref{hypothesis case4}), the directional $p$-value (\ref{directed p-value}) is exactly uniformly distributed.
\label{theoremcase4}
\end{thm}
The proof of Theorem  \ref{theoremcase4} is given  in Appendix \ref{proof of theorem 1}.  Theorem \ref{theoremcase4} only requires $n_i \ge p+2$,  $i \in \{1,\dots,k\}$, for ensuring that the maximum likelihood estimate of the covariance matrix exists with probability one. This assumption is weaker than the condition $p/n_i \to \kappa \in (0,1]$  in \cite{jiang:2013} for the validity of their central limit theorem approximation with large $p$. Moreover, although the numder of groups $k$ is considered here as fixed, simulation results show that the accuracy of the directional test is not affected by the value of $k$ (see Supplementary Material S3).




 \subsection{Testing the equality of several multivariate normal distributions}\label{case3}

Under the same framework introduced in Section \ref{case4}, we are interested in testing whether the multivariate normal distributions in $k$ independent groups are identical, meaning
\begin{eqnarray}
 H_{\psi}: \mu_1 = \dots = \mu_k,\Lambda_1 =\dots = \Lambda_k .
\label{hypothesis case3}
\end{eqnarray}
The empirical  within-groups variance  $A/n$ and the empirical between-groups variance $B/n$ depend on the quantities  $A  = \sum_{i=1}^{k} y_i^T y_i -n_i\bar{y}_i\bar{y}_i^T$ and $ B =  \sum_{i=1}^k n_i \bar{y}_i\bar{y}_i^T -n\bar{y}\bar{y}^T$, such that $A+B =\sum_{i=1}^{k} y_i^T y_i - n\bar{y}\bar{y}^T$,
where  $ \bar{y}= \sum_{i=1}^k n_i \bar{y}_i/n$.
The unconstrained maximum likelihood estimates for all $i \in \{1,\dots,k\}$ are the same as in hypothesis (\ref{hypothesis case4}), while the constrained maximum likelihood estimates are  $\hat{\mu}_0  = \bar{y},\; \hat{\Lambda}_0^{-1} = {(A+B)}/{n}$. In this case, the log-likelihood ratio statistic is
\begin{eqnarray}
W &=& 
 n\log |\hat{\Lambda}_0^{-1}| - \sum_{i=1}^k n_i \log |\hat{\Lambda}_i^{-1}| \nonumber
\label{LRT case3}
\end{eqnarray}
and  asymptotically  has a  $\chi^2_d$ null distribution with $d=\{ p(p+1)/2+p\}(k-1)=p(p+3)(k-1)/2$,  provided that  $p=o(n_i^{1/2})$ for all $i \in \{1,\dots,k\}$.  The analogous condition for the Bartlett correction is  $ p =o(n_i^{2/3})$ for all $i \in \{1,\dots,k\}$ \citep{He2020B}. 
The expression for the modification of the likelihood ratio statistic of \cite{skovgaard:2001} can be found in the Supplementary Material S1.2. 


In order to obtain the directional $p$-value, we find the components of $s_\psi$
\begin{eqnarray}
 -\left\{n_i(\bar{y}_i- \bar{y})^T,\; \frac{n_i}{2} \text{vech}\left( \hat{\Lambda}_0^{-1}  -\frac{ y_i^T y_i}{n_i} + \bar{y}\bar{y}^T \right)^T \right\}^T, \quad i \in \{1,\dots,k\}, \nonumber
\end{eqnarray}
and the $i$-th group contribution to the tilted log-likelihood function $\ell(\varphi;t)=\sum_{i=1}^k$ $ \ell_i(\varphi_i;t)$ with
\begin{eqnarray}
 \ell_i(\varphi_i;t) 
  &=& n_i\xi_i^T \left\{t\bar{y}_i +(1-t)\bar{y}\right\} - \frac{n_i}{2} \tr \left[\Lambda_i\left\{\frac{ty_i^Ty_i}{n_i}+(1-t)\left(\hat{\Lambda}_0^{-1} + \bar{y}\bar{y}^T \right)\right\}\right]\nonumber\\
  && + \frac{n_i}{2} \log |\Lambda_i| -\frac{n_i}{2} \xi_i^T \Lambda_i^{-1}\xi_i. \nonumber
 \label{tilted case3}
\end{eqnarray}
The resulting maximum likelihood estimates from $\ell(\varphi;t)$ are $ \hat{\mu}_i(t) = (1-t)\bar{y} + t\bar{y}_i$ and $ \hat{\Lambda}_i(t)^{-1} = (1-t) \hat{\Lambda}_0^{-1} +t \hat{\Lambda}_i^{-1} +t(1-t) (\bar{y}_i-\bar{y})(\bar{y}_i-\bar{y})^T$.
Hence, the saddlepoint approximation (\ref{saddle L0}) along the line $s(t)$ is
\begin{eqnarray}
h\{s(t);\psi\}
  &=& c\exp\left\{ \sum_{i=1}^k \frac{n_i-p-2}{2} \log |\hat{\Lambda}_i^{-1}(t)| \right\}, \nonumber
  \label{saddle case3}
\end{eqnarray}
where $c$ is a normalizing constant. The value $t_{sup}$ in (\ref{directed p-value}) is the largest $t$ for which each $\hat{\Lambda}_i(t)^{-1}$ is positive definite and has to be found iteratively.  The following theorem gives conditions for the exactness of the directional $p$-value.

\begin{thm}
Assume that $p=p_n$ such that $n_i\ge p+2$ for all $n_i\ge 3, i \in \{1,\dots,k\}$, with $k$ fixed. Then, under the null hypothesis $H_\psi$ (\ref{hypothesis case3}), the directional $p$-value (\ref{directed p-value}) is exactly uniformly distributed. 
\label{theoremcase3}
\end{thm}

The proof of Theorem \ref{theoremcase3} is similar to the one of Theorem \ref{theoremcase4} and is  given in  Appendix \ref{proof of theorem 2}.

\section{Computational aspects}\label{section4:computational aspects}


\subsection{Determination of $t_{sup}$}\label{section:determination of tsup}

The upper bound  $t_{sup}$ of the integrals in formula (\ref{directed p-value})   is the largest value of $t$ such that the maximum likelihood estimate $\hat{\Lambda}^{-1}(t)$ or $\hat{\Lambda}_i^{-1}(t), i \in \{1,\dots,k\}$, is positive definite. Depending on the case, $t_{sup}$ can be found analytically or approximated numerically. For instance, we can derive Lemma \ref{lemmacase1} and Lemma \ref{lemmacase4} to compute $t_{sup}$ analytically for hypotheses (III)--(V) and (I), respectively. In particular, for hypotheses (III)--(V), we have that $t_{sup}=\{1-\nu_{(1)}\}^{-1}$, where $\nu_{(1)}$ is  the smallest eigenvalue of $\hat{\Lambda}_0 \hat{\Lambda}^{-1}$, while for hypothesis (I),  $t_{sup}=\{1-\mathop {\min }\limits_{1 \le i \le k} \nu^i_{(1)} \}^{-1}$, where $\nu_{(1)}^i$ is the smallest eigenvalue of $\hat{\Lambda}_0 \hat{\Lambda}_i^{-1}$.  On the contrary, there is no available closed form for  $t_{sup}$ when testing hypotheses (II) and (VI).  In such cases,  we need to find $t_{sup}$ by searching iteratively values of $t>1$ until matrices  $\hat{\Lambda}^{-1}(t)$ for hypothesis (VI) or $\hat{\Lambda}_i^{-1}(t), i \in \{1,\dots,k\}$ for hypothesis (II) are no longer positive definite.  

\begin{lem}
	The estimator $\hat{\Lambda}^{-1}(t)$ is positive definite if and only if all elements  $1-t+t\nu_l,l \in \{1,\dots,p\}$, are positive, where $\nu_l$ are the eigenvalues of the  matrix $\hat{\Lambda}_0 \hat{\Lambda}^{-1}$. Specifically, $\hat{\Lambda}^{-1}(t)$ is positive definite in $t \in [0,\{1-\nu_{(1)}\}^{-1}]$,  where $\nu_{(1)}$  is the smallest eigenvalue of $\hat{\Lambda}_0 \hat{\Lambda}^{-1}$. 
	\label{lemmacase1}
\end{lem}

The proof of Lemma \ref{lemmacase1} is given in Appendix \ref{proof of lemma 1}.

\begin{lem}
	The estimator $\hat{\Lambda}_i^{-1}(t), i \in \{1,\dots,k\}$, is positive definite if and only if all elements  $1-t+t\nu_l^i,l \in \{1,\dots,p\}$, are positive, where $\nu_l^i$ are the eigenvalues of the  matrix $\hat{\Lambda}_0 \hat{\Lambda}_i^{-1}$. Specifically, $\hat{\Lambda}_i^{-1}(t), i \in \{1,\dots,k\}$,  are all positive definite in   $t \in \left[0,\{1-\mathop {\min }\limits_{1 \le i \le k} \nu^i_{(1)} \}^{-1}\right]$, where $\nu^i_{(1)}$  is the smallest eigenvalue of $\hat{\Lambda}_0 \hat{\Lambda}_i^{-1}$. 
	\label{lemmacase4}
\end{lem}

The proof of Lemma \ref{lemmacase4} is given in Appendix \ref{proof of lemma 2}.


\subsection{ Numerical integration for  the directional $p$-value}
\begin{figure}[t]
	\centering
	\captionsetup{font=footnotesize}
	\subfigure{
		\begin{minipage}[b]{.46\linewidth}
			\centering
			\includegraphics[height = 7cm, width = 8cm]{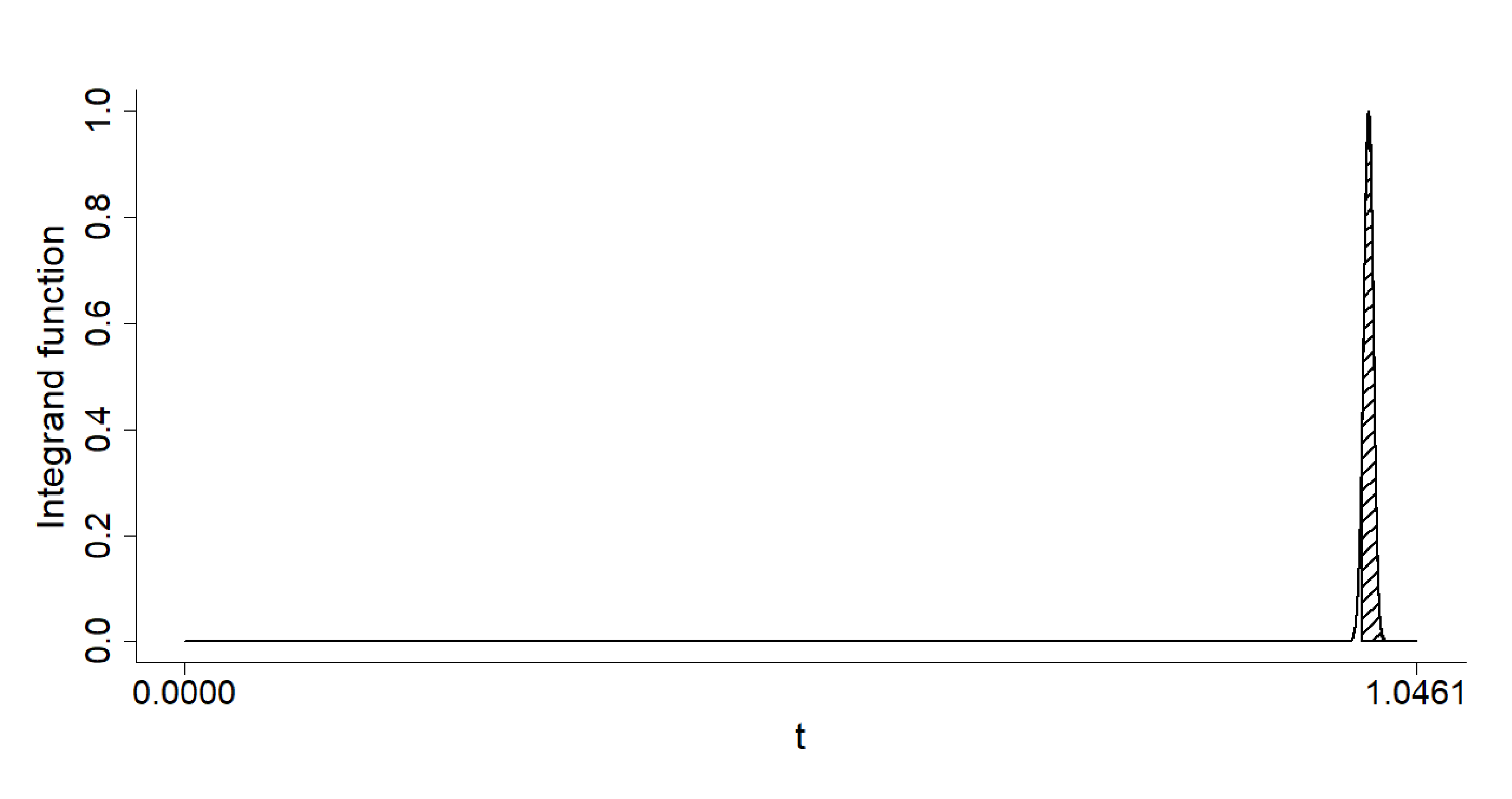}
		\end{minipage}
	}
	\subfigure{
		\begin{minipage}[b]{.46\linewidth}
			\centering
			\includegraphics[height = 7cm, width = 8cm]{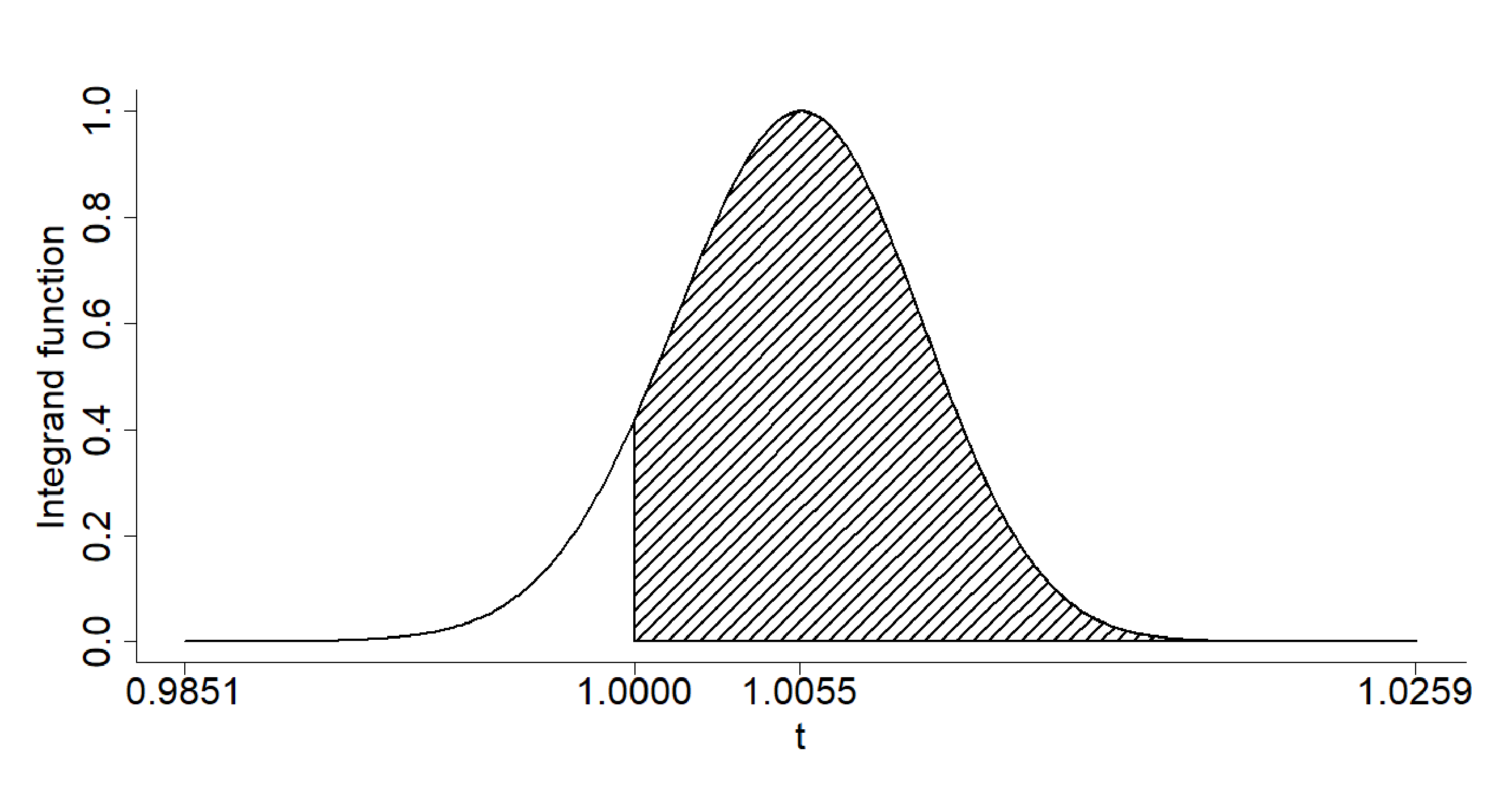}
		\end{minipage}
	}
	\caption{Integrand function $\exp \{\bar{g}(t;\psi) - \bar{g}(\hat{t};\psi)\}$ in the directional $p$-value to test the hypothesis (\ref{hypothesis case4}). The $n=100$  observations are sampled from a  $N_p({0}_p, \bI_p)$ distribution with $p=70$. The left panel refers to the  interval $[0,t_{sup}]=[0,1.046135]$, the right panel to the interval  $[t_{min},t_{max}]=[0.985117, 1.025937]$.}   
	\label{fig: numerical computation }
\end{figure}

Let $g(t;\psi) = t^{d-1}h\{s(t);\psi\} = \exp [(d-1)\log t + \log h\{s(t);\psi\}] = \exp \{\bar{g}(t;\psi)\}$ be the integrand function in formula (\ref{directed p-value}).
In order to account for situations in which $g(t;\psi)$ is numerically too small or too large, we consider  rescaling  $\bar{g}(t;\psi)$  in the interval $[0,t_{sup}]$ using $\bar{g}(\hat{t};\psi)=\mathop { \sup }\limits_{t \in [0,t_{sup}]} \bar{g}(t;\psi) $. The directional $p$-value can then be computed as
\begin{eqnarray}
p(\psi)  = \frac{\int_{1}^{t_{sup}} \exp \{\bar{g}(t;\psi) - \bar{g}(\hat{t};\psi)\} dt}{\int_{0}^{t_{sup}} \exp \{\bar{g}(t;\psi) - \bar{g}(\hat{t};\psi)\} dt}. \nonumber
\end{eqnarray}

Moreover, when the dimension $p$ is large, the integrand function often  concentrates on a very small range, meaning that it is significantly different from zero in a very small interval around $\hat{t}$.  Using the hypothesis problem (\ref{hypothesis case4}) as an illustration, in the left hand  panel of Figure \ref{fig: numerical computation } the integrand function is plotted in the  interval $[0,t_{sup}]$. We can observe that only for very few $t$ values the function is appreciably different  from zero. For a more accurate and efficient numerical integration, we can apply  the Gauss-Hermite quadrature \citep{gausshermite}, and focus on a narrower integration interval $[t_{min},t_{max}]$.  The integrand function curve in such an interval is displayed in the right hand panel of  Figure \ref{fig: numerical computation }. 
Hence, the directional $p$-value  can be well approximated by
\begin{eqnarray}
p(\psi) \;  \doteq  \; \frac{\int_{1}^{t_{max}} \exp \{\bar{g}(t;\psi) - \bar{g}(\hat{t};\psi)\} dt}{\int_{t_{min}}^{t_{max}} \exp \{\bar{g}(t;\psi) - \bar{g}(\hat{t};\psi)\} dt}. 
\label{directed-p adjust}
\end{eqnarray}
Details on the implemention of the Gauss-Hermite quadrature (\ref{directed-p adjust}) for the hypotheses considered in 
 Section \ref{section3:main results} and Section S2 in Supplementary Material are described in the Supplementary Material S1.3.


\section{Simulation studies}\label{section5:simulation study}

\subsection{Setup}\label{simulation setting}

The  performance  of the directional test  for the hypotheses of Section \ref{section3:main results} in the high dimensional multivariate normal framework is here assessed via Monte Carlo simulations based on $100,000$ replications. The exact directional test is compared with the $\chi_d^2$ approximation for the log-likelihood ratio test, its Bartlett correction, two Skovgaard's  modifications \cite{skovgaard:2001}, and with the normal approximation for the test  proposed by Jiang and Yang \cite{jiang:2013}. The  six tests are evaluated in terms of empirical distribution, empirical distribution of the corresponding $p$-values,  estimated size and power. 
 Simulation results for the hypotheses (III)--(VI) are reported in the Supplementary Material S4.

Samples of size $n_i, i \in \{1,\dots, k\}$,  are  generated from the $p$-variate standard normal distribution $N_p(0_p, \bI_p)$ under the null hypothesis.  For each simulation experiment, we show  results for $k=3$, $n_i=100$ for all $i=1,2,3$, and $p/n_i \in \{0.05,0.1,0.3,0.5,0.7,$ $0.9\}$. Additional results for different values of  $n_i$ and $p/n_i$  are reported in the Supplementary Material S5--S6. 
The various simulation setups are detailed below, partly taken from \cite{jiang:2013}.

Hypothesis (I): {testing the equality of covariance matrices in $k$ normal distributions}.     
When evaluating power, four settings are considered for the alternative hypothesis: (1)  $\Lambda_1^{-1} =   \bI_p$, $\Lambda_2^{-1} = 1.21 \bI_p$ and $\Lambda_3^{-1} = 0.81 \bI_p$; (2) $\Lambda_1^{-1} = \bI_p$, $ \Lambda_2^{-1} = \Lambda_3^{-1} = \Lambda_1^{-1} +\delta(pn_i)^{-1/2} \bI_p$; (3)  $\Lambda_1^{-1} = \bI_p$, $\Lambda_2^{-1} = \Lambda_3^{-1} = (1-\rho) \bI_p  + \rho 1_p 1_p^T$ with $\rho = \delta(pn_i)^{-1/2}$; (4) $\Lambda_1^{-1} = \bI_p$, $ \Lambda_2^{-1} = \Lambda_3^{-1} = \diag( \eta,1_{p-1}^T)$ where $\eta \in \mathbb{R}^+$. 

Hypothesis (II): {testing the equality of $k$ multivariate normal distributions}. 
When evaluating power, four settings are considered for the alternative hypothesis: (1) $\mu_1=0_p$, $\mu_2=\mu_3=0.1\cdot {{1}}_p$ and  $\Lambda_1^{-1} = 0.5 {{1}}_p{{1}}_p^T + 0.5 \bI_p$, $\Lambda_2^{-1} = 0.6 {{1}}_p{{1}}_p^T + 0.4 \bI_p$, $\Lambda_3^{-1} = 0.5 {{1}}_p{{1}}_p^T + 0.31 \bI_p$; (2) and (3) $\mu_1=0_p$, $\mu_2 =\mu_3 = \delta(pn_i)^{-1/2} {{1}}_p$;
 (4) $\mu_1=0_p$, $\mu_2 =\mu_3 = \{10 (pn_i)^{-1/2},0_{p-1}^T\}^T$, and the setup of covariance matrices of (2)--(4) as in  Hypothesis (I).

In the Supplementary Material S3 we report the empirical results for Hypotheses (I)--(II) for large  group values of $k \in \{30, 300\}$, which shows that the accuracy of the directional $p$-value does not change.


\subsection{Null distribution}\label{limiting null distribuiton: results}

\setlength{\tabcolsep}{3.5mm}{
	\begin{table}[]
		\centering
		\caption{Empirical probability of Type I error for the directional test (DT), central limit theorem test (CLT), log-likelihood ratio test (LRT), Bartlett correction (BC) and two Skovgaard's modifications \cite{skovgaard:2001} (Sko1 and Sko2, respectively) at the nominal level $\alpha = 0.05$}
		{\begin{tabular}{cccccccc}
				\hline
				Hypothesis &  $p/n_i$ &  DT & CLT & LRT &BC & Sko1 & Sko2\\
				\hline
				(I) & 0.05 & 0.050 & 0.078 & 0.062 & 0.050 & 0.048 & 0.048 \\ 
				&  0.1 & 0.049 & 0.064 & 0.102 & 0.049 & 0.041 & 0.040 \\ 
				&  0.3 & 0.051 & 0.057 & 0.950 & 0.067 & 0.010 & 0.006 \\ 
				&  0.5 & 0.050 & 0.054 & 1.000 & 0.183 & 0.000 & 0.000 \\ 
				&  0.7 & 0.050 & 0.054 & 1.000 & 0.865 & 0.000 & 0.000 \\ 
				&  0.9 & 0.049 & 0.054 & 1.000 & 1.000 & 0.065 & 0.000 \\ 
				(II) & 0.05 & 0.049 & 0.061 & 0.068 & 0.049 & 0.045 & 0.045 \\ 
				&  0.1 & 0.048 & 0.055 & 0.115 & 0.049 & 0.037 & 0.036 \\ 
				&  0.3 & 0.051 & 0.055 & 0.967 & 0.068 & 0.007 & 0.003 \\ 
				&  0.5 & 0.050 & 0.053 & 1.000 & 0.192 & 0.000 & 0.000 \\ 
				&  0.7 & 0.050 & 0.053 & 1.000 & 0.880 & 0.000 & 0.000 \\ 
				&  0.9 & 0.049 & 0.053 & 1.000 & 1.000 & 0.032 & 0.000  \\
				\hline              
			\end{tabular}
		}
		\label{table type I normal:cases}
\end{table}}

The Monte Carlo simulations for the hypotheses (I) and (II)  described  in Section \ref{section3:main results}  are here illustrated. The Type I error at level $\alpha=0.05$ based on the approximate null distribution of the various statistics is  evaluated here.
The empirical distribution of $p$-values for the six tests is examined by comparison with the Uniform$(0,1)$ distribution in the Supplementary Materials S3--S4. The limiting null distribution of the statistics is also compared with their corresponding chi-square or standard normal distribution in the Supplementary Material S3--S4.

Tables \ref{table type I normal:cases} reports the empirical Type I error  at the nominal level $\alpha = 0.05$ under the null hypothesis. 
The directional $p$-value exhibits an excellent performance in terms of the empirical  Type I error, not needing essentially correction  over the different choices of $p$ as suggested by the theory in Section \ref{section3:main results}.  
 In this respect, it is significantly better than that of the central limit theorem test of Jiang and Yang \cite{jiang:2013},  which has a slightly liberal empirical Type I error. In addition,  the four statistics with chi-square approximate distributions  are not  very accurate, and even remarkably unreliable  with increasing $p$. 
 This behavior confirms the results in \cite{He2020B} stating that  the chi-square approximation to the log-likelihood ratio statistic distribution applies if and only if $p=o(n_i^{1/2})$, and that  to its Bartlett corrected version if and only if $p=o(n_i^{2/3})$, $i \in \{1,\dots,k\}$, which are both instances of low dimensional asymptotic regimes.
 There is no analogous theoretical result for Skovgaard's statistics \cite{skovgaard:2001}, yet the numerical evidence suggests an intermediate condition between those of the log-likelihood ratio statistic and its Bartlett correction.
 The performance of the directional test is  stable over all scenarios, outperforming the other methods.  
 Below, we provide more details on the simulation outcomes for the empirical probability of Type I error.

 Table \ref{table type I normal:cases}  from top to bottom displays  the empirical Type I error of  hypotheses (I) and (II), respectively. The directional test is exact up to simulation error for all different choices of $p$. Instead,  the empirical Type I error of the central limit theorem test is slightly larger than the nominal level.
The chi-square approximation to the distributions of the log-likelihood test, Bartlett correction, and Skovgaard's modifications \cite{skovgaard:2001} is accurate only for small $p$, becoming completely unreliable as  $p/n_i$ increases. However, the empirical Type I error for one of Skovgaard's modifications \cite{skovgaard:2001},  somehow surprisingly, improves for the largest value $p/n_i=0.9$.

\begin{figure}[t]
	\centering
	\captionsetup{font=footnotesize}
	\subfigure{
		\begin{minipage}[b]{.3\linewidth}
			\centering
			\includegraphics[scale=0.0735]{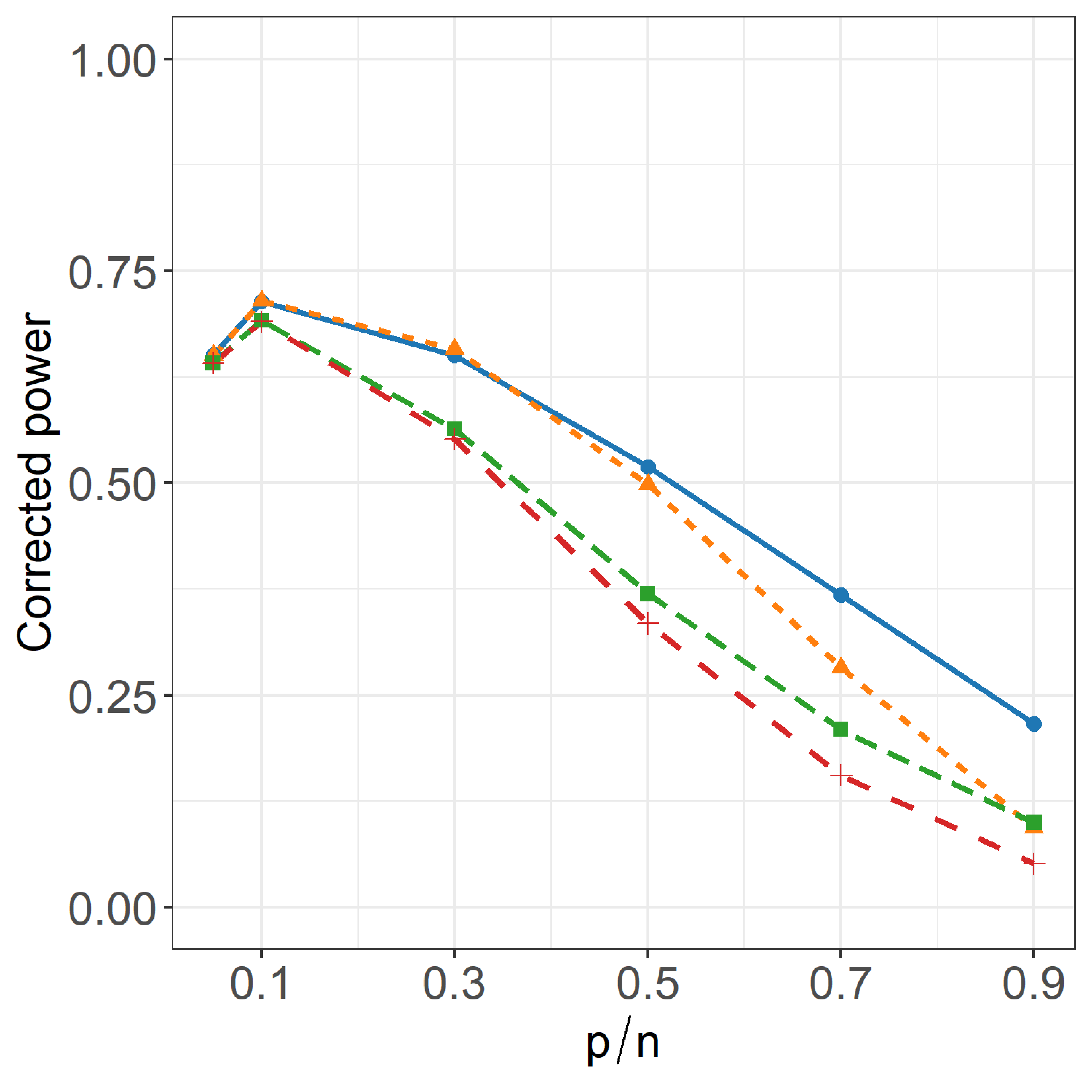}
		\end{minipage}
	}
	\subfigure{
		\begin{minipage}[b]{.3\linewidth}
			\centering
			\includegraphics[scale=0.0735]{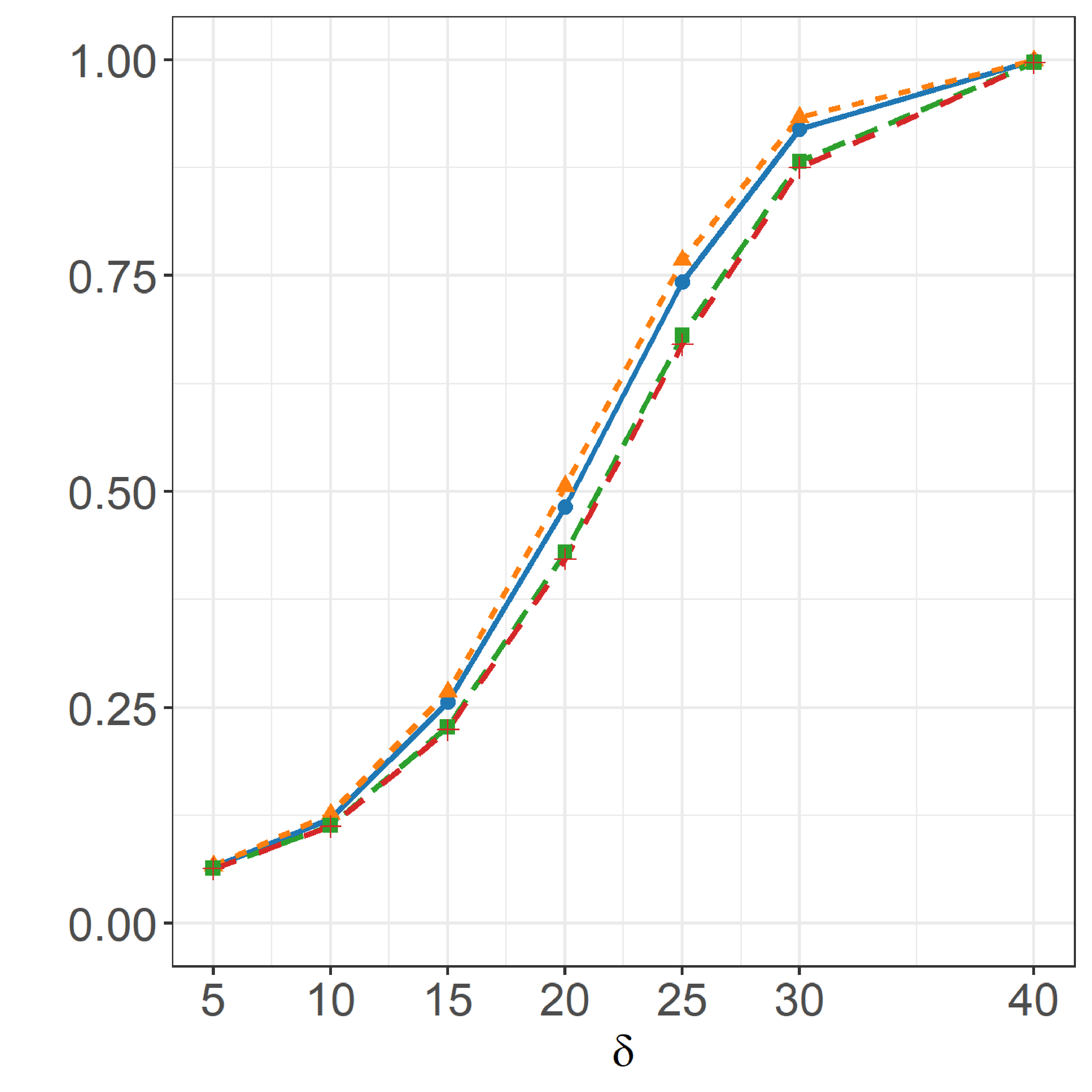}
		\end{minipage}
	}
	\subfigure{
	\begin{minipage}[b]{.3\linewidth}
		\centering
		\includegraphics[scale=0.0735]{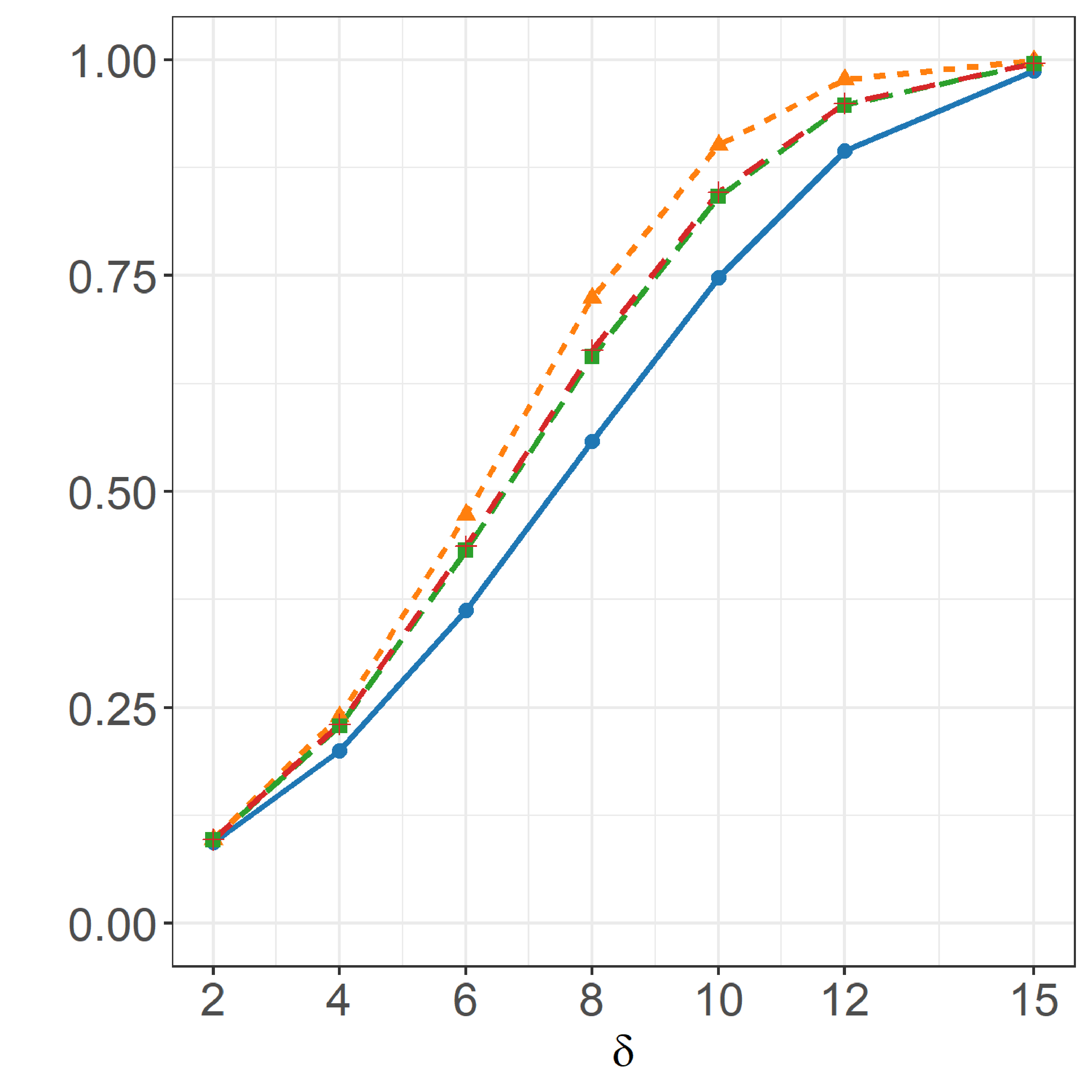}
	\end{minipage}
}
	\\
	\subfigure{
		\begin{minipage}[b]{.3\linewidth}
			\centering
			\includegraphics[scale=0.0735]{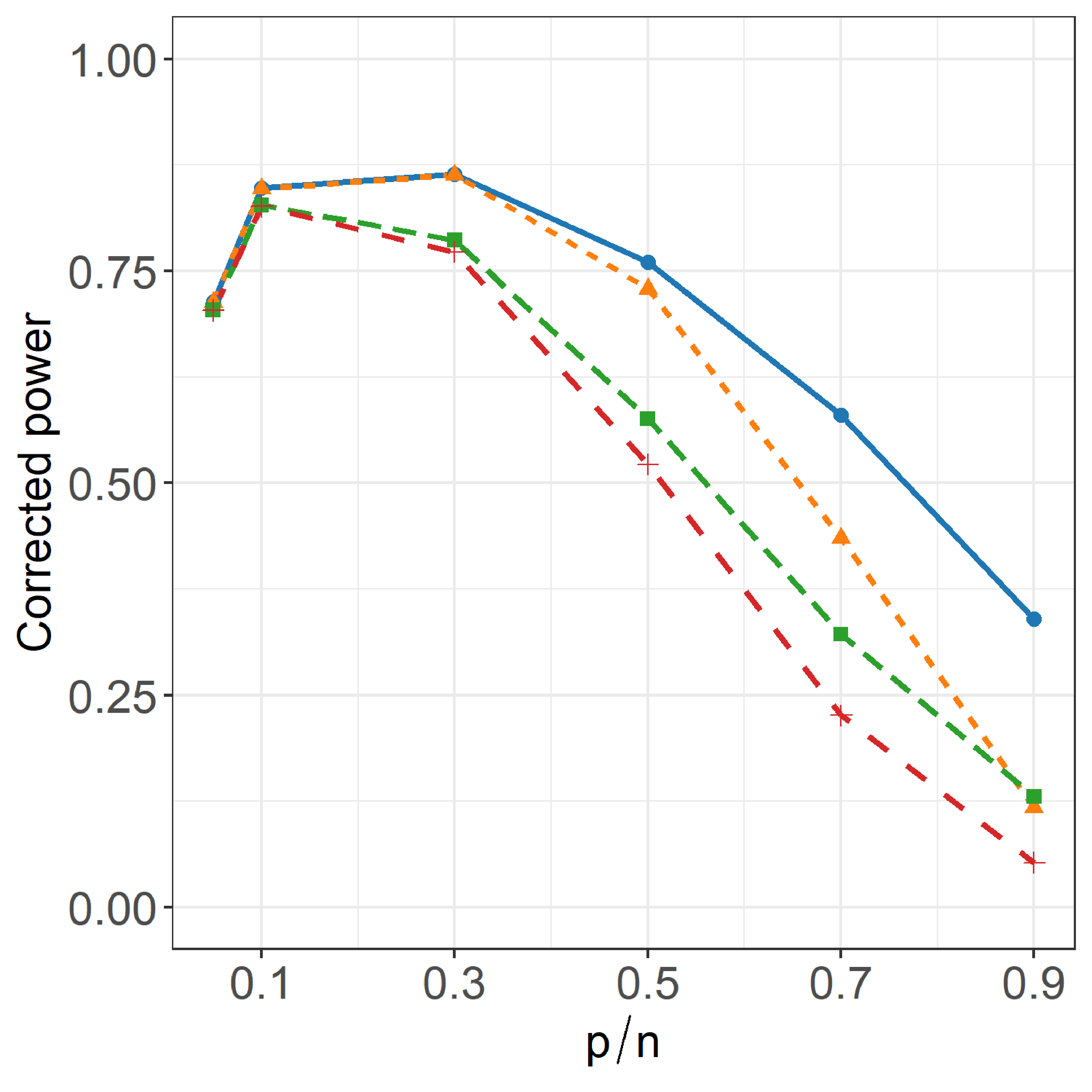}
		\end{minipage}
	}
	\subfigure{
		\begin{minipage}[b]{.3\linewidth}
			\centering
			\includegraphics[scale=0.0735]{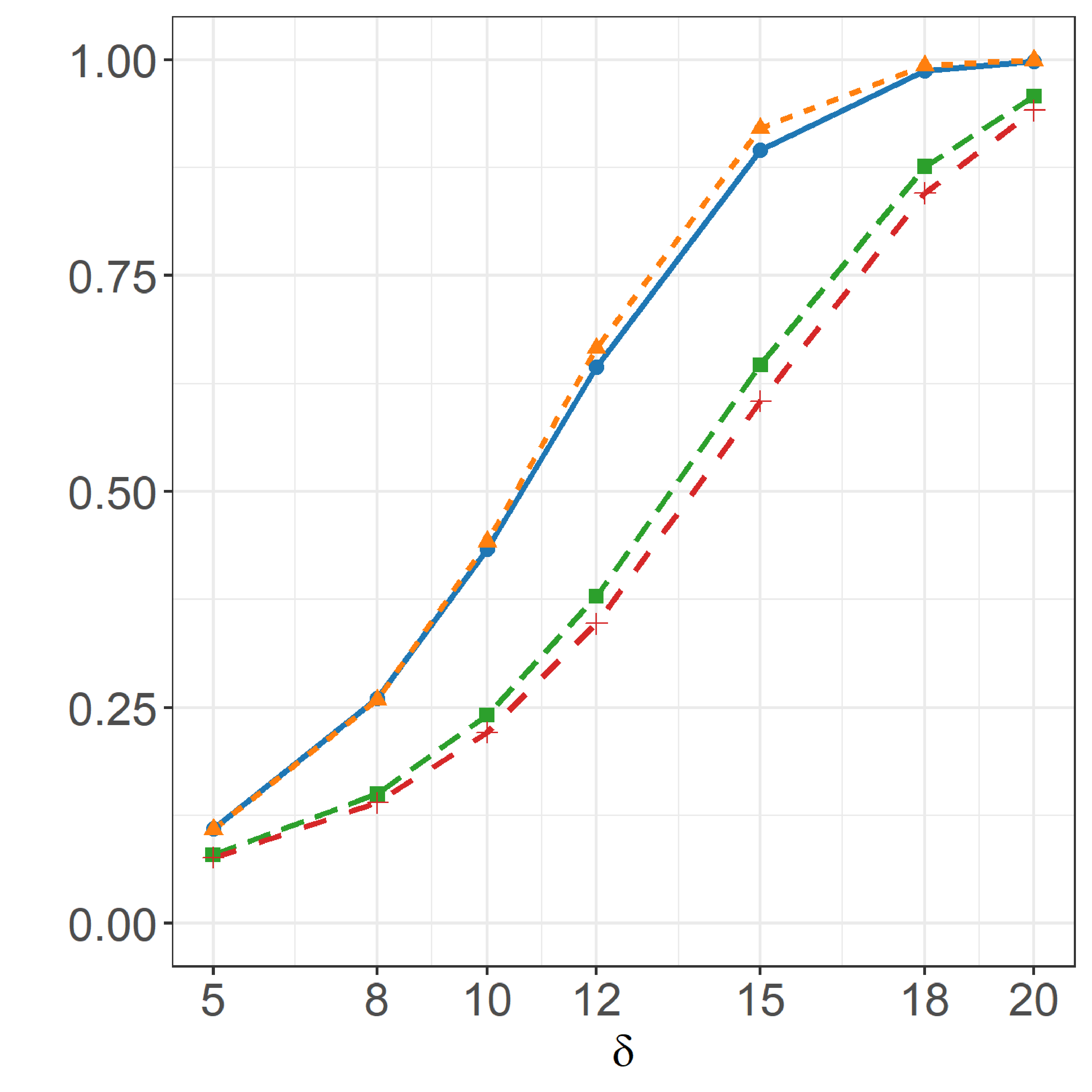}
		\end{minipage}
	}
	\subfigure{
	\begin{minipage}[b]{.3\linewidth}
		\centering
		\includegraphics[scale=0.0735]{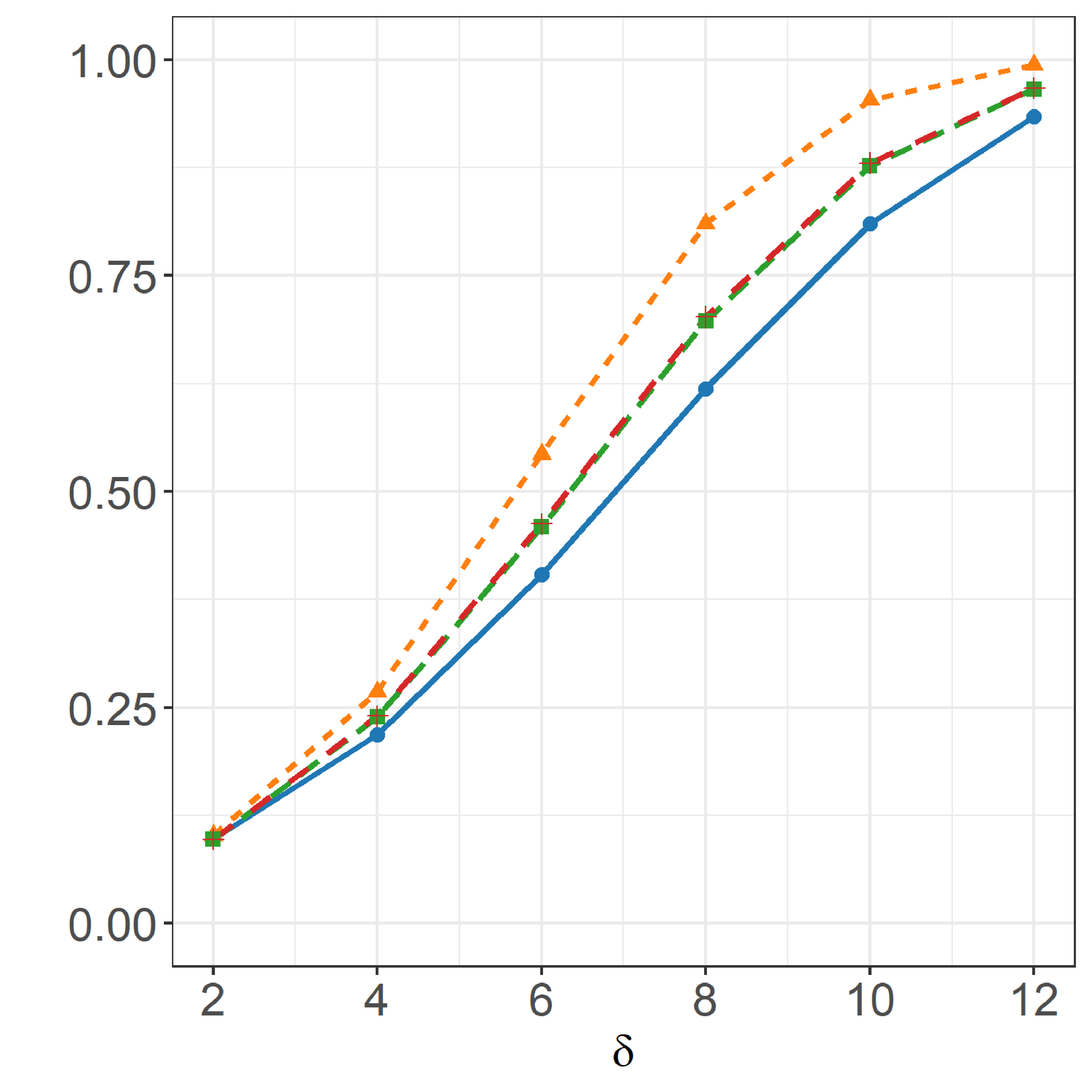}
	\end{minipage}
}
	\caption{Empirical corrected powers of four tests. The solid, dashed, long-dashed, and dot-dashed curves are the empirical power functions of the central limit theorem test, directional test and two Skovgaard's modifications \cite{skovgaard:2001}, respectively. The top and bottom rows correspond to hypotheses (I) and (II), respectively; the left, middle and right columns correspond to alternative hypothesis settings (1), (2) and (3), respectively. The two right-most columns refer to the scenario with $p/n_i=0.3$.}   
	\label{fig:power}
\end{figure}

\subsection{Empirical corrected power}
The power 
of the tests considered for the hypotheses problems in the previous section are here investigated empirically for some alternative settings. In particular, four  possible  choices for $\mu_i$ and $\Lambda_i^{-1}$ under the alternative hypotheses detailed in Section \ref{simulation setting} are studied. The first alternative setting (1) is taken from Jiang and Qi \cite{jiang:2015}, who extended the use of the central limit theorem test developed in \cite{jiang:2013} to cases where $p$ is very close to $n$ (see Section \ref{section6:disscussion} for further details). The second alternative setting (2) deals with
situations where the Frobenius norm between the null and alternative parameters  converges to zero as $n_i$ goes to infinity. The third alternative setting (3) is based on the compound symmetry structure of the covariance matrix with correlation going to zero as $n_i$ diverges, while only one group has the identity structure. The last alternative setting (4) is motivated by Jensen \cite{Jensen:2021} and considers a situation where only one or two elements of the parameter differ between the null and alternative hypotheses. Due to space constraints,  we report here results referred to the corrected power only. Corrected power is based on the corrected Type I error,  which is the $5\%$ quantile of the empirical $p$-values obtained under the null hypothesis, and is reported in the Supplementary Material S3--S4. This allows a fair comparison among the tests, since power is intended with a given significance level. However, it is important to remark that the directional $p$-value is the only approach that does not need a correction for the Type I error, being exact under the null hypothesis.  The central limit theorem, log-likelihood ratio test and Bartlett correction have the same corrected power as they use the same test statistic  $W$ and result in different cutoff values for the corrected Type I error.

 The left-most column of Figure \ref{fig:power} summarizes simulation results for the  hypotheses (I) and (II) and $p/n_i \in \{0.05$, $0.1$, $0.3$, $0.5$, $0.7$, $0.9\}$. The alternative setting for each hypothesis is the same as in \cite{jiang:2015}, where the use of the central limit theorem test was recommended. The power of the directional test across the different ratios $p/n_i$ is always greater than the nominal level $0.05$; it is comparable with the corrected power of the central limit theorem test, log-likelihood ratio test and Bartlett correction when $p$ is moderate, but it is lower otherwise. However, it must be taken into account that the log-likelihood ratio test and Bartlett correction do not  control the Type I error when $p$ is large, therefore their  power is meaningless in such scenarios. Finally, Skovgaard's tests have uniformly the  lowest corrected power.

We also investigate the local power,  i.e. how large $\delta$ in the alternative settings of Section \ref{simulation setting} needs to be so that the power can tend to 1. The middle and right columns of  Figure \ref{fig:power} display the empirical local corrected power 
of the tests for various values of $\delta$ and ratio $p/n_i=0.3$. Under the alternative setting (2), shown in the middle column, the power of the directional test is comparable or slightly superior to the corrected power of the central limit theorem test, and clearly higher than the corrected power of the Skovgaard's modifications. Under the alternative setting (3), shown in the right column, the directional test is the most powerful while the central limit theorem test has the worst power performance even after correction for Type I error.

\begin{figure}[t]
	\centering
	\captionsetup{font=footnotesize}
	\subfigure{
		\begin{minipage}[b]{.3\linewidth}
			\centering
			\includegraphics[scale=0.0735]{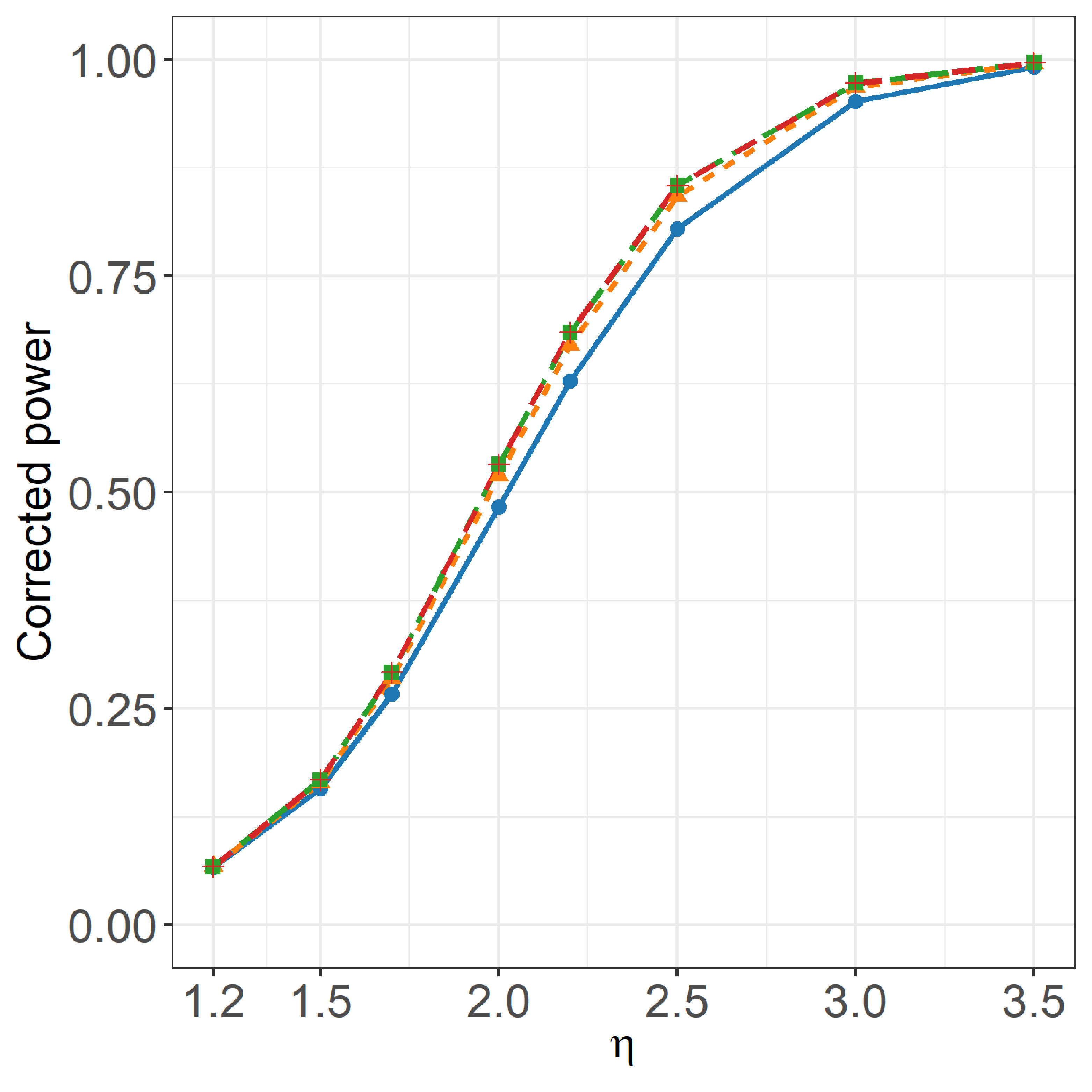}
		\end{minipage}
	}
	\subfigure{
		\begin{minipage}[b]{.3\linewidth}
			\centering
			\includegraphics[scale=0.0735]{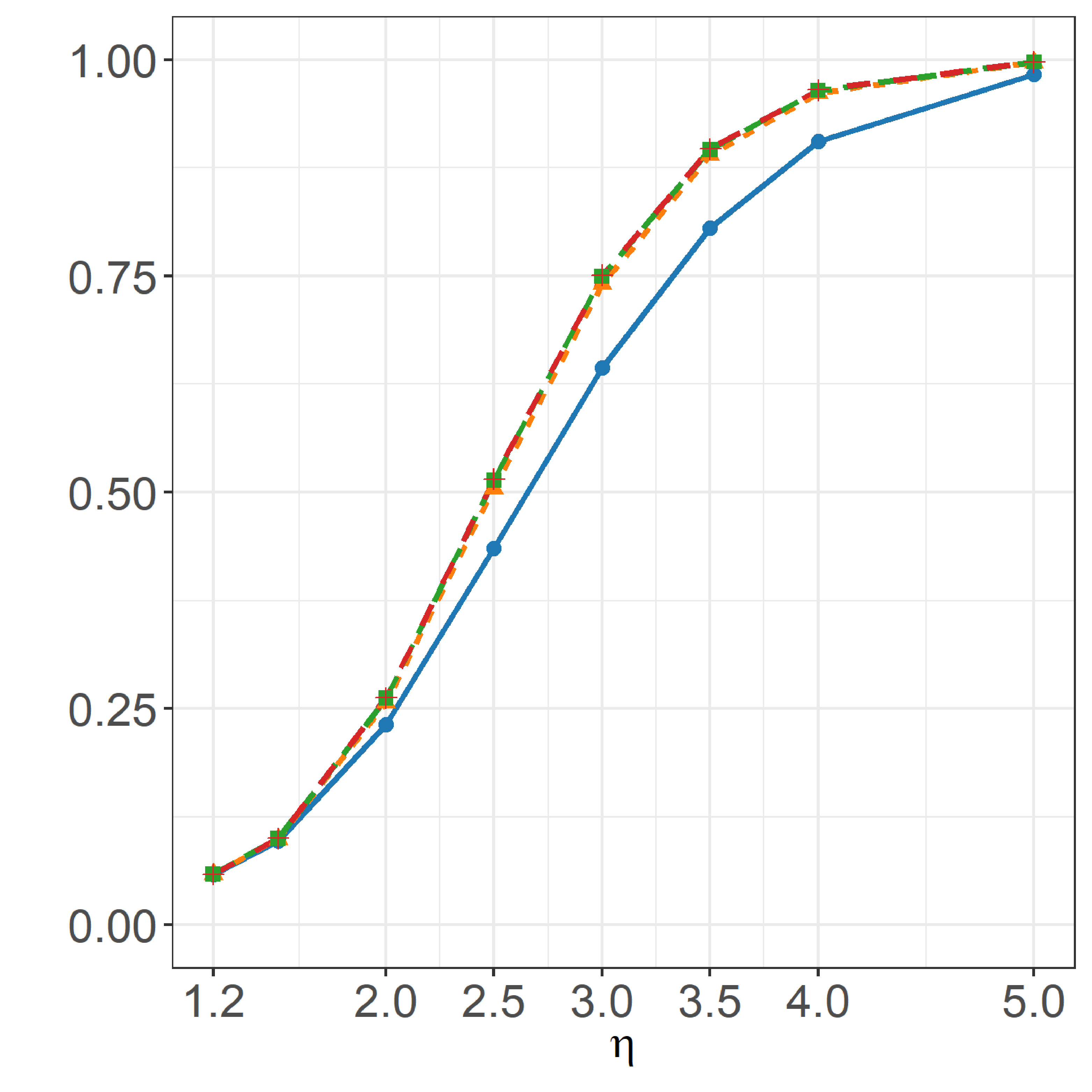}
		\end{minipage}
	}
	\subfigure{
		\begin{minipage}[b]{.3\linewidth}
			\centering
			\includegraphics[scale=0.0735]{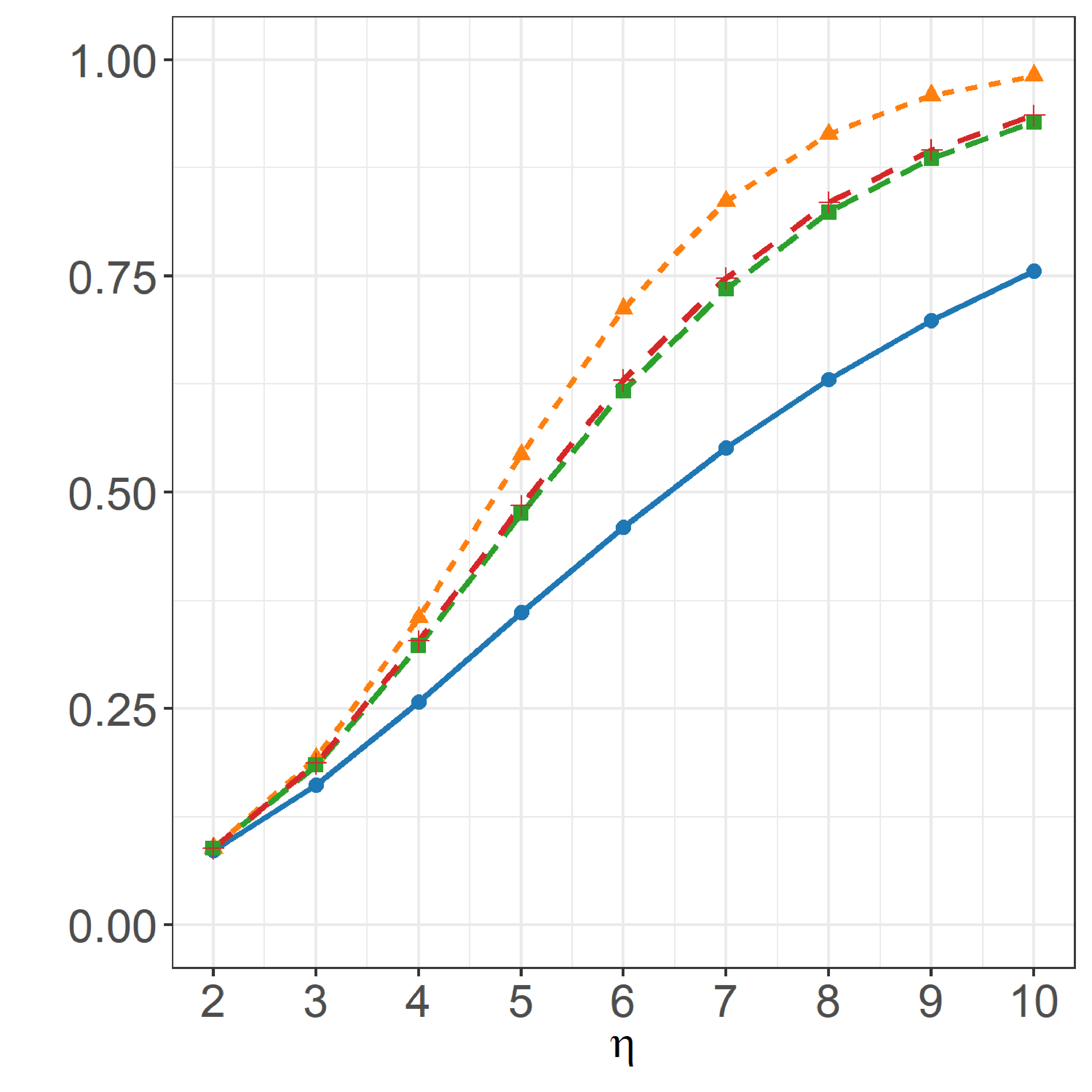}
		\end{minipage}
	}
	\subfigure{
		\begin{minipage}[b]{.3\linewidth}
			\centering
			\includegraphics[scale=0.0735]{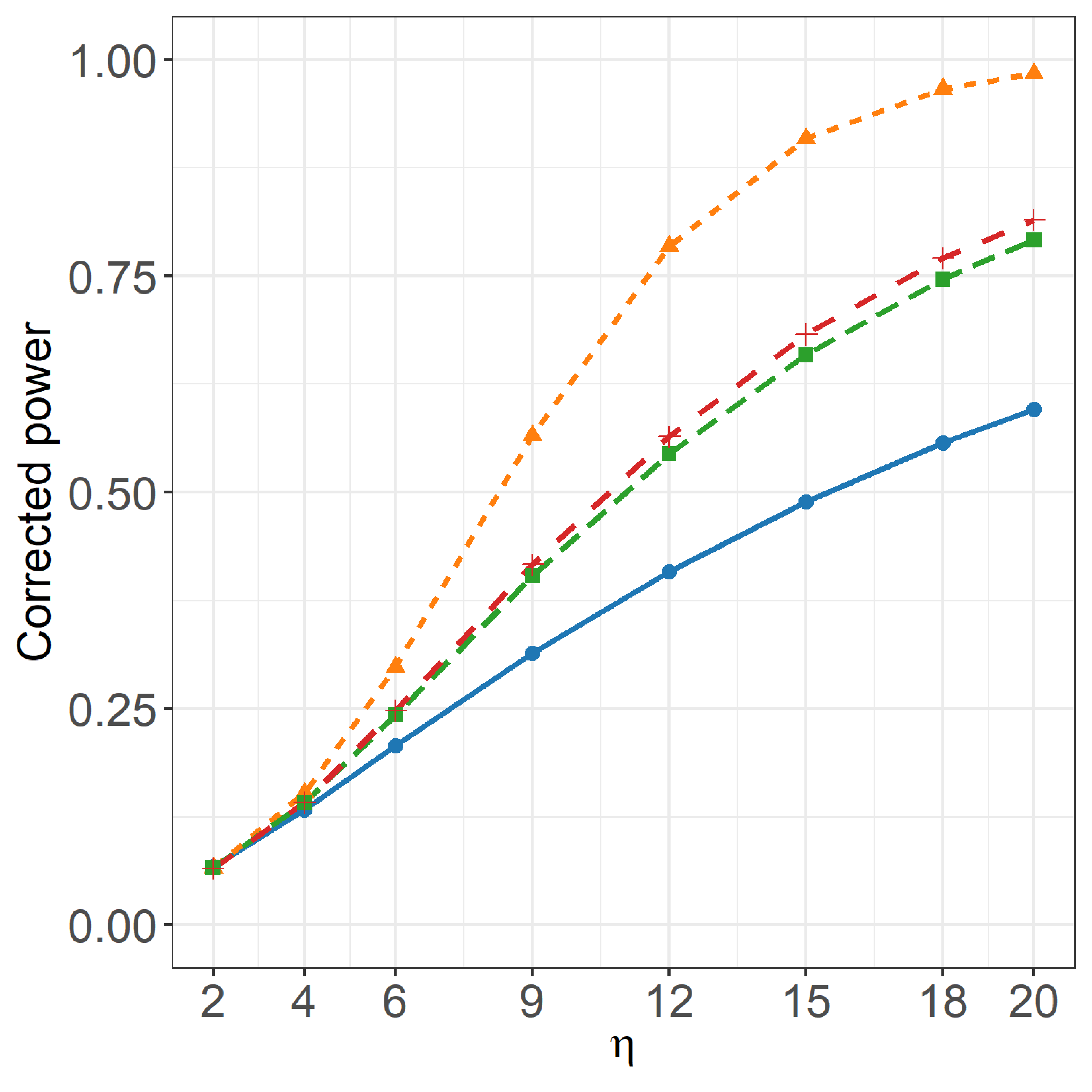}
		\end{minipage}
	}
	\subfigure{
		\begin{minipage}[b]{.3\linewidth}
			\centering
			\includegraphics[scale=0.0735]{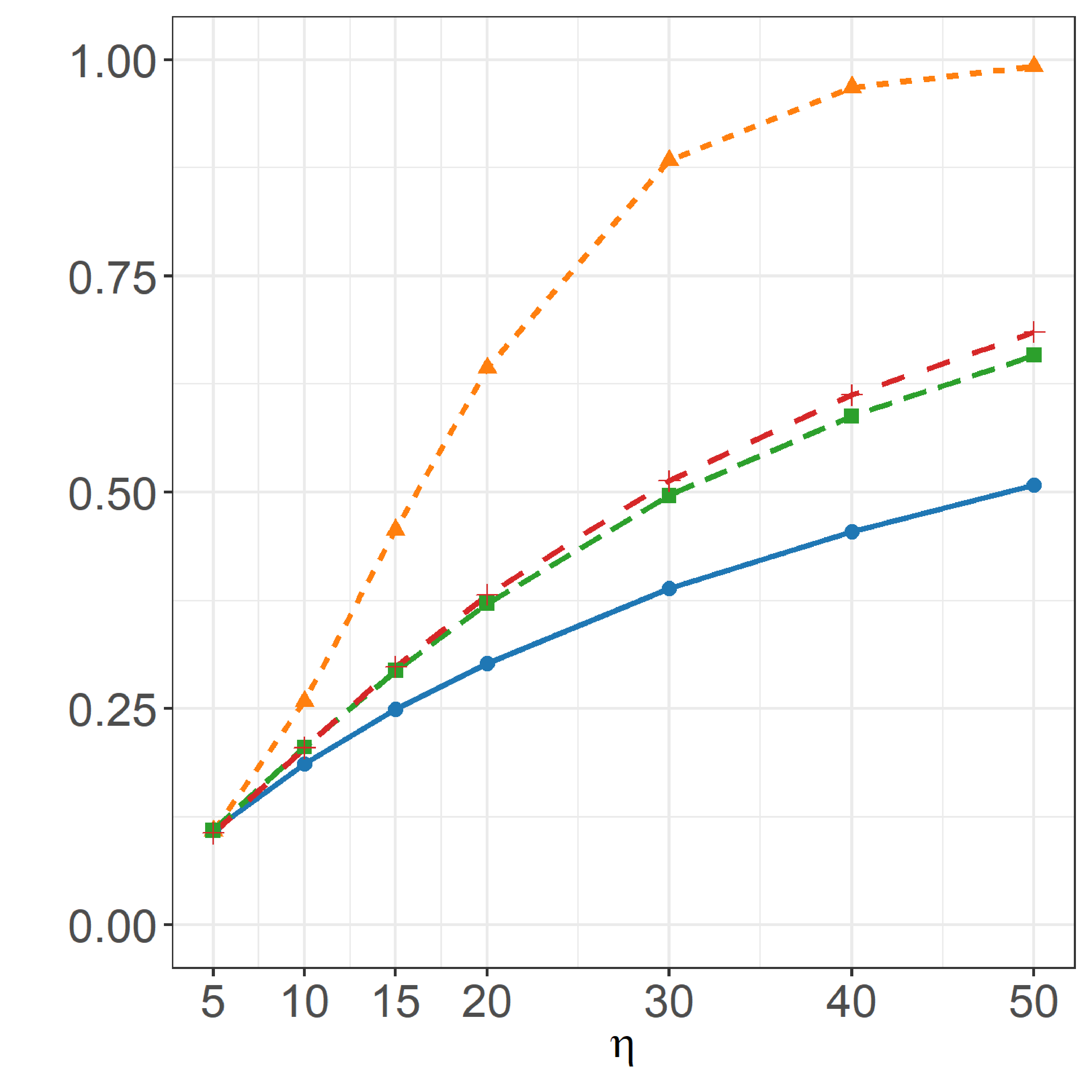}
		\end{minipage}
	}
	\subfigure{
		\begin{minipage}[b]{.3\linewidth}
			\centering
			\includegraphics[scale=0.0735]{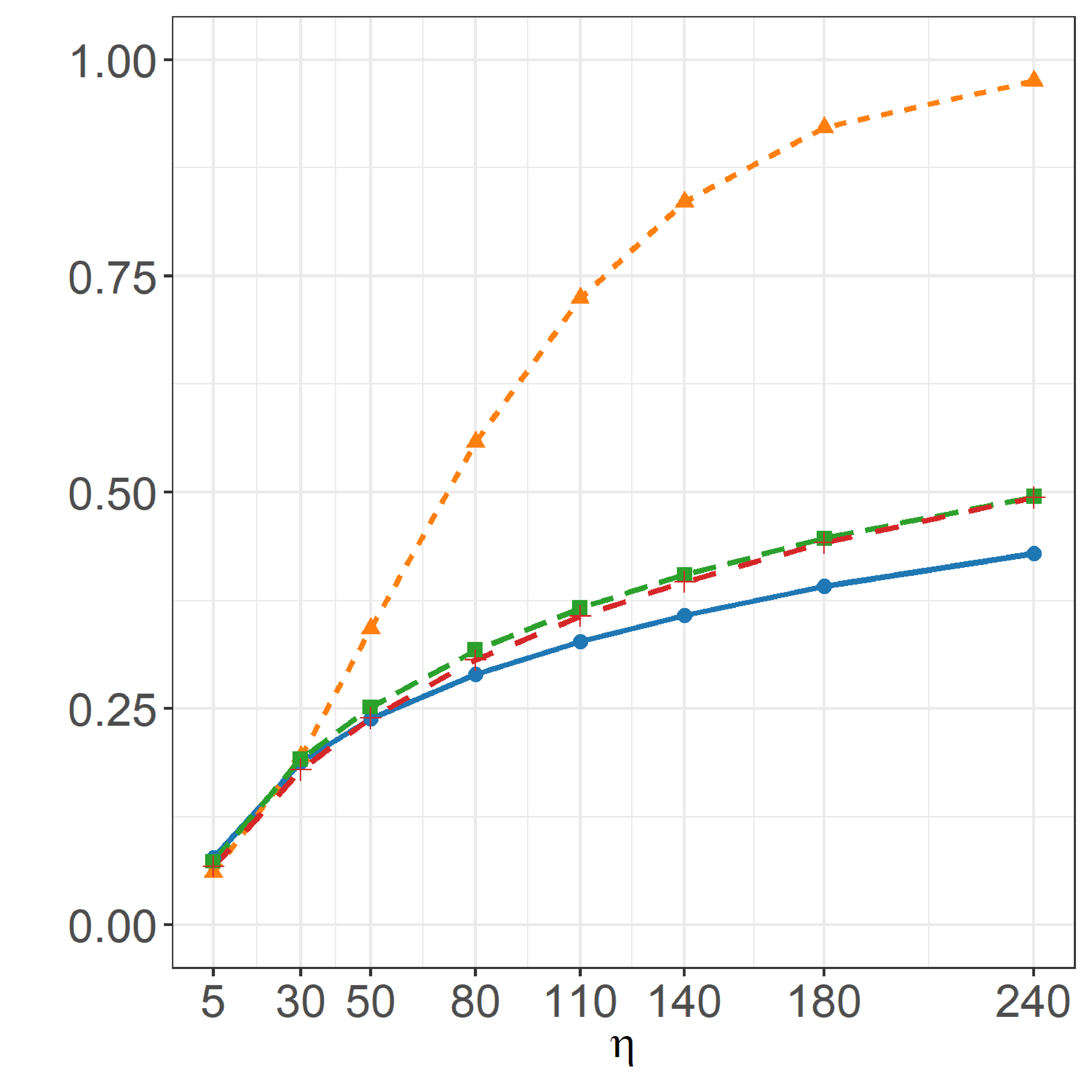}
		\end{minipage}
	}
	\caption{Empirical corrected power of four tests for hypothesis (I) with different values of $\eta$ and $p/n_i$. The solid, dashed, long-dashed, and dot-dashed curves are the empirical power functions of the central limit theorem test, directional test and two Skovgaard's modifications \cite{skovgaard:2001}, respectively.   The alternative setting (4) is given in Section \ref{simulation setting}. The six plots correspond to $p/n_i \in \{0.05, 0.1, 0.3, 0.5, 0.7, 0.9\}$, starting from top left and proceeding by row.}  
	\label{fig:power extreme case4}
\end{figure}

\begin{figure}[t]
	\centering
	\captionsetup{font=footnotesize}
	\subfigure{
		\begin{minipage}[b]{.3\linewidth}
			\centering
			\includegraphics[scale=0.0735]{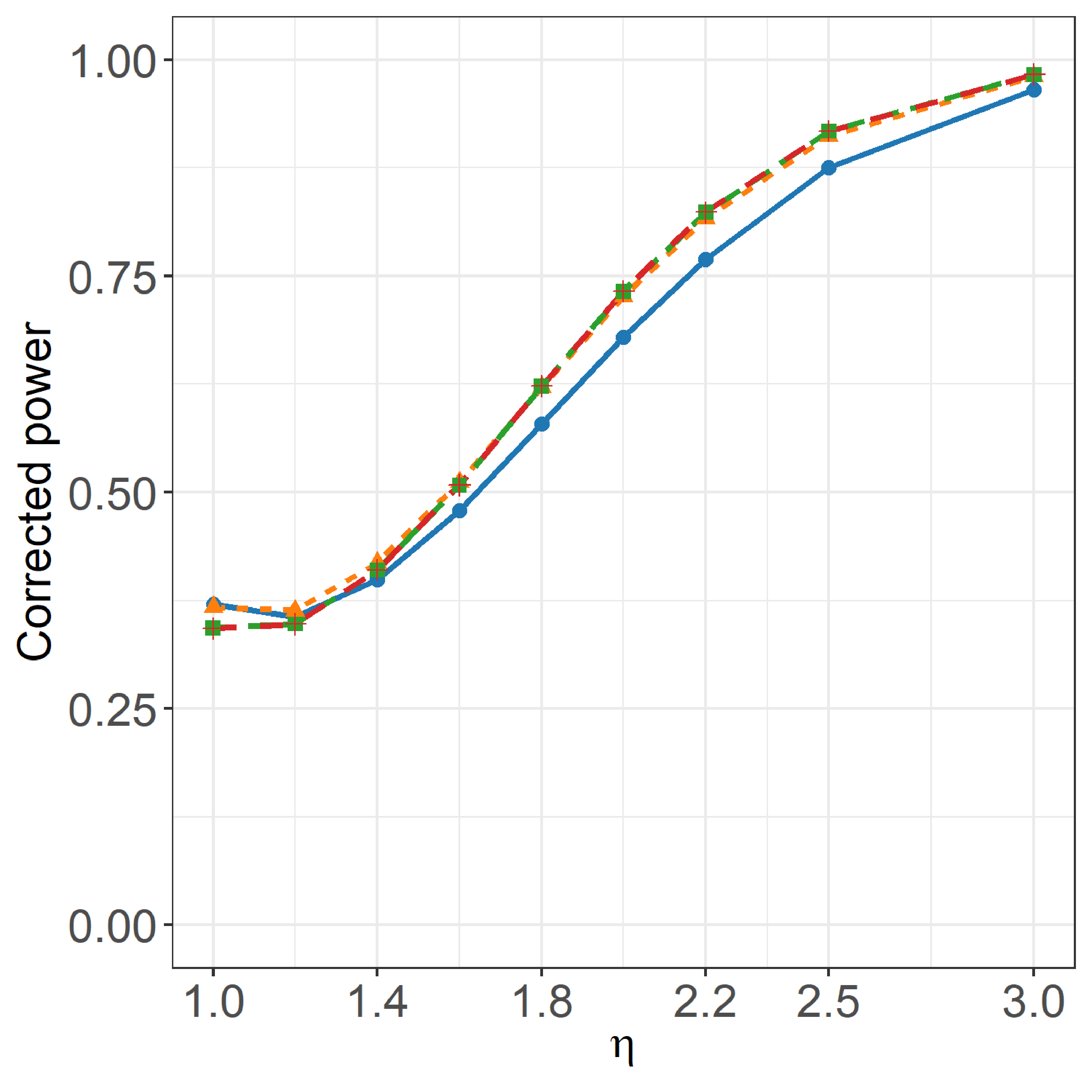}
		\end{minipage}
	}
	\subfigure{
		\begin{minipage}[b]{.3\linewidth}
			\centering
			\includegraphics[scale=0.0735]{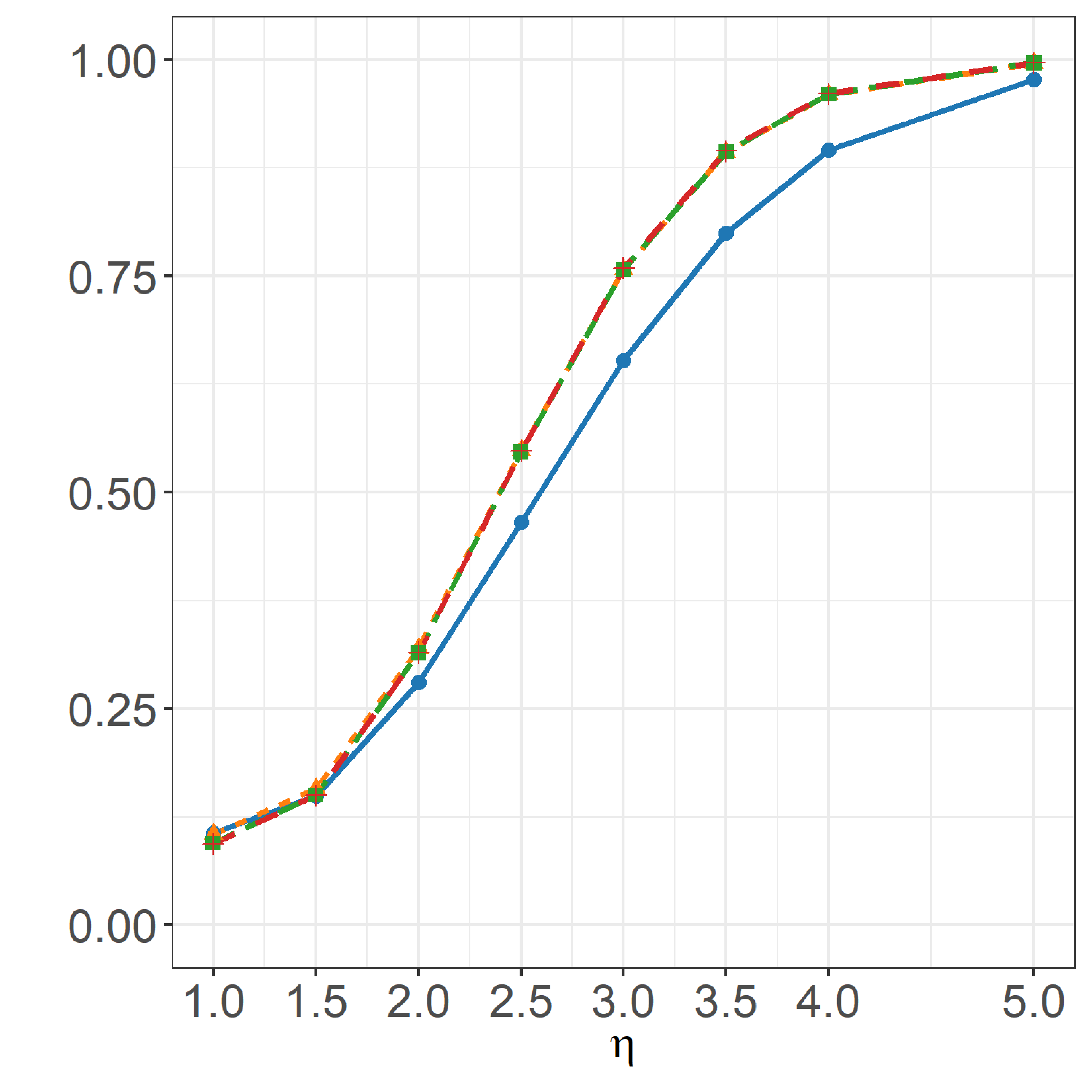}
		\end{minipage}
	}
	\subfigure{
		\begin{minipage}[b]{.3\linewidth}
			\centering
			\includegraphics[scale=0.0735]{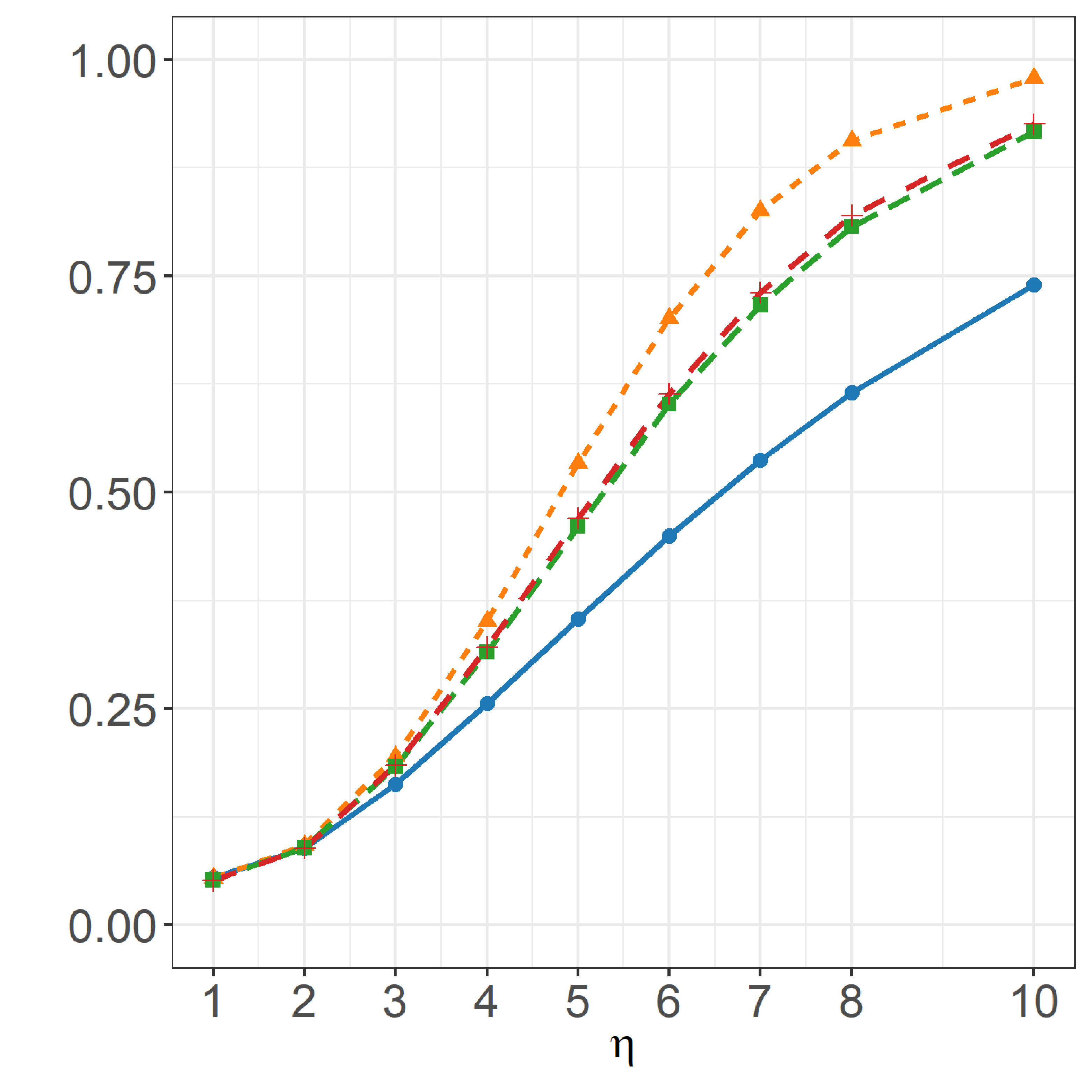}
		\end{minipage}
	}
	\subfigure{
		\begin{minipage}[b]{.3\linewidth}
			\centering
			\includegraphics[scale=0.0735]{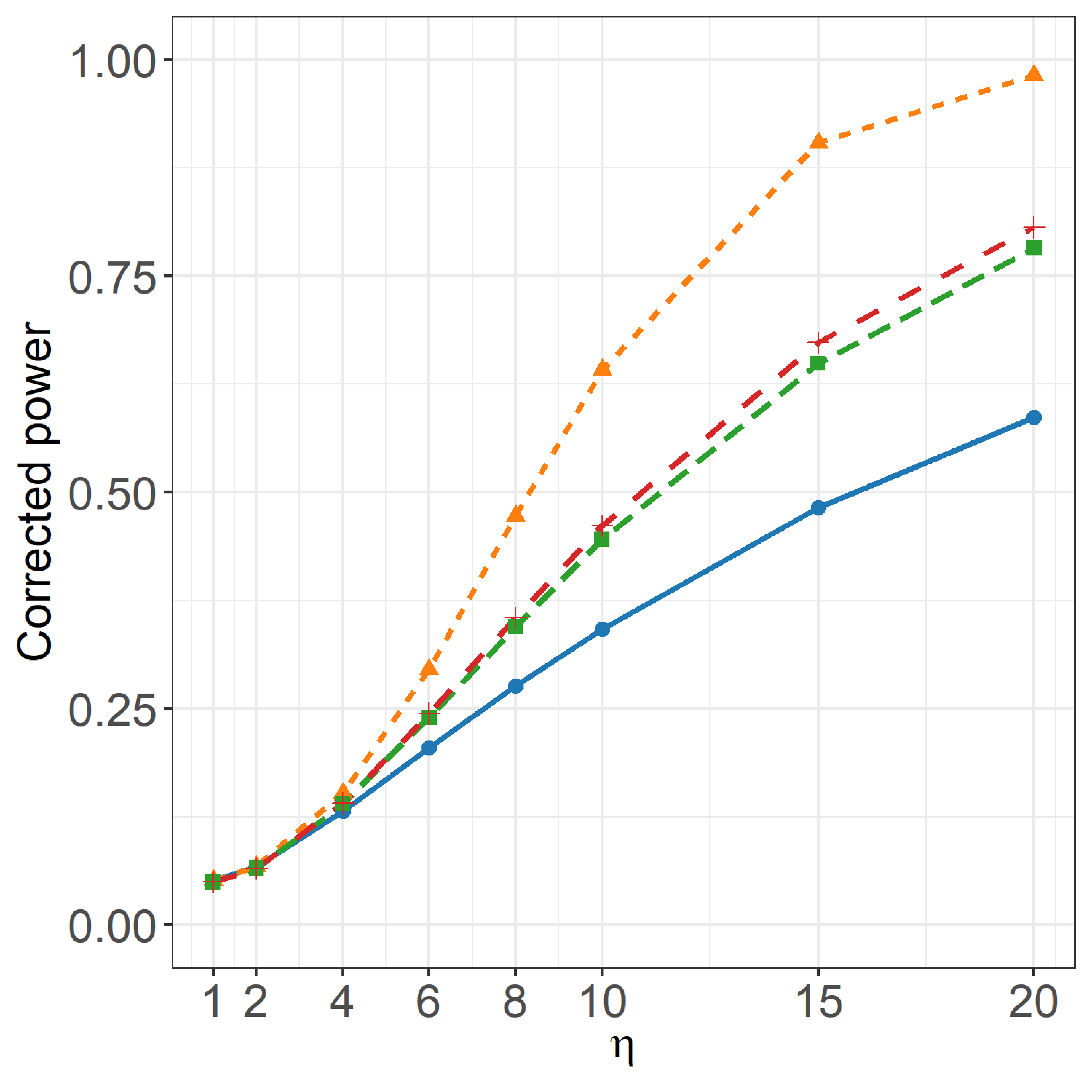}
		\end{minipage}
	}
	\subfigure{
		\begin{minipage}[b]{.3\linewidth}
			\centering
			\includegraphics[scale=0.0735]{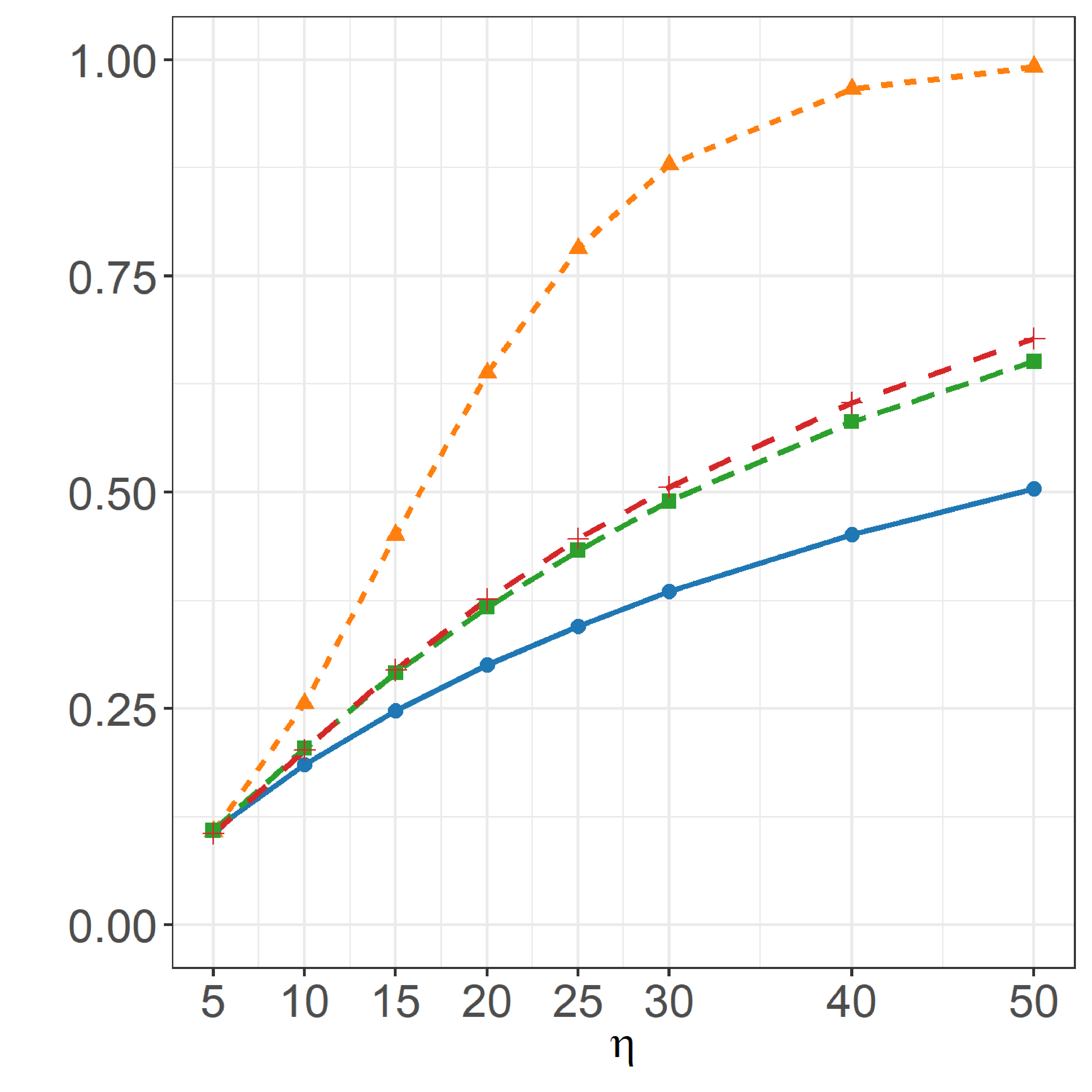}
		\end{minipage}
	}
	\subfigure{
		\begin{minipage}[b]{.3\linewidth}
			\centering
			\includegraphics[scale=0.0735]{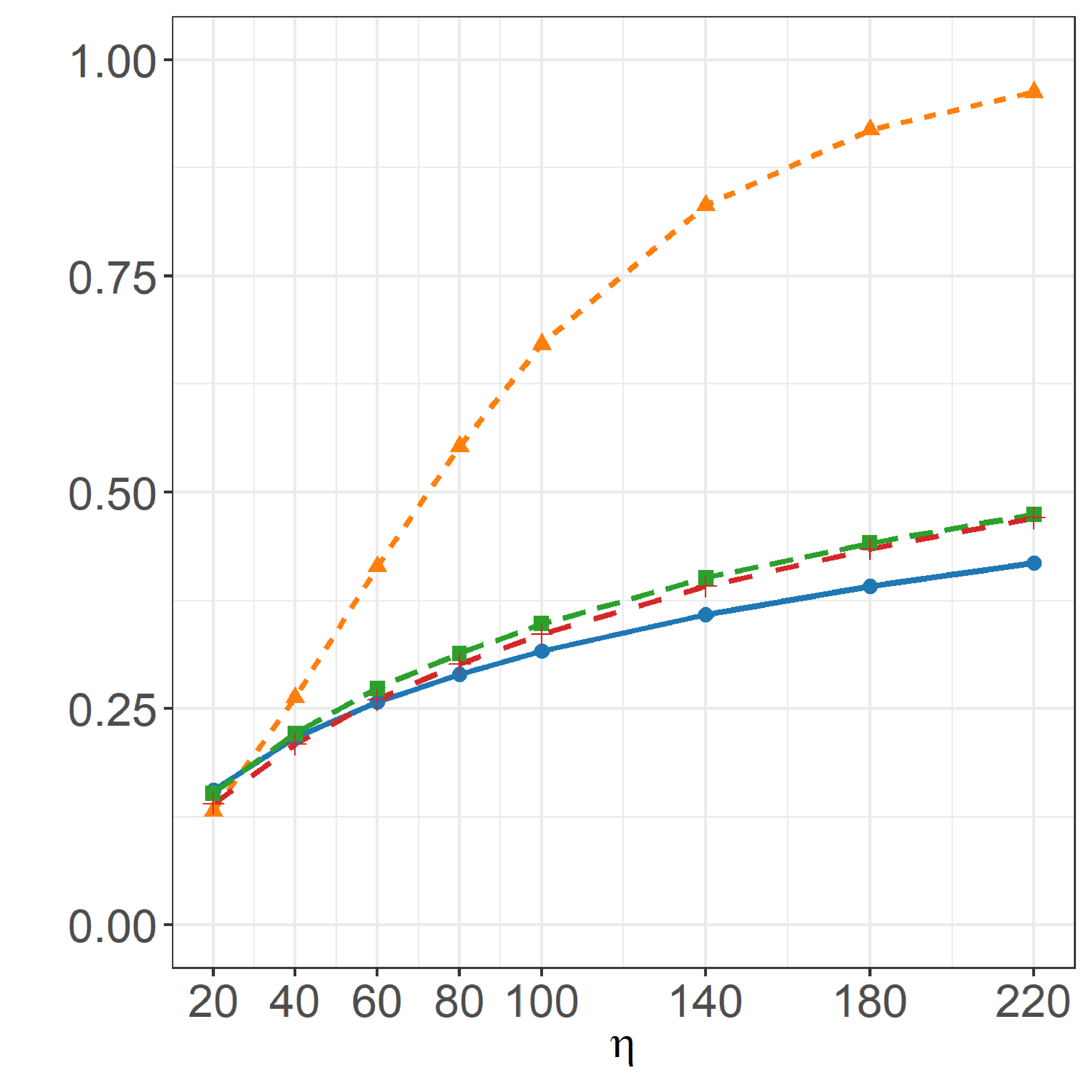}
		\end{minipage}
	}
	\caption{Empirical corrected power of four tests for hypothesis (II) with different values of $\eta$ and $p/n_i$. The solid, dashed, long-dashed, and dot-dashed curves are the empirical power functions of the central limit theorem test, directional test and two Skovgaard's modifications \cite{skovgaard:2001}, respectively.   The  alternative  setting (4) is given in Section \ref{simulation setting}. The six plots correspond to $p/n_i \in \{0.05, 0.1, 0.3, 0.5, 0.7, 0.9\}$, starting from top left and proceeding by row.}  
	\label{fig:power extreme case3}
\end{figure}

Finally, Figures \ref{fig:power extreme case4} and \ref{fig:power extreme case3} analyse the empirical corrected power of the tests for various ratios $p/n_i$ as $\eta$ in Section \ref{simulation setting} varies under the alternative setting (4) \cite{Jensen:2021} of hypotheses (I) and (II), respectively.  The directional test enjoys the best properties, proving to be particularly powerful with respect to its competitors when $p/n_i \geq  0.5$. Even in this case, the corrected power of the central limit theorem test is uniformly lowest.


\section{Discussion}\label{section6:disscussion}

This work examines directional testing for hypotheses on a vector parameter of interest in $p$-variate normal distributions when $n_i$ independent observations are available for the $i$th group $(i = 1,\dots,k)$ in the high dimensional regime with $p/n_i \to \kappa \in (0,1]$ \citep{battey2022some}. The construction of the directional test is based on the saddlepoint approximation to the  density of the canonical sufficient statistic, which is found to be exact provided that each $n_i \ge p+2$.  The numerical results support the  theoretical findings  on the  exact control of Type I error  of the directional approach under these mild conditions. The simulation outcomes show that the directional test outperforms the omnibus tests which look in all directions of the parameter space for alternatives both when $p$ is large and small relatively to $n_i$. Our formal derivations of the exactness of the underlying saddlepoint-type expansions provide also a theoretical ground to previous numerical findings obtained in the high dimensional simulation setting \cite{sartori:2014}.


The six hypotheses testing problems considered here and in the Supplementary Material mainly come from \cite{jiang:2013} and \cite{jiang:2015}. Jiang and Qi \cite{jiang:2015} showed that the central limit theorem test works well when $p$ is very close to $n_i$, assuming that $n_i>p+a$ for some  constant $1 \le a \le 4$. In our Monte Carlo experiments, the central limit theorem result seems inaccurate when the dimension $p$ is small, while the directional test is able to control exactly the Type I error for every value of $p$, provided that $n_i \ge p+2$. The two tests have been compared empirically also in terms of corrected power for some alternative hypotheses. Similarly to the log-likelihood ratio test, the central limit theorem approach is an omnibus test, whereas  the directional test measures the departure from the null hypothesis along the direction determined by the observed data point. In this respect, the latter is not constructed based on any kind of optimality \citep{skogaard:1988} and its marginal power may change  according to the specific alternative setting  \cite{Jensen:2021}. Nevertheless, our empirical results found not only that the power of the directional test does not need any correction for Type I error, but also that it is overall comparable with the corrected power of its main competitor.


The asymptotic theory for the directional test derived in this paper applies to linear exponential family models with hypotheses regarding linear functions of the canonical parameter, as in \cite{sartori:2014}.  Similar results for tests regarding the mean vector and/or covariance matrix that cannot be expressed as hypotheses on linear function of the canonical parameter could be obtained under the more general framework in \cite{sartori:2016}. Further research might focus on 
deriving the properties of  directional inference when $p$ increases with $n$ under the models considered previously by  
\cite{sartori:2016} and 
\cite{sartori:2019Ftest} for fixed $p$ only.

In general, the accuracy of the directional $p$-value stems from the accuracy of the underlying saddlepoint approximation to the conditional density of the canonical sufficient statistic.
{For instance, in the high dimensional regime  the directional test for the one-sample hypothesis on the normal mean vector $H_\psi\!: \mu = \mu_0$ is expected to behave as those seen here, since it was shown equivalent to the Hotelling's $T^2$ statistic \cite{sartori:2019Ftest}. 
In the multiple-sample case, preliminary results reveal that the high dimensional accuracy determined by the exactness of the directional $p$-value for testing the equality of the mean vectors is preserved only when assuming an identical covariance matrix for the $k$ independent groups.

For multivariate continuous distributions there are other instances of exactness of the saddlepoint approximation \cite[Section 10.9]{salvan:1997} where we can expect accuracy comparable with the high dimensional normal case. On the other hand, saddlepoint methods cannot be exact with discrete probability functions.
In  settings where the hypotheses are not linear in the canonical parameter or the saddlepoint approximation is not exact, 
ongoing simulation results and previous works \cite[Section 4.2]{sartori:2016, sartori:2014} suggest that
a  low dimensional asymptotic regime where $p/n \to 0$, typically with $p=O(n^\alpha)$, $0\leq\alpha<1$, might be required for observing the same accuracy of the directional $p$-value found in this paper.}

{Our interest in this work lies exclusively in the high dimensional asymptotic regime because maximum likelihood estimation is generally feasible in such situation, and so is the computation of the directional $p$-value.
That being said, under particular sparsity assumptions \cite[Section 4.4]{battey2022some} it is possible that the maximum likelihood estimator exists even if $p > n$, thus also the directional approach can be adopted in the ultra-high dimensional regime. 
As an example, consider hypothesis (III) in the Supplementary Material S2, testing the sphericity of the concentration matrix. The maximum likelihood estimate exists as long as $n$ is larger than the maximal clique size of the corresponding graph \cite{buhl1993existence}; hence, if the concentration matrix is assumed sparse enough the directional inference can still be applied \citep{di2021accurate}.}


\vspace*{20pt}

\appendix

\setcounter{equation}{0}
\renewcommand{\theequation}{A\arabic{equation}}

\setcounter{section}{0}
\renewcommand{\thesection}{Appendix A}

\setcounter{subsection}{0}
\renewcommand{\thesubsection}{A.\arabic{subsection}}

\section*{Appendix} \label{appendix}

\subsection{ Proof of Theorem \ref{theoremcase4}}\label{proof of theorem 1}

\begin{proof}
Suppose $y_{ij} \sim N_p(\mu_i, \Lambda^{-1}_i)$, $i \in \{1,\dots,k\}$, $j \in \{1,\dots,n_i\}$. For each $i$-th group, the random variables $y_{ij}$ are independent. Let $\bar{y}_i = n_i^{-1} {{1}}_{n_i}^T y_i$, $\hat{\Lambda}^{-1}_i = n_i^{-1}y_i^Ty_i - \bar{y}_i\bar{y}_i^T$. 
In this case, $\hat{\Lambda}_i^{-1} \sim W_p(n_i-1,n_i^{-1}\Lambda_i^{-1})$, for $i \in \{1,\dots, k\}$. Due to the groups  independence, the joint distribution of $\hat{\Lambda}_i$, $i \in \{1,\dots,k\}$,  is the product of Wishart densities $ \prod_{i=1}^k    f(\hat{\Lambda}_i^{-1};\Lambda_i^{-1})$ with 
\begin{eqnarray}
  f(\hat{\Lambda}_i^{-1};\Lambda_i^{-1}) &=& \left(\frac{n_i}{2}\right)^{\frac{p(n_i-1)}{2}} \Gamma_p \left(\frac{n_i-1}{2}\right)^{-1}\nonumber\\
  && \times |\Lambda_i^{-1}|^{-\frac{n_i-1}{2}} \text{etr}\left(-\frac{n_i}{2} \Lambda \hat{\Lambda}_i^{-1}\right) |\hat{\Lambda}_i^{-1}|^{\frac{n_i-p-2}{2}}.
  \label{proof:wishart distribution}
\end{eqnarray}
The  log-likelihood for  the canonical parameter  $\varphi$ under the multivariate normal distribution  is
\begin{eqnarray}
  \ell(\varphi;s) &=& \sum_{i=1}^k \frac{n_i}{2} \log |\Lambda_i| - \frac{1}{2} \tr(\Lambda_i y_i^Ty_i) + n\bar{y}_i^T \xi_i -\frac{n_i}{2} \xi_i^T  \Lambda_i^{-1} \xi_i.\nonumber
\end{eqnarray}
In order to assess the exactness of the saddlepoint approximation, it is convenient to express the log-likelihood function for $\varphi$ as
\begin{eqnarray}
  \ell(\varphi;s) = \sum_{i=1}^k {-\frac{n_i}{2}} \log|\Lambda_i^{-1}| -\frac{n_i}{2} \tr (\Lambda_i \hat{\Lambda}_i^{-1}) -\frac{n_i}{2}(\bar{y}_i -\Lambda_i^{-1}\xi_i )^T \Lambda_i (\bar{y}_i-\Lambda_i^{-1}\xi_i ). 
  \label{proof:loglikelihood case4}
\end{eqnarray}
The  maximum likelihood estimate $\hat{\varphi}$ has components $  \{\hat{\xi}^T_i , \text{vech}(\hat{\Lambda}_i)^T\}^T = \{\bar{y}_i^T\hat{\Lambda}_i ,$ $\text{vech}(\hat{\Lambda}_i)^T\}^T$  and the constrained  maximum likelihood estimate $\hat{\varphi}_ {\psi} $ has components $\{\bar{y}_i^T{\hat{\Lambda}}_0 , \text{vech}({\hat{\Lambda}}_0)^T\}^T$, $i \in \{1,\dots,k\}$. 
Evaluating (\ref{proof:loglikelihood case4}) at the unconstrained and constrained maximum likelihood estimates for $\varphi$, the correponding log-likelihood at $\hat{\varphi}$ and $\hat{\varphi}_{\psi}$ are $\ell(\hat{\varphi};s) = 2^{-1}\sum_{i=1}^k -{n_i} \log|\hat{\Lambda}_i^{-1}| - {n_ip}$ and $\ell(\hat{\varphi}_\psi;s) = 2^{-1}\sum_{i=1}^k -{n_i} $ $\log|{\hat{\Lambda}_0}^{-1}| - {n_i} \tr(\hat{\Lambda}_0\hat{\Lambda}_i^{-1})$, respectively. Then, under the null hypothesis $H_\psi$, 
and using the fact that  $|J_{\varphi \varphi}(\hat{\varphi})|$ is propotional to $\prod_{i=1}^k |\hat{\Lambda}_i^{-1}|^{p+2}$ (see Supplementary Material S1.1), the saddlepoint approximation (\ref{saddle L0}) is
\begin{eqnarray}
  h(s;\psi)
  &=&\prod_{i=1}^k  c_i(\psi) \; |\hat{\Lambda}_0^{-1}|^{-\frac{n_i-1}{2}} \exp\left\{-\frac{n_i}{2} \tr (\hat{\Lambda}_0 \hat{\Lambda}_i^{-1})\right\}  |\hat{\Lambda}_i^{-1}|^{\frac{n_i-p-2}{2}}.
  \label{proof:saddle case4}
\end{eqnarray}
Formula (\ref{proof:saddle case4}) is the exact joint distribution of $\hat{\Lambda}_1^{-1}, \dots, \hat{\Lambda}_k^{-1}$, i.e. a product of Wishart densities with parameters $(n_i-1, \hat{\Lambda}_0^{-1} )$ given in (\ref{proof:wishart distribution}), if $\hat{\Lambda}_0^{-1}$ is considered as fixed. In particular, we have $c_i(\psi) = c_i = \left({n_i}/{2}\right)^{p(n_i-1)/2} \Gamma_p\left\{({n_i-1})/{2}\right\}^{-1}$. It is indeed correct to fix $\hat{\Lambda}_0^{-1}$ because when considering the saddlepoint approximation density along the line $s(t)$, by construction the constrained maximum likelihood estimates of $\Lambda_i^{-1}$ is fixed and equal to the observed value $\hat{\Lambda}_0^{-1}$. 

When we consider the density of $s(t)$, we just need to replace  $\hat{\Lambda}_i^{-1}$ in (\ref{proof:saddle case4}) with $\hat{\Lambda}_i^{-1}(t)$, i.e. the value which maximizes $\ell\left\{\varphi;s(t)\right\}$. Then, given that  $\hat{\Lambda}_i^{-1}(t)=(1-t)\hat{\Lambda}_0^{-1} + t\hat{\Lambda}_i^{-1}$ and the groups are independent, under $H_{\psi}$ we have 
\begin{eqnarray}
  h\{s(t);\psi\} &=& \prod_{i=1}^k c_{i}  \; |\hat{\Lambda}_0^{-1}|^{-\frac{n_i-1}{2}} \exp\left[-\frac{n_i}{2} \tr\{\hat{\Lambda}_0 \hat{\Lambda}_i(t)^{-1}\}\right] |\hat{\Lambda}_i(t)^{-1}|^{\frac{n_i-p-2}{2}} \nonumber \\
  & \propto&  \exp\left\{ \sum_{i=1}^k \frac{n_i-p-2}{2} \log |\hat{\Lambda}_i^{-1}(t)| \right\}, \nonumber
\end{eqnarray}
where we have used the equality  $n^{-1} \tr(\sum_{i=1}^{k} n_i  \hat{\Lambda}_0 \hat{\Lambda}^{-1}_i )=p$ with $n = \sum_{i=1}^{k}n_i$.
 Since  the saddlepoint approximation $ h\{s(t);\psi\}$  is exact, apart from the normalizing constant, the integral in the denominator of the directional $p$-value (\ref{directed p-value})  is just the normalizing constant of the conditional distribution of $||s||$ given the direction $s/||s||$.  Therefore, the directional $p$-value is the exact probability of $||s|| > ||s^0||$ given the direction $s/||s||$ under the null hypothesis, and  is thus exactly uniformly distributed. 
\end{proof}

\subsection{ Proof of Theorem \ref{theoremcase3}}\label{proof of theorem 2}

\begin{proof}
	We know that $\bar{y}_i \sim N_p(\mu_i, n_i^{-1} \Lambda_i^{-1})$ and $\hat{\Lambda}_i^{-1} \sim W_p(n_i-1,n_i^{-1}\Lambda_i^{-1})$, $i \in \{1,\dots, k\}$. In addition, $\bar{y}_i$ and $\hat{\Lambda}_i$ are independent \citep[][Section\ 10.8]{muirhead1982}, thus the joint distribution of $\bar{y}_i$ and $\hat{\Lambda}_i$ takes the form $\prod_{i=1}^k   f(\bar{y}_i;\mu_i,\Lambda_i^{-1}) f(\hat{\Lambda}_i^{-1};\Lambda_i^{-1})$
	with 
	\begin{eqnarray}
	f(\bar{y}_i;\mu_i,\Lambda_i^{-1}) &=& (2\pi)^{-\frac{p}{2}} |\Lambda_i^{-1}|^{-\frac{1}{2}} \exp \left\{-\frac{n_i}{2} (\bar{y}_i -\mu_i) ^T \Lambda_i (\bar{y}_i - \mu_i)\right\},  \nonumber\\
	f(\hat{\Lambda}_i^{-1};\Lambda_i^{-1}) &=& (n_i/2)^{\frac{p(n_i-1)}{2}} \Gamma_p \left(\frac{n_i-1}{2}\right)^{-1}\nonumber\\
	&& \times |\Lambda_i^{-1}|^{-\frac{n_i-1}{2}} \text{etr}\left(-\frac{n_i}{2} \Lambda \hat{\Lambda}_i^{-1}\right) |\hat{\Lambda}_i^{-1}|^{\frac{n_i-p-2}{2}}. \nonumber
	\end{eqnarray}
	Similarly to the proof of Theorem \ref{theoremcase4}, we can easily obtain 
	the saddlepoint approximation to the density of the sufficient statistic $s$ as
	\begin{eqnarray}
	h(s;\psi)
	&=& \prod_{i=1}^k  c_{i1} \; |\hat{\Lambda}_0^{-1}|^{-\frac{1}{2}} \exp \left\{ -\frac{n_i}{2} (\bar{y}_i - \hat{\mu}_0)^T \hat{\Lambda}_0 (\bar{y}_i-\hat{\mu}_0)\right\} \nonumber \\
	&& \quad \;\; \times  c_{i2} \; |\hat{\Lambda}_0^{-1}|^{-\frac{n_i-1}{2}} \exp\left\{-\frac{n_i}{2} \tr (\hat{\Lambda}_0 \hat{\Lambda}_i^{-1})\right\}  |\hat{\Lambda}_i^{-1}|^{\frac{n_i-p-2}{2}}.
	\label{proof:saddle case3}
	\end{eqnarray}
	Expression (\ref{proof:saddle case3}) equals the exact joint distribution of $\bar{y}_1, \dots, \bar{y}_k$ and $\hat{\Lambda}_1^{-1},\dots , \hat{\Lambda}_k^{-1}$ with $c_{i1} = (2\pi)^{-p/2}$, $c_{i2} = \left(\frac{n_i}{2}\right)^{p(n_i-1)/2} \Gamma_p\left(\frac{n_i-1}{2}\right)^{-1}$ and with fixed $\hat{\mu}_0$ and $\hat{\Lambda}_0^{-1}$. It is indeed correct to consider $\hat{\mu}_0$ and $\hat{\Lambda}_0^{-1}$ as fixed since the constrained maximum likelihood estimate is fixed and equal to the observed value when considering  the saddlepoint approximation along the line $s(t)$ under $H_\psi$.  In such case we have  $\hat{\mu}_{i}(t)=(1-t) \hat{\mu}_0+t\bar{y}_i$ and $\hat{\Lambda}_i(t)^{-1} = (1-t)\hat{\Lambda}_0^{-1} + t\hat{\Lambda}_i^{-1}+t(1-t)(\bar{y}-\bar{y}_i)(\bar{y}-\bar{y}_i)^T$ where $\hat{\mu}_0=\bar{y}$ and $\hat{\Lambda}_0^{-1} = n^{-1}(A+B)$ (see Section \ref{section3:main results} for more details). Then, the saddlepoint approximation for the distribution of $s(t)$ under $H_{\psi}$ follows from (\ref{proof:saddle case3}), and is equal to
	\begin{eqnarray}
	h\{s(t);\psi\} &=& \prod_{i=1}^k c_{1i} \; |\hat{\Lambda}^{-1}_0|^{-\frac{1}{2}} \exp \left[ -\frac{n_i}{2} \{\hat{\mu}_{i}(t) - \hat{\mu}_0\}^T \hat{\Lambda}_0 \{\hat{\mu}_{i}(t)-\hat{\mu}_0\}\right] \nonumber \\
	&& \quad \; \; \times  c_{i2} \; |\hat{\Lambda}_0^{-1}|^{-\frac{n_i-1}{2}} \exp\left[-\frac{n_i}{2} \tr \{\hat{\Lambda}_0 \hat{\Lambda}_i(t)^{-1}\}\right]  |\hat{\Lambda}_i(t)^{-1}|^{\frac{n_i-p-2}{2}} \nonumber \\
	&\propto& \exp\left\{\sum_{i=1}^k \frac{n_i-p-2}{2} \log|\hat{\Lambda}_i(t)^{-1}| \right\}. \nonumber
	\end{eqnarray}
	The remaining part of the proof is similar to that for hypothesis (I). It follows then  that the directional $p$-value is exactly uniformly distributed under the null hypothesis $H_\psi$.
\end{proof}

 \subsection{ Proof of Lemma \ref{lemmacase1}} \label{proof of lemma 1}
 

 \begin{proof}
 	If $t \in [0,1]$ the result is straightforward,  because a convex combination of positive definite matrices is positive definite. Indeed, for all $x \in \mathbb{R}^p$, $x \ne 0$, $ x^T  \hat{\Lambda}^{-1}(t) x = (1-t) x^T \hat{\Lambda}^{-1}_0 x + t x^T \hat{\Lambda}^{-1} x >0$  since $1-t \ge 0$ and $t\ge0$.
 	Let us focus on the case  $t>1$. Consider a square root $B_0$ of $\hat{\Lambda}^{-1}_0$ such that $\hat{\Lambda}^{-1}_0=B_0B_0^T=B_0^TB_0$, which always exists if $\hat{\Lambda}^{-1}_0$ is positive definite. Hence, the estimator $\hat{\Lambda}^{-1}(t) = (1-t)\hat{\Lambda}_0^{-1} + t\hat{\Lambda}^{-1}$ can be  rewritten as
 	\begin{eqnarray}
 	\hat{\Lambda}^{-1}(t) = B_0^T \left\{(1-t)\bI_p + t (B_0^T)^{-1} \hat{\Lambda}^{-1}B_0^{-1}\right\}B_0. \nonumber
 	\end{eqnarray}
 	The matrix $(B_0^T)^{-1} \hat{\Lambda}^{-1}B_0^{-1}$ is  symmetric since $\hat{\Lambda}^{-1}$ is symmetric.
 	Moreover, according to the eigen decomposition \citep[][Theorem 1.13]{magnus1999}, there exists an orthogonal $p\times p$ matrix $P$ whose columns are eigenvectors of $(B_0^T)^{-1} \hat{\Lambda}^{-1}B_0^{-1}$ and a diagonal matrix $Q$ whose diagonal elements are the eigenvalues of $(B_0^T)^{-1} \hat{\Lambda}^{-1}B_0^{-1}$, such that $(B_0^T)^{-1} \hat{\Lambda}^{-1}B_0^{-1} = P Q P^T$. Therefore, we have $\hat{\Lambda}^{-1}(t) = B_0^T P \left\{(1-t)\bI_p \right.$ $\left.+ t Q \right\}P^TB_0 $.
 	Lemma \ref{lemmacase1} can then be proved through the following three steps.

 	\textit{Step} 1: checking that  $\hat{\Lambda}^{-1}(t)$ is positive definite is equivalent to checking that $(1-t)\bI_p + tQ$ is positive definite.
 	Indeed, for all $x \in \mathbb{R}^p$, $x\ne0$, then
 	\begin{eqnarray}
 	x^T \hat{\Lambda}^{-1}(t) x &=& x^T B_0^T P \left\{(1-t)\bI_p + t Q \right\} P^T B_0 x \nonumber \\
 	&= & \tilde{x}^T \left\{(1-t)\bI_p + t Q \right\} \tilde{x} >0, \nonumber
 	\end{eqnarray}
 	where $\tilde{x}=P^TB_0 x$, with $\tilde{x}\ne0$ if  $x\ne0$.

 	\textit{Step} 2: checking that $(1-t)\bI_p + tQ$ is positive definite is equivalent to checking that all  elements of the diagonal matrix $(1-t)\bI_p + tQ = \diag(1-t+t\nu_l)$ are positive, where $\nu_l$, $l \in \{1,\dots,p\}$,  are the eigenvalues of the matrix $(B_0^T)^{-1} \hat{\Lambda}^{-1}B_0^{-1}$. 
 	We now need to find out the largest $t$ such that $1-t +t \nu_l >0,  l \in \{1,\dots,p\} $:
 	\begin{itemize}
 		\item if $1-\nu_{(1)}>0$, where $\nu_{(1)}$ is the smallest eigenvalue of $(B_0^T)^{-1} \hat{\Lambda}^{-1}B_0^{-1}$, then $t<\frac{1}{1-\nu_l} \le \frac{1}{1-\nu_{(1)}}$;
 		\item if $1-\nu_{(1)} \le 0$, then $t>\frac{1}{1-\nu_l}$ as $\frac{1}{1-\nu_l}<0$, and this condition holds true  $\forall \; t \in \mathbb{R}^+$.
 	\end{itemize}
 	\textit{Step} 3: The last step consists of  checking that  the eigenvalues $\nu_1,\dots,\nu_p$ of $(B_0^T)^{-1}\hat{\Lambda}^{-1}$ $B_0^{-1}$ are the same as those of $\hat{\Lambda}_0\hat{\Lambda}^{-1}$, which is equivalent to show that the matrices $(B_0^T)^{-1}\hat{\Lambda}^{-1}B_0^{-1}$ and $\hat{\Lambda}_0\hat{\Lambda}^{-1}$ are similar. 
 	In addition, $\hat{\Lambda}_0^{-1} = B_0^T B_0$, given the invertible matrix $B_0$ such that 
 	\begin{eqnarray}
 	B_0^{-1} (B_0^T)^{-1}\hat{\Lambda}^{-1}B_0^{-1} B_0 =   B_0^{-1} (B_0^T)^{-1}\hat{\Lambda}^{-1} = \hat{\Lambda}_0\hat{\Lambda}^{-1}. \nonumber
 	\end{eqnarray}
 	According to matrix similarity,   $(B_0^T)^{-1}\hat{\Lambda}^{-1}B_0^{-1}$ and $\hat{\Lambda}_0\hat{\Lambda}^{-1}$ are similar and therefore have the same eigenvalues.

 	Finally, since $\hat{\Lambda}_0\hat{\Lambda}^{-1}$ is positive definite and $\tr( \hat{\Lambda}_0\hat{\Lambda}^{-1})=p$, the smallest eigenvalue $\nu_{(1)}$ must be lower than 1. Therefore,  $\hat{\Lambda}^{-1}(t)$ is positive definite in $t \in [0,\{1-\nu_{(1)}\}^{-1}]$.
 \end{proof}


 \subsection{ Proof of Lemma \ref{lemmacase4}} \label{proof of lemma 2}
 
 
 
 \begin{proof}
 	Based on the proof of Lemma \ref{lemmacase1}, it is easy to show that $\hat{\Lambda}^{-1}_i(t)$ for all $i \in \{1,\dots,k\}$, is positive definite if and only if all elements $1-t+t\nu_l^i>0$, where $\nu_l^i, l \in \{1,\dots,p\}$, are the eigenvalues of the matrix $\hat{\Lambda}_0\hat{\Lambda}^{-1}_i$, $i \in \{1,\dots, k\}$.  
 	Since $\hat{\Lambda}_0\hat{\Lambda}^{-1}_i$ are positive definite and $\tr(\hat{\Lambda}_0\hat{\Lambda}^{-1}_i) = p$ for all $ i \in \{1,\dots,k\}$, there exists at least one of  the $\nu_{(1)}^i$  lower than 1,
 	where  $\nu_{(1)}^i$ denotes the smallest eigenvalue of $\hat{\Lambda}_0\hat{\Lambda}^{-1}_i$. 
 	In this respect,  $\hat{\Lambda}^{-1}_i(t)$, $\forall i \in \{1, \dots, k\}$, are positive definite in $t \in \left[0, \{1- \mathop {\min }\limits_{1 \le i \le k} \nu_{(1)}^i\}^{-1}\right]$.
 \end{proof}

\bibliographystyle{main}
\bibliography{main}

\begin{thebibliography}{38}
\expandafter\ifx\csname natexlab\endcsname\relax\def\natexlab#1{#1}\fi

\bibitem[{Anderson(2003)}]{anderson1958}
\textsc{Anderson, T.~W.} (2003).
\newblock \textit{An Introduction to Multivariate Statistical Analysis}.
\newblock Wiley, 3rd ed.

\bibitem[{Barndorff-Nielsen(1986)}]{barndorff1986}
\textsc{Barndorff-Nielsen, O.} (1986).
\newblock Inference on full or partial parameters based on the standardized
  signed log likelihood ratio.
\newblock \textit{Biometrika} \textbf{73}, 307--322.

\bibitem[{Bartlett(1937)}]{bartlett:1937}
\textsc{Bartlett, M.} (1937).
\newblock Properties of sufficiency and statistical tests.
\newblock \textit{Proc. Roy. Soc. London Ser. A} \textbf{160}, 268--282.

\bibitem[{Battey \& Cox(2022)}]{battey2022some}
\textsc{Battey, H.} \& \textsc{Cox, D.} (2022).
\newblock Some perspectives on inference in high dimensions.
\newblock \textit{Statistical Science} \textbf{37}, 110--122.

\bibitem[{Buhl(1993)}]{buhl1993existence}
\textsc{Buhl, S.~L.} (1993).
\newblock On the existence of maximum likelihood estimators for graphical
  gaussian models.
\newblock \textit{Scandinavian Journal of Statistics} , 263--270.

\bibitem[{Cheah et~al.(1994)Cheah, Fraser \& Reid}]{fraser1994}
\textsc{Cheah, P.~K.}, \textsc{Fraser, D. A.~S.} \& \textsc{Reid, N.} (1994).
\newblock Multiparameter testing in exponential models: Third order
  approximations from likelihood.
\newblock \textit{Biometrika} \textbf{81}, 271--278.

\bibitem[{Cordeiro \& Cribari-Neto(2014)}]{cordeiro:2014}
\textsc{Cordeiro, G.} \& \textsc{Cribari-Neto, F.} (2014).
\newblock \textit{An Introduction to Bartlett Correction and Bias Reduction}.
\newblock Springer-Verlag.

\bibitem[{Davison(2003)}]{davison:2003}
\textsc{Davison, A.~C.} (2003).
\newblock \textit{Statistical Models}.
\newblock Cambridge University Press.

\bibitem[{Davison et~al.(2014)Davison, Fraser, Reid \& Sartori}]{sartori:2014}
\textsc{Davison, A.~C.}, \textsc{Fraser, D. A.~S.}, \textsc{Reid, N.} \&
  \textsc{Sartori, N.} (2014).
\newblock Accurate directional inference for vector parameters in linear
  exponential families.
\newblock \textit{J. Amer. Statist. Assoc.} \textbf{109}, 302--314.

\bibitem[{Dette \& Dörnemann(2020)}]{dette2020}
\textsc{Dette, H.} \& \textsc{Dörnemann, N.} (2020).
\newblock Likelihood ratio tests for many groups in high dimensions.
\newblock \textit{Journal of Multivariate Analysis} \textbf{178}, 104605.

\bibitem[{Di~Caterina et~al.(2021)Di~Caterina, Reid \&
  Sartori}]{di2021accurate}
\textsc{Di~Caterina, C.}, \textsc{Reid, N.} \& \textsc{Sartori, N.} (2021).
\newblock Accurate directional inference in gaussian graphical models.
\newblock \textit{arXiv preprint arXiv:2103.15394} .

\bibitem[{Eriksen(1996)}]{Eriksen1996}
\textsc{Eriksen, P.~S.} (1996).
\newblock Tests in covariance selection models.
\newblock \textit{Scand. J. Stat.} \textbf{23}, 275--284.

\bibitem[{Fan et~al.(2019)Fan, Demirkaya \& Lv}]{fan2019nonuniformity}
\textsc{Fan, Y.}, \textsc{Demirkaya, E.} \& \textsc{Lv, J.} (2019).
\newblock Nonuniformity of $p$-values can occur early in diverging dimensions.
\newblock \textit{J. Mach. Learn. Res.} \textbf{20}, 77--1.

\bibitem[{Fraser \& Massam(1985)}]{fraser:1985}
\textsc{Fraser, D. A.~S.} \& \textsc{Massam, H.} (1985).
\newblock Conical tests: Observed levels of significance and confidence
  regions.
\newblock \textit{Statistische Hefte} \textbf{26}, 1--17.

\bibitem[{Fraser et~al.(2016)Fraser, Reid \& Sartori}]{sartori:2016}
\textsc{Fraser, D. A.~S.}, \textsc{Reid, N.} \& \textsc{Sartori, N.} (2016).
\newblock {Accurate directional inference for vector parameters}.
\newblock \textit{Biometrika} \textbf{103}, 625--639.

\bibitem[{Guo \& Qi(2021)}]{guo2021asymptotic}
\textsc{Guo, W.} \& \textsc{Qi, Y.} (2021).
\newblock Asymptotic distributions for likelihood ratio tests for the equality
  of covariance matrices.
\newblock \textit{arXiv} , 2110.02384.

\bibitem[{He et~al.(2021)He, Meng, Zeng \& Xu}]{He2020B}
\textsc{He, Y.}, \textsc{Meng, B.}, \textsc{Zeng, Z.} \& \textsc{Xu, G.}
  (2021).
\newblock {On the phase transition of Wilks’ phenomenon}.
\newblock \textit{Biometrika} \textbf{108}, 741--748.

\bibitem[{Højsgaard(2020)}]{gRim}
\textsc{Højsgaard, S.} (2020).
\newblock \textit{gRim: Graphical interaction models}.
\newblock R package version 0.2.5.

\bibitem[{Jensen(2021)}]{Jensen:2021}
\textsc{Jensen, J.~L.} (2021).
\newblock On the use of saddlepoint approximations in high dimensional
  inference.
\newblock \textit{Sankhya A} \textbf{83}, 379--392.

\bibitem[{Jiang \& Qi(2015)}]{jiang:2015}
\textsc{Jiang, T.} \& \textsc{Qi, Y.} (2015).
\newblock Likelihood ratio tests for high-dimensional normal distributions.
\newblock \textit{Scand. J. Stat.} \textbf{42}, 988--1009.

\bibitem[{Jiang \& Yang(2013)}]{jiang:2013}
\textsc{Jiang, T.} \& \textsc{Yang, F.} (2013).
\newblock Central limit theorems for classical likelihood ratio tests for
  high-dimensional normal distributions.
\newblock \textit{Ann. Statist.} \textbf{41}, 2029--2074.

\bibitem[{Lauritzen(1996)}]{Lauritzen:1996}
\textsc{Lauritzen, S.~L.} (1996).
\newblock \textit{Graphical Models}.
\newblock Oxford University Press.

\bibitem[{Lawley(1956)}]{lawley1956}
\textsc{Lawley, D.~N.} (1956).
\newblock {A general method for approximating to the distribution of likelihood
  ratio criteria}.
\newblock \textit{Biometrika} \textbf{43}, 295--303.

\bibitem[{Liu \& Pierce(1994)}]{gausshermite}
\textsc{Liu, Q.} \& \textsc{Pierce, D.~A.} (1994).
\newblock A note on {G}auss-{H}ermite quadrature.
\newblock \textit{Biometrika} \textbf{81}, 624--629.

\bibitem[{Magnus \& Neudecker(1999)}]{magnus1999}
\textsc{Magnus, J.} \& \textsc{Neudecker, H.} (1999).
\newblock \textit{Matrix Differential Calculus with Applications in Statistics
  and Econometrics}.
\newblock Wiley, 3rd ed.

\bibitem[{Marchetti et~al.(2020)Marchetti, Drton \& Sadeghi}]{ggm}
\textsc{Marchetti, G.~M.}, \textsc{Drton, M.} \& \textsc{Sadeghi, K.} (2020).
\newblock \textit{ggm: Graphical markov models with mixed graphs}.
\newblock R package version 0.2.5.

\bibitem[{Mccormack et~al.(2019)Mccormack, Reid, Sartori \&
  Theivendran}]{sartori:2019Ftest}
\textsc{Mccormack, A.}, \textsc{Reid, N.}, \textsc{Sartori, N.} \&
  \textsc{Theivendran, S.~A.} (2019).
\newblock A directional look at {$F$}-tests.
\newblock \textit{Canad. J. Statist.} \textbf{47}, 619--627.

\bibitem[{McCullagh(2018)}]{mccullage:2018}
\textsc{McCullagh, P.} (2018).
\newblock \textit{Tensor Methods in Statistics}.
\newblock Dover Publications, 2nd ed.

\bibitem[{Muirhead(2009)}]{muirhead1982}
\textsc{Muirhead, R.~J.} (2009).
\newblock \textit{Aspects of Multivariate Statistical Theory}.
\newblock Wiley.

\bibitem[{Pace \& Salvan(1997)}]{salvan:1997}
\textsc{Pace, L.} \& \textsc{Salvan, A.} (1997).
\newblock \textit{{Principles of Statistical Inference from a Neo-Fisherian
  Perspective}}.
\newblock World Scientific Press.

\bibitem[{{R Core Team}(2020)}]{Rcode}
\textsc{{R Core Team}} (2020).
\newblock \textit{R: A Language and Environment for Statistical Computing}.
\newblock R Foundation for Statistical Computing, Vienna, Austria.

\bibitem[{Severini(2001)}]{Severini:2000}
\textsc{Severini, T.~A.} (2001).
\newblock \textit{Likelihood Methods in Statistics}.
\newblock Oxford University Press.

\bibitem[{Skovgaard(1988)}]{skogaard:1988}
\textsc{Skovgaard, I.} (1988).
\newblock Saddlepoint expansions for directional test probabilities.
\newblock \textit{J. R. Stat. Soc. Ser. B. Stat. Methodol.} \textbf{50},
  269--280.

\bibitem[{Skovgaard(2001)}]{skovgaard:2001}
\textsc{Skovgaard, I.} (2001).
\newblock Likelihood asymptotics.
\newblock \textit{Scand. J. Stat.} \textbf{28}, 3--32.

\bibitem[{Sur \& Candès(2019)}]{sur:2019}
\textsc{Sur, P.} \& \textsc{Candès, E.~J.} (2019).
\newblock A modern maximum-likelihood theory for high-dimensional logistic
  regression.
\newblock \textit{Proc. Natl. Acad. Sci.} \textbf{116}, 14516--14525.

\bibitem[{Sur et~al.(2019)Sur, Chen \& Cand{\`e}s}]{sur2019likelihood}
\textsc{Sur, P.}, \textsc{Chen, Y.} \& \textsc{Cand{\`e}s, E.~J.} (2019).
\newblock The likelihood ratio test in high-dimensional logistic regression is
  asymptotically a rescaled chi-square.
\newblock \textit{Probab. Theory Related Fields} \textbf{175}, 487--558.

\bibitem[{Tang \& Reid(2020)}]{tang2020modified}
\textsc{Tang, Y.} \& \textsc{Reid, N.} (2020).
\newblock Modified likelihood root in high dimensions.
\newblock \textit{J. R. Stat. Soc. Ser. B. Stat. Methodol.} \textbf{82},
  1349--1369.

\bibitem[{Tang \& Reid(2021)}]{tang2021laplace}
\textsc{Tang, Y.} \& \textsc{Reid, N.} (2021).
\newblock Laplace and saddlepoint approximations in high dimensions.
\newblock \textit{arXiv} , 2107.10885.

\end{thebibliography}

\newpage
\setcounter{page}{1}
\setcounter{section}{0}
\setcounter{figure}{0}
\setcounter{table}{0}

\renewcommand{\thefigure}{S\arabic{figure}}
\renewcommand{\thetable}{S\arabic{table}}
\renewcommand{\theequation}{S\arabic{equation}}
\renewcommand{\thesection}{S\arabic{section}}
\renewcommand{\thesubsection}{S\arabic{section}.\arabic{subsection}}

\section*{Supplementary material to directional testing for high-dimensional multivariate normal distributions}

	\begin{abstract}
	In this supplementary material,  some auxiliary computational results are available in Section \ref{S1}.  Section \ref{B1:determinant} gives the determinant of the observed information matrix in multivariate normal distribution. Section \ref{B2:gamma}  derives the quantities needed in Skovgaard's modifications proposed by \cite{skovgaard:2001}. Section \ref{B3:tmintmax} details the computation of the directional $p$-value. We also investigate some theoretical results in Section 3 of the paper. In particular, the theoretical outcomes for hypotheses (III)-(V) and (VI) are detailed in Sections \ref{case1} and \ref{case5}, respectively. In all hypotheses problems, we prove the exact uniformity of the directional $p$-value under the null hypothesis.  
  In Section \ref{S3:simulation (I)-(II)}--\ref{S6:various sample sizes}, we report results of additional simulation studies. 
	\end{abstract}


\section{Auxiliary computational results}\label{S1}
\subsection{Determinant of the observed information matrix in multivariate normal distributions}\label{B1:determinant}

In Section 2.3 of the paper,  the observed information matrix with respect to the canonical parameters  $\xi$ and $\zeta = \text{vech} (\Lambda)$ is
\begin{eqnarray}
J_{\varphi \varphi}(\varphi) &=& \left[\begin{array}{cc}
J_{\xi\xi}(\varphi)  & J_{\xi\zeta}(\varphi) \\
J_{\zeta\xi}(\varphi) &  J_{\zeta\zeta}(\varphi)
\end{array}\right]\nonumber\\
&=&  \left[\begin{array}{cc}
n\Lambda^{-1}   & \quad -n(\xi^T \Lambda^{-1} \otimes \Lambda^{-1} ) D_p \\
-n D_p^T (\Lambda^{-1} \xi \otimes \Lambda^{-1})    & \quad  \frac{n}{2} D_p^T \{ \Lambda^{-1} (\bI_p + 2\xi \xi^T \Lambda^{-1}) \otimes \Lambda^{-1}\} D_p
\end{array} \right].\nonumber
\end{eqnarray}
In order to obtain the saddlepoint approximation to the conditional distribution of the  canonical sufficient statistic and the modification term of Skovgaard's  statistics \cite{skovgaard:2001}, we compute the determinant of $ J_{\varphi \varphi}(\varphi)$ as follows:
\begin{eqnarray}
\left| J_{\varphi \varphi}(\varphi)\right| =  \left|   n\Lambda^{-1}  \right| \cdot \left|  C_2  \right| = n^p  \left|   \Lambda^{-1}  \right| n^{\frac{p(p+1)}{2}} 2^{-p}  \left|   \Lambda  \right|^{-p-1} = n^{\frac{p(p+3)}{2}} 2^{-p} \left|   \Lambda  \right|^{-p-2}, \nonumber
\end{eqnarray}
where $C_2$ can be found as
\begin{eqnarray}
C_2 &=& \frac{n}{2} D_p^T \{ \Lambda^{-1} (\bI_p + 2\xi \xi^T \Lambda^{-1}) \otimes \Lambda^{-1}\} D_p - n D_p^T (\Lambda^{-1} \xi \otimes \Lambda^{-1})\nonumber\\
&&  \cdot n^{-1} \Lambda \cdot n(\xi^T \Lambda^{-1} \otimes \Lambda^{-1} ) D_p\nonumber \\
&=&  \frac{n}{2} D_p^T ( \Lambda^{-1}  \otimes \Lambda^{-1}) D_p  +  n D_p^T ( \Lambda^{-1} \xi \xi^T \Lambda^{-1}) \otimes \Lambda^{-1})D_p \nonumber\\
&& - n D_p^T (\Lambda^{-1} \xi \otimes \Lambda^{-1}) \cdot  (1 \otimes \Lambda) \cdot (\xi^T \Lambda^{-1} \otimes \Lambda^{-1} ) D_p\nonumber \\
&=&  \frac{n}{2} D_p^T ( \Lambda^{-1}  \otimes \Lambda^{-1}) D_p.\nonumber 
\end{eqnarray}
According to Theorem 3.14 of \cite{magnus1999}, we then have $\left|C_2\right|=n^{\frac{p(p+1)}{2}} 2^{-p} \left|\Lambda\right|^{-p-1}$.

\subsection{Quantities needed in Skovgaard's modifications (2.4) of the paper}\label{B2:gamma}
Skovgaard's  modified likelihood ratio statistics $W^{*}$ and $W^{**}$ \cite{skovgaard:2001} for hypothesis (I) in the paper can be computed explicitly based on the formula for the  correction factor $\gamma$ defined in (2.4) of the paper, i.e.
\begin{eqnarray}
\gamma 
&=& \frac{ \left[\sum_{i=1}^k {n_i}/{2} \left\{ \tr (\hat{\Lambda}_i^{-1}\hat{\Lambda}_0 \hat{\Lambda}_i^{-1} \hat{\Lambda}_0) -p \right \}\right] ^{d/2} \prod_{i=1}^k |\hat{\Lambda}_0^{-1}\hat{\Lambda}_i|^\frac{p+2}{2}}{ \left\{\sum_{i=1}^k n_i\log |\hat{\Lambda}_i\hat{\Lambda}_0^{-1}|\right\}^{d/2-1} \sum_{i=1}^k {n_i}/{2}  \left\{\tr( \hat{\Lambda}_i\hat{\Lambda}_0^{-1})  -p\right\} }. \nonumber
\end{eqnarray}

The quantities required in  the correction factor $\gamma$ of the modified likelihood ratio test proposed by \cite{skovgaard:2001} for testing the hypothesis (II) in the paper are 
\begin{eqnarray}
(\hat{\varphi}-\hat{\varphi}_\psi)^T (s-s_\psi)
&=& \frac{1}{2} \sum_{i=1}^k n_i  (\bar{y}_i-\bar{y})^T\hat{\Lambda}_i(\bar{y}_i-\bar{y})+ \frac{1}{2} \sum_{i=1}^k n_i \tr(\hat{\Lambda}_i \hat{\Lambda}_0^{-1})\nonumber\\
&&   -\frac{1}{2}np; \nonumber\\
(s-s_\psi)^T J_{\varphi\varphi}(\hat{\varphi}_\psi)^{-1} (s-s_\psi) &=&  \sum_{i=1}^k n_i (1+2\bar{y}\hat{\Lambda}_0\bar{y}) (\bar{y}_i -\bar{y})^T \hat{\Lambda}_0 (\bar{y}_i -\bar{y}) \nonumber \\
&& -2 \sum_{i=1}^k n_i (\bar{y}_i -\bar{y})^T \hat{\Lambda}_0 \left(\frac{y_i^Ty_i}{n_i}-\hat{\Lambda}_0^{-1} -\bar{y}\bar{y}^T \right) \hat{\Lambda}_0\bar{y} \nonumber \\
&&  + \frac{1}{2} \sum_{i=1}^k n_i \tr \left\{\left(\frac{y_i^Ty_i}{n_i} -\bar{y}\bar{y}^T\right) \hat{\Lambda}_0 \left(\frac{y_i^Ty_i}{n_i}  -\bar{y}\bar{y}^T\right) \hat{\Lambda}_0 \right\} \nonumber\\
&& -\frac{1}{2}np;\nonumber\\
\left\{\frac{|{J}_{\varphi \varphi}(\hat{\varphi}_\psi)|}{|{J}_{\varphi \varphi}(\hat{\varphi})|}\right\}^{1/2}  &=& |\hat{\Lambda}_0^{-1}|^\frac{p+2}{2}  \prod_{i=1}^k |\hat{\Lambda}_i|^\frac{p+2}{2}. \nonumber
\end{eqnarray}


 \subsection{Details on the computation of the directional $p$-value}\label{B3:tmintmax}

The directional $p$-value in formula (2.7) of the paper  is obtained via one-dimensional integration.  The integrand is a function of the determinant $|\hat{\Lambda}^{-1}(t)|$, and  the computational cost increases with  the dimension $p$ of the square matrix $\hat{\Lambda}^{-1}(t)$.  Based on  Lemmas 4.1--4.2 in the paper, and  according to the Jordan decomposition \citep[][Theorem 1.14]{magnus1999}, the determinant of $\hat{\Lambda}^{-1}(t)$, defined in hypothesis \ref{case1} - \ref{case5}, is such that
\begin{eqnarray}
|\hat{\Lambda}^{-1}(t)| 
=|\hat{\Lambda}_0^{-1} | \prod_{l=1}^p  (1-t+t \nu_l), \nonumber
\label{simplify:determinant case1}
\end{eqnarray}
where $\nu_l, l \in \{1,\dots,p\}$, are the eigenvalues of the matrix $\hat{\Lambda}_0 \hat{\Lambda}^{-1}$. Since the matrix $\hat{\Lambda}_0 \hat{\Lambda}^{-1}$ is constant in $t$, the eigenvalues  $\nu_l$ can be calculated only  once. Even the determinant $|\hat{\Lambda}_0^{-1} |$ does not depend on $t$, so it can be neglected when calculating the two one-dimensional integrals for the directional $p$-value in formula (2.7) of the paper.  These expedients enable to greatly speed up the computation of the $p$-value and increase the accuracy. It is however convenient to include a multiplicative constant in the integrand function in order to improve numerical stability, so that we use
\begin{eqnarray}
g(t;\psi) &\propto&  \exp \left\{(d-1)\log t + 0.5 (n-p-2) \sum _{l=1}^p \log \left (1-t+t \nu_l\right) \right\}. \nonumber
\label{simplify:integrand function case1}
\end{eqnarray}

Similarly, the determinant of $\hat{\Lambda}_i^{-1}(t),i \in \{1,\dots,k\}$, defined in hypothesis (I) can be simplified to $|\hat{\Lambda}_0^{-1} | \prod_{l=1}^p  (1-t+t \nu_{l}^i)$, where $\nu_{l}^i$ are the eigenvalues of the matrix $\hat{\Lambda}_0 \hat{\Lambda}_i^{-1}$, $i \in \{1,\dots,k\}$. The integrand function in this case can be rewritten as 
\begin{eqnarray}
g(t;\psi)
&\propto& \exp \left\{(d-1)\log t + \sum_{i=1}^k \sum_{l=1}^p 0.5 (n_i - p -2) \log \left(1-t+t\nu_{l}^i\right) \right\}.\nonumber
\end{eqnarray}
On the other hand, the  determinant  of the matrix $\hat{\Lambda}^{-1}(t)$ in the hypotheses problems (II) and  (VI)  cannot be easily simplified, since it is a quadratic function of $t$.

The 5-sigma narrower integration interval $[t_{min}, t_{max}]$ can be computed with $t_{min} = \max\{0, \hat{t} - 5 j(\hat{t};\psi)^{-1}\}$ and $t_{max} = \min\{\hat{t} + 5 j(\hat{t};\psi)^{-1}, t_{sup}\} $, where $j(\hat{t};\psi) = - \left.{\partial^2 \bar{g}(t;\psi)}/{\partial t^2}  \right|_{t=\hat{t}}$.  The second-order derivative of the integrand function for each of the hypotheses is:
\begin{itemize}
	\item hypothesis (I): 
	\begin{eqnarray}
	\frac{\partial^2 \bar{g}(t;\psi)}{\partial t^2}   &=& -(d-1)t^{-2} \nonumber\\
	&& - \sum_{i=1}^k \sum_{l=1}^p 0.5 (n_i - p -2) (1-t+t\nu_{l}^i)^{-2} (-1+\nu_{l}^i)^2, \nonumber
	\end{eqnarray}
	where $\nu_{l}^i,l \in \{1,\dots,p\}$, are the eigenvalues of the matrix $\hat{\Lambda}_0\hat{\Lambda}_i^{-1}$ for the $i$-th group;
	\item  hypothesis (II): the estimator equals $\hat{\Lambda}_i^{-1}(t) = (1-t) \hat{\Lambda}_0^{-1} + t \hat{\Lambda}_i^{-1} + t(1-t) B_i$, where $B_i = (\bar{y}_i - \bar{y})(\bar{y}_i-\bar{y})^T$. Then we find
	\begin{eqnarray}
	\frac{\partial^2 \bar{g}(t;\psi) }{\partial t^2} 
	&=& -(d-1)t^{-2} \nonumber \\
	&& - \sum_{i=1}^k 0.5 (n_i - p -2)\tr \left\{\hat{\Lambda}_i(t) \frac{\partial \hat{\Lambda}_i^{-1} (t)}{\partial t} \hat{\Lambda}_i(t)\frac{\partial \hat{\Lambda}_i^{-1} (t)}{\partial t}\right\}\nonumber\\
	&&  - \sum_{i=1}^k 0.5 (n_i - p -2)\tr \left\{2\hat{\Lambda}_i(t)B_i \right\}, \nonumber
	\end{eqnarray}
	where ${\partial \hat{\Lambda}_i^{-1} (t)}/{\partial t} = -\hat{\Lambda}_0^{-1} +\hat{\Lambda}_i^{-1} + (1-2t)B_i$.
		\item  hypothesis (III)-(V)
	\begin{eqnarray}
	\frac{\partial^2 \bar{g}(t;\psi)}{\partial t^2}   &=& -(d-1)t^{-2}  -  0.5 (n - p -2)\sum_{l=1}^p  (1-t+t\nu_{l})^{-2} (1-\nu_{l})^2, \nonumber
	\end{eqnarray}
	where $\nu_{l}, l \in \{1,\dots,p\}$, are the eigenvalues of the matrix $\hat{\Lambda}_0 \hat{\Lambda}^{-1}$;
	
	\item hypothesis (VI):
	the estimator of $\Lambda^{-1}$ is
	$\hat{\Lambda}^{-1}(t) = (1-t) \bI_p + t \hat{\Lambda}^{-1} + t(1-t) \bar{y}\bar{y}^T$. Then
	\begin{eqnarray}
	\frac{\partial^2  \bar{g}(t;\psi)}{\partial t^2}
	&=& -(d-1)t^{-2} \nonumber\\
	&& - 0.5 (n - p -2) \tr \left\{\hat{\Lambda}(t) \frac{\partial \hat{\Lambda}^{-1} (t)}{\partial t} \hat{\Lambda}(t)\frac{\partial \hat{\Lambda}^{-1} (t)}{\partial t}  + 2\hat{\Lambda}(t)\bar{y}\bar{y}^T \right\}, \nonumber
	\end{eqnarray}
	where  ${\partial \hat{\Lambda}^{-1} (t)}/{\partial t} = -\bI_p +\hat{\Lambda}^{-1} + (1-2t)\bar{y}\bar{y}^T$.

\end{itemize}




%
%
\section{Directional test for one-sample hypotheses}\label{one-sample cases}
%
%
%

\subsection{Testing conditional independence} \label{case1}
As in \citep[][Section 5.3]{sartori:2014}, we first focus on testing  the conditional independence of normal random components, meaning that some off-diagonal  elements of the concentration matrix $\Lambda$ are equal to zero.
Indeed, a zero entry in the concentration matrix $\Lambda$ implies the conditional independence of the two corresponding  variables, given the others.  The hypothesis can be formulated as  
\begin{eqnarray}
H_{\psi}: \Lambda = \Lambda_0,
\label{hypothesis:condition independence}
\end{eqnarray}
where the matrix $\Lambda_0$ is unknown but with some off-diagonal elements equal to zero.  The  maximum likelihood estimates of parameters $\mu$ and $\Lambda^{-1}$ are those given in Section \ref{section:multivariate normal distribution}. The constrained maximum likelihood estimate of $\mu$ coincides with the unconstrained one, while that  of  $\Lambda^{-1}$ is denoted by $\hat{\Lambda}^{-1}_0$, and  is typically obtained numerically. For instance, the functions \texttt{fitConGraph} and \texttt{cmod} in the R  \citep{Rcode} packages \texttt{ggm}  \citep{ggm} and \texttt{gRim} \citep{gRim}, respectively, can be used.   The log-likelihood ratio test statistic simplifies to
\begin{eqnarray}
W = - n \log |\hat{\Lambda}^{-1}\hat{\Lambda}_0|,
\label{LRT case1}
\end{eqnarray}
and follows approximately a  $\chi^2_d$ distribution, with $d$ equal to the difference  between the number of free parameters under the alternative and  null hypotheses. However, when $p$ is large relative to the sample size $n$, the chi-square approximation to the distribution of $W$  fails \citep{jiang:2013,jiang:2015,He2020B}. In fact, assuming that $p$ depends on $n$, i.e. $p=p_n$, in some particular cases, such as the hypotheses (\ref{hypothesis: case1}), (\ref{hypothesis:case2}) and (\ref{hypothesis:case6}) below, the chi-square approximation to the log-likelihood ratio statistic  works if and only if $p =o(n^{1/2})$. The analogous condition for its Bartlett correction is $p=o(n^{2/3})$ \cite{He2020B}. 

 Skovgaard's modifications \cite{skovgaard:2001} can be computed easily \citep[][(22)]{sartori:2014}, using
\begin{eqnarray}
\gamma = \frac{ \left[{n}/{2} \left\{ \tr (\hat{\Lambda}^{-1}\hat{\Lambda}_0 \hat{\Lambda}^{-1} \hat{\Lambda}_0) -p \right \}\right] ^{d/2} |\hat{\Lambda}_0^{-1}\hat{\Lambda}|^\frac{p+2}{2}}{ \left\{n\log |\hat{\Lambda}\hat{\Lambda}_0^{-1}|\right\}^{d/2-1} {n}/{2}  \left\{ \tr(\hat{\Lambda}\hat{\Lambda}_0^{-1})   -p\right\} }.\nonumber
\label{skovgaard gamma}
\end{eqnarray}
Also $W^*$ and $W^{**}$ have approximate $\chi^2_d$ distributions, when $p$ is fixed and $n \to \infty$.


In order to  obtain  the directional test in this setting,  we derive  the expected value ${s_{\psi}}$ under $H_\psi$  as
\begin{eqnarray}
s_{\psi} = -\ell_\varphi (\hat{\varphi}_\psi)=-\left\{{0}_p^T,\; \frac{n}{2} \text{vech} \left( \hat{\Lambda}_0^{-1}  - \hat{\Lambda}^{-1} \right)^T \right\}^T.\nonumber
\end{eqnarray}
The tilted log likelihood  for the canonical parameter $\varphi = \{\xi^T, \text{vech} (\Lambda)^T\}^T$ along the line $s(t)=(1-t)s_\psi$ takes the form
\begin{eqnarray}
\ell(\varphi;t) 
&=& n \xi^T \bar{y} -\frac{n}{2}\tr \left[\Lambda \left\{\frac{y^Ty}{n}+(1-t)\left(\hat{\Lambda}_0^{-1}-\hat{\Lambda}^{-1} \right)\right\}\right] \nonumber\\
&&  + \frac{n}{2} \log |\Lambda| -\frac{n}{2} \xi^T \Lambda^{-1}\xi,  \nonumber
\end{eqnarray}
and its maximization  yields $\hat{\varphi}(t) = \left[\hat{\xi}(t)^T,\text{vech} \{ \hat{\Lambda}(t)\}^T \right]^T$, with $\hat{\xi}(t)=\hat{\Lambda}(t) \hat{\mu}$ and $\hat{\Lambda}^{-1}(t) = (1-t)\hat{\Lambda}_0^{-1} + t\hat{\Lambda}^{-1}$.  Hence, the saddlepoint approximation (2.5) in the paper to the density of $s$ along the line $s(t)$ is 
\begin{eqnarray}
h\{s(t);\psi\} 
&=& c\exp\left\{  \frac{n-p-2}{2} \log |\hat{\Lambda}^{-1}(t)| \right\},
\label{saddle case1}
\end{eqnarray}
where $c$ is a normalizing constant. 

The value $t_{sup}$ in (2.7) of the paper is the largest $t$ for which $\hat{\Lambda}(t)^{-1}$ is positive definite and is equal to $\{1-  \nu_{(1)} \}^{-1}$ where $\nu_{(1)}$ is the smallest eigenvalue of $\hat{\Lambda}_0 \hat{\Lambda}^{-1}$ (see Lemma 4.1 in Section 4 of the paper).

 Thanks to the exactness of the saddlepoint approximation (2.5) in the paper to the density of $s$, it is possible to prove the following theorem. 


\begin{thm}
	Assume that $p=p_n$ such that $n\ge p+2$ for all $n\ge 3$. Then, under the null hypothesis $H_\psi$ in (\ref{hypothesis:condition independence}), the directional $p$-value (2.7) in the paper is exactly uniformly distributed. 
	\label{theoremcase1}
\end{thm}


\begin{proof}
	Suppose $y=[y_1 \cdots y_n]^T$ with $y_i \sim N_p(\mu, \Lambda^{-1})$, $i \in \{1,\dots,n\}$. Let $\bar{y} = n^{-1} {{1}}_n^T y$ and $A = y^Ty -n \bar{y}\bar{y}^T$. It is well known that $A$ has a Wishart distribution, i.e. $A \sim W_p(n-1,\Lambda^{-1})$ with density function  \citep[][Theorem 3.2.1]{muirhead1982}
	\begin{eqnarray}
	f(A;\Lambda^{-1}) &=&   2^{-\frac{p(n-1)}{2}} \Gamma_p \left(\frac{n-1}{2}\right)^{-1} |\Lambda^{-1}|^{-\frac{n-1}{2}} \text{etr}\left(-\frac{1}{2} \Lambda A\right) |A|^{\frac{n-p-2}{2}},\nonumber
	\end{eqnarray}
	where $\text{etr}(\cdot) = \exp\{\tr(\cdot)\}$ with $\tr(\cdot)$ the trace operator and $\Gamma_p(\cdot)$ is a multivariate gamma function. 
	
	Moreover, the canonical parameter for the multivariate normal distribution is $\varphi = \{\xi^T,\text{vech}(\Lambda)^T\}^T = \{\mu^T\Lambda ,\text{vech}( \Lambda)^T\}^T$, with corresponding log-likelihood 
	\begin{eqnarray}
	l(\varphi;s) ={-\frac{n}{2}} \log|\Lambda^{-1}| -\frac{n}{2} \tr ( \Lambda \hat{\Lambda}^{-1}) -\frac{n}{2}(\bar{y} -\Lambda^{-1}\xi )^T \Lambda (\bar{y}-\Lambda^{-1}\xi ), \nonumber
	\label{proof:loglikegroup1}
	\end{eqnarray}
	where $\hat{\Lambda}^{-1} = n^{-1} A$. The maximum likelihood estimate  of $\varphi$ is  $\hat{\varphi} = \{\hat{\xi}^T , \text{vech}(\hat{\Lambda})^T \}^T$ $ = \{ \bar{y}^T\hat{\Lambda}, \text{vech}(\hat{\Lambda})^T\}^T$
	and the constrained  maximum likelihood estimate is $  \hat{\varphi}_{\psi} = \{\hat{\xi}^T_{\psi} , \text{vech}( \hat{\Lambda}_0)^T\}^T$ $= \{\bar{y}^T\hat{\Lambda}_0 ,\text{vech}( \hat{\Lambda}_0)^T\}^T$.
	The corresponding  log-likelihoods at $\hat{\varphi}$ and $\hat{\varphi}_{\psi}$ are $\ell(\hat{\varphi};s) = {-{n}/{2}} \log|\hat{\Lambda}^{-1}| -{(np)}/{2}$ and $\ell(\hat{\varphi}_{\psi};s)
	= {-{n}/{2}} \log|\hat{\Lambda}_0^{-1}| -{n}/{2} \tr (\hat{\Lambda}_0 \hat{\Lambda}^{-1})$, respectively. Then,  under the null hypothesis $H_{\psi}$, the main factor of the saddlepoint approximation of $s$  takes the form
	\begin{eqnarray}
	\exp\left\{\ell(\hat{\varphi}_{\psi};s) - \ell(\hat{\varphi};s) \right\} &=& \exp\left\{-\frac{n}{2} \log \frac{|\hat{\Lambda}_0^{-1}|}{|\hat{\Lambda}^{-1}|} -\frac{n}{2} \tr \left(\hat{\Lambda}_0 \hat{\Lambda}^{-1} -\bI_p\right)\right\} \nonumber \\
	&=& |\hat{\Lambda}_0^{-1}|^{-\frac{n}{2}} |\hat{\Lambda}^{-1}|^{\frac{n}{2}} \exp\left\{-\frac{n}{2} \tr (\hat{\Lambda}_0 \hat{\Lambda}^{-1}) \right\}  \exp \left(\frac{np}{2}\right).\nonumber
	\end{eqnarray}
	As $|J_{\varphi\varphi}(\hat{\varphi})| \propto |\hat{\Lambda}^{-1}|^{p+2}$ (see Supplementary Material S1.1), the saddlepoint approximation  results 
	\begin{eqnarray}
	h(s;\psi) 
	&=& c \; |\hat{\Lambda}_0^{-1}|^{-\frac{n}{2}} \exp\left\{-\frac{n}{2} \tr (\hat{\Lambda}_0 \hat{\Lambda}^{-1})\right\}  |\hat{\Lambda}^{-1}|^{\frac{n-p-2}{2}}.
	\label{proof saddle case1}
	\end{eqnarray}
	Suppose that the normalizing constant $c = \left({n}/{2}\right)^{p(n-1)/2} \Gamma_p\left\{{(n-1)}/{2}\right\}^{-1} |\hat{\Lambda}_0|$,  
	the saddlepoint approximation (\ref{proof saddle case1}) is the exact conditional distribution of $\hat{\Lambda}^{-1}$ given $\hat{\Lambda}_0^{-1}$, which is a Wishart random variable with parameters $(n-1, \hat{\Lambda}_0^{-1})$ under $H_\psi$. Here we used the result  $\tr(\hat{\Lambda}_0\hat{\Lambda}^{-1}) =p$ \citep[][page 278]{Eriksen1996}.
	The saddlepoint approximation density along the line  $s(t)$ replaces $\hat{\Lambda}^{-1}$ with $\hat{\Lambda}^{-1}(t)$, which maximizes $\ell\{\varphi;s(t)\}$.
	Since $\hat{\Lambda}^{-1}(t)=(1-t)\hat{\Lambda}_0^{-1} + t\hat{\Lambda}^{-1} $, the exact saddlepoint approximation for the distribution of $s(t)$  under $H_{\psi}$ is simply obtained from (\ref{proof saddle case1}) as
	\begin{eqnarray}
	h\{s(t);\psi\} &=&   c \; |\hat{\Lambda}_0^{-1}|^{-\frac{n}{2}} \exp\left[-\frac{n}{2} \tr\{ \hat{\Lambda}_0 \hat{\Lambda}(t)^{-1}\}\right] |\hat{\Lambda}(t)^{-1}|^{\frac{n-p-2}{2}} \nonumber \\
	&\propto& 
	\exp\left\{ \frac{n-p-2}{2} \log |\hat{\Lambda}(t)^{-1}| \right\}. 
	\label{proof:saddleL0 case1}
	\end{eqnarray}
	In the last step we used $\tr\{\hat{\Lambda}_0 \hat{\Lambda}^{-1}(t)\}$ $=\tr \left[ \hat{\Lambda}_0 \left\{(1-t)\hat{\Lambda}_0^{-1} + t\hat{\Lambda}^{-1} \right\} \right]=\tr\left\{(1-\right.$  $t)\bI_p\left.+ t\hat{\Lambda}_0\hat{\Lambda}^{-1} \right\}$ $= (1-t)p + tp = p$. 
	Since  the saddlepoint approximation $ h\{s(t);\psi\}$ (\ref{proof:saddleL0 case1}) is exact, apart from the normalizing constant, the integral in the denominator of the directional $p$-value (2.7) in the paper  is just the normalizing constant of the conditional distribution of $||s||$ given the direction $s/||s||$.  Therefore, the directional $p$-value is the exact probability of $||s|| > ||s^0||$ given the direction $s/||s||$ under the null hypothesis, and  is thus exactly uniformly distributed. 
\end{proof}

Theorem \ref{theoremcase1} only requires $n\ge p+2$ for ensuring that the maximum likelihood estimate of the covariance matrix exists with probability one. This assumption is also weaker than the condition $p/n \to \kappa \in (0,1]$  in \cite{jiang:2013} for the validity of their central limit theorem approximation with large $p$.

Following \cite{jiang:2013}, we consider three specific  hypotheses of the form (\ref{hypothesis:condition independence}) of potential interest. In each case we give details on how to find the  constrained maximum likelihood estimate under the null hypothesis $H_{\psi}$. First, we focus on testing whether the covariance matrix of a  multivariate normal distribution is proportional to the identity matrix. The hypothesis of interest is then
\begin{eqnarray}
H_{\psi}: \Lambda^{-1} = \sigma^2 \bI_p,
\label{hypothesis: case1}
\end{eqnarray}
where $\sigma^2$ is unspecified and $\bI_p$ denotes the $p\times p$ identity matrix. This corresponds to checking if the covariance  matrix is diagonal with equal elements.  In this particular case, it is equivalent to testing that  the concentration matrix $\Lambda$ is proportional to the identity matrix, i.e. $H_{\psi}: \Lambda = \sigma^{-2} \bI_p$ with unspecified $\sigma^{-2}$. The constrained maximum likelihood estimate of $ \Lambda^{-1}$ is  $\hat{\Lambda}^{-1}_0 = { \tr(y^Ty/n - \bar{y}\bar{y}^T)} / {p} \; \bI_p$, and  can be used in the formula for the log-likelihood statistic (\ref{LRT case1}),  Skovgaard's statistics  (2.3) in the paper  and the directional $p$-value (2.7) in the paper.

The second hypothesis concerns block-independence in the multivariate normal distribution. For $k \ge 2$, let $p_1,\dots,p_k$ be positive integers with $p=\sum_{j=1}^k p_j$. Then the covariance matrix can be expressed as  $\Lambda^{-1}=\left(\Lambda^{-1}_{ij}\right)_{p\times p}$,  where $\Lambda^{-1}_{ij}$ is a $p_i \times p_j$ sub-matrix for all $1 \le i ,j \le k$.  We are interested in testing the null hypothesis
\begin{eqnarray}
H_{\psi}: \Lambda^{-1}_{ij} = {0}, \quad  \; 1\le i<j \le k,
\label{hypothesis:case2}
\end{eqnarray}
which is equivalent to testing that the off-diagonal matrix $\Lambda_{ij}=0$, $ 1\le i<j \le k$. The constrained maximum likelihood estimate of $\Lambda^{-1}$ is then
\begin{eqnarray}
\hat{\Lambda}^{-1}_0 = \left(\begin{array}{cccc}
\hat{\Lambda}^{-1}_{11} & 0 & \cdots & 0 \\
0 &  \hat{\Lambda}^{-1}_{22} & \cdots & 0 \\
\vdots &  & & \vdots \\
0 & 0 & \cdots &  \hat{\Lambda}^{-1}_{kk}
\end{array}\right),\nonumber
\end{eqnarray}
where $ \hat{\Lambda}^{-1}_{ii}= [\hat{\Lambda}^{-1}]_{p_ip_i}, i \in \{1,\dots,k\}$, with $[\hat{\Lambda}^{-1}]_{p_ip_i}$ denoting   the $p_i \times p_i$ diagonal sub-matrix of $\hat{\Lambda}^{-1}$.

Lastly, we consider the complete-independence hypothesis.  Let $\bR=(r_{ij})_{p\times p}$ be the correlation matrix corresponding to $\Lambda^{-1}$.   The hypothesis of interest here is 
\begin{eqnarray}
H_{\psi}: \bR = \bI_p,
\label{hypothesis:case6}
\end{eqnarray}
which implies the covariance  (or concentration) matrix to be  diagonal.
The constrained maximum likelihood estimate of $\Lambda^{-1}$ equals $\hat{\Lambda}^{-1}_0 = \diag (\hat{\Lambda}^{-1})\; \bI_p$, where $\diag (\hat{\Lambda}^{-1})$ is a $p$-vector with entries equal to the diagonal elements of the unconstrained maximum likelihood estimate $\hat{\Lambda}^{-1}$. Simulation studies for the hypothesis (\ref{hypothesis: case1}) - (\ref{hypothesis:case6}) are reported in Section \ref{S4:simulation (III)-(VI)}.


\subsection{Testing a specific multivariate normal distribution } \label{case5}

Under the same framework introduced in Section \ref{case1}, let us consider the null hypothesis
\begin{eqnarray}
H_{\psi}: \mu = \mu_0\; \text{and}\; \Lambda = \Lambda_0 , \nonumber
\end{eqnarray}
where $\mu_0 \in \mathbb{R}^p$  and the  $p \times p$ non-singular matrix $\Lambda_0$   are completely specified. Considering the standardized data $\tilde{y}_i = \Lambda_0^{1/2} (y_i - \mu_0)$, which has distribution $N_p(\tilde{\mu},\tilde{\Lambda}^{-1})$ with $\tilde{\mu} = \Lambda_0^{1/2}(\mu - \mu_0)$ and $\tilde{\Lambda}^{-1} = \Lambda_0^{1/2} \Lambda^{-1} \Lambda_0^{1/2}$, where $\Lambda_0^{1/2}$ is a square root matrix of $\Lambda_0$, the above hypothesis is equivalent to 
\begin{eqnarray}
H_{\psi}: \tilde{\mu} = {0}_p\; \text{and}\; \tilde{\Lambda} = \bI_p  .
\label{hypothesis case5}
\end{eqnarray}
Under $H_\psi$ we also have $ \tilde{\Lambda}^{-1}= \bI_p$, since joint marginal independence is equivalent to  joint conditional independence.
The log-likelihood statistic in the canonical parameterization is
\begin{eqnarray}
W &=& (n-1)\log |\hat{\Lambda}|+\tr(y^Ty)-np. \nonumber
\label{LRT}
\end{eqnarray}
It can be shown that $W$ follows a $\chi^2_d$  asymptotic null distribution with $d= p(p+1)/2+p=p(p+3)/2$ if and only if $p=o(n^{1/2})$. The analogous condition for Bartlett correction of $W$ is  $ p =o(n^{2/3})$ \citep{He2020B}.

The correction factor  $\gamma$ in the modified log-likelihood ratio test proposed by \cite{skovgaard:2001} is
\begin{eqnarray}
\gamma = \frac{ \left\{{n}/{2} \tr \left({y^Tyy^Ty}/{n^2}\right)-n \tr(\hat{\Lambda}^{-1})  +{np}/{2} \right\} ^{d/2} |\hat{\Lambda}|^{(p+2)/{2}}} { \left\{n\log |\hat{\Lambda}|+\tr(y^Ty)-np\right\}^{d/2-1} {n}/{2} \tr \left(\hat{\Lambda} \bar{y} \bar{y}^T + \hat{\Lambda} +{y^Ty}/{n}  -2 \bI_p\right) }. \nonumber
\end{eqnarray}

Moreover, for the directional test, the expected value of $s$ under $H_\psi$ is
\begin{eqnarray}
s_{\psi} 
&=& -\left\{n\bar{y}^T,\;  \frac{n}{2}\text{vech}\left( \bI_p-\frac{y^T y}{n}\right)^T \right\}^T. \nonumber
\end{eqnarray}
The tilted log-likelihood along the line $s(t)=(1-t)s_\psi$ can then be expressed as
\begin{eqnarray}
\ell(\varphi;t)&=&n\xi^T\{\bar{y}-(1-t)\bar{y}\}- \frac{n}{2}\tr \left[ \Lambda \left\{\frac{y^T y}{n} +(1-t) \left(\bI_p-\frac{y^T y}{n}\right)\right\}\right] \nonumber\\
&& +\frac{n}{2}\log |\Lambda| -\frac{n}{2}\xi^T\Lambda^{-1}\xi, \nonumber
\end{eqnarray}
with corresponding  maximum likelihood estimate  $\hat{\varphi}(t) = \left[\hat{\mu}(t)^T\hat{\Lambda}(t) ,\text{vech}\left\{ \hat{\Lambda}(t)\right\}^T\right]^T$,  with $\hat{\mu}(t) = t\bar{y}$ and $\hat{\Lambda}^{-1}(t) = (1-t) \bI_p + t\hat{\Lambda}^{-1} +t(1-t) \bar{y}\bar{y}^T$. 
Finally, the saddlepoint approximation (2.5) in the paper to the density of $s$ along the line $s(t)$ is
\begin{eqnarray}
h\{s(t);\psi\} 
&=& c \exp\left\{\frac{n-p-2}{2} \log |\hat{\Lambda}^{-1}(t)|+\frac{nt}{2} \tr \left(\bI_p - \frac{y^T y}{n}\right)  \right\}, \nonumber
\label{saddle case5}
\end{eqnarray}
where $c$ is a normalizing constant. The value $t_{sup}$ in (2.7) of the paper is the largest $t$ for which $\hat{\Lambda}(t)^{-1}$ is positive definite and is found iteratively. Since the saddlepoint approximation to the density of $s$ is exact also in this case, then we can derive the following theorem to explain why the directional $p$-value is exact.

\begin{thm}
	Assume that $p=p_n$ such that $n\ge p+2$ for all $n\ge 3$. Then, under the null hypothesis $H_{\psi}$ (\ref{hypothesis case5}), the directional $p$-value (2.7) of the paper is exactly uniformly distributed.
	\label{theoremcase5}
\end{thm}



\begin{proof}
	Since $\bar{y} \sim N_p(\mu, n^{-1} \Lambda^{-1})$ and $A \sim W_p(n-1,\Lambda^{-1})$  are independent  \citep[][Theorem 3.1.2]{muirhead1982},  the joint density of $\bar{y}$ and $A$ is
	\begin{eqnarray}
	f(\bar{y},A;\mu,\Lambda^{-1}) &=&   (2\pi)^{-\frac{p}{2}} |\Lambda^{-1}|^{-\frac{1}{2}} \exp \left\{-\frac{n}{2} (\bar{y} -\mu) ^T \Lambda (\bar{y} - \mu)\right\}  \nonumber\\
	&& \times 2^{-\frac{p(n-1)}{2}} \Gamma_p \left(\frac{n-1}{2}\right)^{-1} |\Lambda^{-1}|^{-\frac{n-1}{2}} \text{etr}\left(-\frac{1}{2} \Lambda A\right) |A|^{\frac{n-p-2}{2}}.\nonumber
	\end{eqnarray}
	Similarly to the proof of Theorem \ref{theoremcase1}, under the null hypothesis $H_\psi$, and using the fact that $|J_{\varphi \varphi}(\hat{\varphi})|$ is proportional to $|\hat{\Lambda}^{-1}|^{p+2}$  (see Supplementary Material S1.1), 
	the saddlepoint approximation is
	\begin{eqnarray}
	h(s;\psi) 
	&=& c_1 \;  |\Lambda_0^{-1}|^{-\frac{1}{2}} \exp \left\{ -\frac{n}{2} (\bar{y} - \mu_0)^T \Lambda_0 (\bar{y}-\mu_0)\right\} \nonumber \\
	&& \times  c_2 \; |\Lambda^{-1}_0|^{-\frac{n-1}{2}} \exp\left\{-\frac{n}{2} \tr (\Lambda_0 \hat{\Lambda}^{-1})\right\}  |\hat{\Lambda}^{-1}|^{\frac{n-p-2}{2}},
	\label{proof saddle case5}
	\end{eqnarray}
	with $c_1 = (2\pi)^{-p/2}$ and $c_2 = \left({n}/{2}\right)^{p(n-1)/2} \Gamma_p\left\{({n-1})/{2}\right\}^{-1}$. Expression (\ref{proof saddle case5}) is the exact joint distribution of $\bar{y}$ and $\hat{\Lambda}^{-1}$. 
	For obtaining the directional test, we just replace $\bar{y}$ and  $\hat{\Lambda}^{-1}$ with the maximizers of $\ell\{\varphi;s(t)\}$, namely  $\hat{\mu}(t) = t\bar{y}$ and $\hat{\Lambda}^{-1}(t) = (1-t)\bI_p + t\hat{\Lambda}^{-1}+t(1-t)\bar{y}\bar{y}^T$, respectively. Then, the saddlepoint approximation for the distribution of  $s(t)$ under  $H_{\psi}$ can be obtained from (\ref{proof saddle case5}) as
	\begin{eqnarray}
	h\{s(t);\psi\} &=& c_1 \; |{\Lambda}^{-1}_0|^{-\frac{1}{2}} \exp \left[ -\frac{n}{2} \{\hat{\mu}(t) - {\mu}_0\}^T {\Lambda}_0 \{\hat{\mu}(t)-{\mu}_0\}\right] \nonumber \\
	&& \times  c_2 \; |{\Lambda}_0^{-1}|^{-\frac{n-1}{2}} \exp\left[-\frac{n}{2} \tr \{ {\Lambda}_0 \hat{\Lambda}(t)^{-1}\}\right] |\hat{\Lambda}(t)^{-1}|^{\frac{n-p-2}{2}} \nonumber \\
	&\propto&  \exp\left\{\frac{nt}{2} \tr \left(\bI_p - \frac{y^Ty}{n}\right) + {\frac{n-p-2}{2}} \log |\hat{\Lambda}(t)^{-1}| \right\}, \nonumber
	\end{eqnarray}
	where we  use  ${\mu}_0 = {{0}}_p$ and ${\Lambda}^{-1}_0=\bI_p$. 
	The directional $p$-value is then exactly uniformly distributed, since the saddlepoint approximation $h\{s(t);\psi\}$ is exact.
\end{proof}


\section{Numerical results in multiple-sample hypotheses (I) and (II)}\label{S3:simulation (I)-(II)}

This section investigates the Monte Carlo simulations referred to the hypotheses which are testing  the equality of covariance matrices in $k$ independent groups and equality of multivariate normal distributions in $k$ independent groups. 
We reports the empirical distribution of $p$-values, the limiting null distribution, estimated size and power of the directional test, the central limit theorem test  proposed by Jiang and Yang \cite{jiang:2013}, the log-likelihood ratio test ($W$), its Bartlett correction ($W_{BC}$) and the tests $W^*$ and $W^{**}$ proposed by Skovgaard \cite{skovgaard:2001}. The  numerical results are based on $100,000$ replications. The empirical distribution of $p$-values for the six tests is examined by comparison with the Uniform$(0,1)$ distribution.
We also obtain the simulated null distribution curve for  the statistics $W$, $W_{BC}$,  $W^*$ and $W^{**}$, and compare them with their corresponding theoretical chi-square distribution.  Similarly, we compare the simulated null distribution of the central limit  proposed by Jiang and Yang \cite{jiang:2013} with its asymptotic standard normal distribution.  In order to compare the directional approach  with its main competitor, the central limit theorem test, we transform the  directional $p$-value via  the quantile function of a $N(0,1)$. If the directional $p$-value is exactly uniformly distributed, this transformation should be exacted $N(0,1)$.



\begin{figure}[H]
	\centering
	\captionsetup{font=footnotesize}
	\includegraphics[scale=0.08]{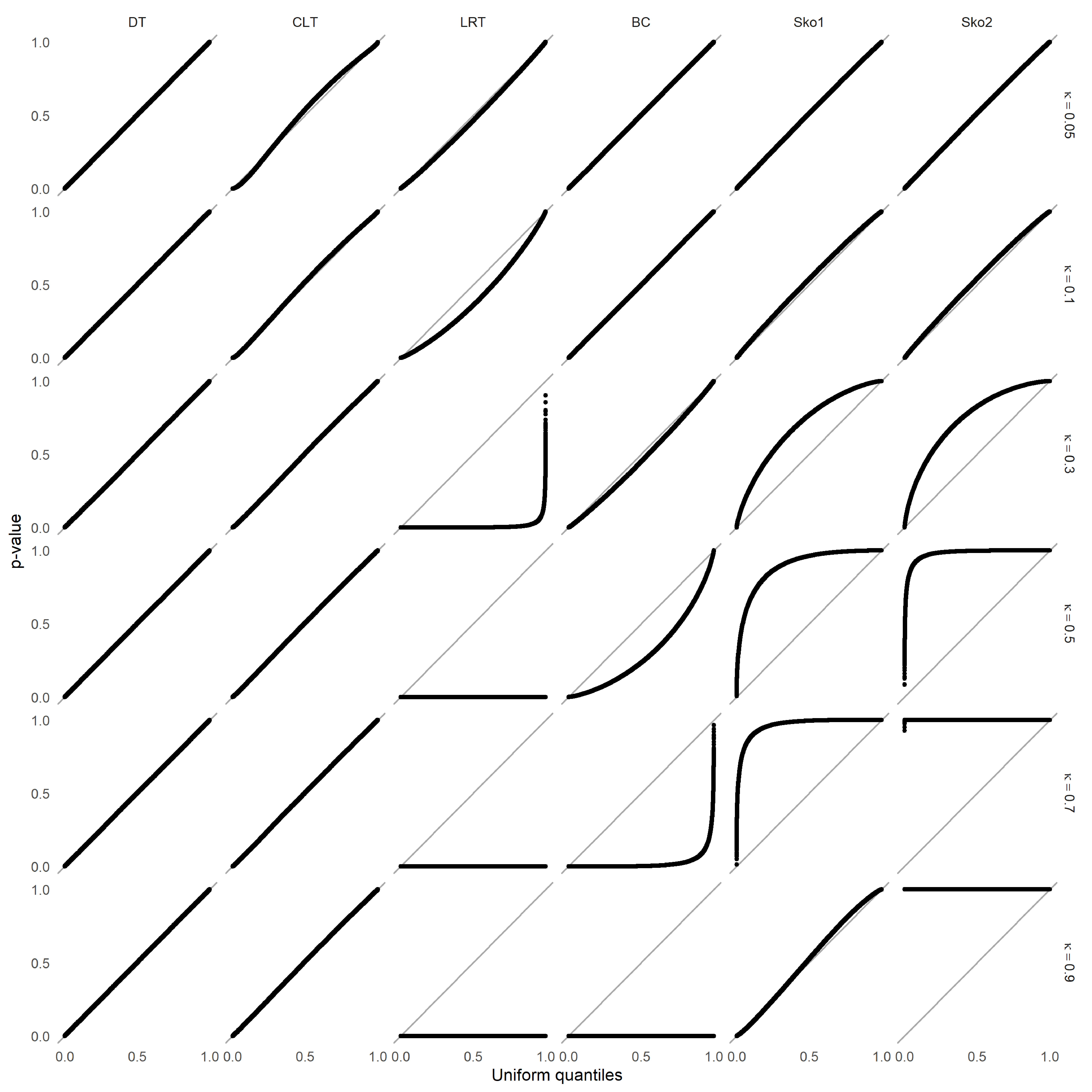}
	\caption{Hypothesis (I). Empirical null  distribution of $p$-values for various ratio $p/n_i$ of the directional test (DT), the central limit theorem test (CLT), log-likelihood ratio test (LRT), Bartlett correction (BC),  and two Skovgaard's modifications \cite{skovgaard:2001} (Sko1 and Sko2, respectively), compared
		with the $U(0,1)$ given by the gray diagonal.}
	\label{SMsimulation:pvalue case4}
\end{figure}

\begin{figure}[H]
	\centering
	\captionsetup{font=footnotesize}
	\includegraphics[scale=0.08]{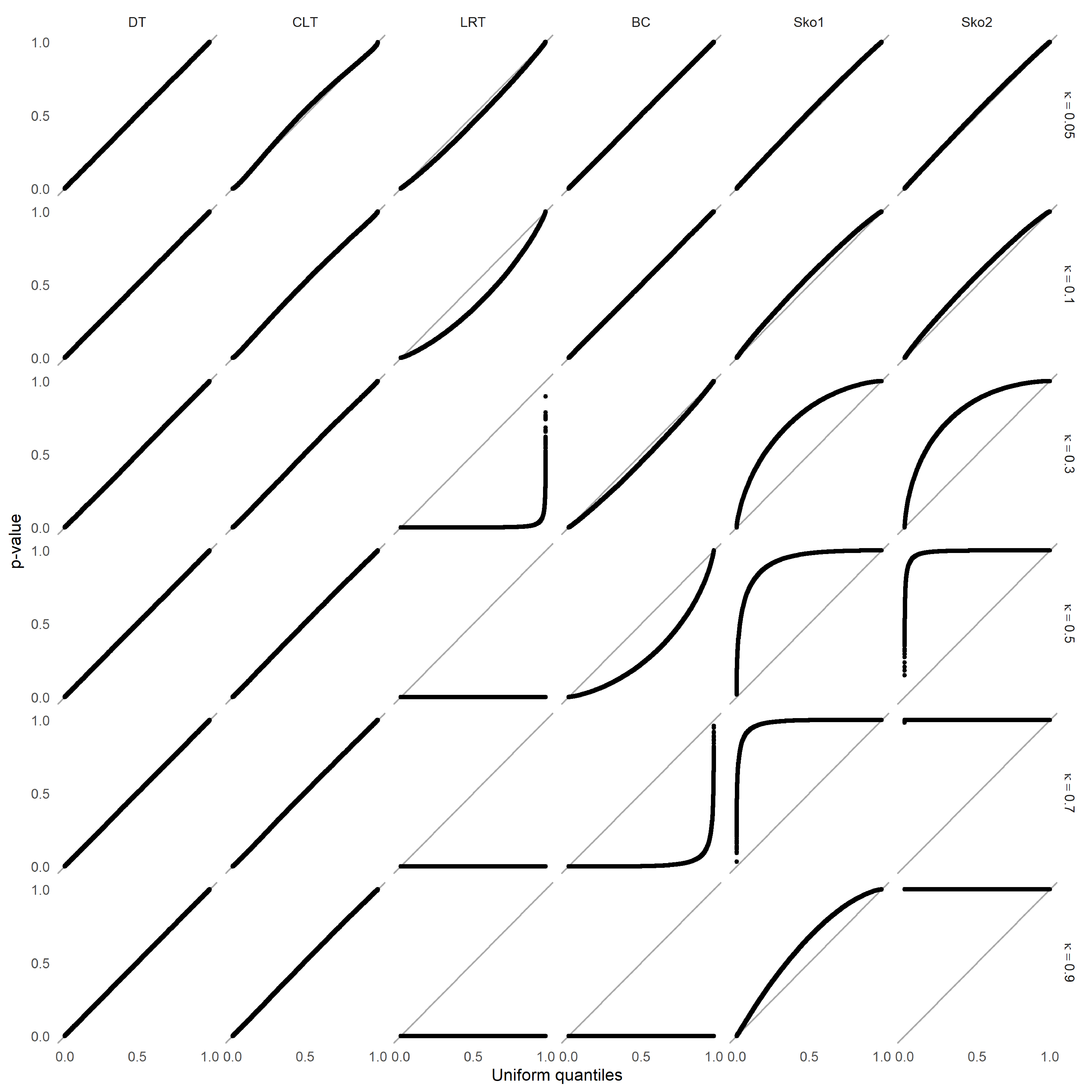}
	\caption{Hypothesis (II). Empirical null  distribution of $p$-values for various ratio $p/n_i$ of the directional test (DT), the central limit theorem test (CLT), log-likelihood ratio test (LRT), Bartlett correction (BC),  and two Skovgaard's modifications \cite{skovgaard:2001} (Sko1 and Sko2, respectively), compared
		with the $U(0,1)$ given by the gray diagonal.}
	\label{SMsimulation:pvalue case3}
\end{figure}

\begin{figure}[H]
	\centering
	\captionsetup{font=footnotesize}
	\includegraphics[scale=0.08]{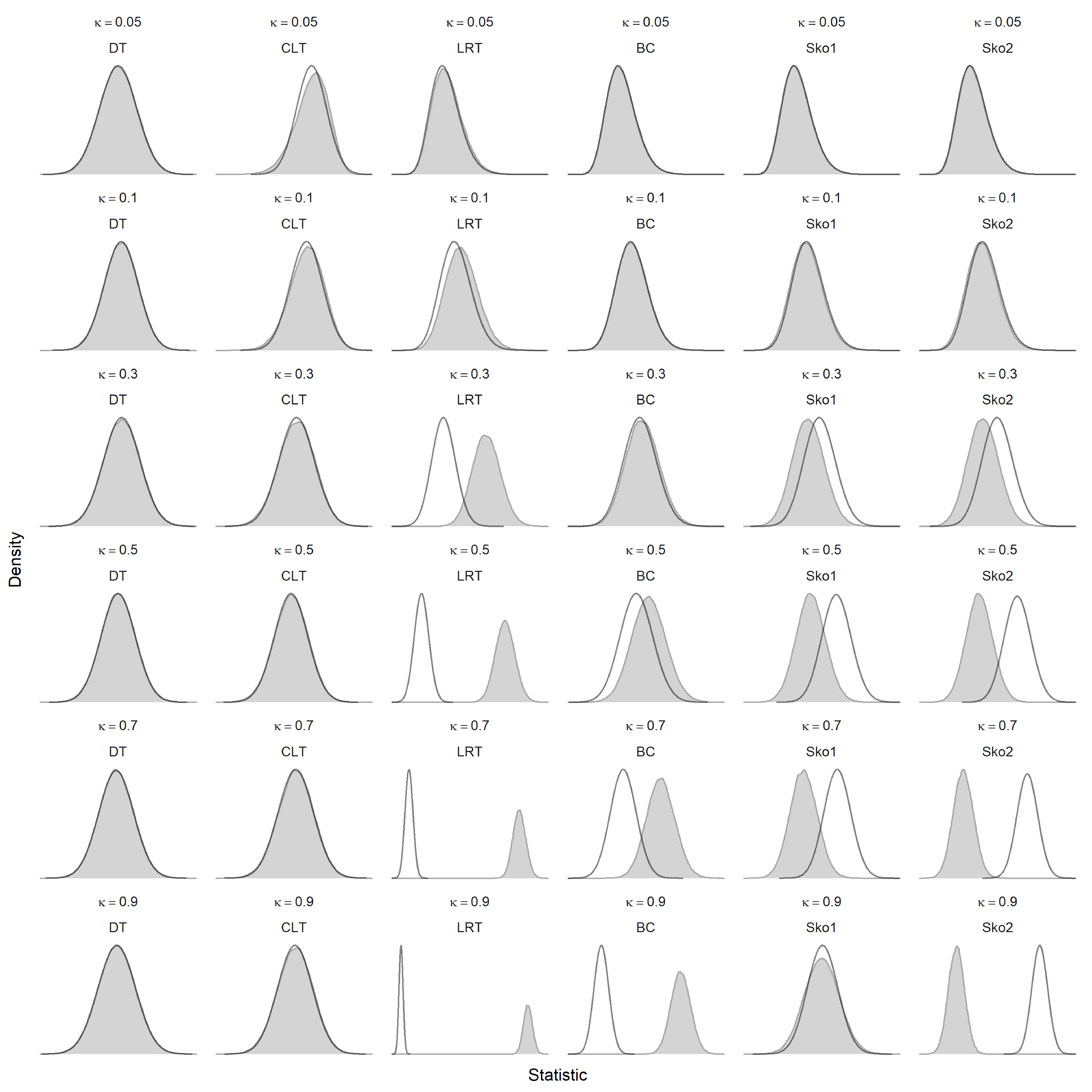}
	\caption{Hypothesis (I). Comparison between   empirical (gray) and theoretical (black) null  distributions of  the directional test (DT), the central limit theorem test (CLT), log-likelihood ratio test (LRT), Bartlett correction (BC),  and two Skovgaard's modifications \cite{skovgaard:2001} (Sko1 and Sko2, respectively)  for  various values of $p/n_i$, $i \in \{1,\dots,k\}$.}
	\label{SMsimulation:density null case4}
\end{figure}

\begin{figure}[H]
	\centering
	\captionsetup{font=footnotesize}
	\includegraphics[scale=0.08]{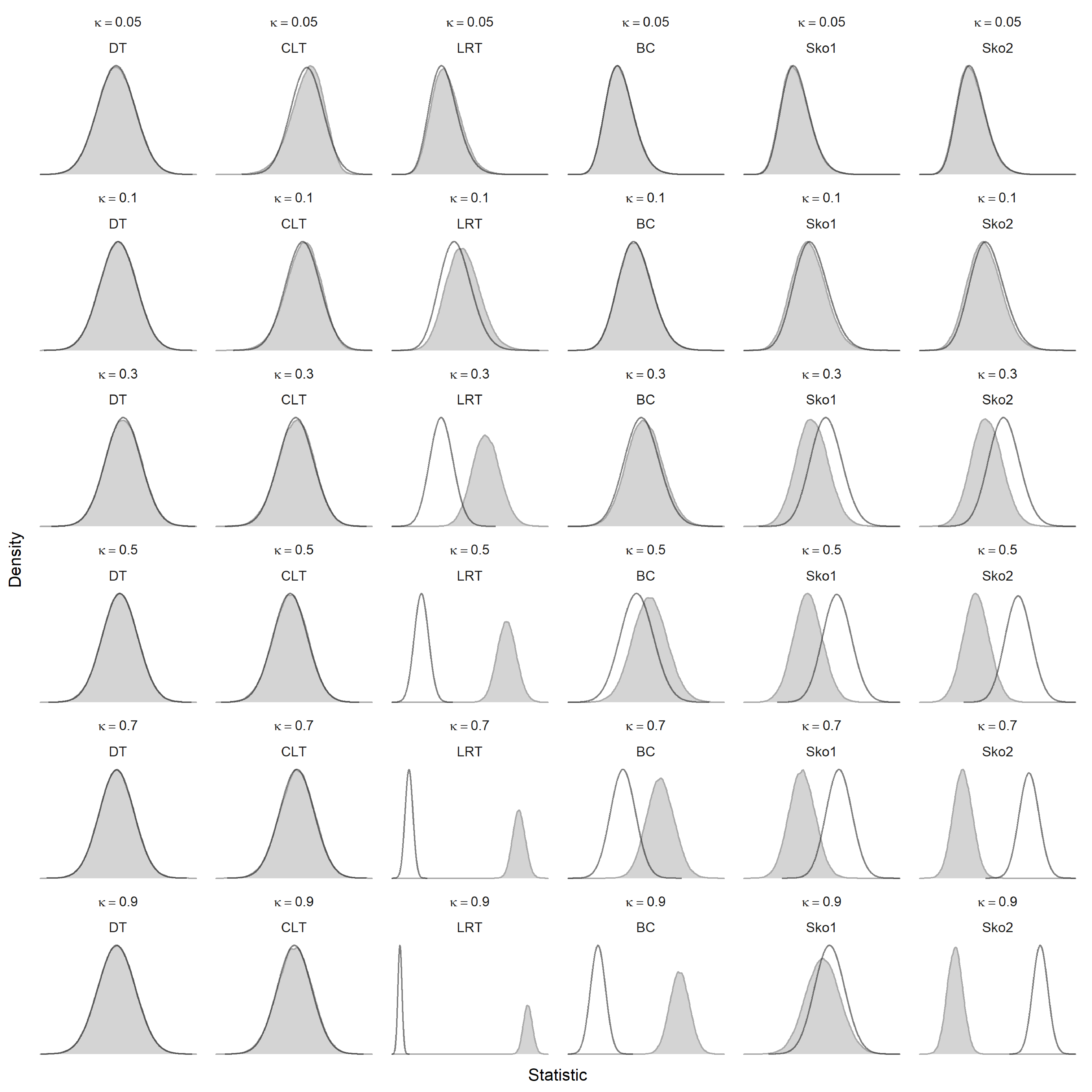}
	\caption{Hypothesis (II). Comparison between   empirical (gray) and theoretical (black) null  distributions of the directional test (DT), the central limit theorem test (CLT), log-likelihood ratio test (LRT), Bartlett correction (BC),  and two Skovgaard's modifications \cite{skovgaard:2001} (Sko1 and Sko2, respectively) for  various values of $ p/n_i$, $i \in \{1,\dots,k\}$.}
	\label{SMsimulation:density null case3}
\end{figure}

\begin{figure}[H]
	\centering
	\captionsetup{font=footnotesize}
	\subfigure{
		\begin{minipage}[b]{.35\linewidth}
			\centering
			\includegraphics[scale=0.085]{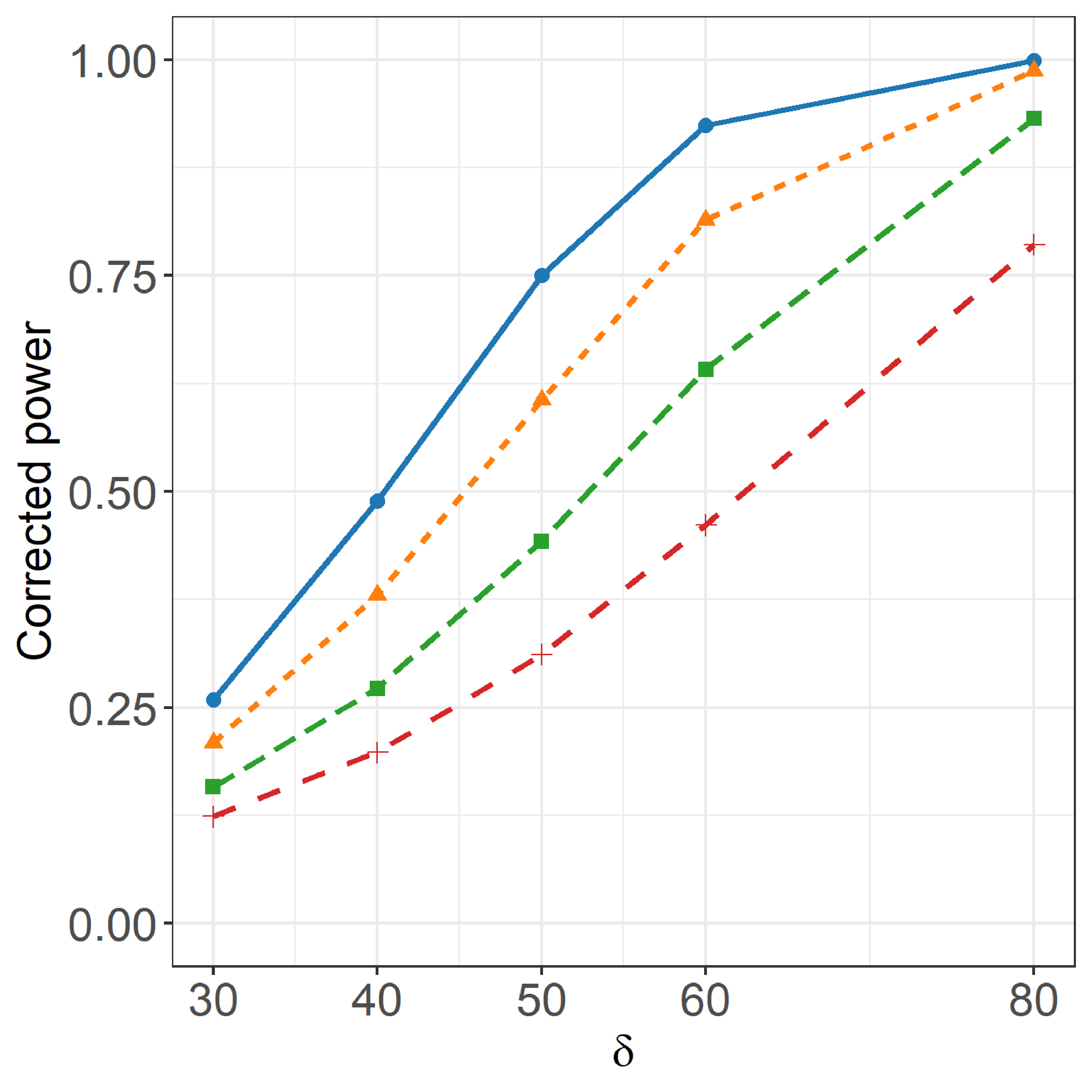}
		\end{minipage}
	}
	\subfigure{
		\begin{minipage}[b]{.35\linewidth}
			\centering
			\includegraphics[scale=0.085]{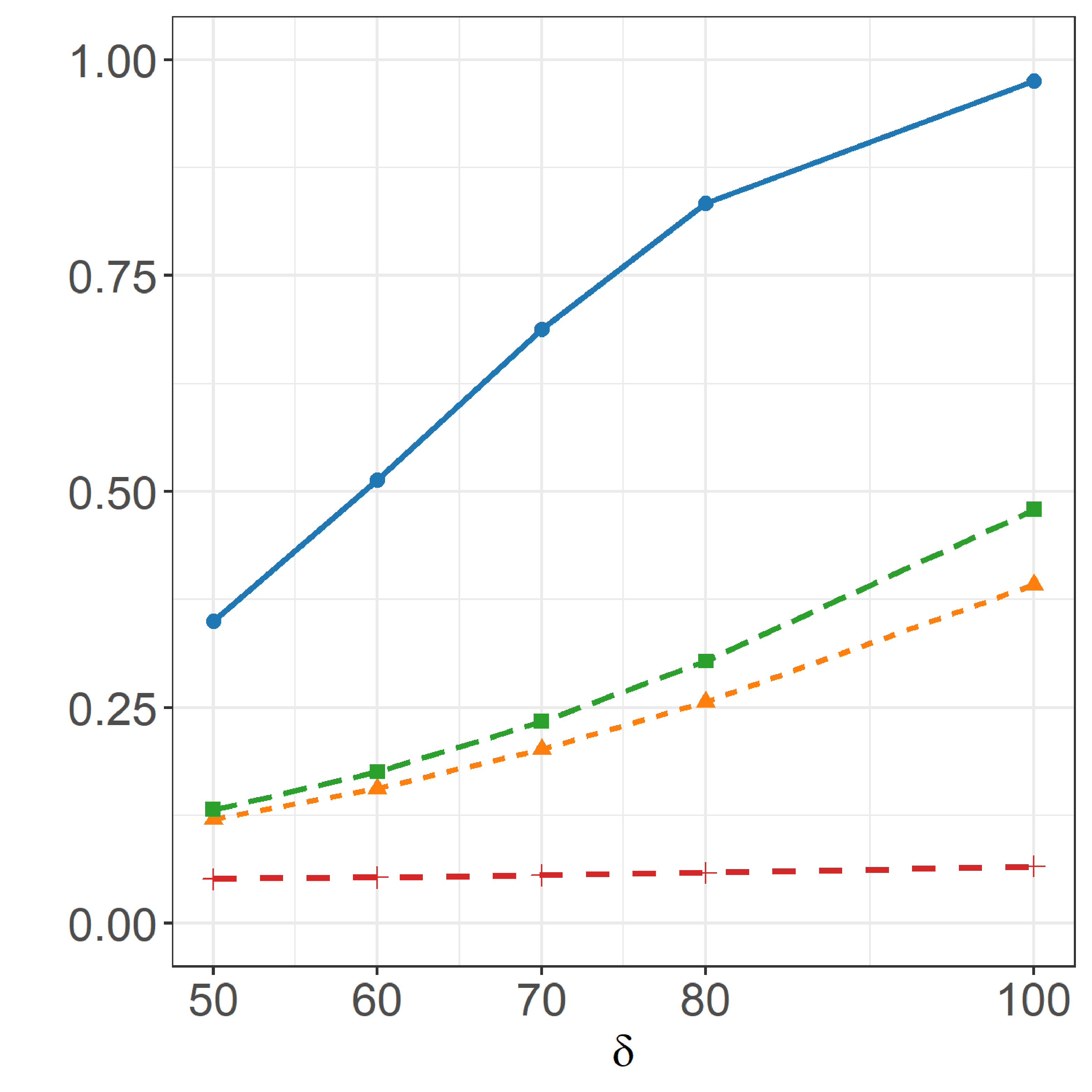}
		\end{minipage}
	}
	\\
	\subfigure{
		\begin{minipage}[b]{.35\linewidth}
			\centering
			\includegraphics[scale=0.085]{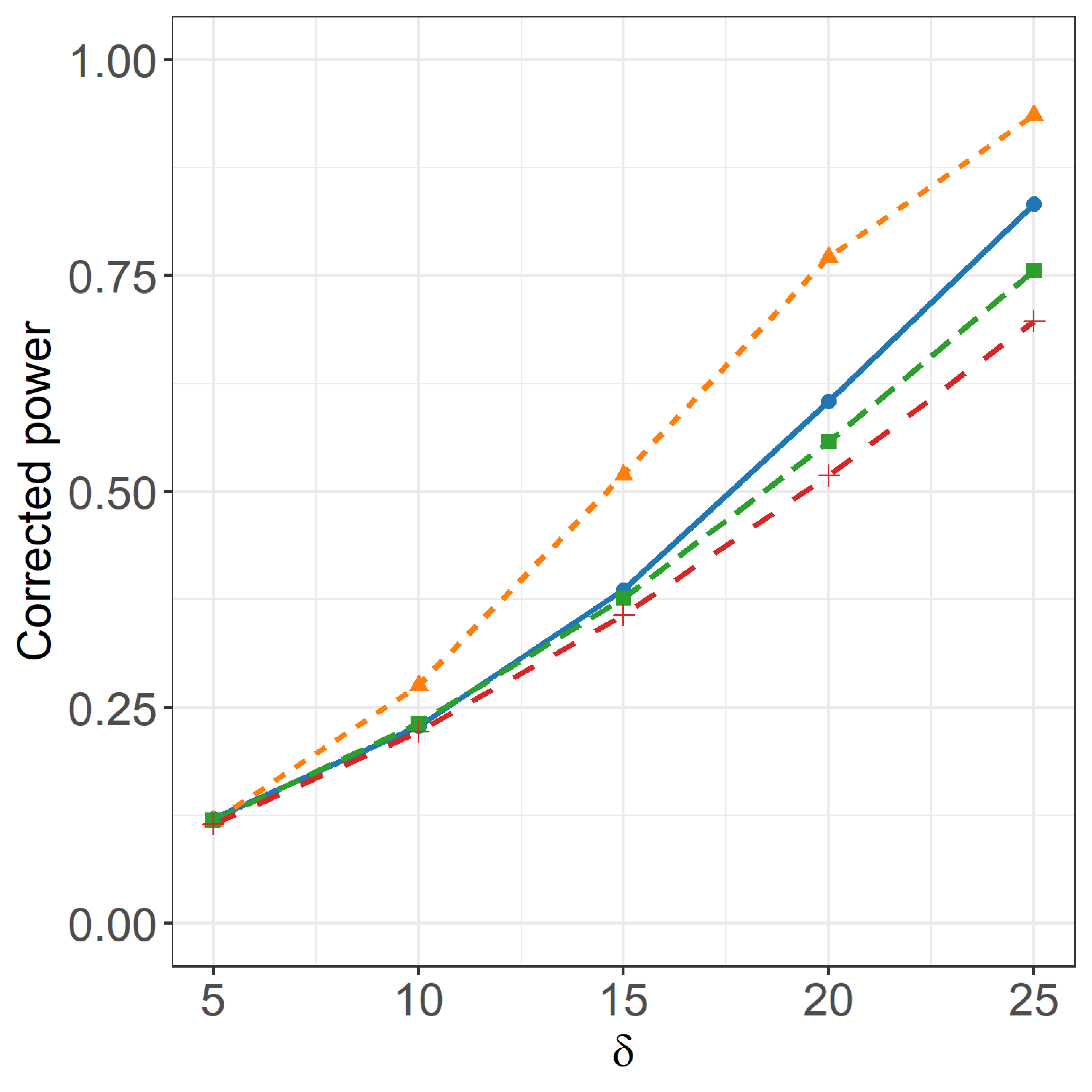}
		\end{minipage}
	}
	\subfigure{
		\begin{minipage}[b]{.35\linewidth}
			\centering
			\includegraphics[scale=0.085]{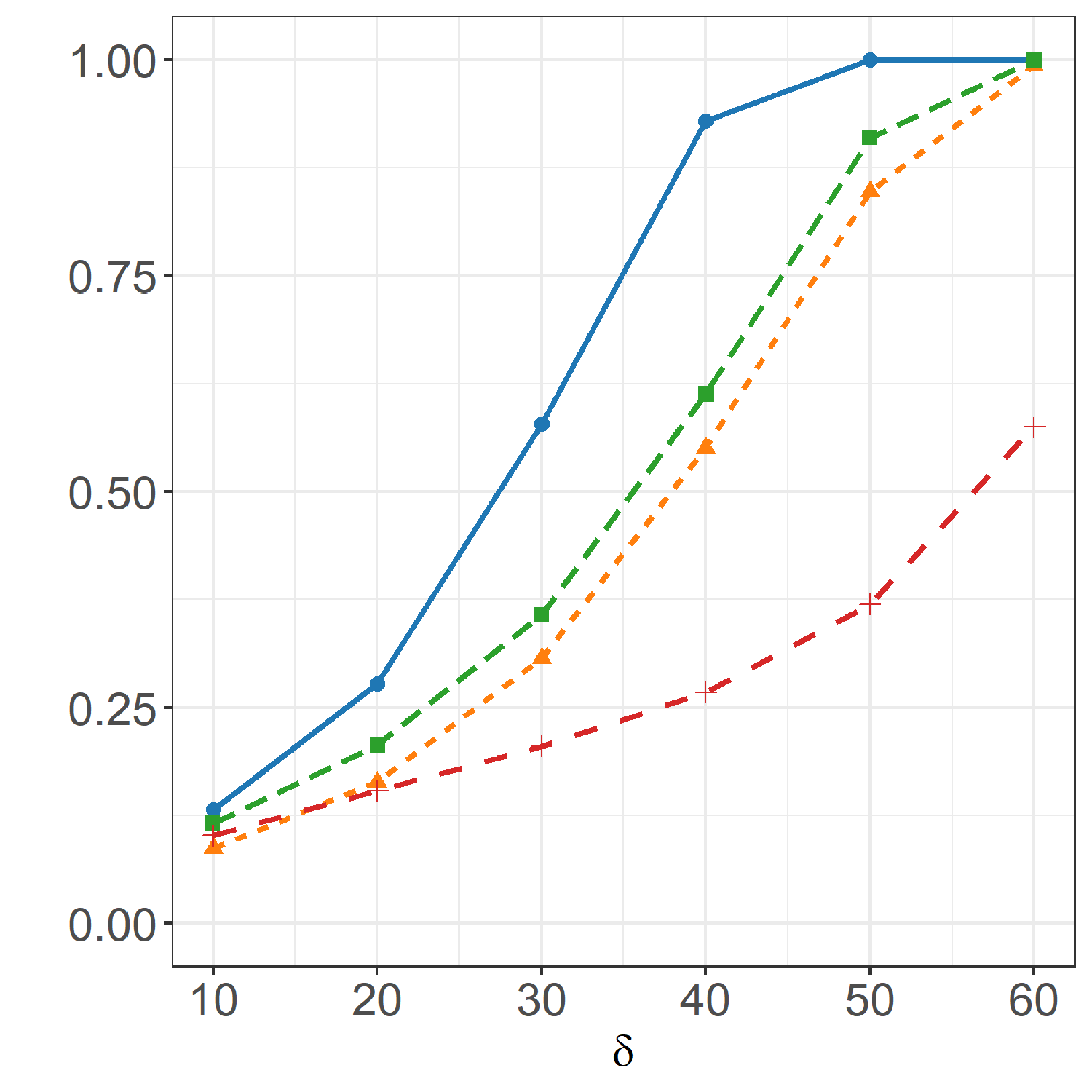}
		\end{minipage}
	}
	\caption{Empirical local corrected powers of four tests. The solid, dashed, longdashed, and dotdashed curves are the empirical power functions of the central limit theorem test, directional test and two Skovgaard's modifications \cite{skovgaard:2001}, respectively. The top and bottom rows correspond to alternative hypothesis settings (2) and (3) of hypothesis (I), respectively; the left and right columns correspond to the ratio $p/n_i = 0.7$ and $0.9$, respectively.}   
	\label{SMfig:power0.7 case4}
\end{figure}

\begin{figure}[H]
	\centering
	\captionsetup{font=footnotesize}
	\subfigure{
		\begin{minipage}[b]{.35\linewidth}
			\centering
			\includegraphics[scale=0.085]{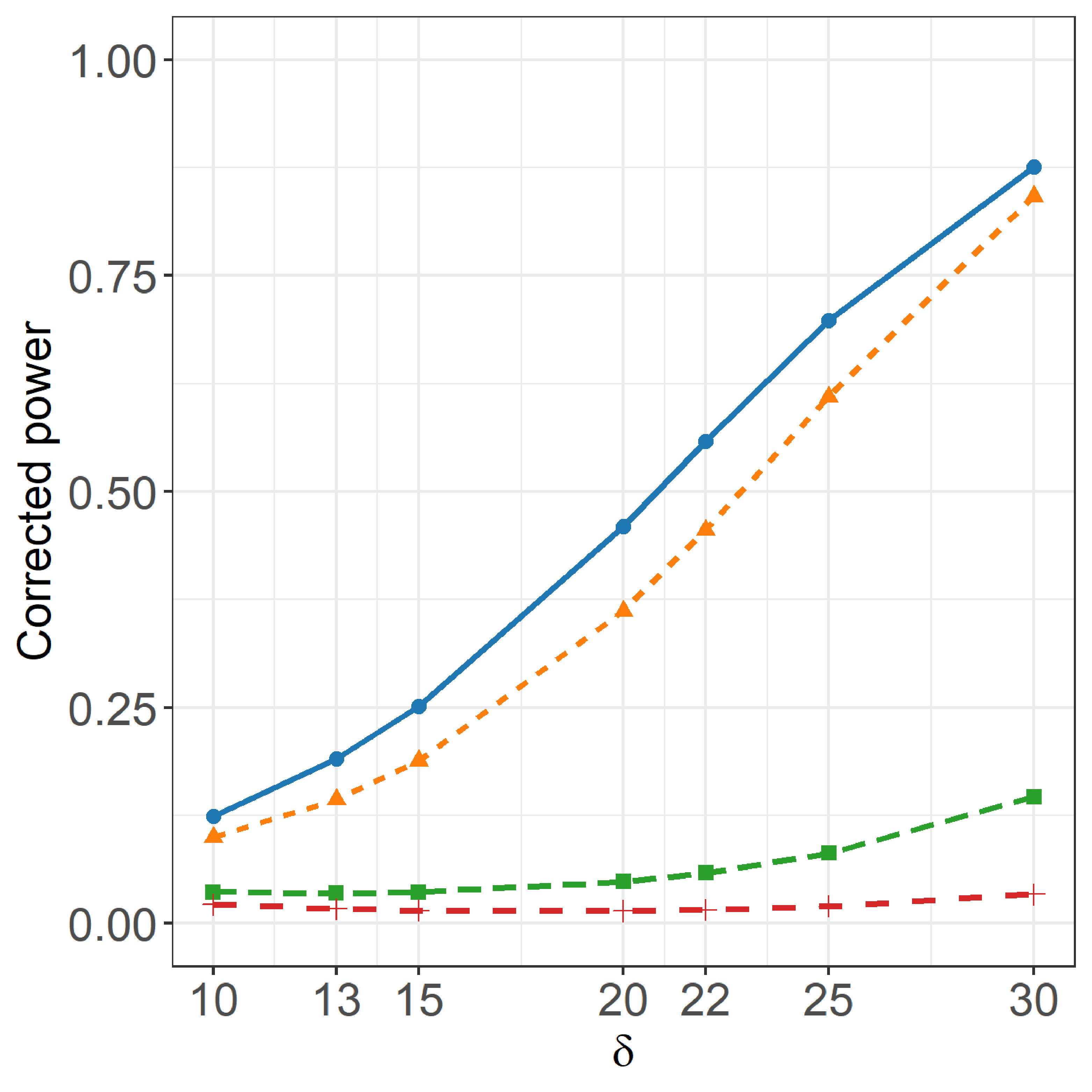}
		\end{minipage}
	}
	\subfigure{
		\begin{minipage}[b]{.35\linewidth}
			\centering
			\includegraphics[scale=0.085]{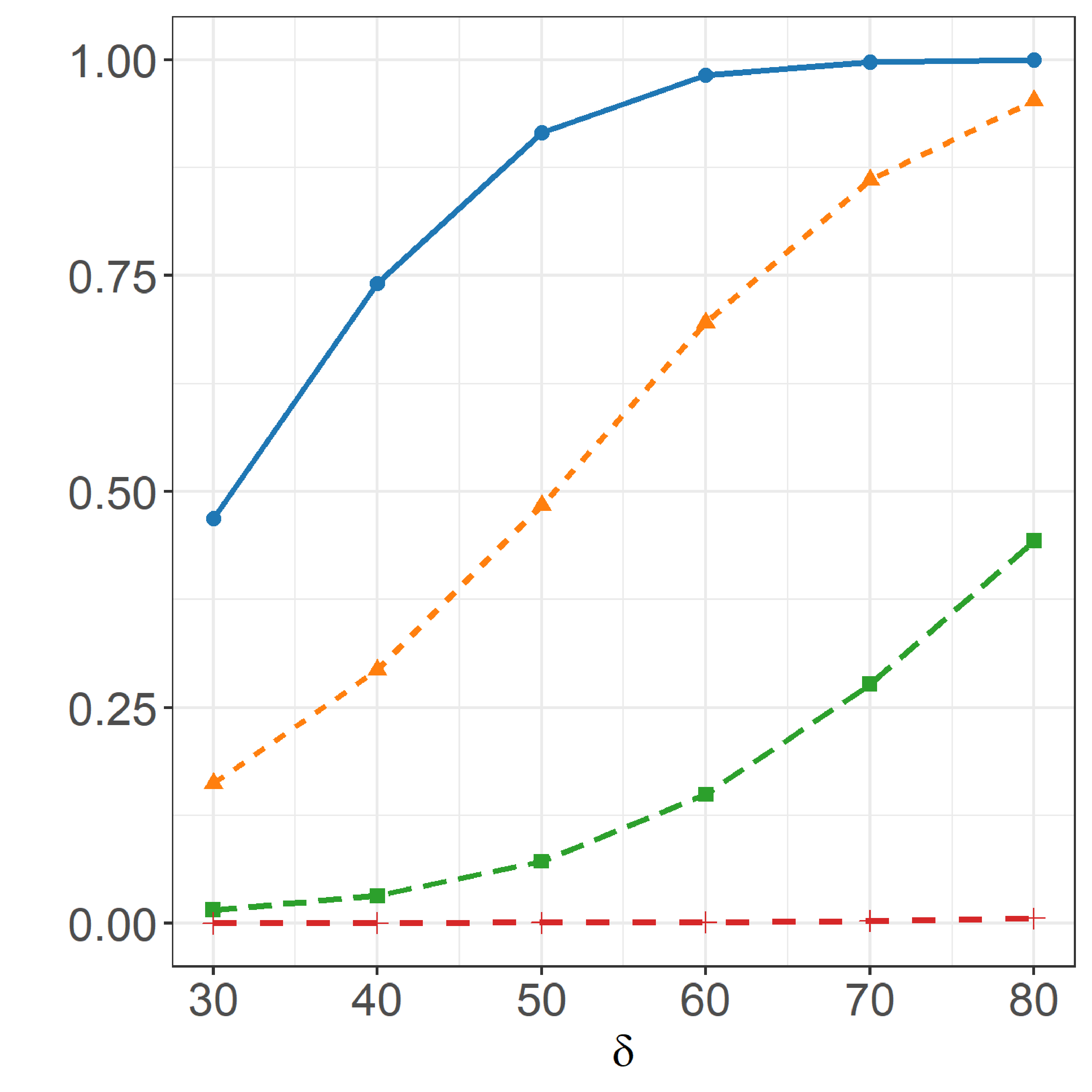}
		\end{minipage}
	}
	\\
	\subfigure{
		\begin{minipage}[b]{.35\linewidth}
			\centering
			\includegraphics[scale=0.085]{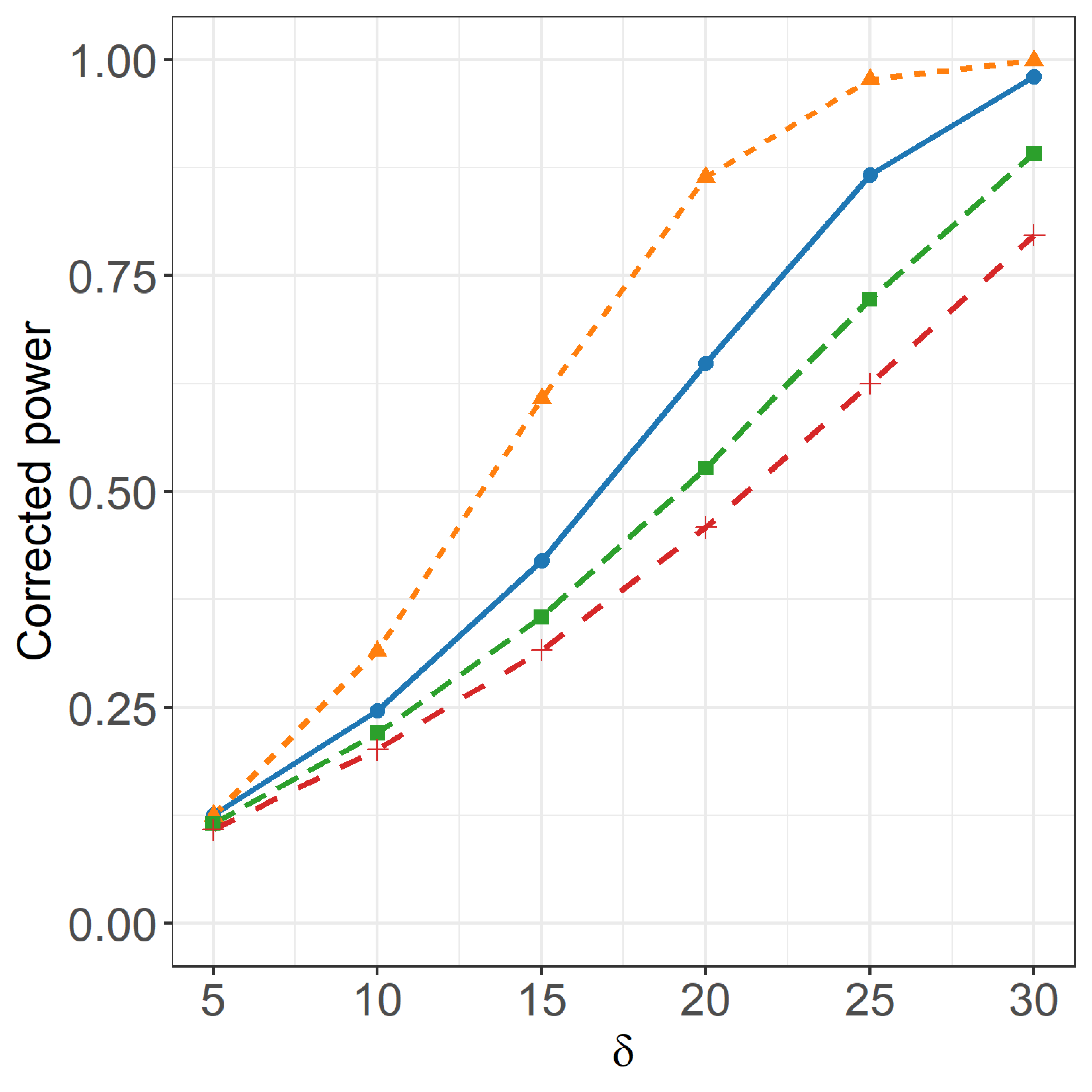}
		\end{minipage}
	}
	\subfigure{
		\begin{minipage}[b]{.35\linewidth}
			\centering
			\includegraphics[scale=0.085]{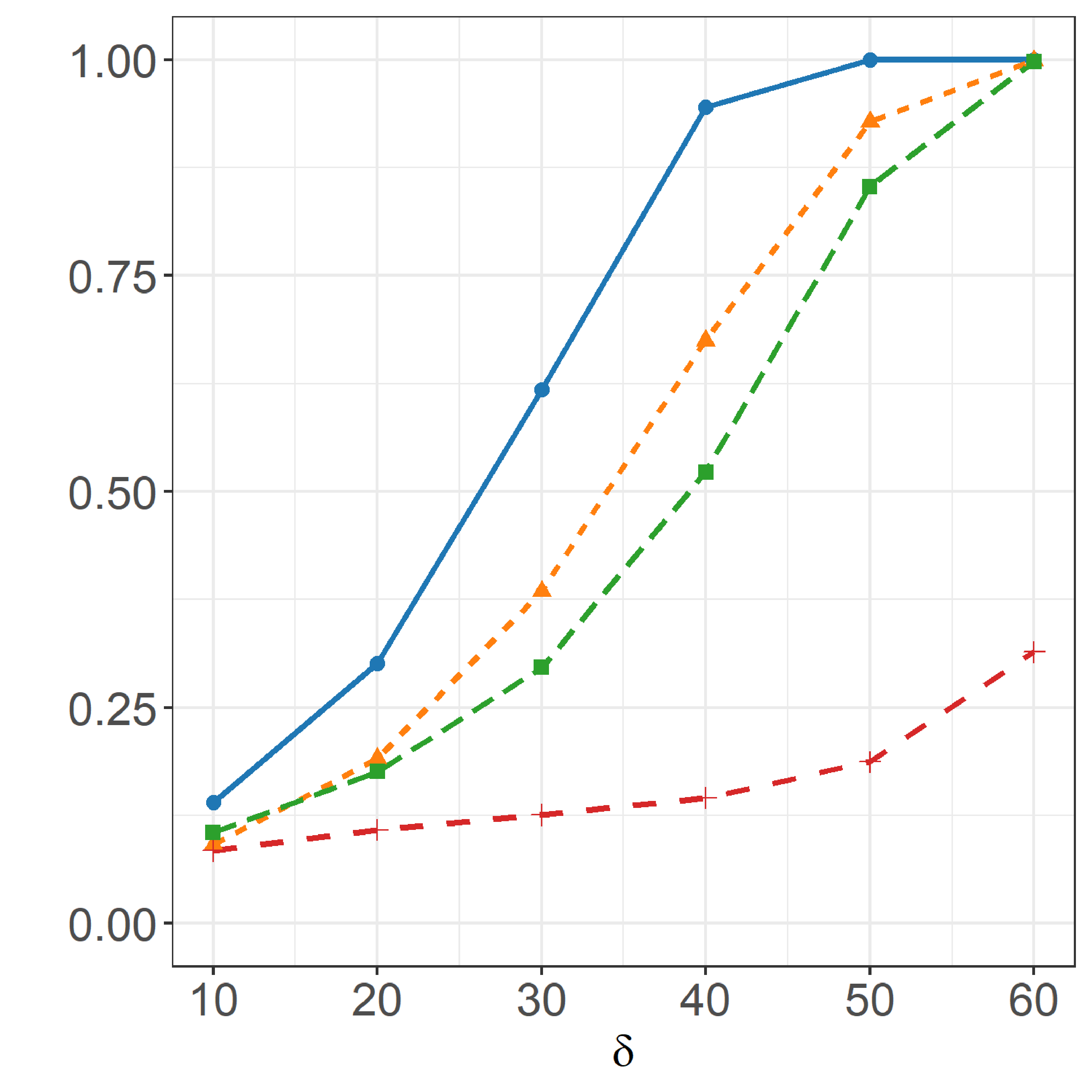}
		\end{minipage}
	}
	\caption{Empirical local corrected powers of four tests. The solid, dashed, longdashed, and dotdashed curves are the empirical power functions of the central limit theorem test, directional test and two Skovgaard's modifications \cite{skovgaard:2001}, respectively. The top and bottom rows correspond to alternative hypothesis settings (2) and (3) of hypothesis (II), respectively; the left and right columns correspond to the ratio $p/n_i = 0.7$ and $0.9$, respectively.}   
	\label{SMfig:power0.7 case3}
\end{figure}

Table \ref{SMsimulation: table type I normal:cases} displays the corrected Type I error,  which is the $5\%$-quantile of the empirical $p$-values obtained under the null hypothesis.
Figures \ref{SMsimulation:pvalue case4}--\ref{SMsimulation:pvalue case3}  confirm the theoretical findings on the directional $p$-value: it is exact, to Monte Carlo accuracy, in every simulation setting.  The  chi-square approximation for the distribution  of $W$, $W_{BC}$,  $W^*$ and $W^{**}$  are sufficiently reliable only when $p$ is small. At the same time, the slight location bias of the normal approximation to the central limit theorem test distribution is apparent. On the other hand, the behavior of $W$, $W_{BC}$, $W^*$ and $W^{**}$  become less accurate  as $p$ gets moderate or large relative to $n$, while the central limit theorem test improves as $p$ increases, as suggested by its theoretical derivation. 
The similar behavior for the six statistics is also in terms of the simulated null distribution (see Figures \ref{SMsimulation:density null case4}--\ref{SMsimulation:density null case3}).   The directional test still outperforms the competitors in each simulation setting. 
Figures \ref{SMfig:power0.7 case4}--\ref{SMfig:power0.7 case3} report  the local corrected power examined in the alternative hypothesis settings (2) and (3).  
With the ratio $p/n_i$ increasing, the central limit theorem is more powerful than that of others. 
But for the lower ratio, the power of directional test in the alternative setting outperforms. 

\begin{figure}[H]
	\centering
	\captionsetup{font=footnotesize}
	\includegraphics[scale=0.08]{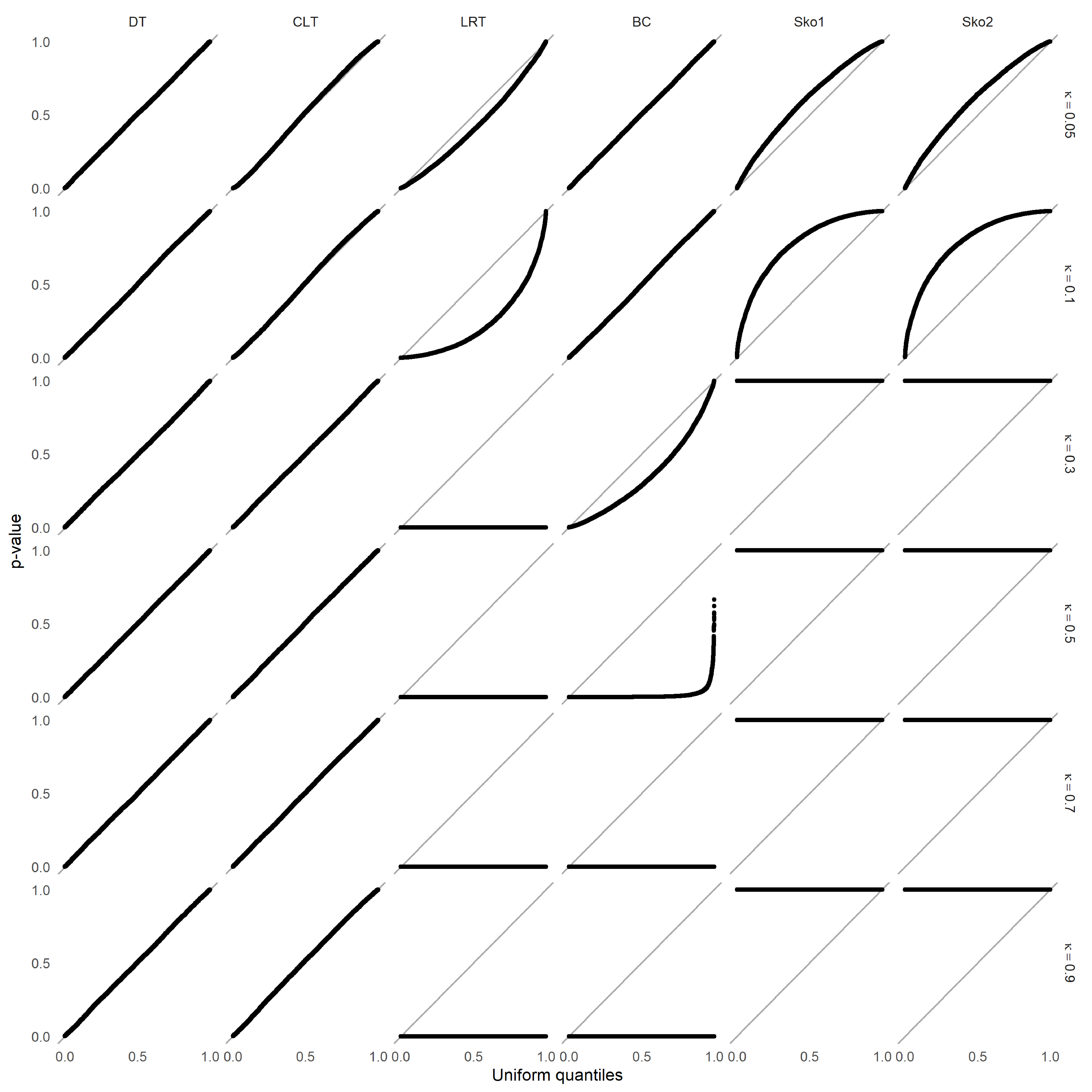}
	\caption{Hypothesis (I). Empirical null  distribution of $p$-values for various ratio $p/n_i$, $i \in \{1,\dots,k\}$, of the directional test (DT), the central limit theorem test (CLT), log-likelihood ratio test (LRT), Bartlett correction (BC),  and two Skovgaard's modifications \cite{skovgaard:2001} (Sko1 and Sko2, respectively) with $k=30$, compared
		with the $U(0,1)$ given by the gray diagonal.}
	\label{SMsimulation:pvalue large k30 case4}
\end{figure}
\begin{figure}[H]
	\centering
	\captionsetup{font=footnotesize}
	\includegraphics[scale=0.08]{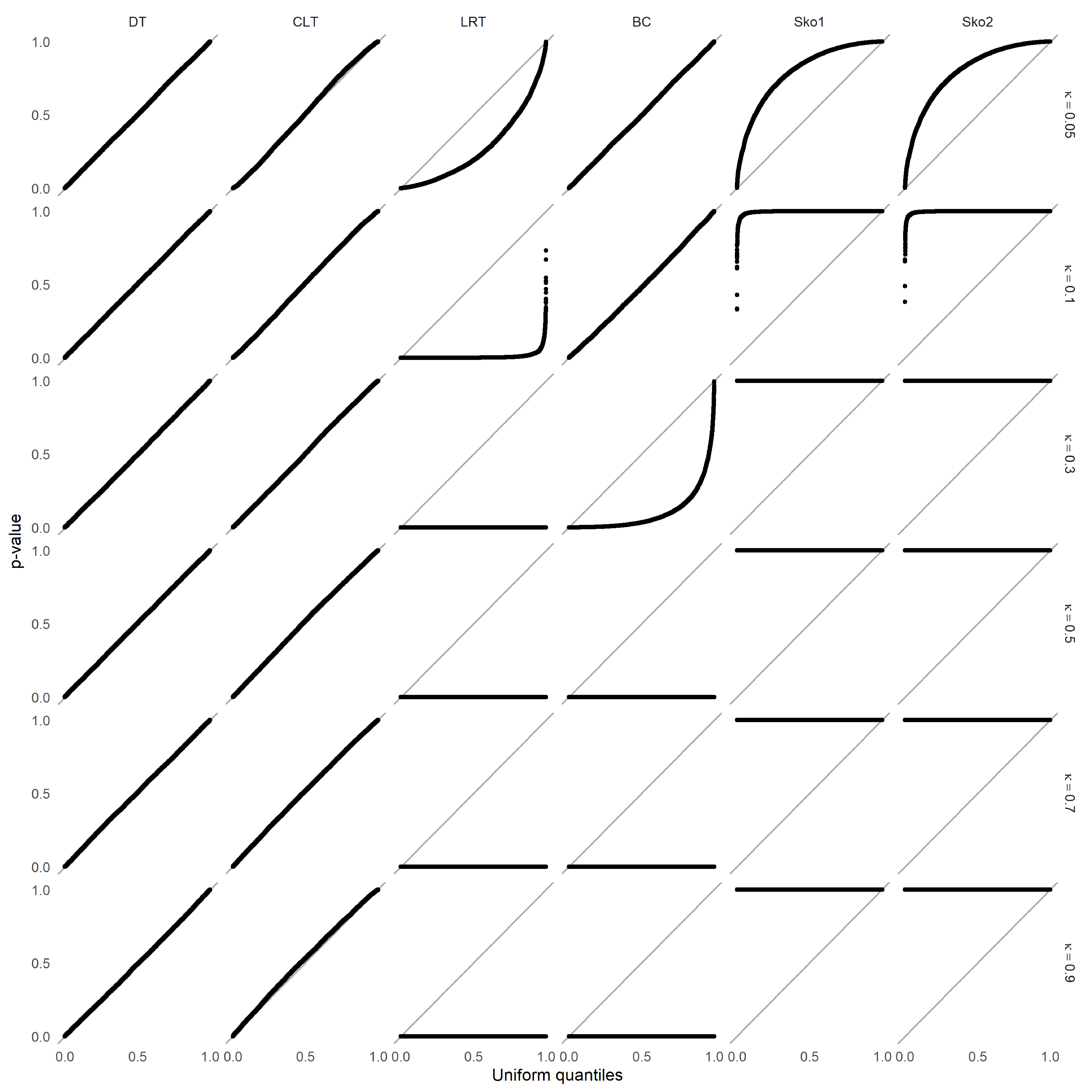}
	\caption{Hypothesis (I). Empirical null  distribution of $p$-values for various ratio $p/n_i$, $i \in \{1,\dots,k\}$, of the directional test (DT), the central limit theorem test (CLT), log-likelihood ratio test (LRT), Bartlett correction (BC),  and two Skovgaard's modifications \cite{skovgaard:2001} (Sko1 and Sko2, respectively) with $k=300$, compared
		with the $U(0,1)$ given by the gray diagonal.}
	\label{SMsimulation:pvalue large k300 case4}
\end{figure}
\begin{figure}[H]
	\centering
	\captionsetup{font=footnotesize}
	\includegraphics[scale=0.08]{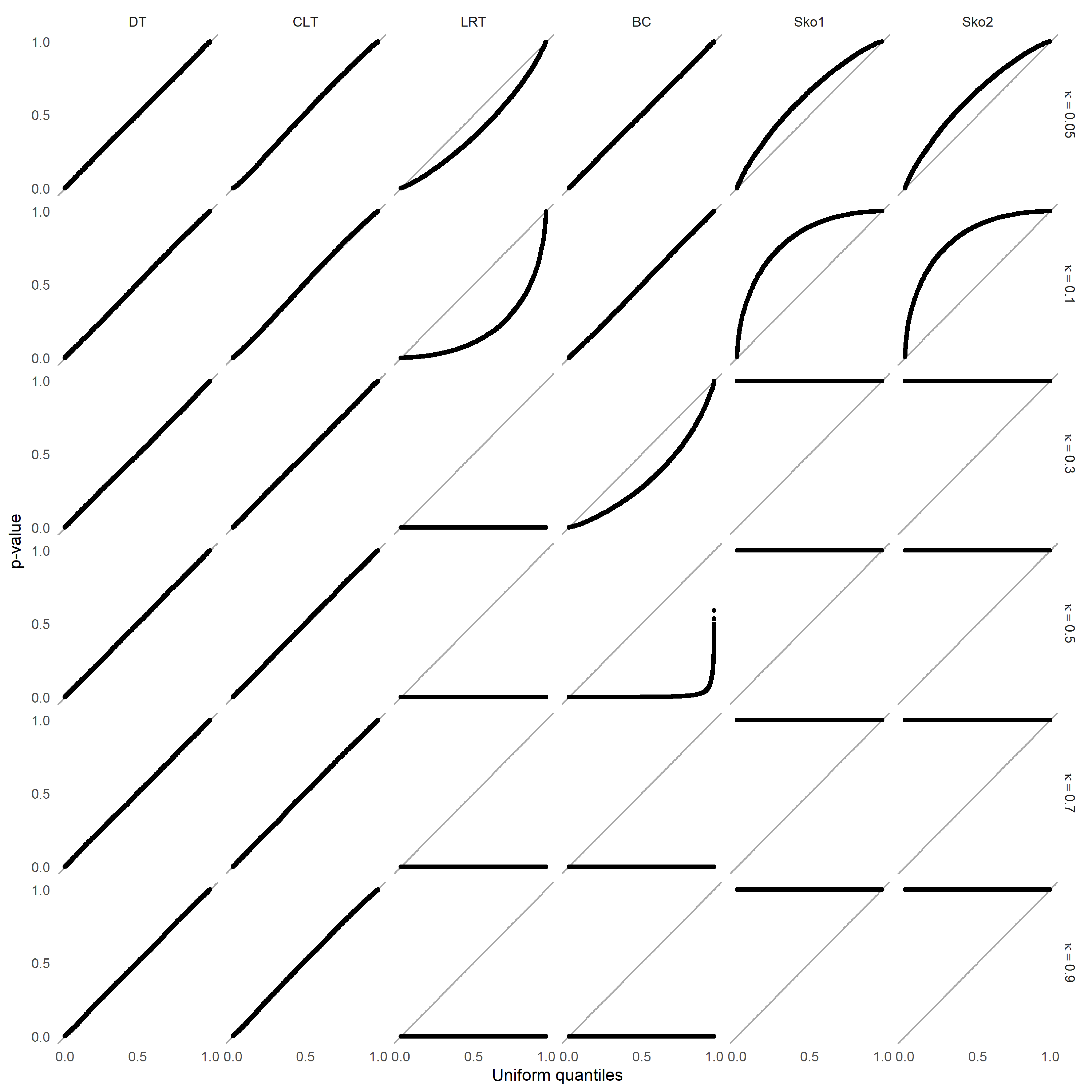}
	\caption{Hypothesis (II). Empirical null  distribution of $p$-values for various ratio $p/n_i$, $i \in \{1,\dots,k\}$, of the directional test (DT), the central limit theorem test (CLT), log-likelihood ratio test (LRT), Bartlett correction (BC),  and two Skovgaard's modifications \cite{skovgaard:2001} (Sko1 and Sko2, respectively) with $k=30$, compared
		with the $U(0,1)$ given by the gray diagonal.}
	\label{SMsimulation:pvalue large k30 case3}
\end{figure}
\begin{figure}[H]
	\centering
	\captionsetup{font=footnotesize}
	\includegraphics[scale=0.08]{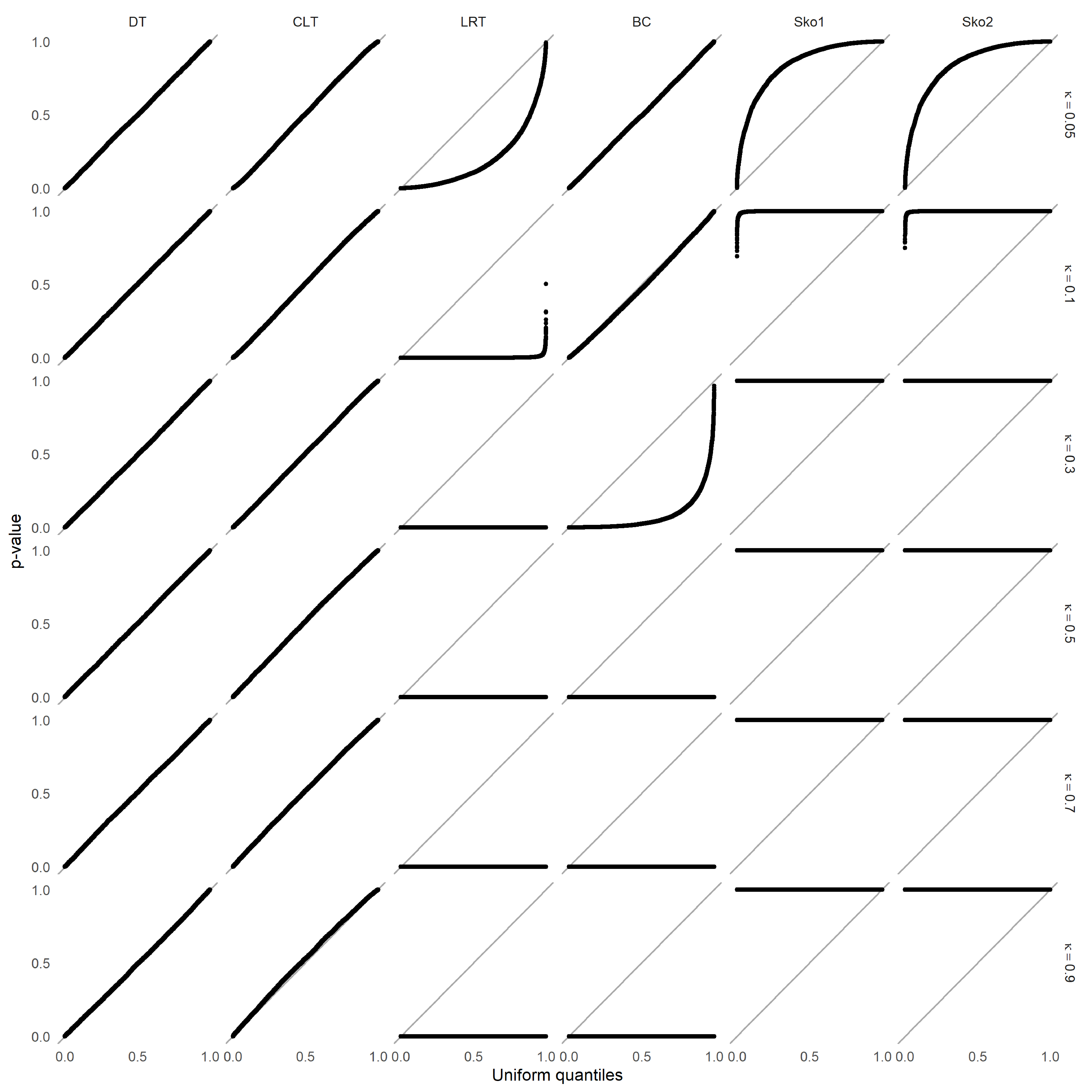}
	\caption{Hypothesis (II). Empirical null  distribution of $p$-values for various ratio $p/n_i$, $i \in \{1,\dots,k\}$, of the directional test (DT), the central limit theorem test (CLT), log-likelihood ratio test (LRT), Bartlett correction (BC),  and two Skovgaard's modifications \cite{skovgaard:2001} (Sko1 and Sko2, respectively) with $k=300$, compared
		with the $U(0,1)$ given by the gray diagonal.}
	\label{SMsimulation:pvalue large k300 case3}
\end{figure}

We also investigate the empirical distribution of $p$-values  for larger values of the number of groups $k \in \{30, 300\}$ in hypotheses (I) and (II).  The numerical results are based on 10,000 replications. Figures \ref{SMsimulation:pvalue large k30 case4}--\ref{SMsimulation:pvalue large k300 case3} show that the directional $p$-value maintains extreme accuracy, apart from simulation errors, while the main competitor proposed by Jiang and Yang \cite{jiang:2013} becomes slightly less accurate as $k$ increases \citep[see also][]{dette2020,guo2021asymptotic}.


\section{Numerical results in one-sample hypotheses (III)--(VI)}\label{S4:simulation (III)-(VI)}

This section investigates the performance of the hypotheses which are testing  the sphericity of the covariance matrix,  block-independence,  complete-independence and specific multivariate normal distribution.
The limiting null distribution, empirical  $p$-value,  corrected Type I error and power are evaluated based on $100,000$ replications. 
We set $n=100$ and $p/n=0.05,0.1,0.3,0.5,0.7,0.9$. Random samples of size $n$  are generated from the standard normal distribution $N_p(0_p,\bI_p)$.  The various simulation setups are details below.

Hypothesis (III):  testing the sphericity of the covariance matrix, namely $H_{\psi}:\Lambda^{-1} = \sigma^2_0\bI_p$, with unspecified $\sigma^2_0$. We set $\sigma^2_0=1$ under the null hypothesis and focus on  three different forms of $\Lambda^{-1}$ under the alternative: (1) $\Lambda^{-1} = \diag(1.69,\dots,1.69,1,\dots,1)$, where the number of diagonal entries equal to $1.69$ is  $\left\lceil {p/2} \right\rceil$; (2)  $\Lambda^{-1} = \diag ({{1}}_p + \delta n^{-1/2} \bu )$ with $\bu = \{ (2/p)^{1/2},\dots,$ $(2/p)^{1/2},0,\dots,0\}^T$ such that $||\bu||=1$,  where the first $\left\lceil {p/2} \right\rceil$ diagonal elements are equal to $(2/p)^{1/2}$; (3) $\Lambda^{-1} = \diag \{(1 + \eta) \; {1}_{p-1}^T, 1\}$, where $\eta \in \mathbb{R}^+$. 

Hypothesis (IV):  testing the block-independence of the normal components as in (\ref{hypothesis:case2}). Under the null hypothesis we set $k=3$,  so that the normal random vector is divided into three sub-vectors with respective dimensions $p_1$, $p_2$ and $p_3$ satisfying $p_1:p_2:p_3 = 2:2:1$. Different alternative settings are considered: (1) compound symmetric covariance structure: $\Lambda^{-1} = 0.15 {{1}}_p{{1}}_p^T + 0.85 \bI_p$; (2) $\Lambda^{-1} = \eta  {{1}}_p{{1}}_p^T + (1-\eta) \bI_p$, where $\eta = \delta \{p(p-1)n\}^{-1/2}$ and  $\delta$ is chosen so that $\eta \in (0,1)$; (3)  $\Lambda^{-1} =(\lambda^{ij})_{p\times p}$  where $\lambda^{ij}=1$ for $i=j$, $\lambda^{1(p_1+1)}=\lambda^{(p_1+1)1}=\eta$ with $\eta \in (0,1)$, and $\lambda^{ij}=0$, otherwise. 

Hypothesis (V): testing complete-independence as in (\ref{hypothesis:case6}). (1) $\Lambda^{-1} = \left(\lambda^{ij}\right)_{p\times p}$ where $\lambda^{ij} = 1$ for $i=j$, $\lambda^{ij} = 0.1$ for $0<|i-j|\le 3$, and $\lambda^{ij} = 0$ for $|i-j|>3$; (2) $\Lambda^{-1} = \left(\lambda^{ij}\right)_{p\times p}$ where $\lambda^{ij} = 1$ for $i=j$, $\lambda^{ij} = \delta n^{-1/2} \bu $ for $0<|i-j|\le 3$ and $\bu=1/\surd{\# \{0 < |i - j|\le 3\}}$ where $\# \{\cdot\}$ counts the number of occurrences, and $\lambda^{ij} = 0$ for $|i-j|>3$; (3) $\Lambda^{-1} = \left(\lambda^{ij}\right)_{p\times p}$ where $\lambda^{ij} = 1$ for $i=j$, $\lambda^{12}= \lambda^{21}= \eta$ with $\eta \in (0,1)$, and $\lambda^{ij} = 0$ otherwise.

Hypothesis (VI): testing specified values for the mean vector and covariance matrix. Three alternative setups are: (1)  $\mu = (0.1,\dots,0.1,0,\dots,0)^T$, where the number of entries equal to $0.1$ is $\left\lceil {p/2} \right\rceil $; (2) $\mu = \{\delta (2/pn)^{1/2} ,\dots,\delta (2/pn)^{1/2},$ $0,\dots,0\}^T$ where the number of entries equal to $\delta (2/pn)^{1/2}$ is  $\left\lceil {p/2} \right\rceil $. The (1) and (2) alternative setups for $\Lambda^{-1}$ are as  in Hypothesis (V); (3) $\mu$ is the same as (1)  and $\Lambda^{-1} = \diag( 1-\eta,1,\dots,1)$ with $\eta \in (0,1)$.

\begin{figure}[H]
	\centering
	\captionsetup{font=footnotesize}
	\includegraphics[scale=0.08]{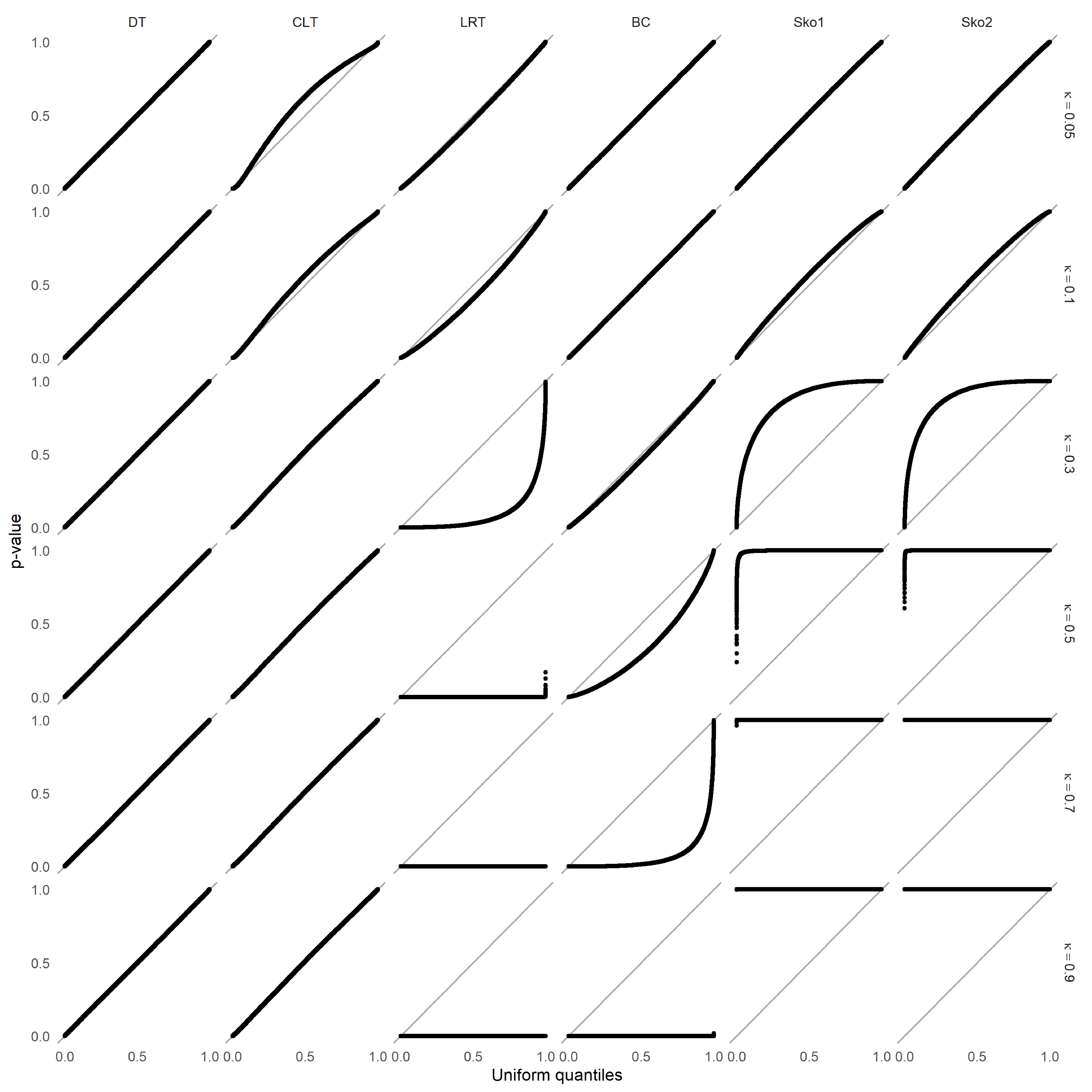}
	\caption{Hypothesis (III).  Empirical null  distribution of $p$-values for various ratios $p/n$  of the directional test (DT), the central limit theorem test (CLT), log-likelihood ratio test (LRT), Bartlett correction (BC),  and two Skovgaard's modifications \cite{skovgaard:2001} (Sko1 and Sko2, respectively),  compared
		with the $U(0,1)$ given by the gray diagonal.}
	\label{SMsimulation:pvalue case1}
\end{figure}
\begin{figure}[H]
	\centering
	\captionsetup{font=footnotesize}
	\includegraphics[scale=0.08]{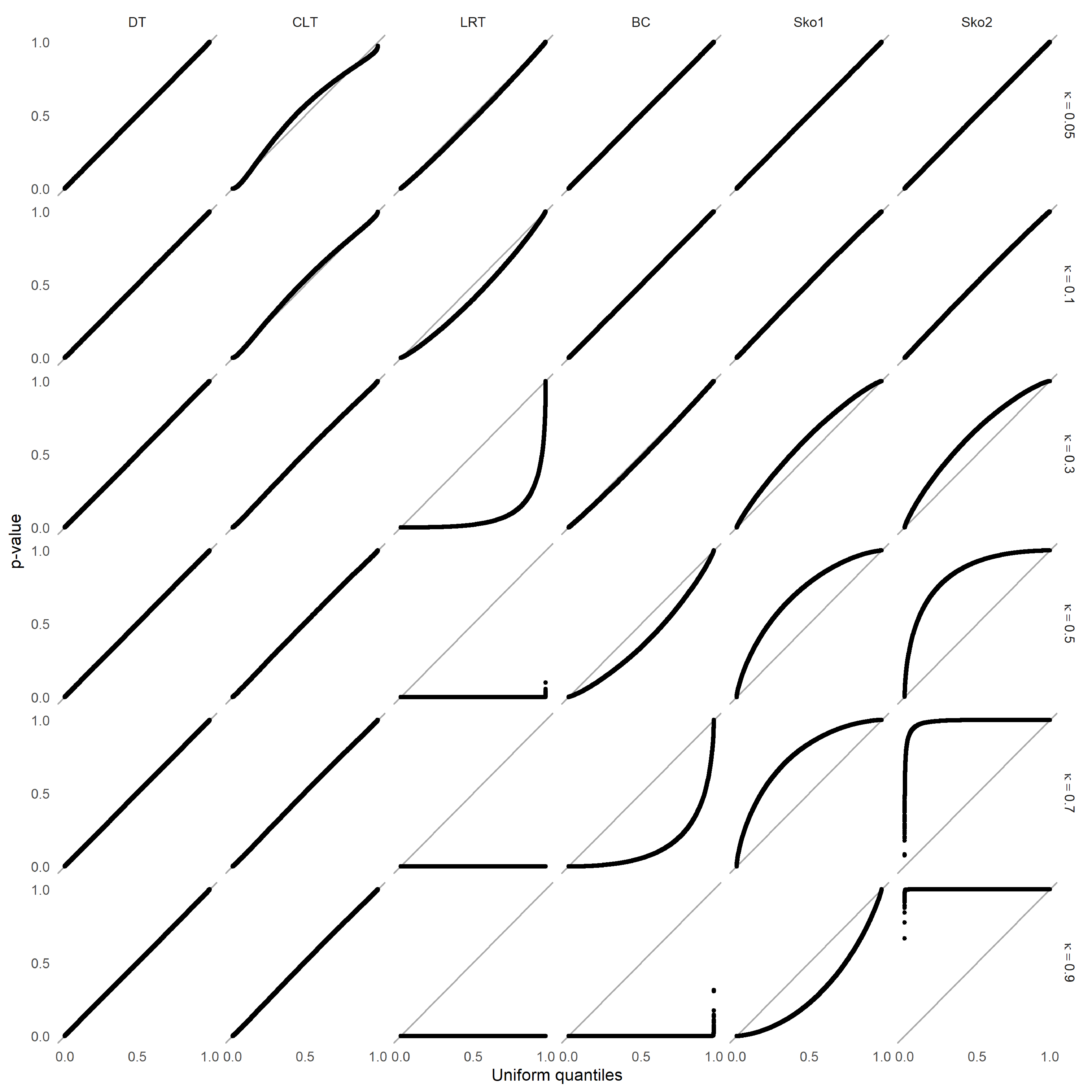}
	\caption{Hypothesis (IV). Empirical null  distribution of $p$-values for various ratio $p/n$  of the directional test (DT), the central limit theorem test (CLT), log-likelihood ratio test (LRT), Bartlett correction (BC),  and two Skovgaard's modifications \cite{skovgaard:2001} (Sko1 and Sko2, respectively),  compared
		with the $U(0,1)$ given by the gray diagonal.}
	\label{SMsimulation:pvalue case2}
\end{figure}
\begin{figure}[H]
	\centering
	\captionsetup{font=footnotesize}
	\includegraphics[scale=0.08]{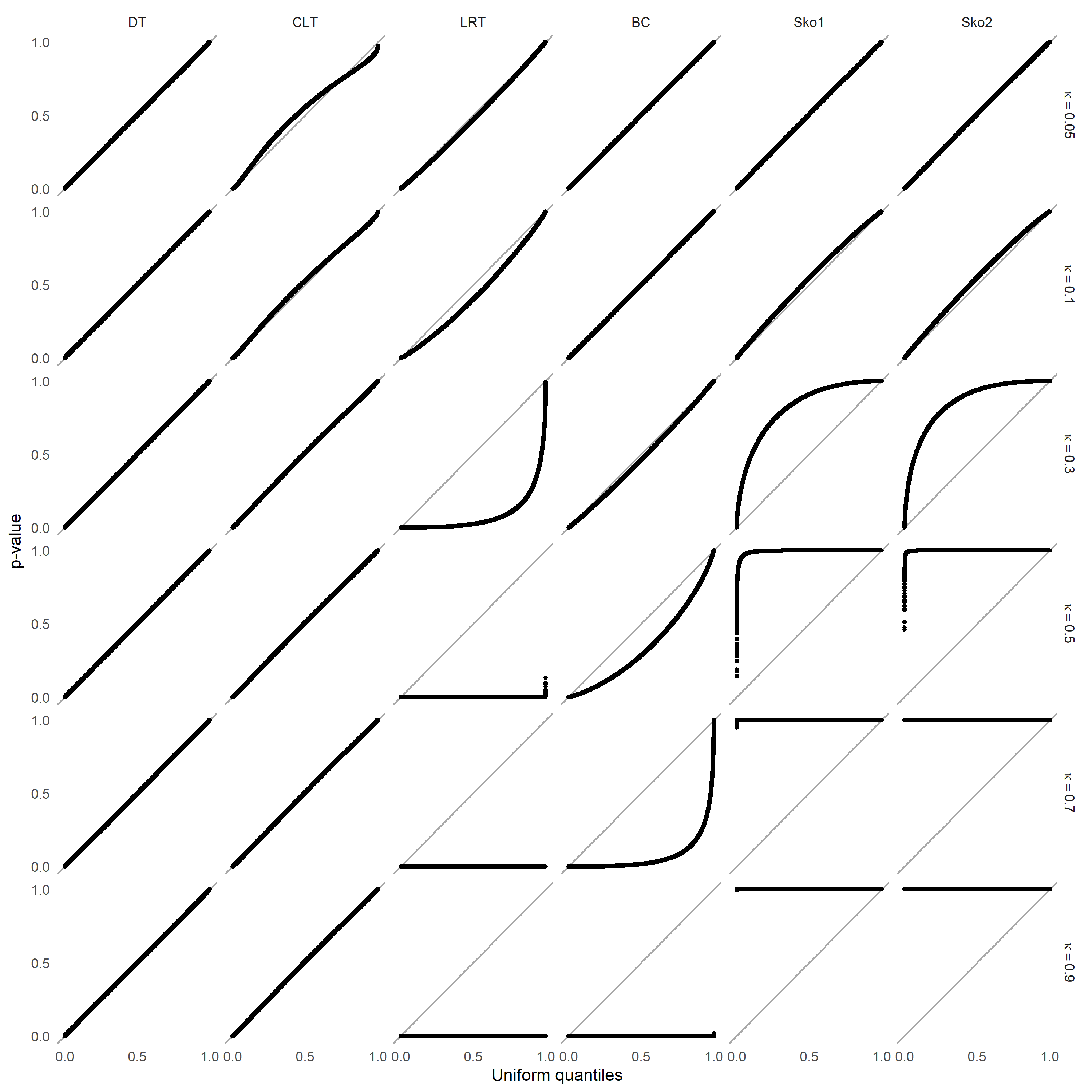}
	\caption{Hypothesis (V). Empirical null  distribution of $p$-values for various ratio $p/n$ of the directional test (DT), the central limit theorem test (CLT), log-likelihood ratio test (LRT), Bartlett correction (BC),  and two Skovgaard's modifications \cite{skovgaard:2001} (Sko1 and Sko2, respectively), compared
		with the $U(0,1)$ given by the gray diagonal.}
	\label{SMsimulation:pvalue case6}
\end{figure}
\begin{figure}[H]
	\centering
	\captionsetup{font=footnotesize}
	\includegraphics[scale=0.08]{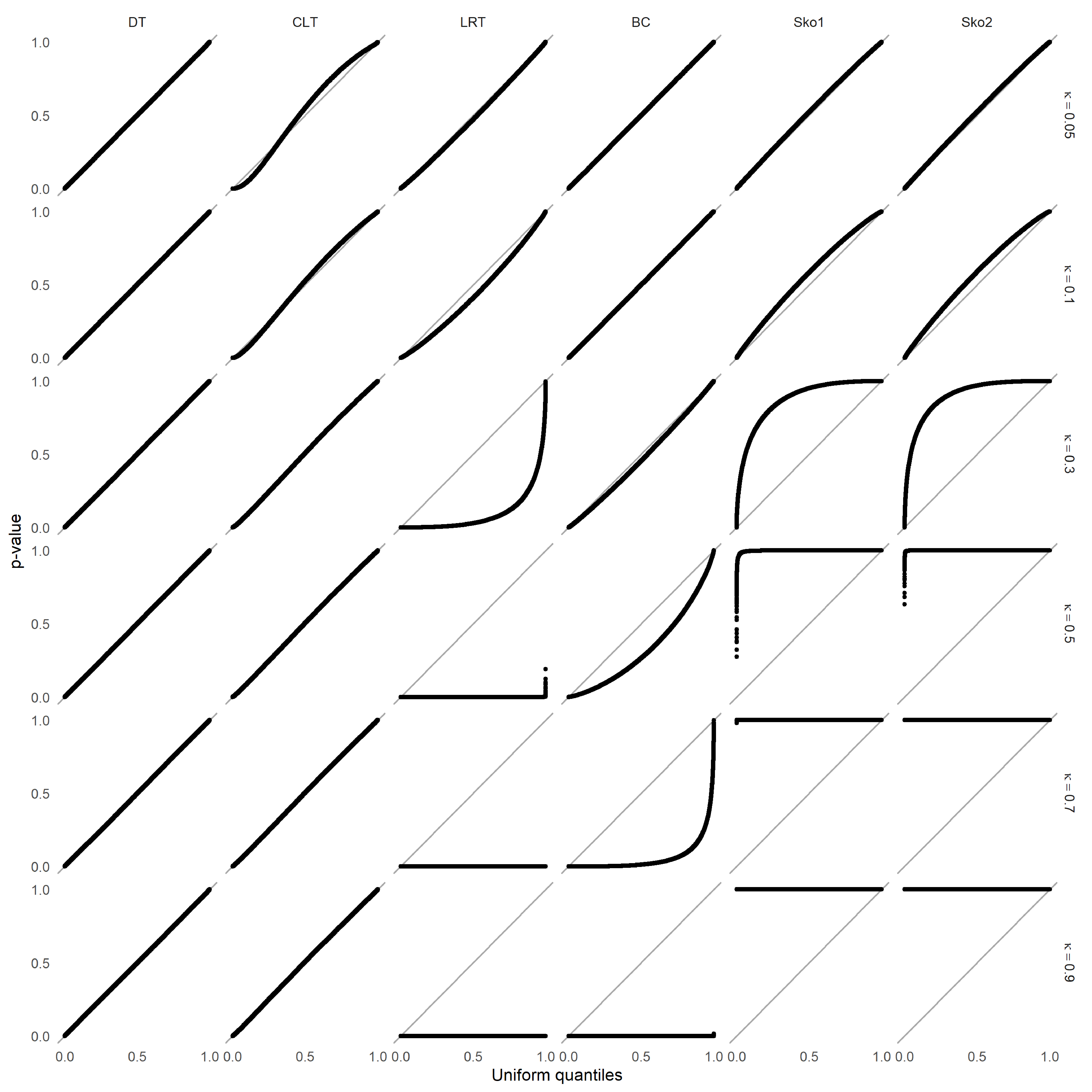}
	\caption{Hypothesis (VI). Empirical null  distribution of $p$-values for various ratio $p/n$ of the directional test (DT), the central limit theorem test (CLT), log-likelihood ratio test (LRT), Bartlett correction (BC),  and two Skovgaard's modifications \cite{skovgaard:2001} (Sko1 and Sko2, respectively), compared
		with the $U(0,1)$ given by the gray diagonal.}
	\label{SMsimulation:pvalue case5}
\end{figure}
\begin{figure}[H]
	\centering
	\captionsetup{font=footnotesize}
	\includegraphics[scale=0.08]{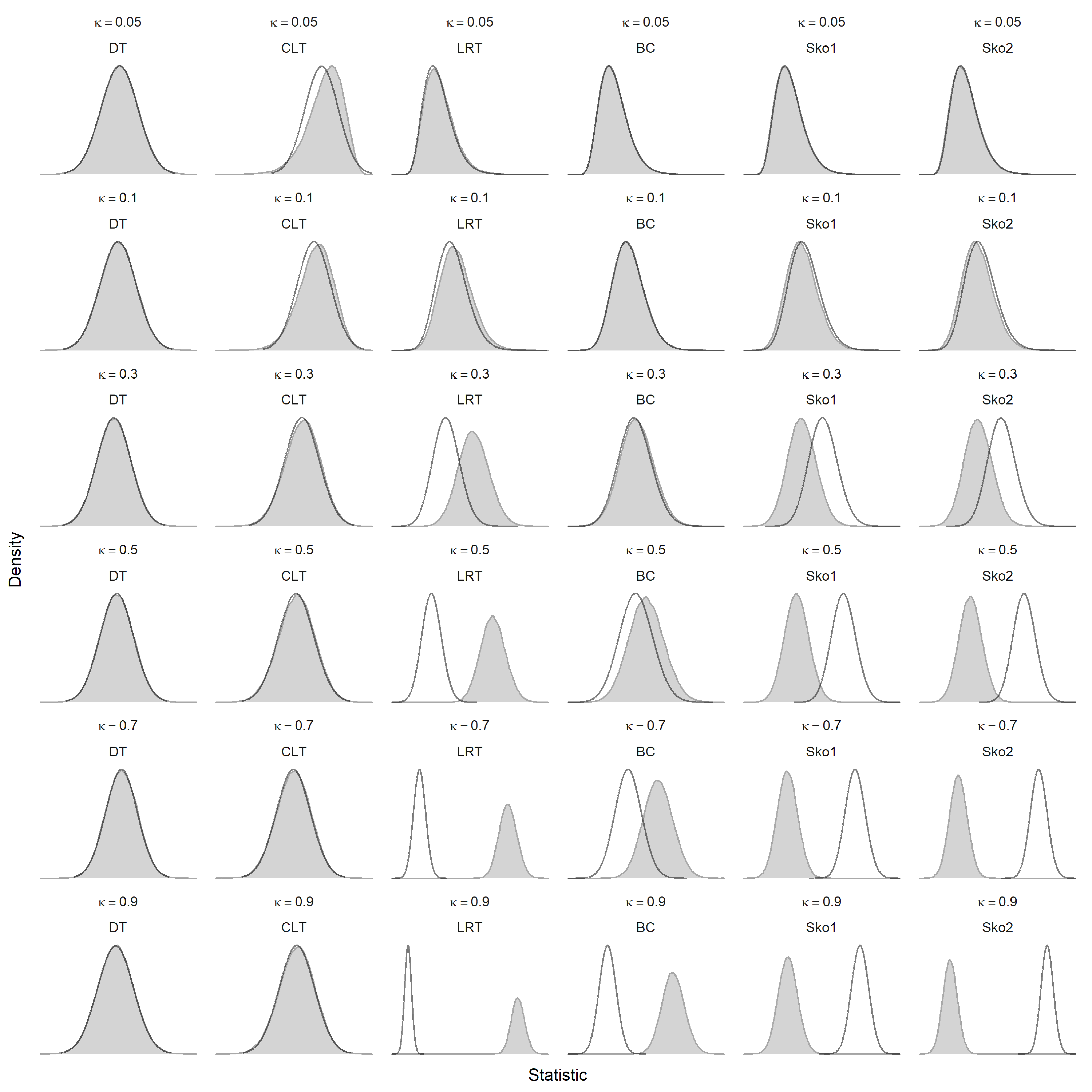}
	\caption{Hypothesis (III). Comparison between   empirical (gray) and theoretical (black) null  distributions of the directional test (DT), the central limit theorem test (CLT), log-likelihood ratio test (LRT), Bartlett correction (BC),  and two Skovgaard's modifications \cite{skovgaard:2001} (Sko1 and Sko2, respectively),  for  various values of $p/n$.}
	\label{SMsimulation:density null case1}
\end{figure}
\begin{figure}[H]
	\centering
	\captionsetup{font=footnotesize}
	\includegraphics[scale=0.08]{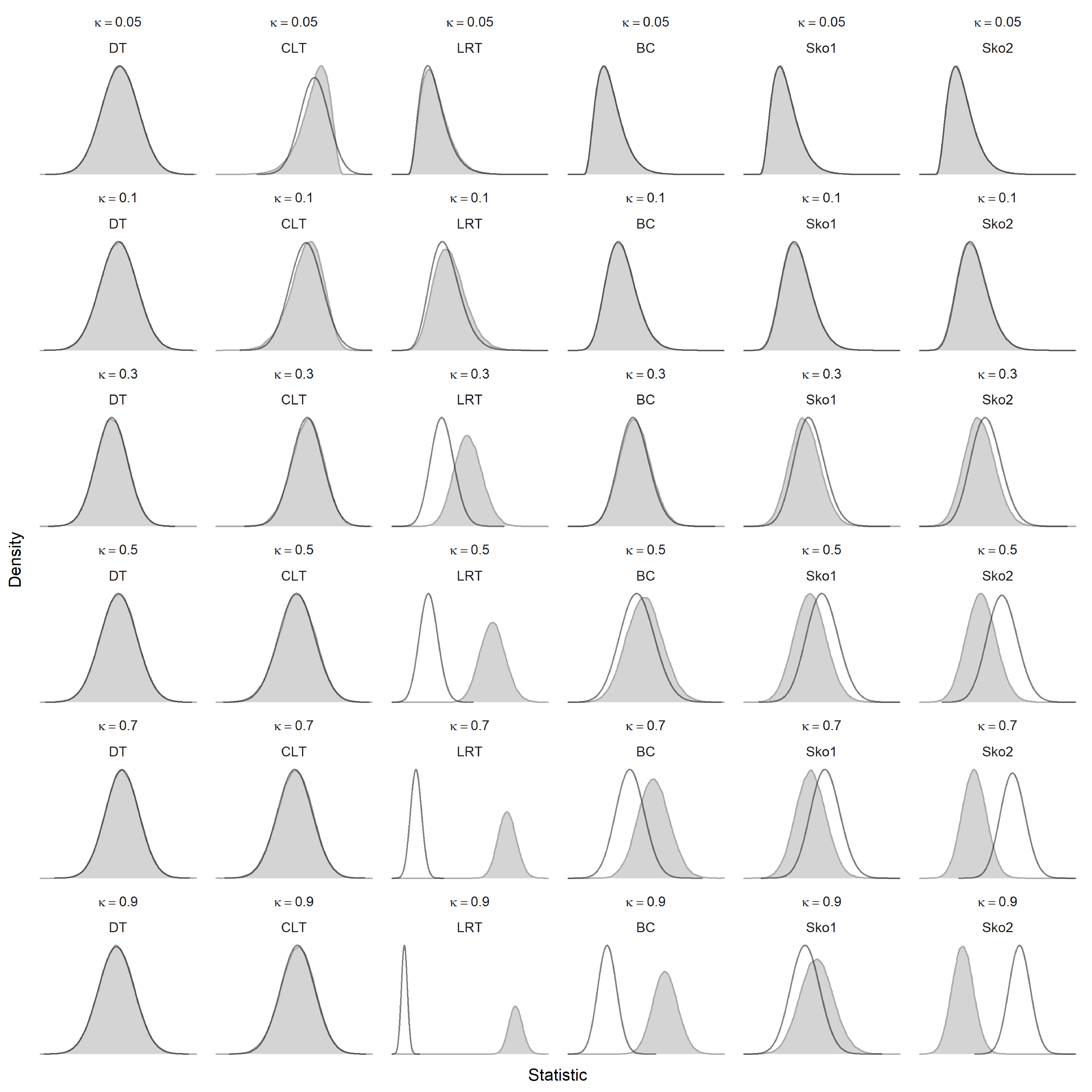}
	\caption{Hypothesis (IV). Comparison between   empirical (gray) and theoretical (black) null  distributions of the directional test (DT), the central limit theorem test (CLT), log-likelihood ratio test (LRT), Bartlett correction (BC),  and two Skovgaard's modifications \cite{skovgaard:2001} (Sko1 and Sko2, respectively),  for  various values of $p/n$.}
	\label{SMsimulation:density null case2}
\end{figure}
\begin{figure}[H]
	\centering
	\captionsetup{font=footnotesize}
	\includegraphics[scale=0.08]{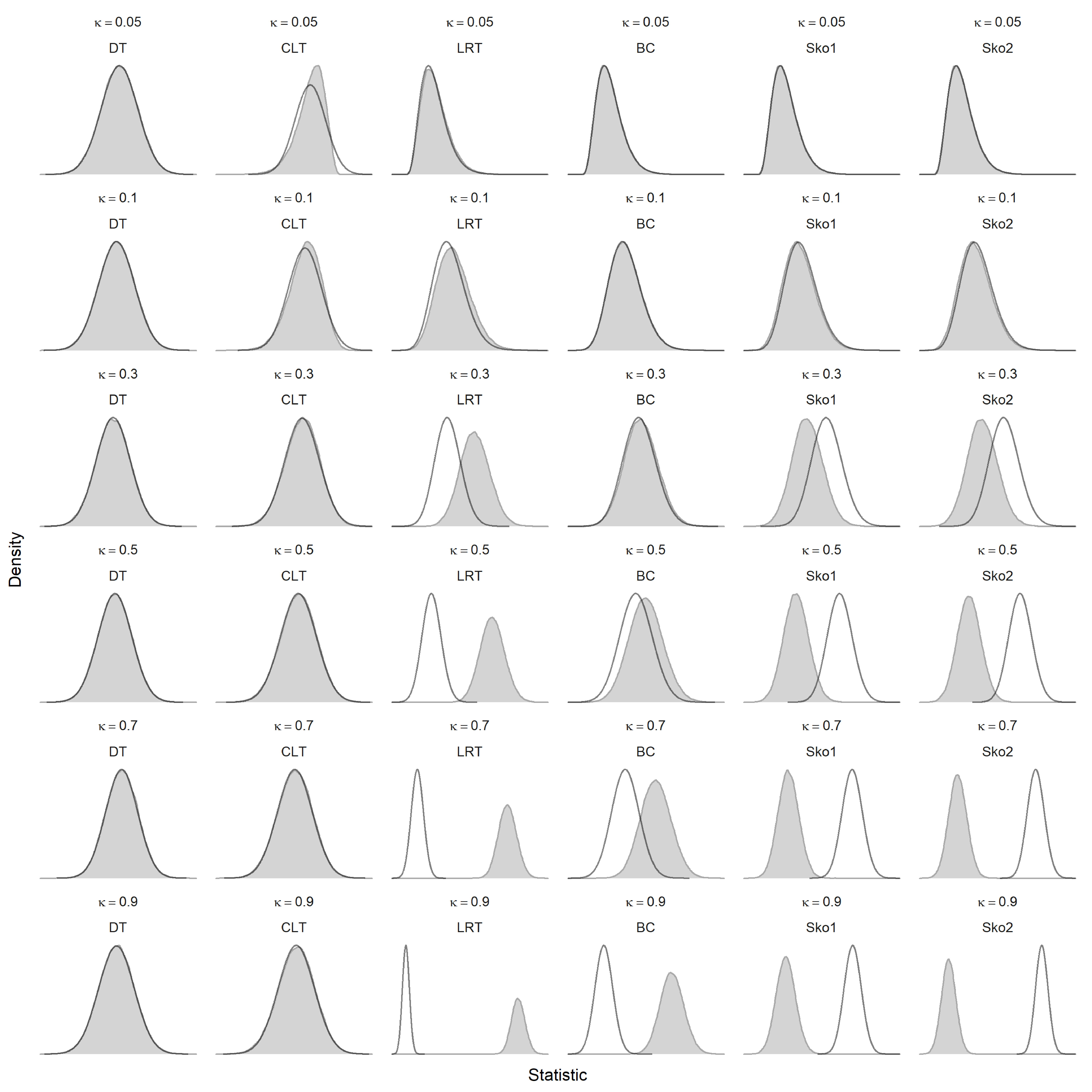}
	\caption{Hypothesis (V). Comparison between   empirical (gray) and theoretical (black) null  distributions of the directional test (DT), the central limit theorem test (CLT), log-likelihood ratio test (LRT), Bartlett correction (BC),  and two Skovgaard's modifications \cite{skovgaard:2001} (Sko1 and Sko2, respectively),  for  various values of $p/n$.}
	\label{SMsimulation:density null case6}
\end{figure}
\begin{figure}[H]
	\centering
	\captionsetup{font=footnotesize}
	\includegraphics[scale=0.08]{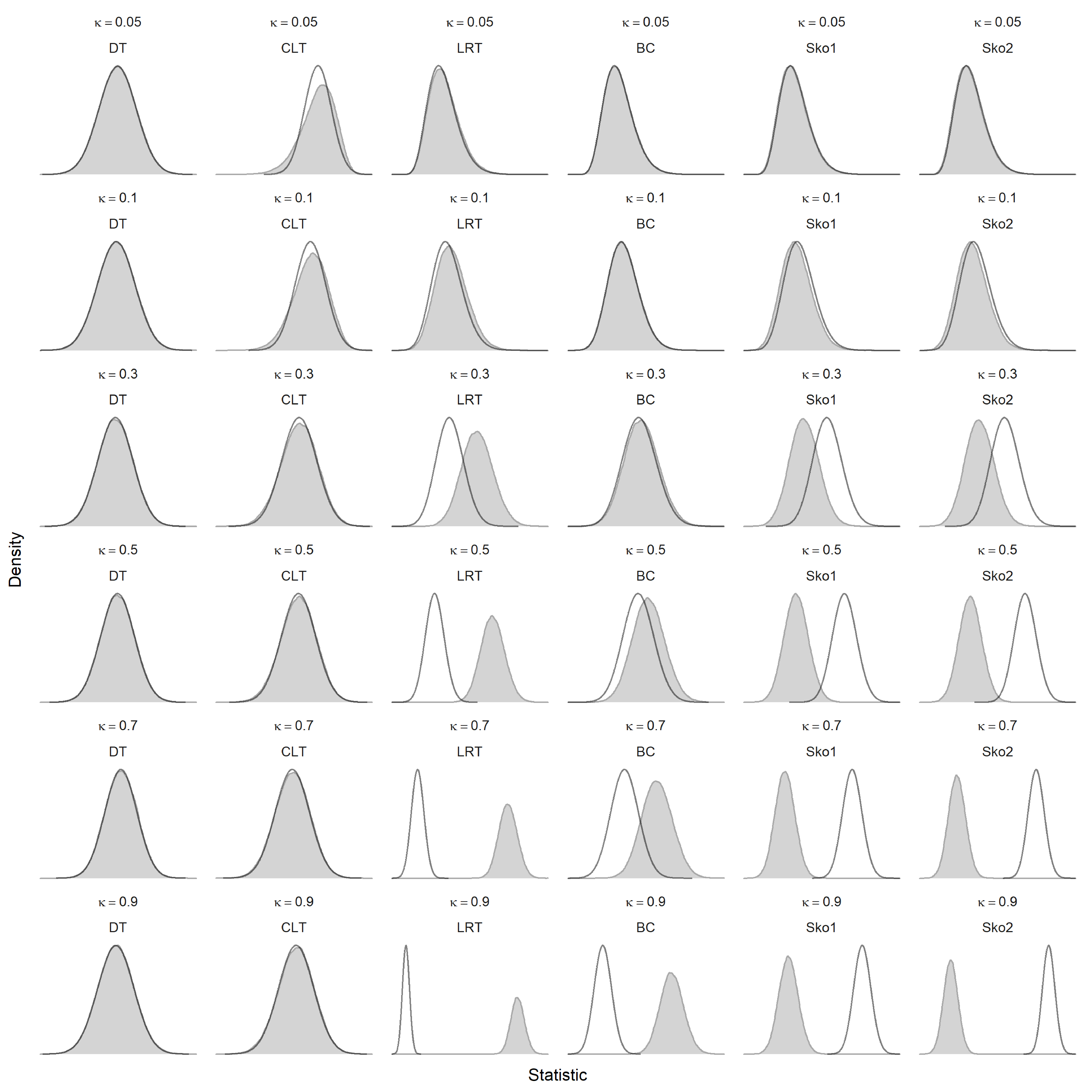}
	\caption{Hypothesis (VI). Comparison between   empirical (gray) and theoretical (black) null  distributions of the directional test (DT), the central limit theorem test (CLT), log-likelihood ratio test (LRT), Bartlett correction (BC),  and two Skovgaard's modifications \cite{skovgaard:2001} (Sko1 and Sko2, respectively),  for  various values of  $p/n$.}
	\label{SMsimulation:density null case5}
\end{figure}

 Figures \ref{SMsimulation:pvalue case1}--\ref{SMsimulation:pvalue case5} display the empirical null distribution of the $p$-values.  Figures \ref{SMsimulation:density null case1} - \ref{SMsimulation:density null case5} show the empirical  null distribution of the statistics.   Tables \ref{SM:table type I normal:case1}--\ref{SM:table type I normal:case5}  report the empirical and corrected probability of Type I error at the nominal level $\alpha=0.05$. The directional test  performs very well across all different values of $p$, while the central limit theorem test is still less accurate when $p$ is small, in particular if $p=5$. When $p=5,10$, $W_{BC}$ proves accurate, with  estimated and corrected Type I error very close to the nominal level.   However,  $W$, $W_{BC}$, $W^*$ and $W^{**}$  break down when $p$ is large.  Below, we provide more comments on the simulation outcomes for the empirical probability of Type I error.
\begin{itemize}
	\item[i)]  Figures \ref{SMsimulation:pvalue case1}, \ref{SMsimulation:density null case1} and Table \ref{SM:table type I normal:case1}. The directional test exhibits an extremely  precise  empirical Type I error in all examined frameworks.   When $p=5,10$, the chi-square approximations to $W$, $W^*$, $W^{**}$ and,  especially, $W_{BC}$ are more accurate than the normal approximation to the central limit theorem test.  However, as $p$ increases,   the four chi-square approximations fail because of the large location and scale bias, whereas the central limit theorem test  has good properties. In particular, for the corrected Type I error at level $0.05$,   $W_{BC}$ is particularly good at controlling  the Type I error when  $p=5,10$.
	\item[ii)]  Figures  \ref{SMsimulation:pvalue case2}, \ref{SMsimulation:density null case2} and Table \ref{SM:table type I normal:case2}. The directional test  performs very well across the different values of $p$, while the central limit theorem test is still less accurate when $p$ is small, in particular if $p=5$. When $p=5,10$, $W_{BC}$ proves accurate, with  estimated and corrected Type I error very close to the nominal level.   However, the  tests $W$, $W_{BC}$, $W^*$ and $W^{**}$ fail, when $p$ is large. 
	\item[iii)]  Figures  \ref{SMsimulation:pvalue case6}, \ref{SMsimulation:density null case6} and Table \ref{SM:table type I normal:case6} indicate that the directional test has the best performance in terms of $p$-value and the limiting null distribution under the null hypothesis.  Moreover, the tests $W$, $W_{BC}$, $W^*$ and $W^{**}$  result in a reasonable empirical Type I error only for small values of $p$, such as $p=5,10$.
	\item[iv)] Figures  \ref{SMsimulation:pvalue case5}, \ref{SMsimulation:density null case5} and  Table \ref{SM:table type I normal:case5}. As in the previous cases, the  directional test and  central limit theorem test perform similarly for general values of $p$, but the directional test is superior when $p$ is small.   Indeed, the central limit theorem test exhibits large bias for the corrected Type I error for small $p$, more evidently than in the other cases.   The tests  $W$, $W_{BC}$, $W^*$ and $W^{**}$ behave similarly as before, getting more and more inaccurate as $p$ increases.
\end{itemize}



\begin{figure}[H]
	\centering
	\captionsetup{font=footnotesize}
	\subfigure{
		\begin{minipage}[b]{.3\linewidth}
			\centering
			\includegraphics[scale=0.0735]{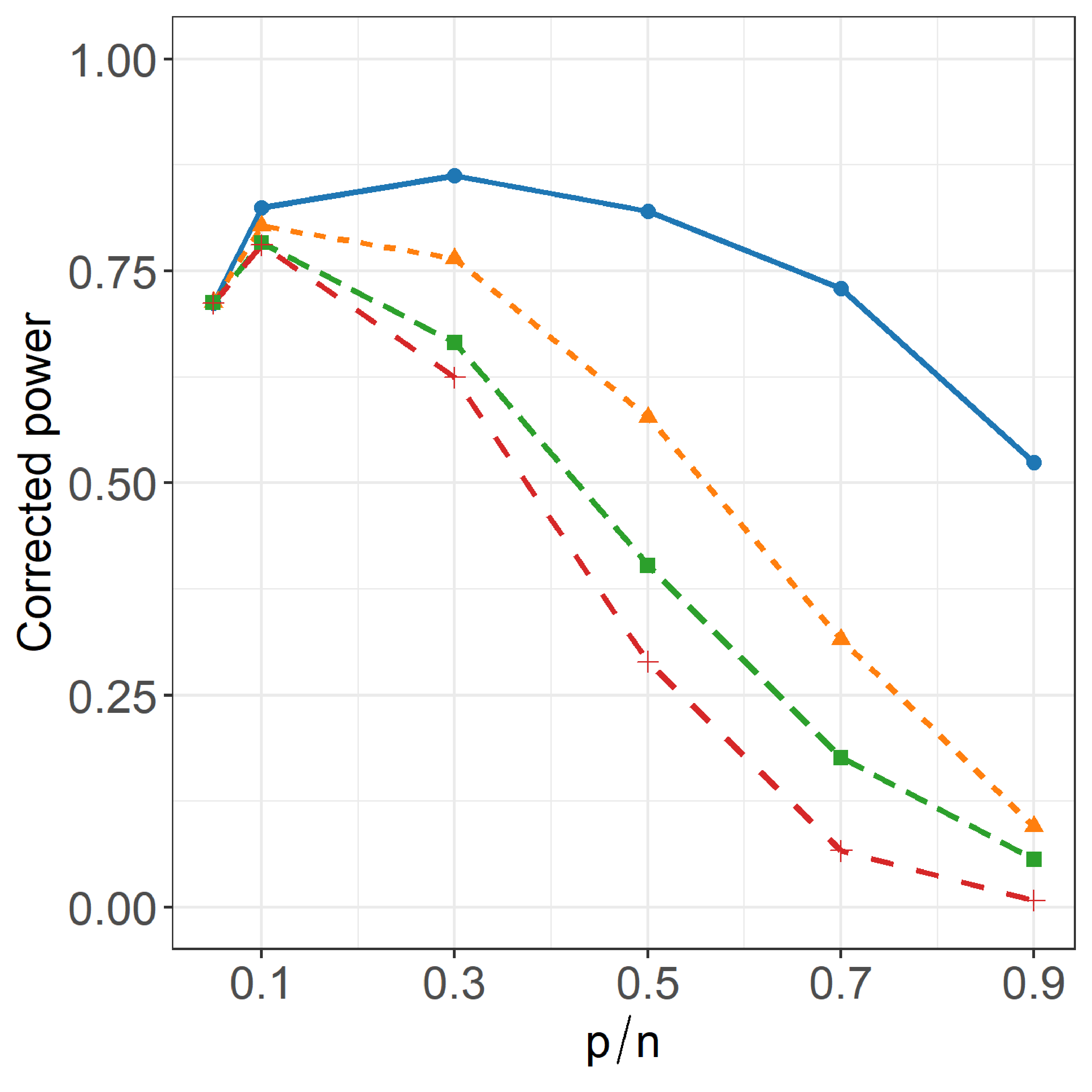}
		\end{minipage}
	}
	\subfigure{
		\begin{minipage}[b]{.3\linewidth}
			\centering
			\includegraphics[scale=0.0735]{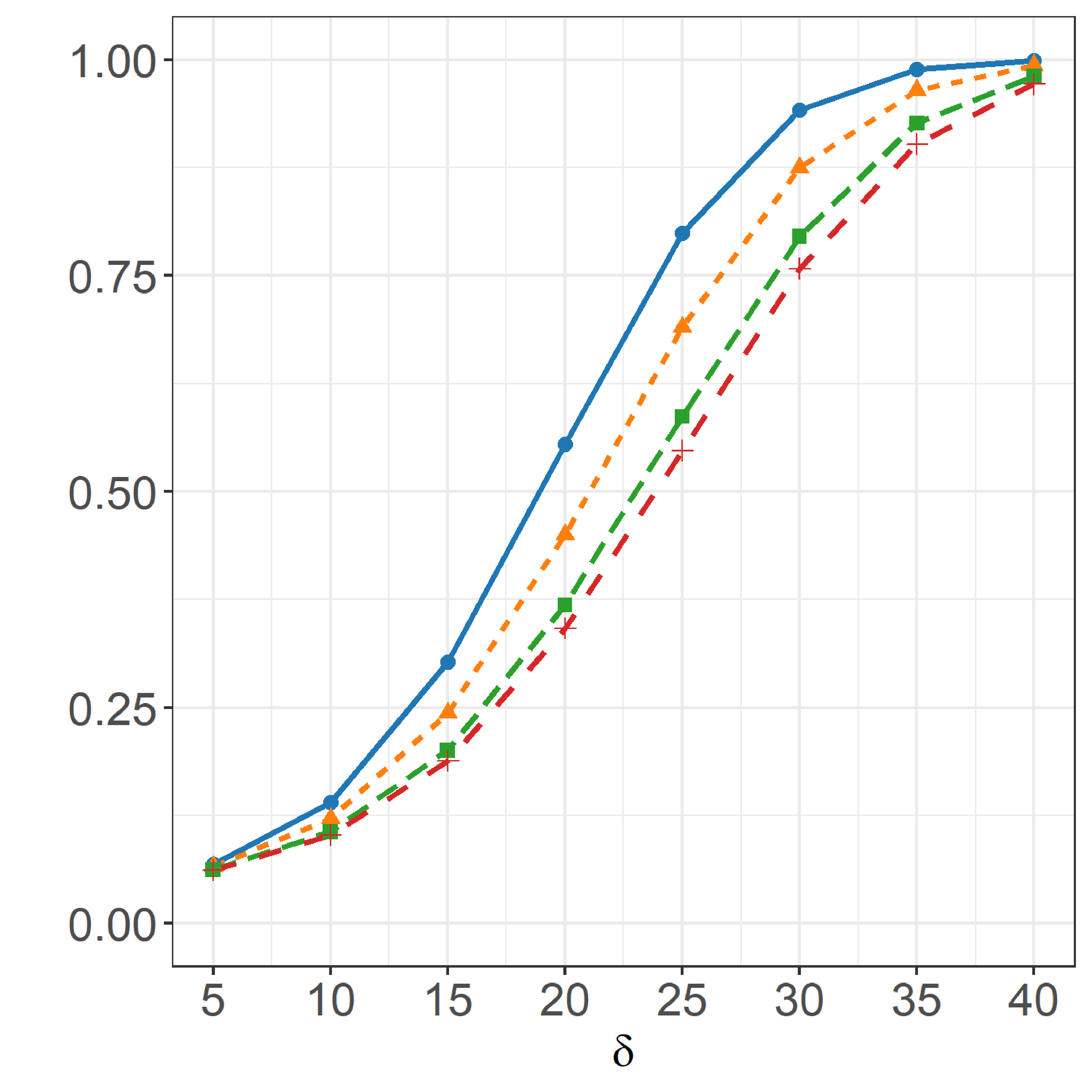}
		\end{minipage}
	}
	\subfigure{
	\begin{minipage}[b]{.3\linewidth}
		\centering
		\includegraphics[scale=0.0735]{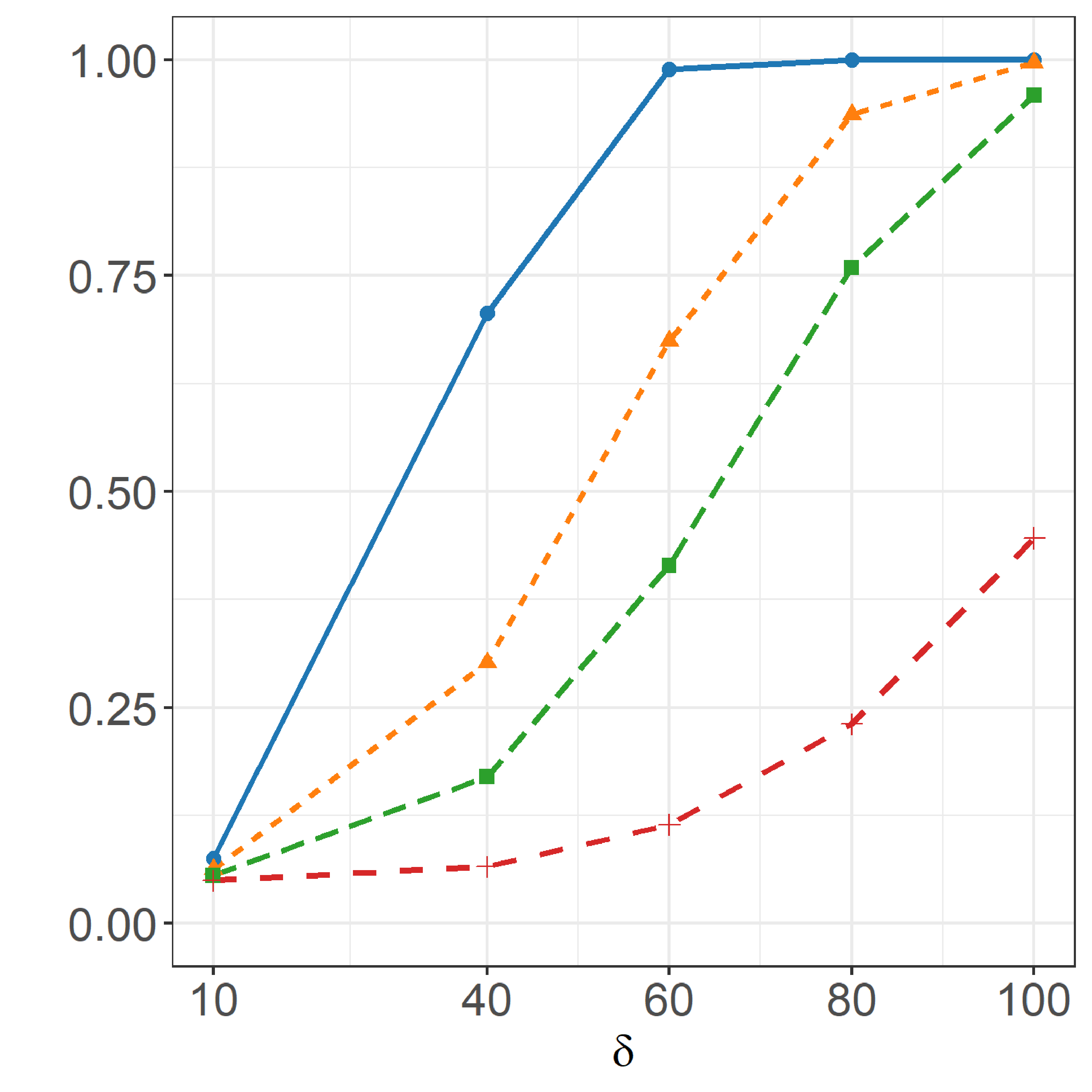}
	\end{minipage}
}
	\\
	\subfigure{
		\begin{minipage}[b]{.3\linewidth}
			\centering
			\includegraphics[scale=0.0735]{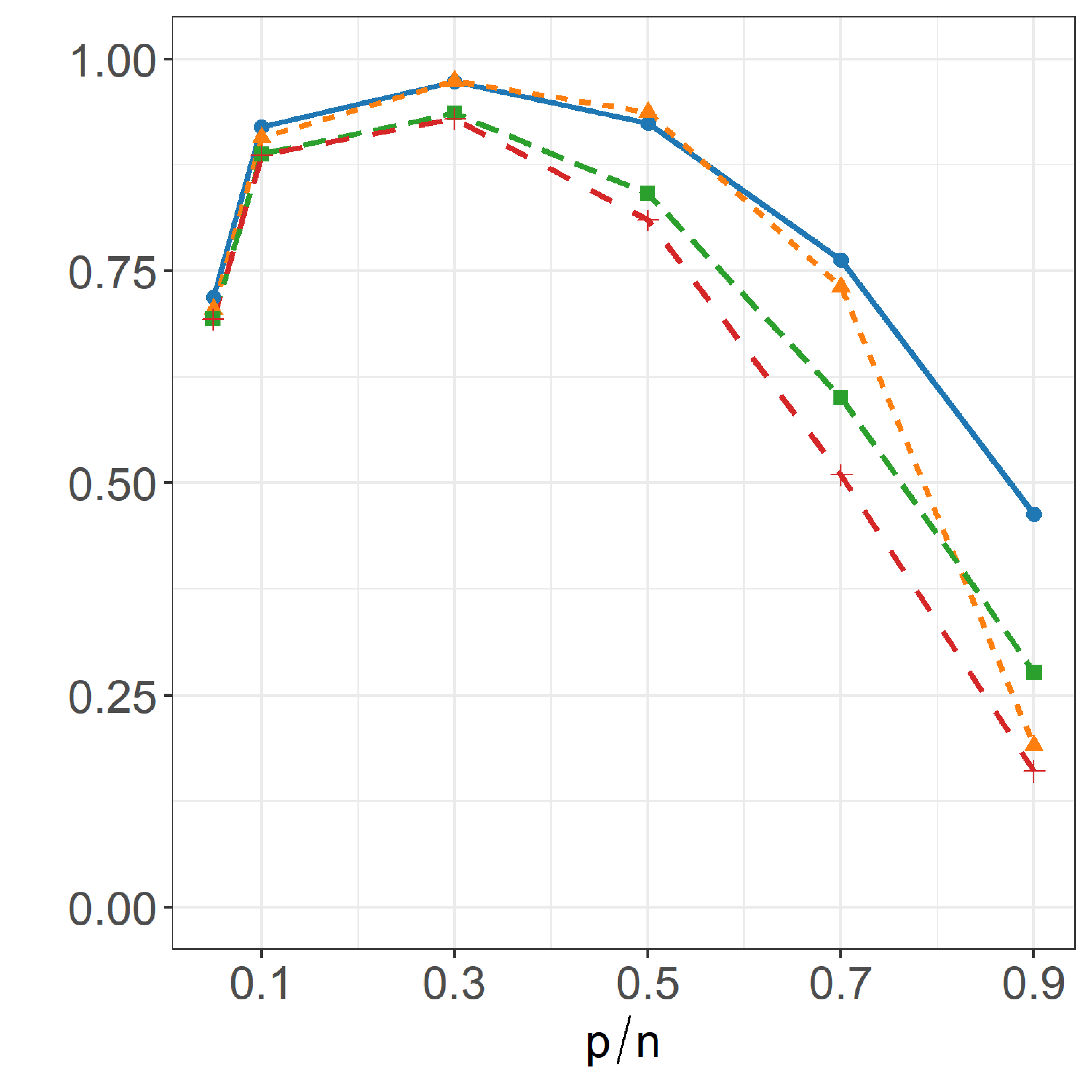}
		\end{minipage}
	}
	\subfigure{
		\begin{minipage}[b]{.3\linewidth}
			\centering
			\includegraphics[scale=0.0735]{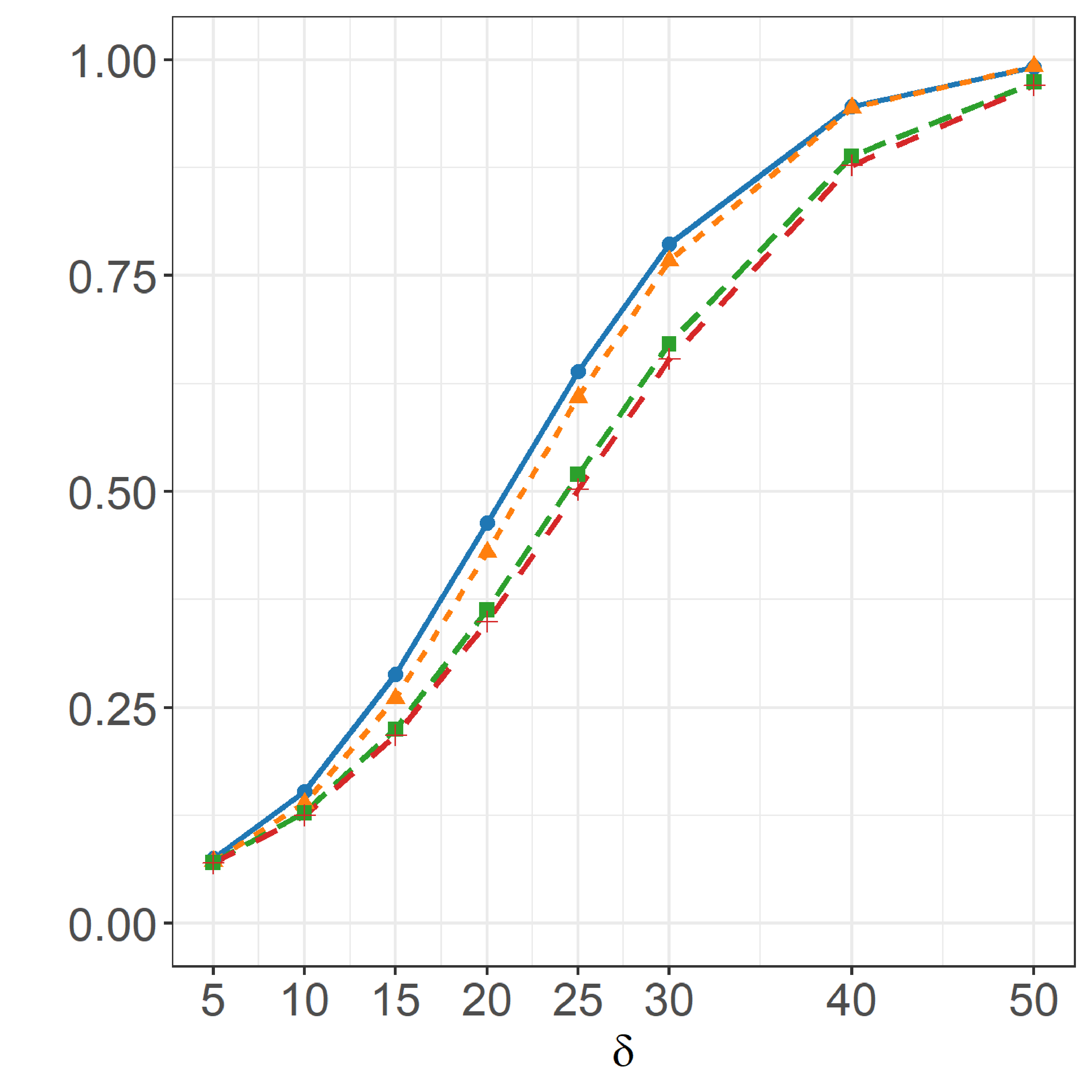}
		\end{minipage}
	}
	\subfigure{
	\begin{minipage}[b]{.3\linewidth}
		\centering
		\includegraphics[scale=0.0735]{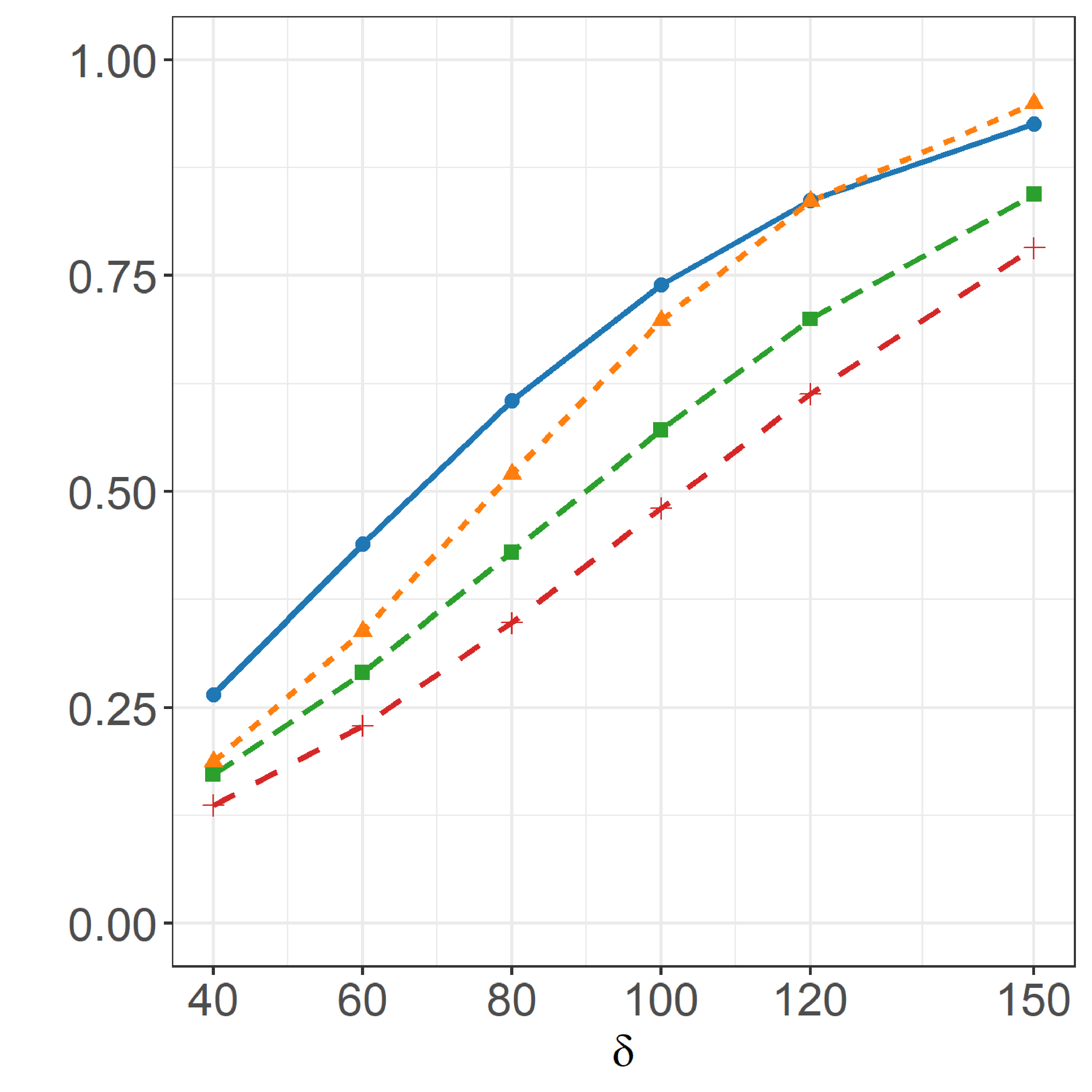}
	\end{minipage}
}
	\\
	\subfigure{
		\begin{minipage}[b]{.3\linewidth}
			\centering
			\includegraphics[scale=0.0735]{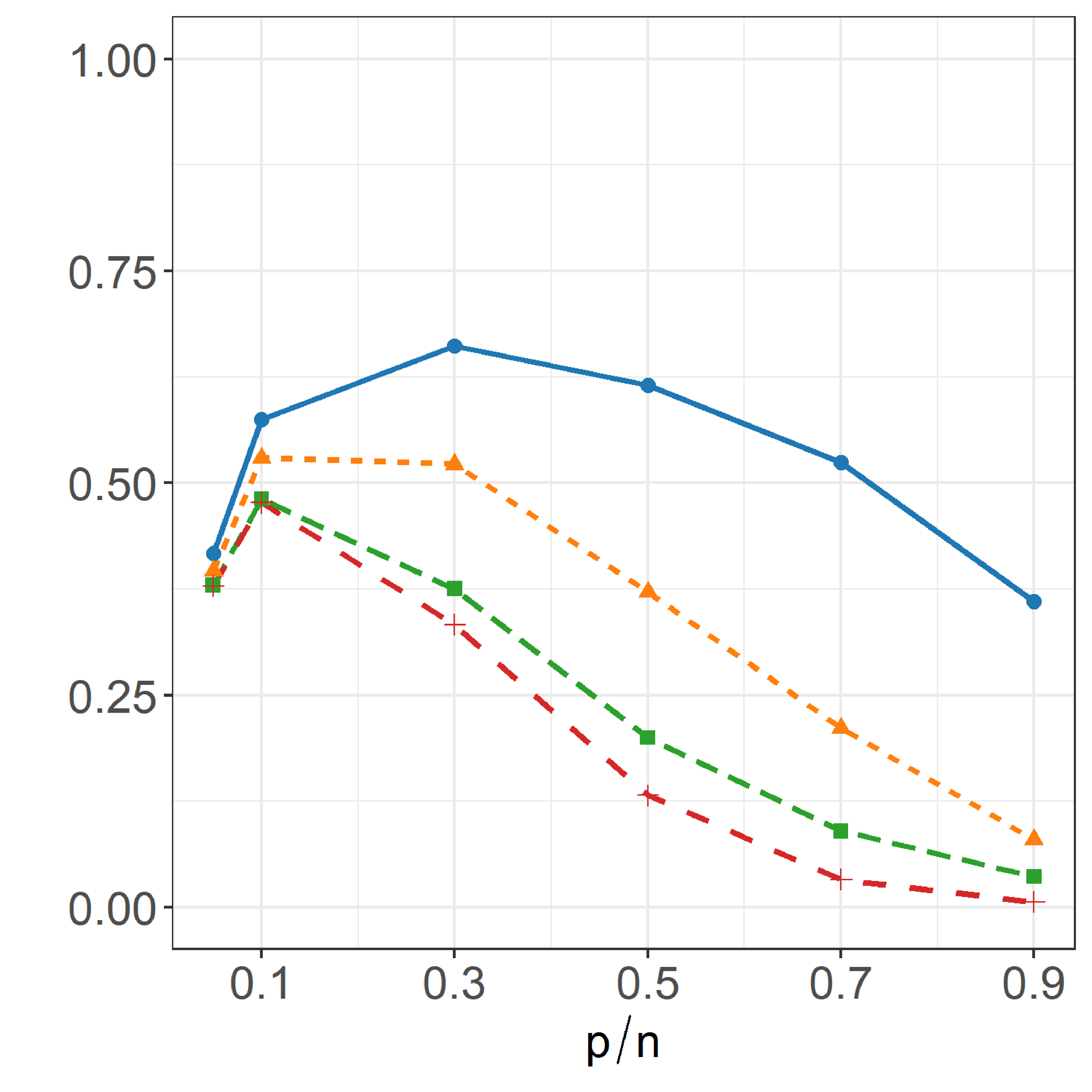}
		\end{minipage}
	}
	\subfigure{
		\begin{minipage}[b]{.3\linewidth}
			\centering
			\includegraphics[scale=0.0735]{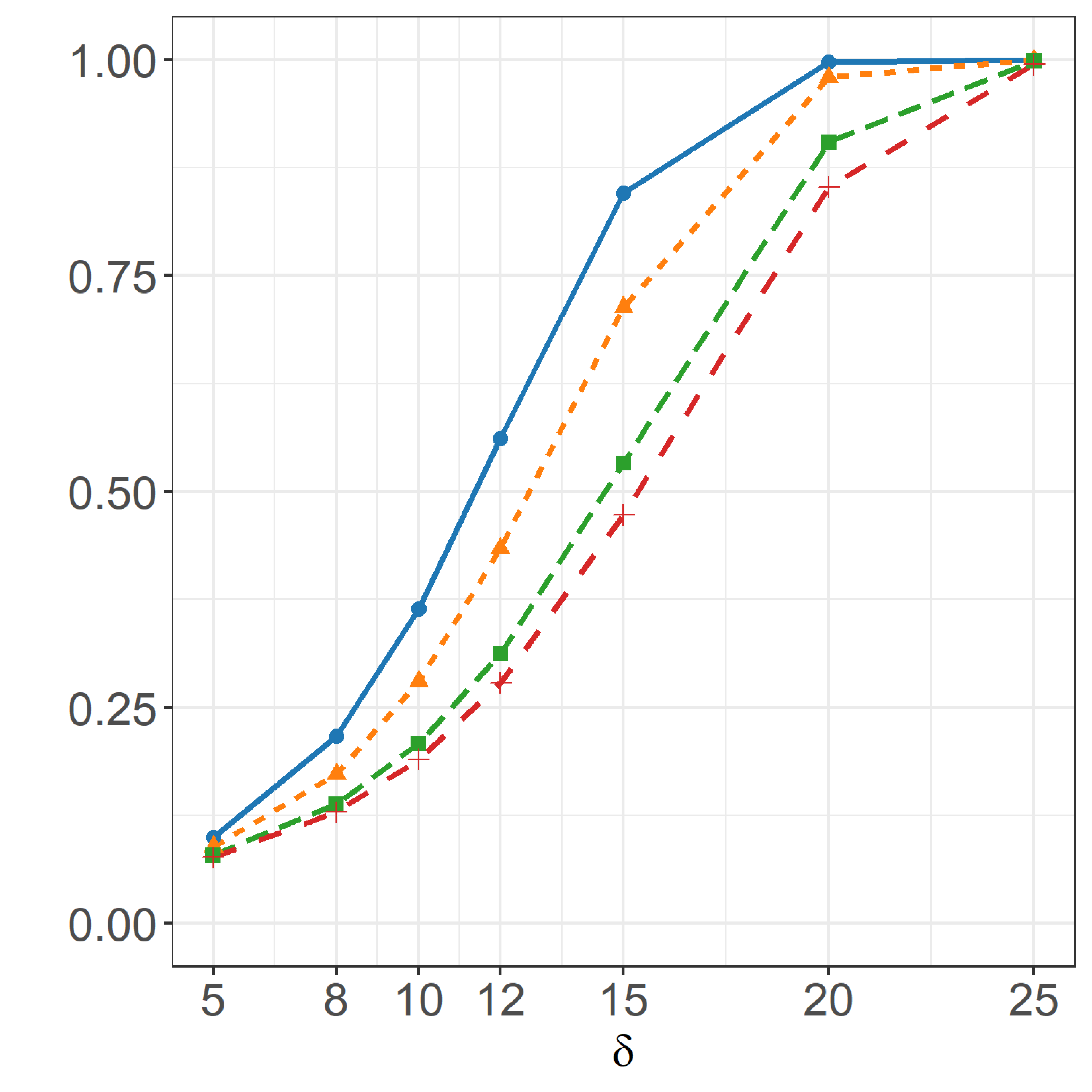}
		\end{minipage}
	}
	\subfigure{
	\begin{minipage}[b]{.3\linewidth}
		\centering
		\includegraphics[scale=0.0735]{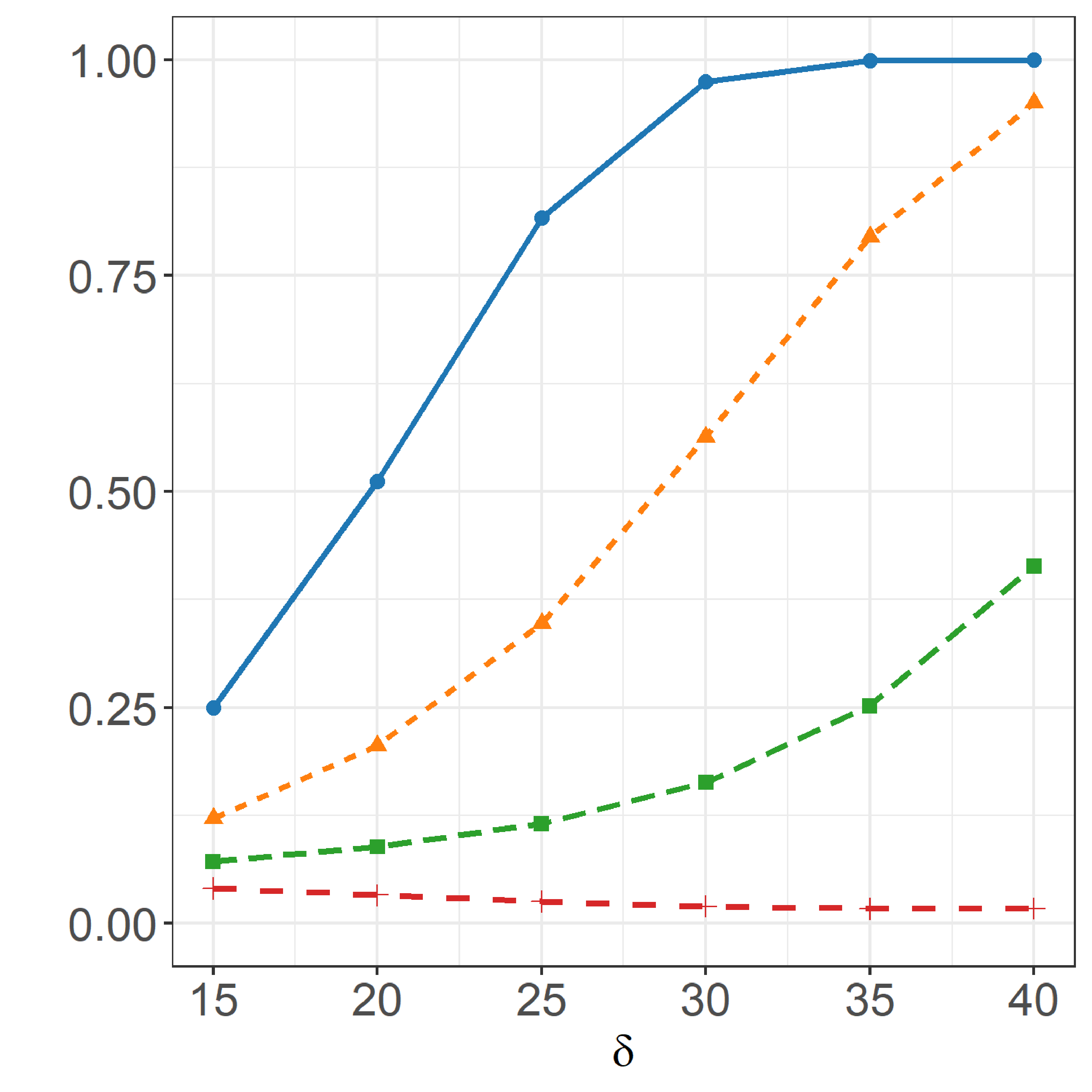}
	\end{minipage}
}
	\\
	\subfigure{
		\begin{minipage}[b]{.3\linewidth}
			\centering
			\includegraphics[scale=0.0735]{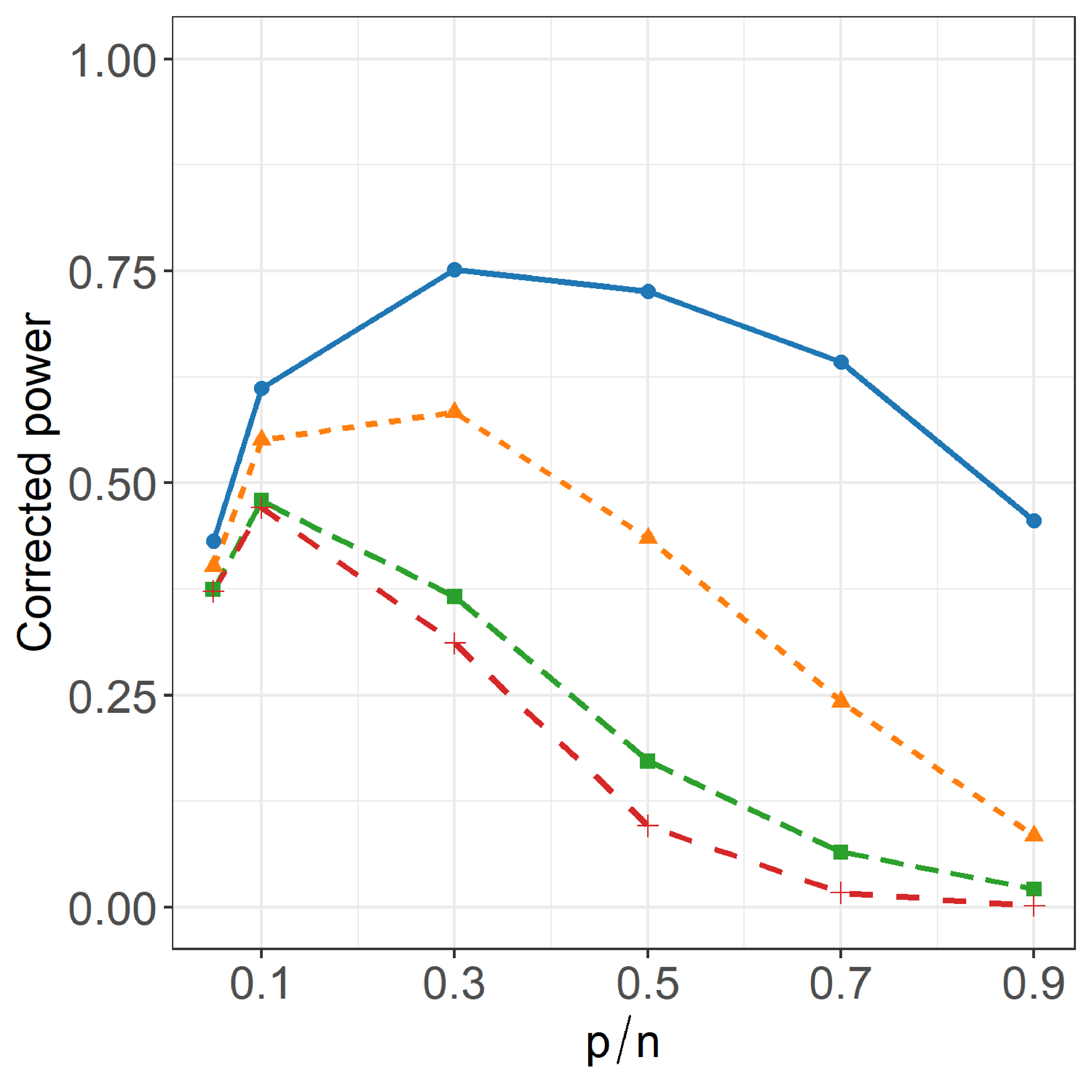}
		\end{minipage}
	}
	\subfigure{
		\begin{minipage}[b]{.3\linewidth}
			\centering
			\includegraphics[scale=0.0735]{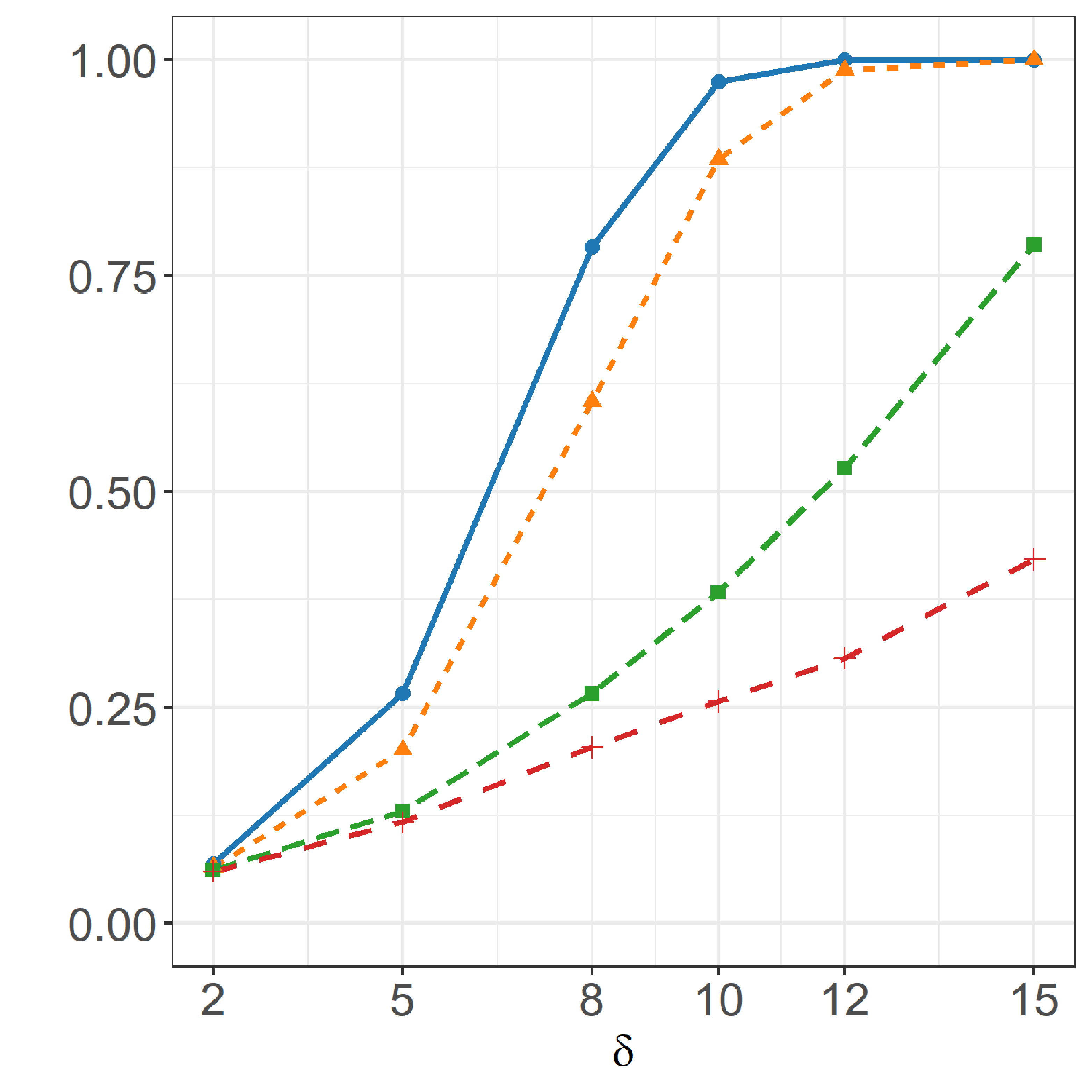}
		\end{minipage}
	}
	\subfigure{
	\begin{minipage}[b]{.3\linewidth}
		\centering
		\includegraphics[scale=0.0735]{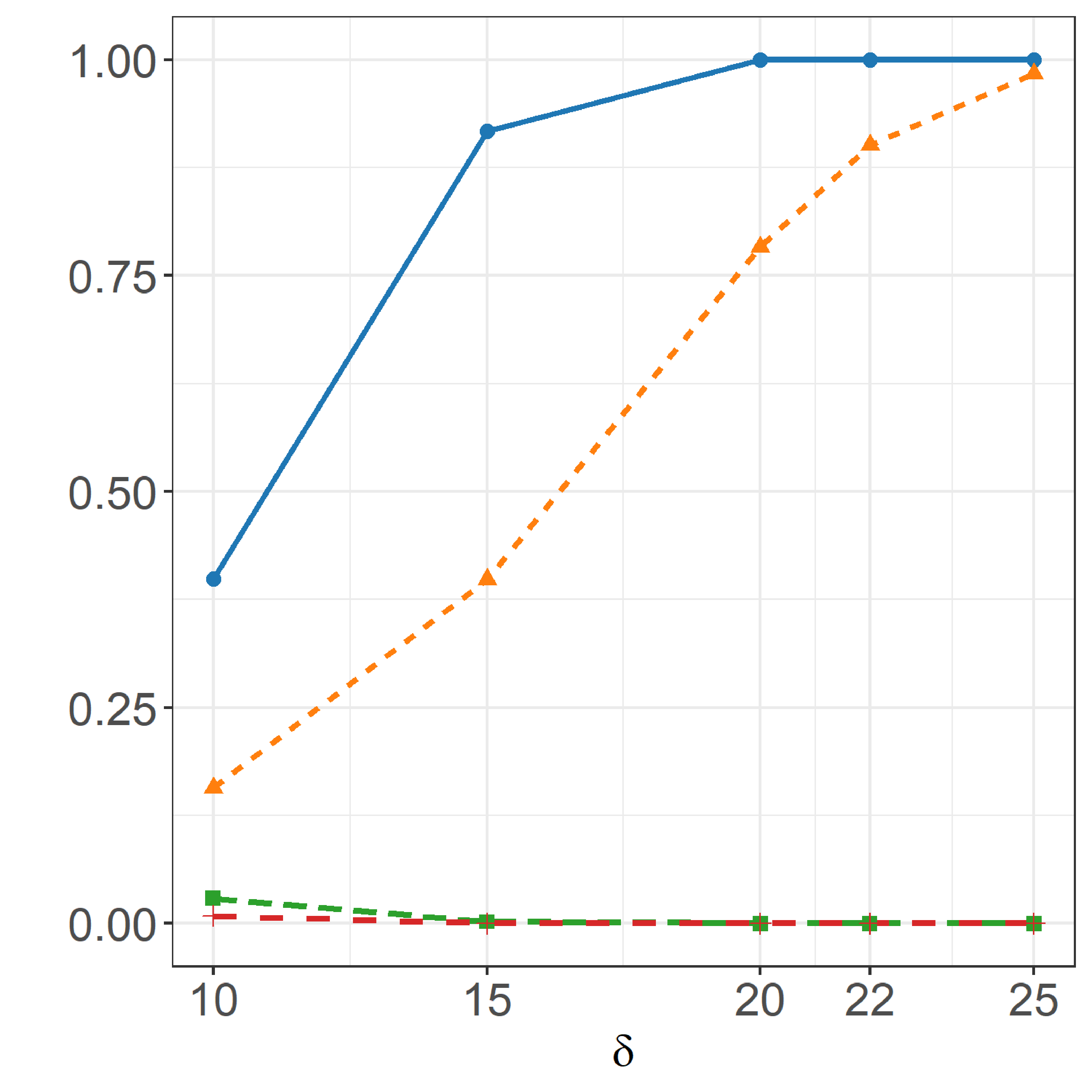}
	\end{minipage}
}
	\caption{Empirical local powers of four tests. The solid, dashed, longdashed, and dotdashed curves are the empirical power functions of the central limit theorem test, directional test and two Skovgaard’s modifications, respectively.   The first, second, third and fourth rows correspond to  hypotheses (III) - (VI), respectively;  the left column corresponds to the alternative hypothesis setting. The middle and right columns correspond to the alternative hypothesis settings 2 for the ratio $0.3$ and $0.7$ in Section \ref{S4:simulation (III)-(VI)}, respectively.}   
	\label{SMfigure:power}
\end{figure}

\begin{figure}[H]
	\centering
	\captionsetup{font=footnotesize}
	\subfigure{
		\begin{minipage}[b]{.3\linewidth}
			\centering
			\includegraphics[scale=0.0735]{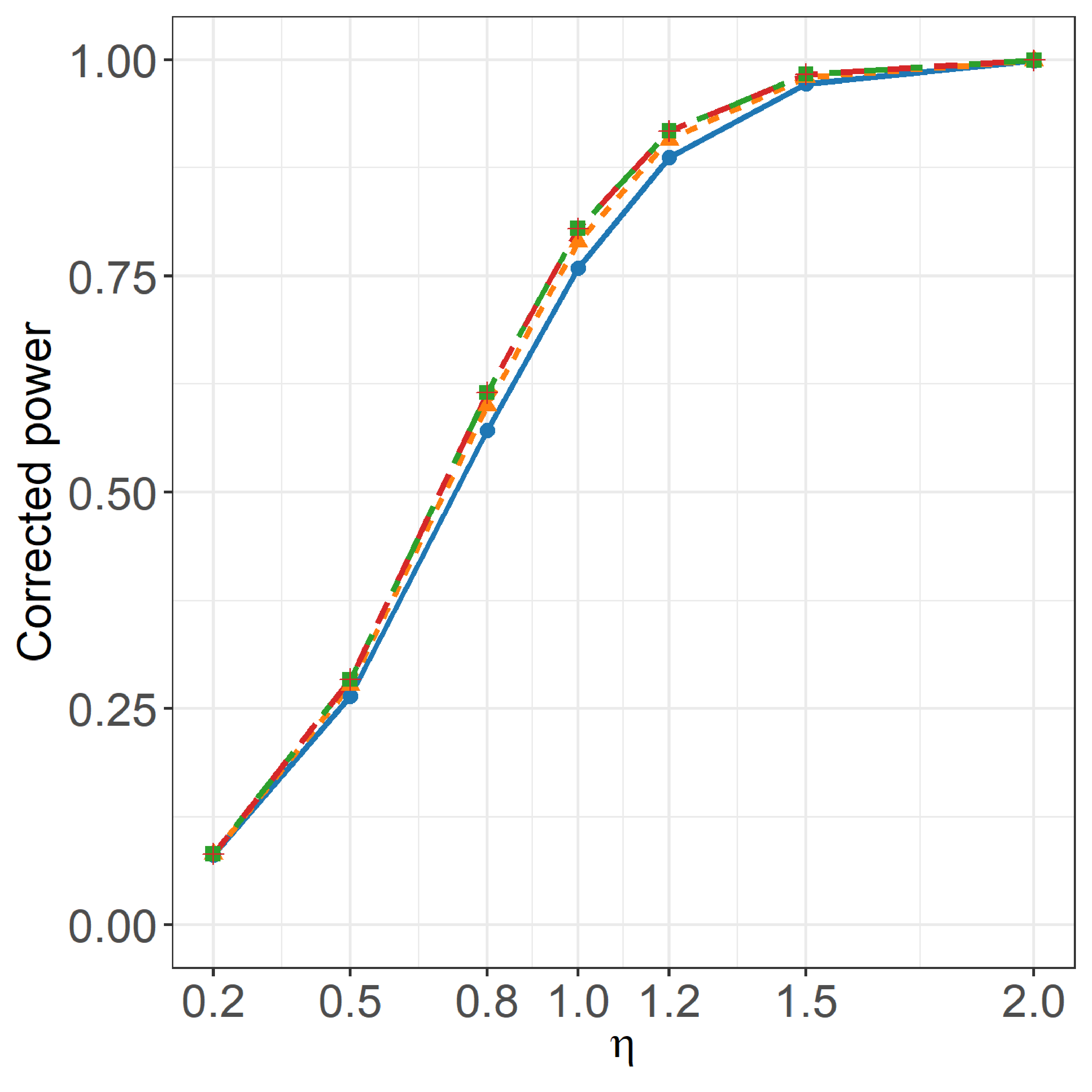}
		\end{minipage}
	}
	\subfigure{
		\begin{minipage}[b]{.3\linewidth}
			\centering
			\includegraphics[scale=0.0735]{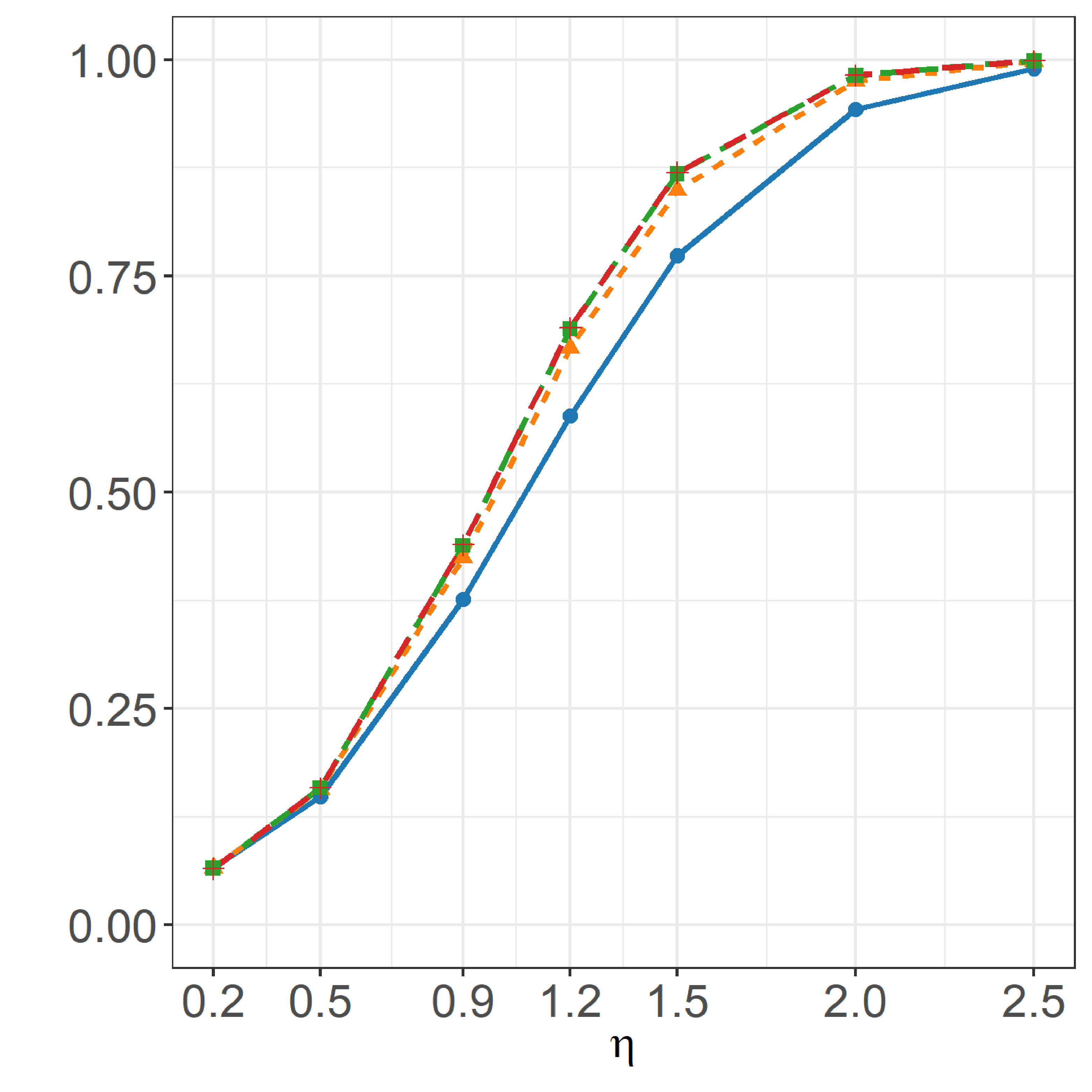}
		\end{minipage}
	}
	\subfigure{
		\begin{minipage}[b]{.3\linewidth}
			\centering
			\includegraphics[scale=0.0735]{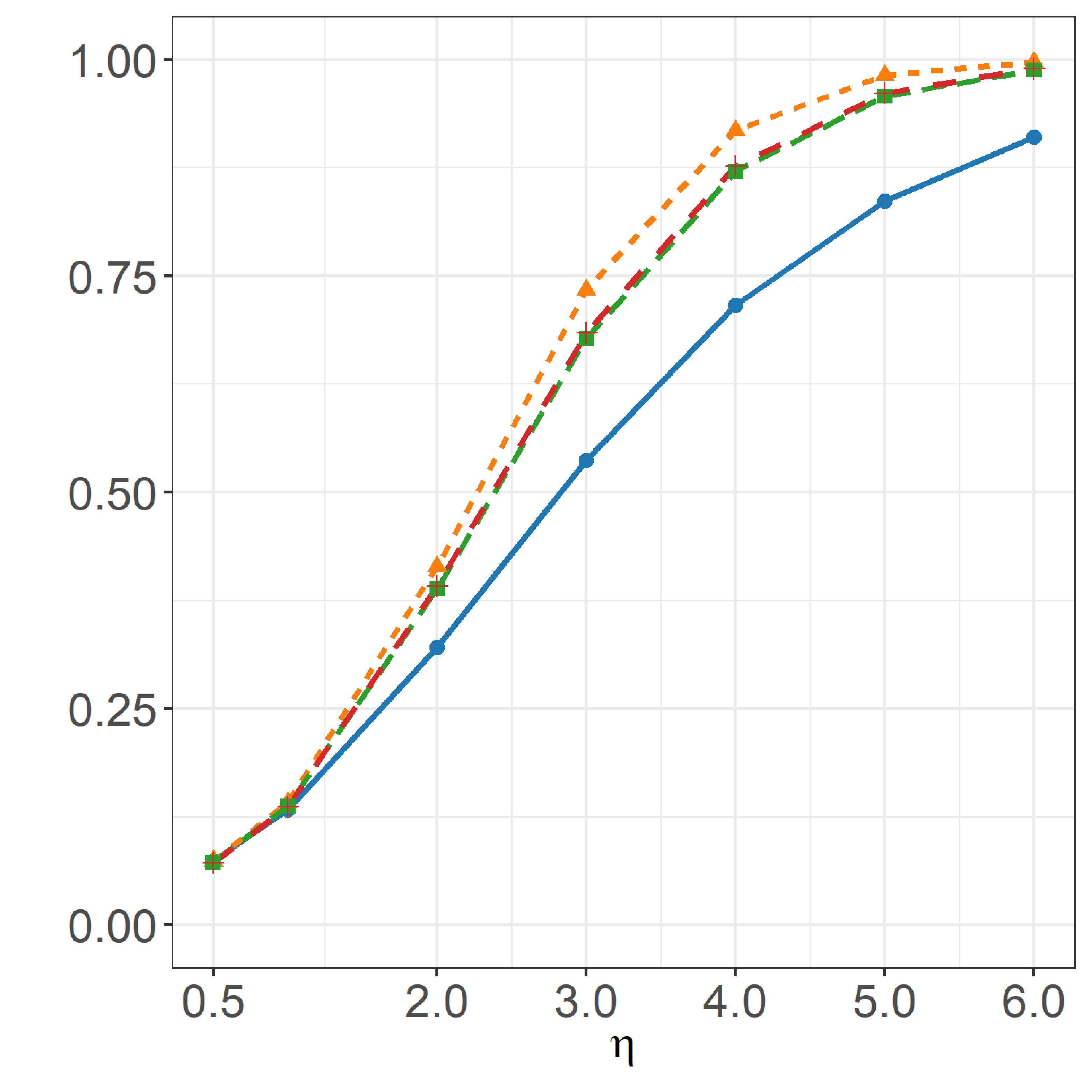}
		\end{minipage}
	}
	\subfigure{
		\begin{minipage}[b]{.3\linewidth}
			\centering
			\includegraphics[scale=0.0735]{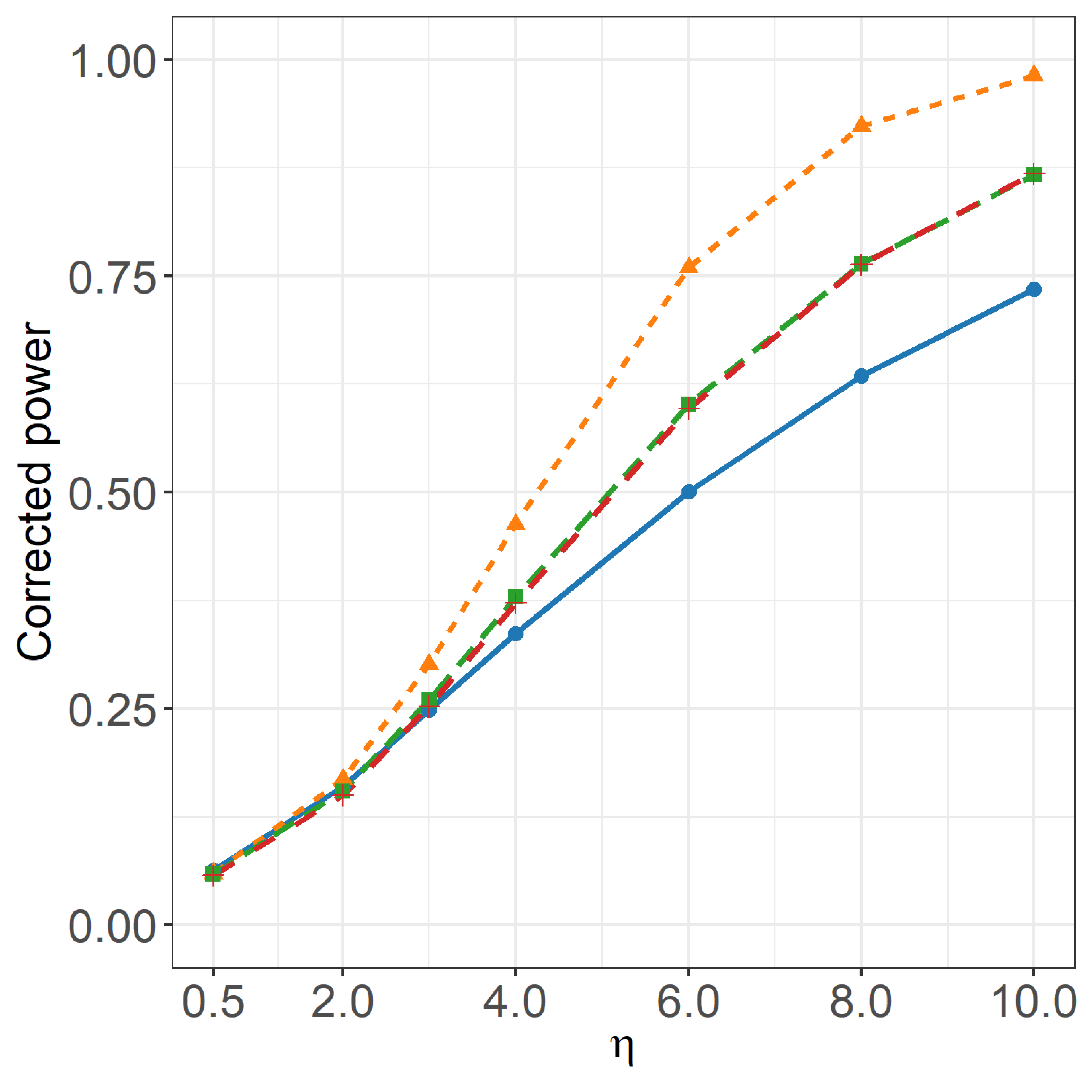}
		\end{minipage}
	}
	\subfigure{
		\begin{minipage}[b]{.3\linewidth}
			\centering
			\includegraphics[scale=0.0735]{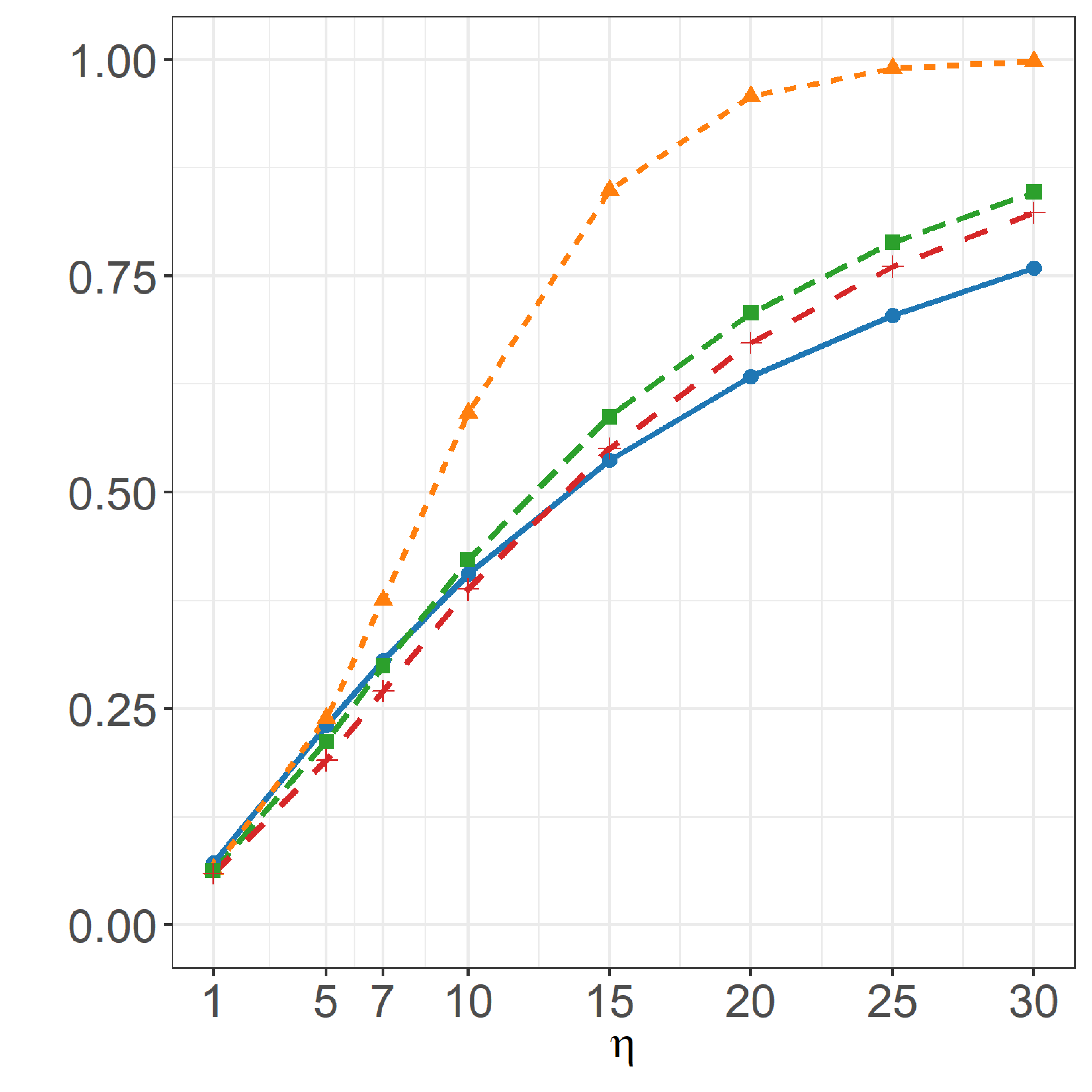}
		\end{minipage}
	}
	\subfigure{
		\begin{minipage}[b]{.3\linewidth}
			\centering
			\includegraphics[scale=0.0735]{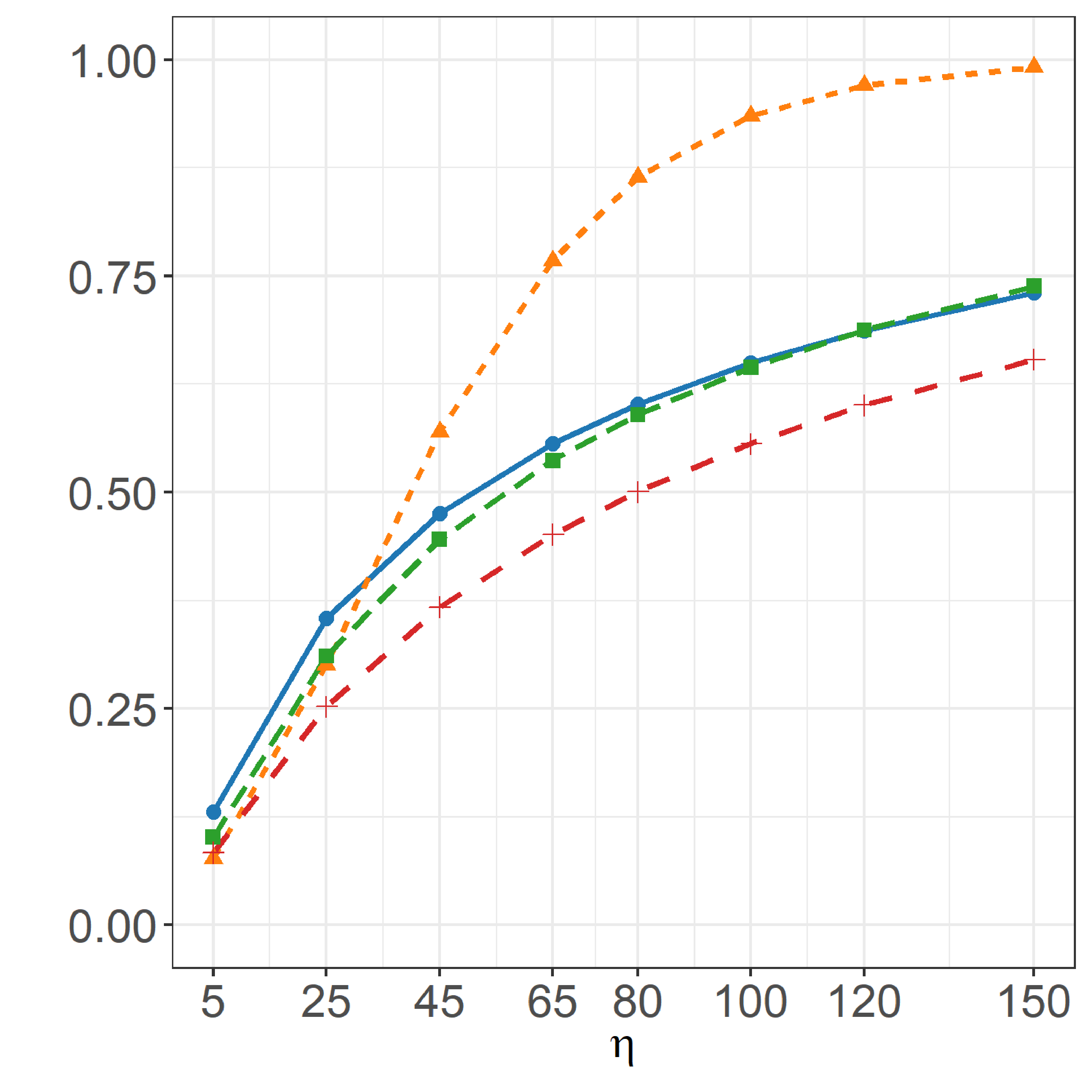}
		\end{minipage}
	}
	\caption{ Hypothesis (III). Empirical local power of four tests with different values of $\eta$ and $p/n$. The solid, dashed, longdashed, and dotdashed curves are the empirical power functions of the central limit theorem test, directional test and two Skovgaard’s modifications, respectively.   The extreme alternative hypothesis setting (3) is given in Section \ref{S4:simulation (III)-(VI)}. The six plots correspond to $p/n \in \{0.05, 0.1, 0.3, 0.5, 0.7, 0.9\}$, starting from top left and proceeding by row.}  
	\label{SMfig:power extreme case1}
\end{figure}

\begin{figure}[H]
	\centering
	\captionsetup{font=footnotesize}
	\subfigure{
		\begin{minipage}[b]{.3\linewidth}
			\centering
			\includegraphics[scale=0.0735]{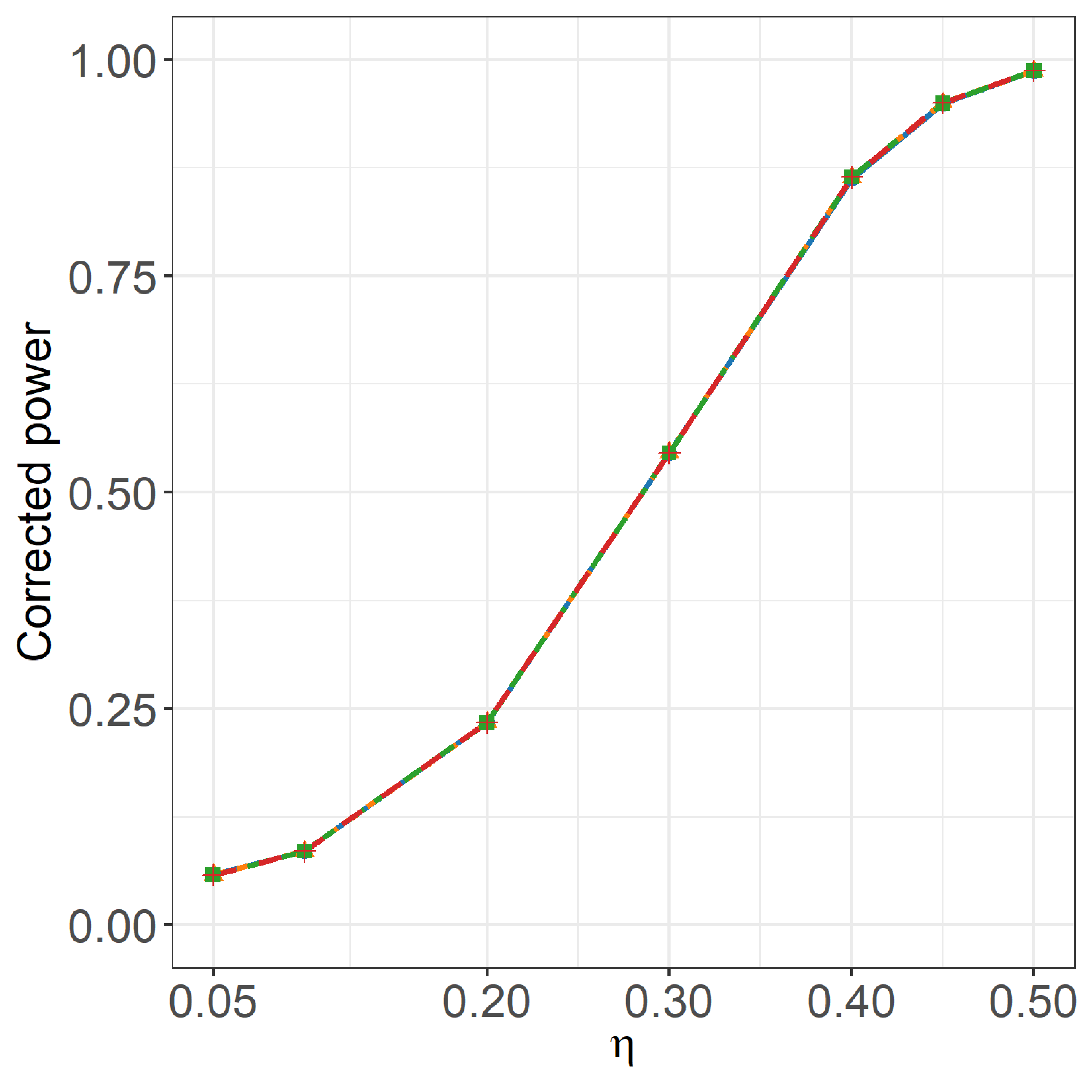}
		\end{minipage}
	}
	\subfigure{
		\begin{minipage}[b]{.3\linewidth}
			\centering
			\includegraphics[scale=0.0735]{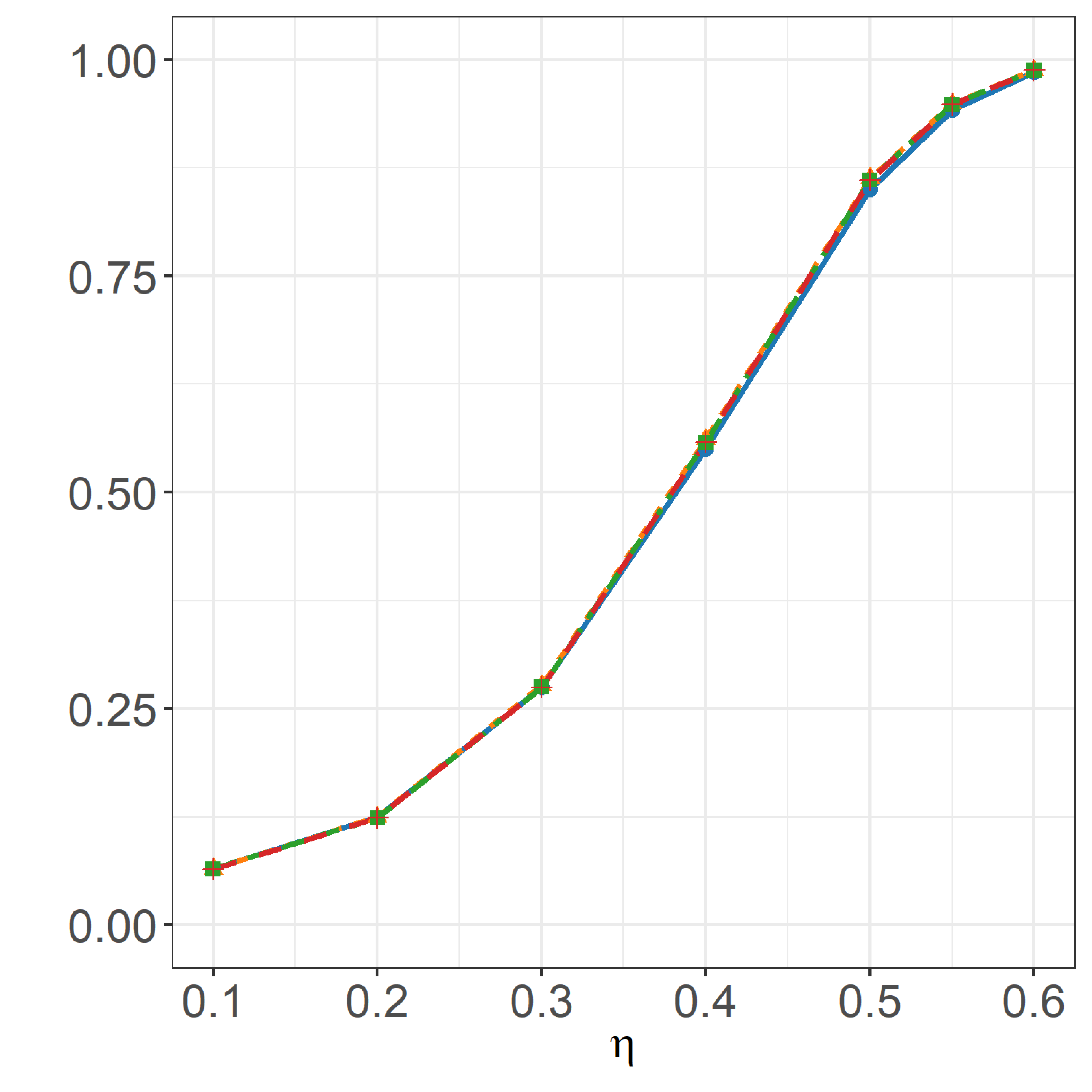}
		\end{minipage}
	}
	\subfigure{
		\begin{minipage}[b]{.3\linewidth}
			\centering
			\includegraphics[scale=0.0735]{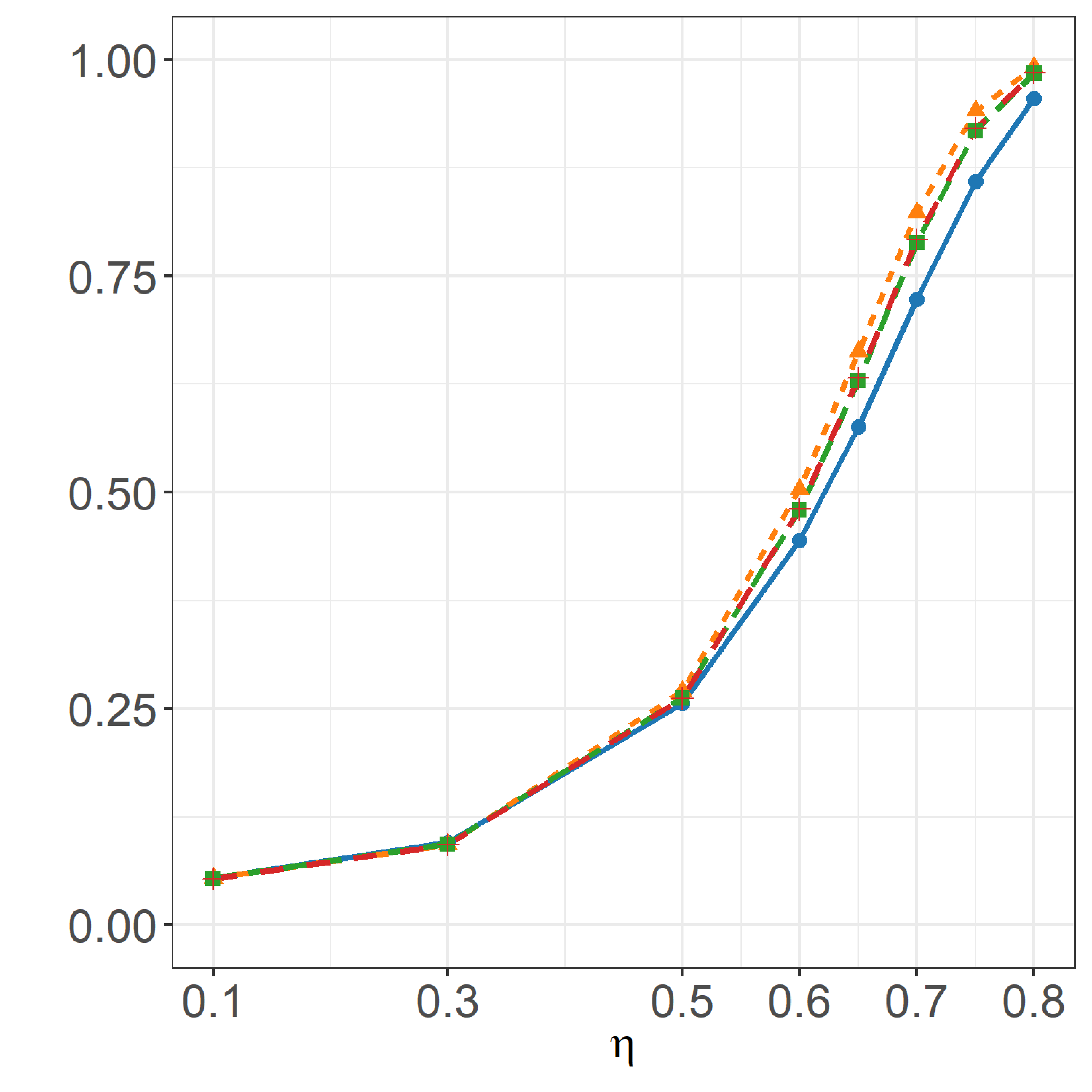}
		\end{minipage}
	}
	\subfigure{
		\begin{minipage}[b]{.3\linewidth}
			\centering
			\includegraphics[scale=0.0735]{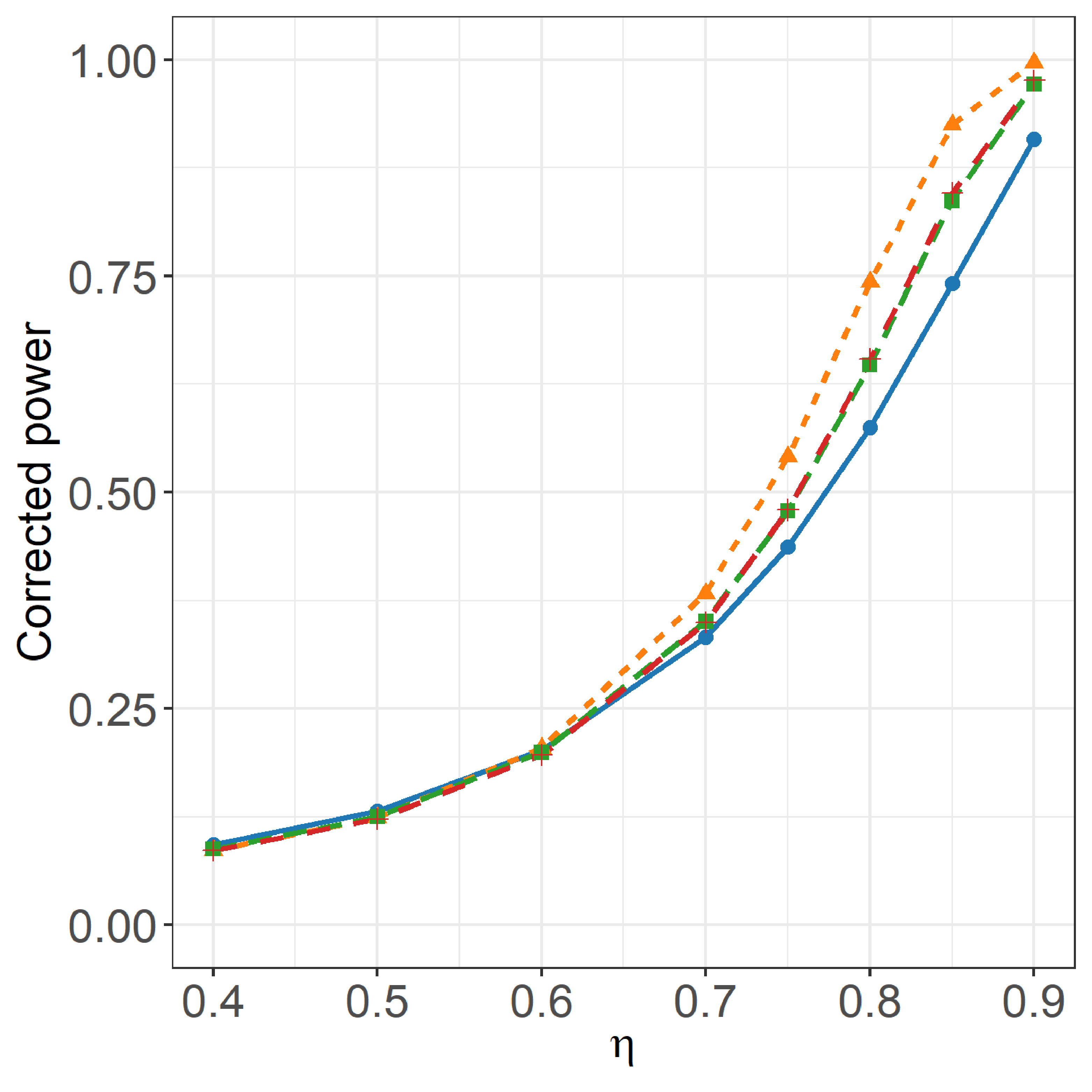}
		\end{minipage}
	}
	\subfigure{
		\begin{minipage}[b]{.3\linewidth}
			\centering
			\includegraphics[scale=0.0735]{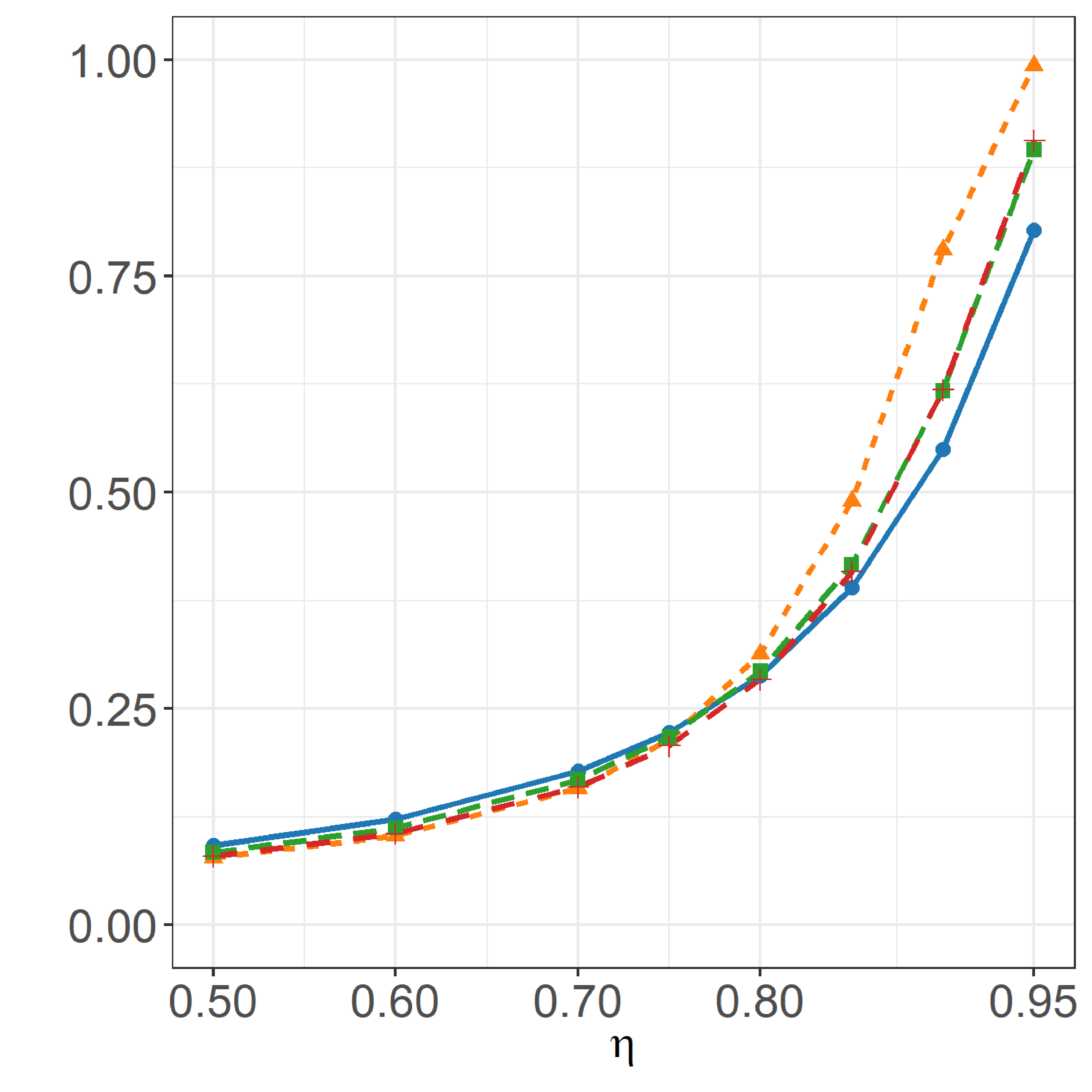}
		\end{minipage}
	}
	\subfigure{
		\begin{minipage}[b]{.3\linewidth}
			\centering
			\includegraphics[scale=0.0735]{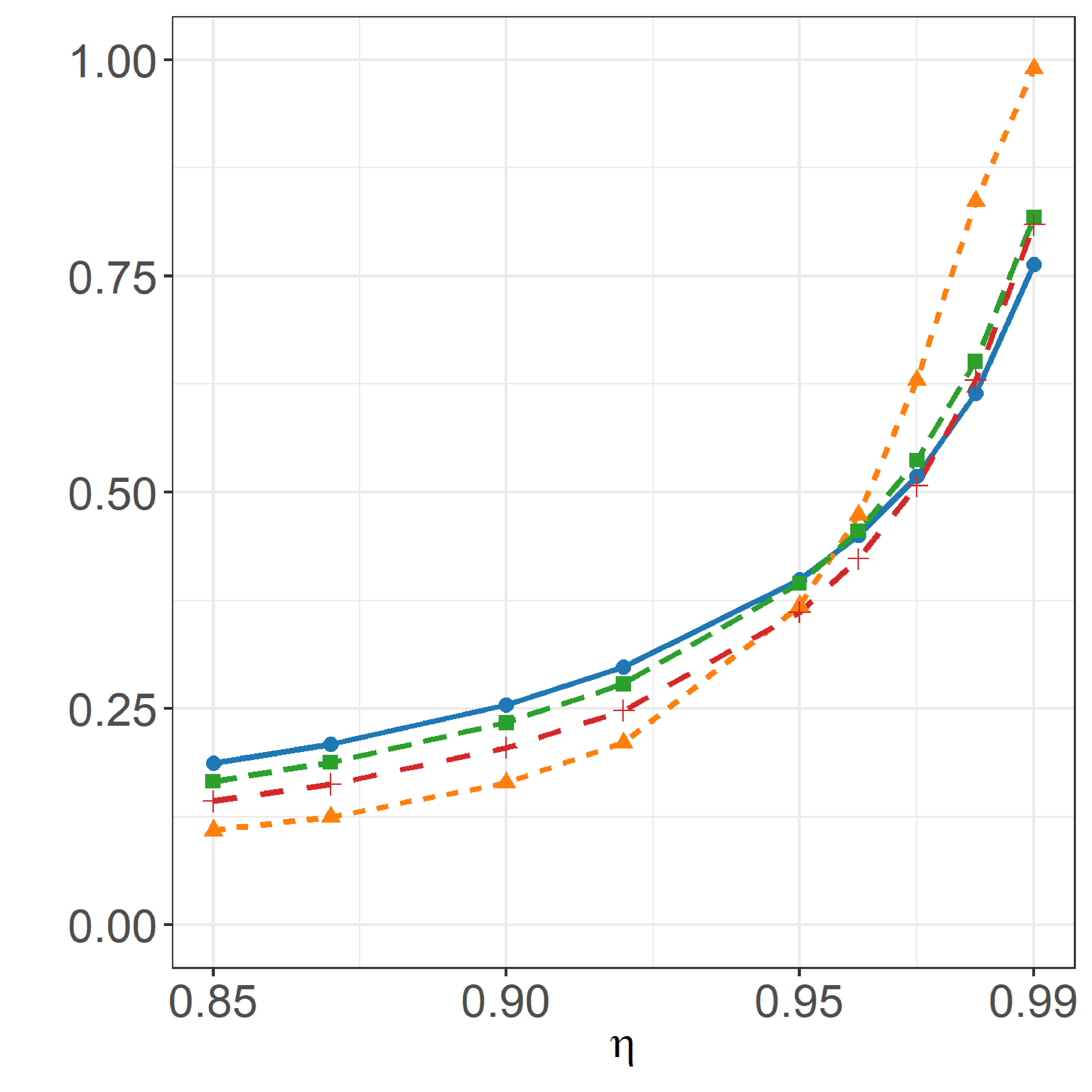}
		\end{minipage}
	}
	\caption{Hypothesis (IV). Empirical local power of four tests with different values of $\eta$ and $p/n$. The solid, dashed, longdashed, and dotdashed curves are the empirical power functions of the central limit theorem test, directional test and two Skovgaard’s modifications \cite{skovgaard:2001}, respectively.   The extreme alternative setting is given in Section \ref{S4:simulation (III)-(VI)}. The six plots correspond to $p/n \in \{0.05, 0.1, 0.3, 0.5, 0.7, 0.9\}$, starting from top left and proceeding by row. }  
	\label{SMfig:power extreme case2}
\end{figure}

\begin{figure}[H]
	\centering
	\captionsetup{font=footnotesize}
	\subfigure{
		\begin{minipage}[b]{.3\linewidth}
			\centering
			\includegraphics[scale=0.0735]{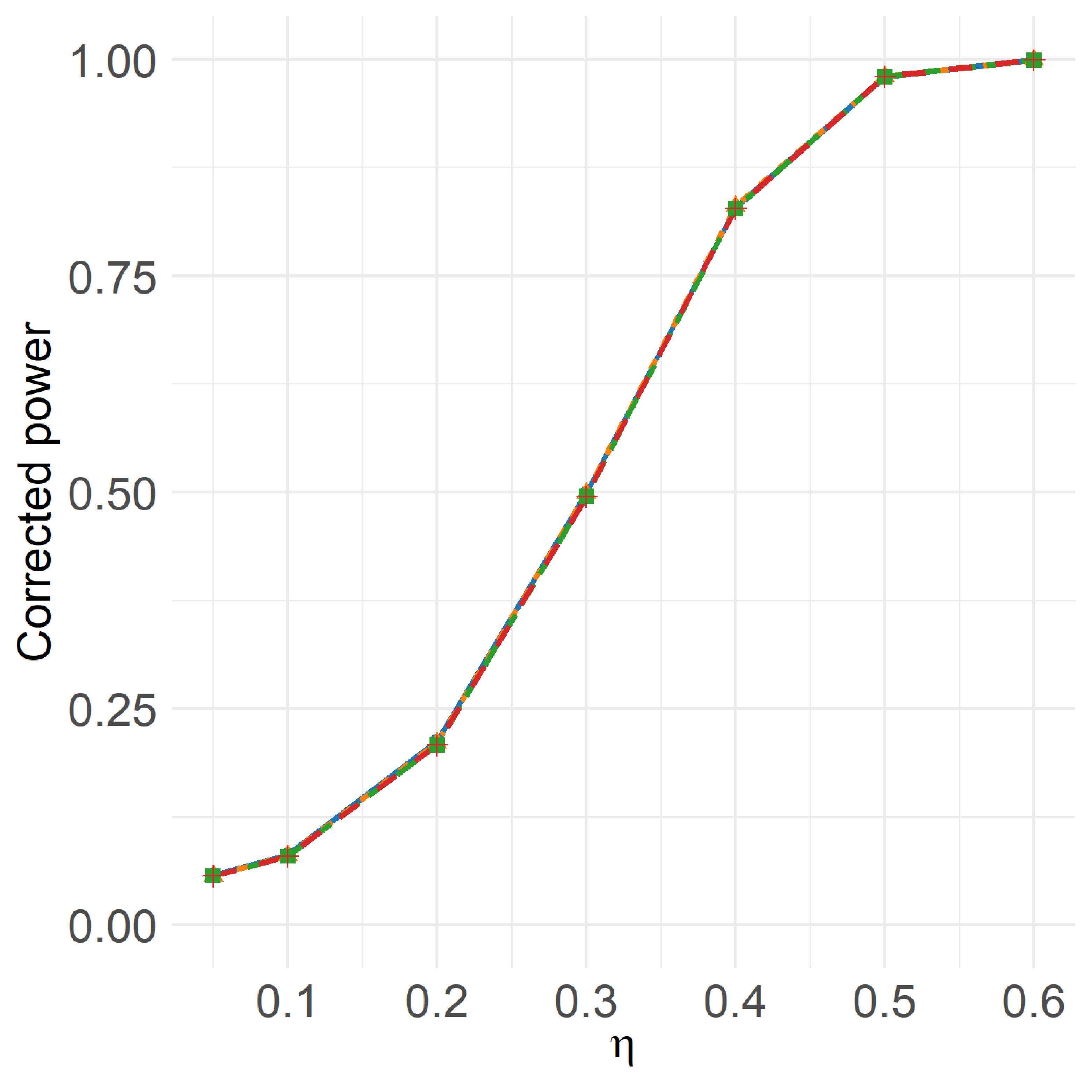}
		\end{minipage}
	}
	\subfigure{
		\begin{minipage}[b]{.3\linewidth}
			\centering
			\includegraphics[scale=0.0735]{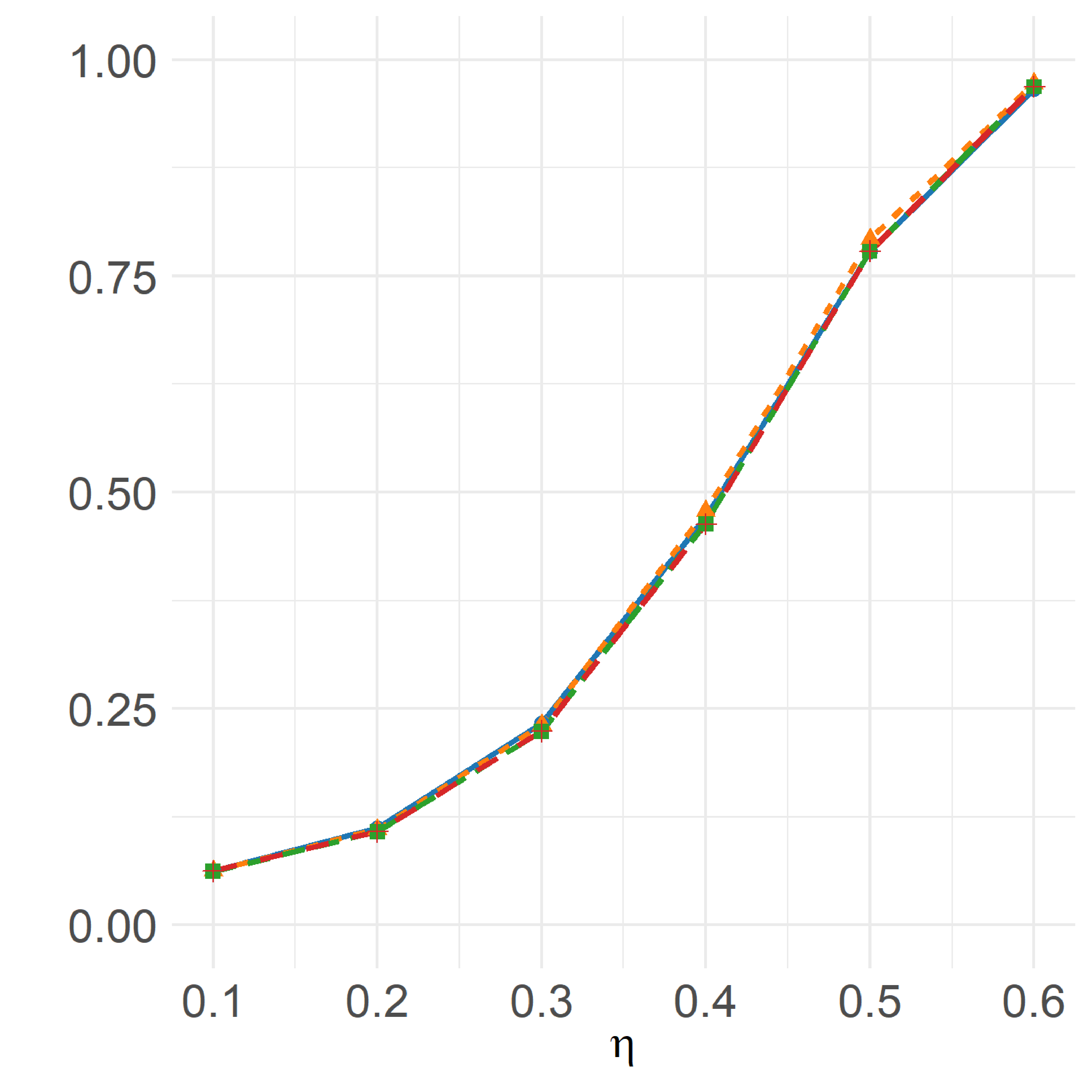}
		\end{minipage}
	}
	\subfigure{
		\begin{minipage}[b]{.3\linewidth}
			\centering
			\includegraphics[scale=0.0735]{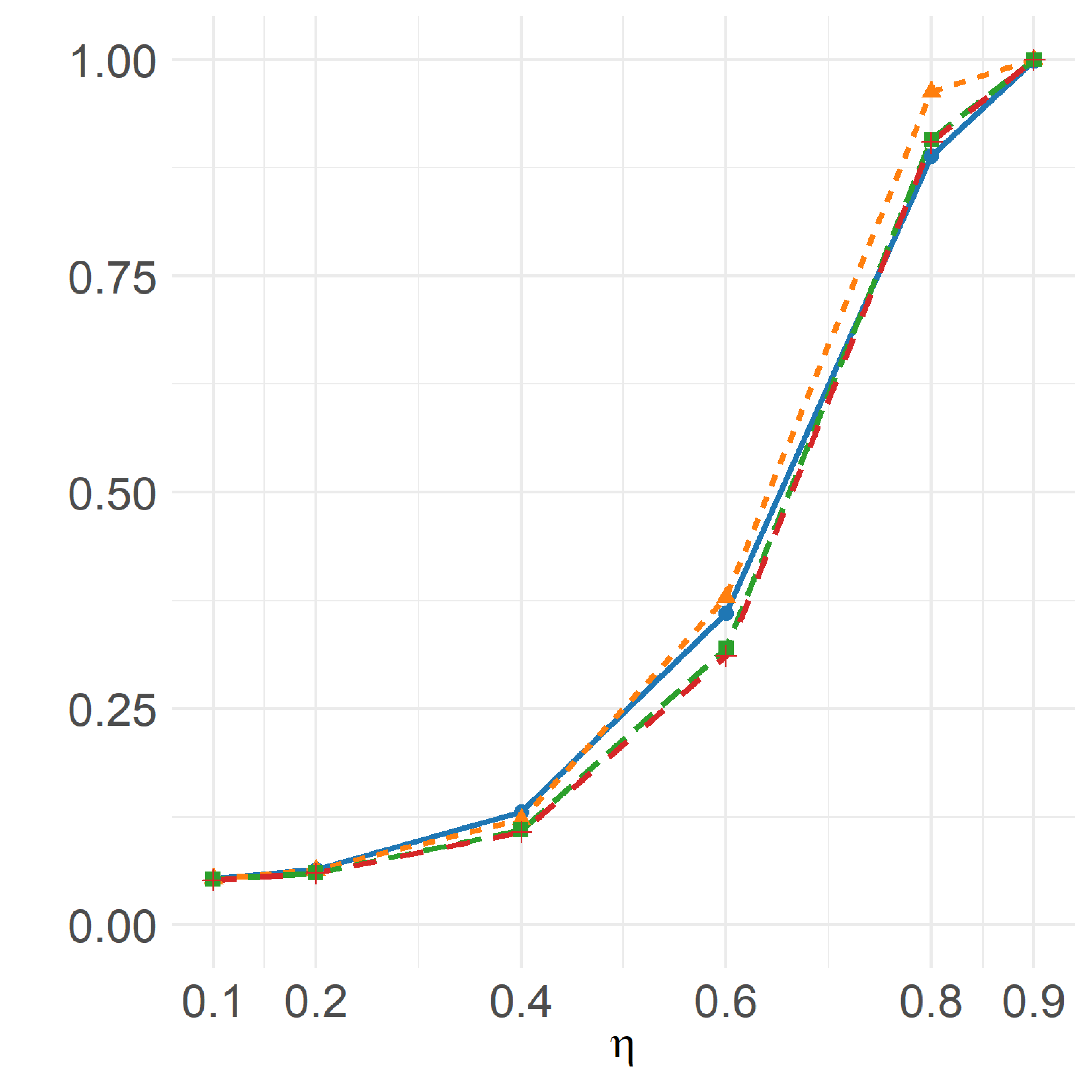}
		\end{minipage}
	}
	\subfigure{
		\begin{minipage}[b]{.3\linewidth}
			\centering
			\includegraphics[scale=0.0735]{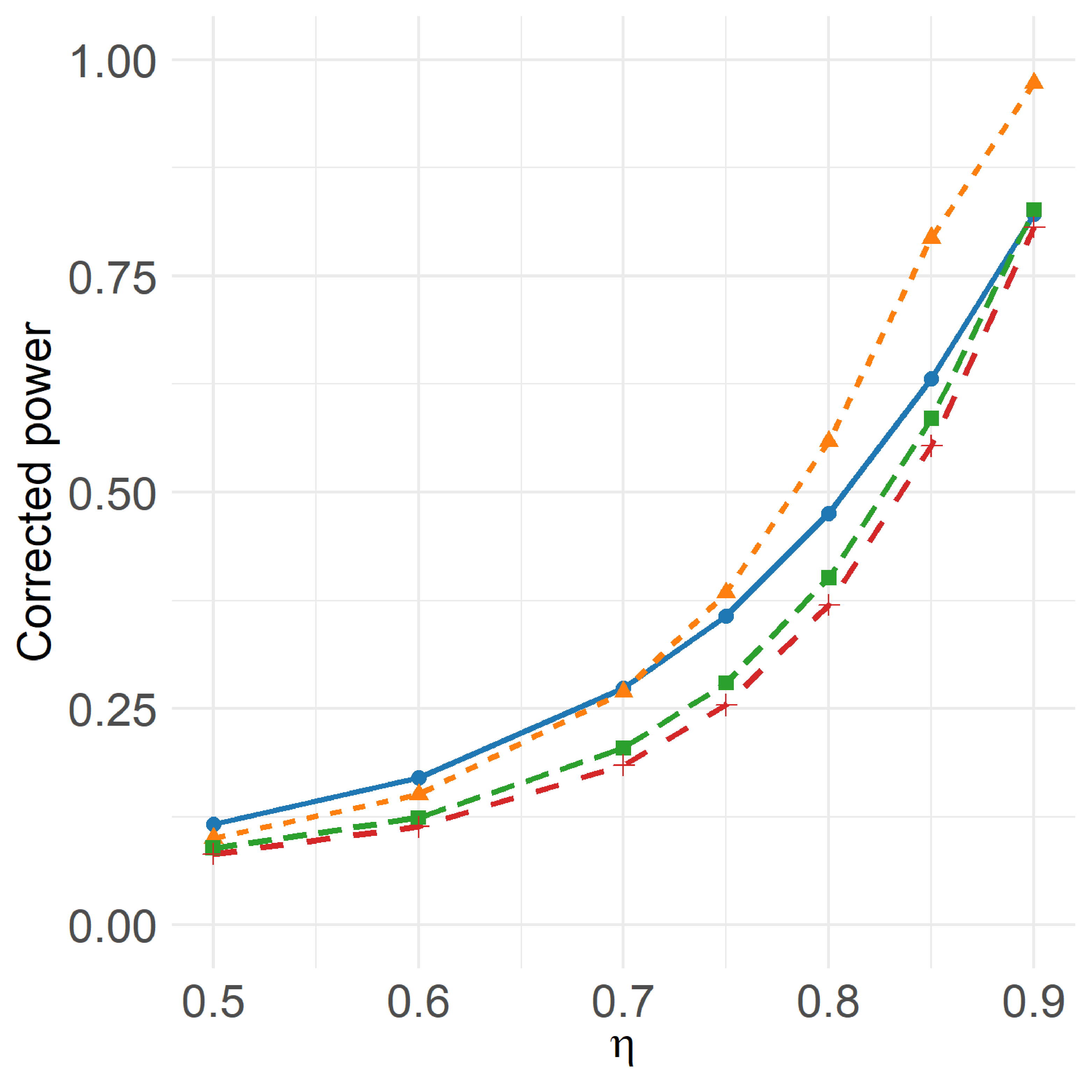}
		\end{minipage}
	}
	\subfigure{
		\begin{minipage}[b]{.3\linewidth}
			\centering
			\includegraphics[scale=0.0735]{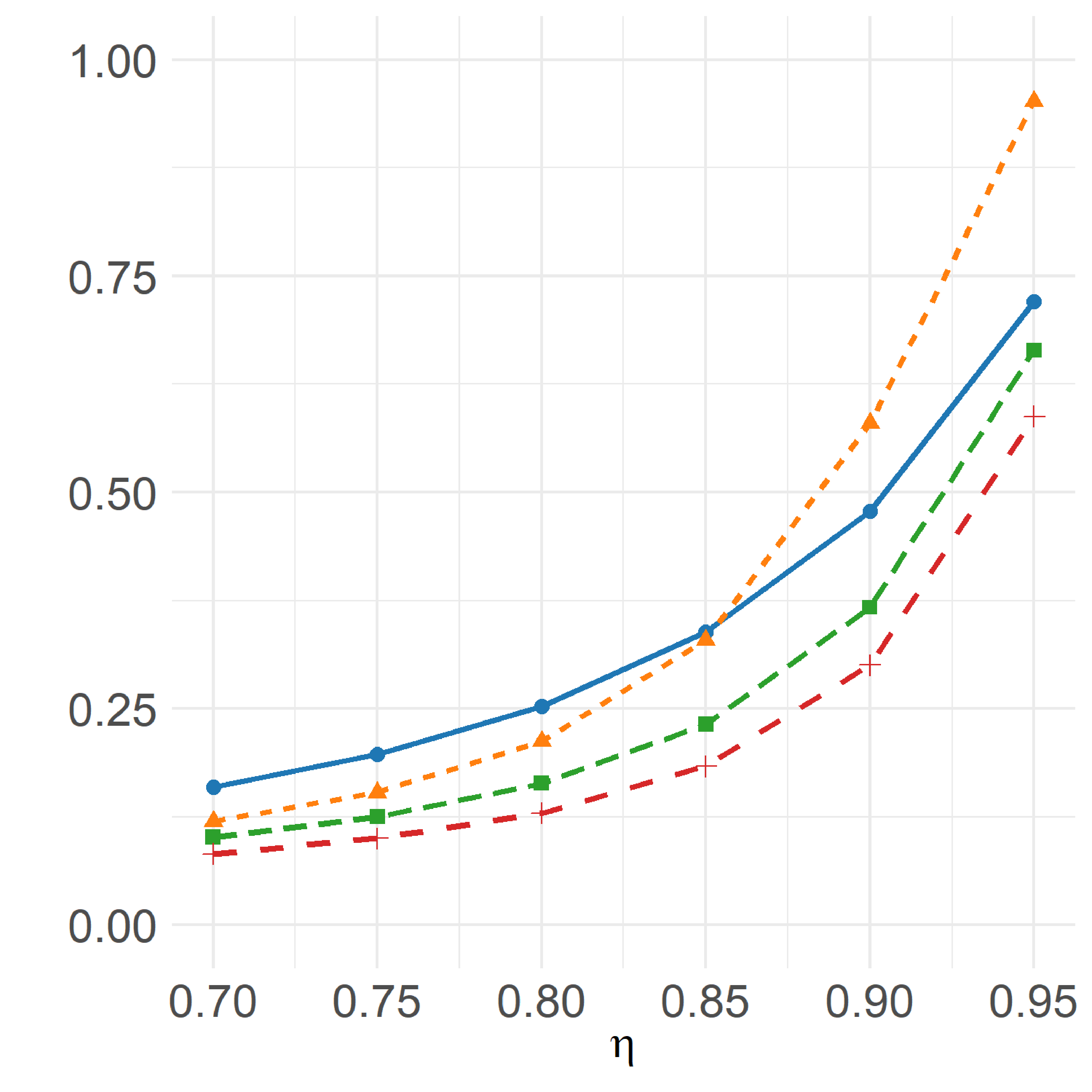}
		\end{minipage}
	}
	\subfigure{
		\begin{minipage}[b]{.3\linewidth}
			\centering
			\includegraphics[scale=0.0735]{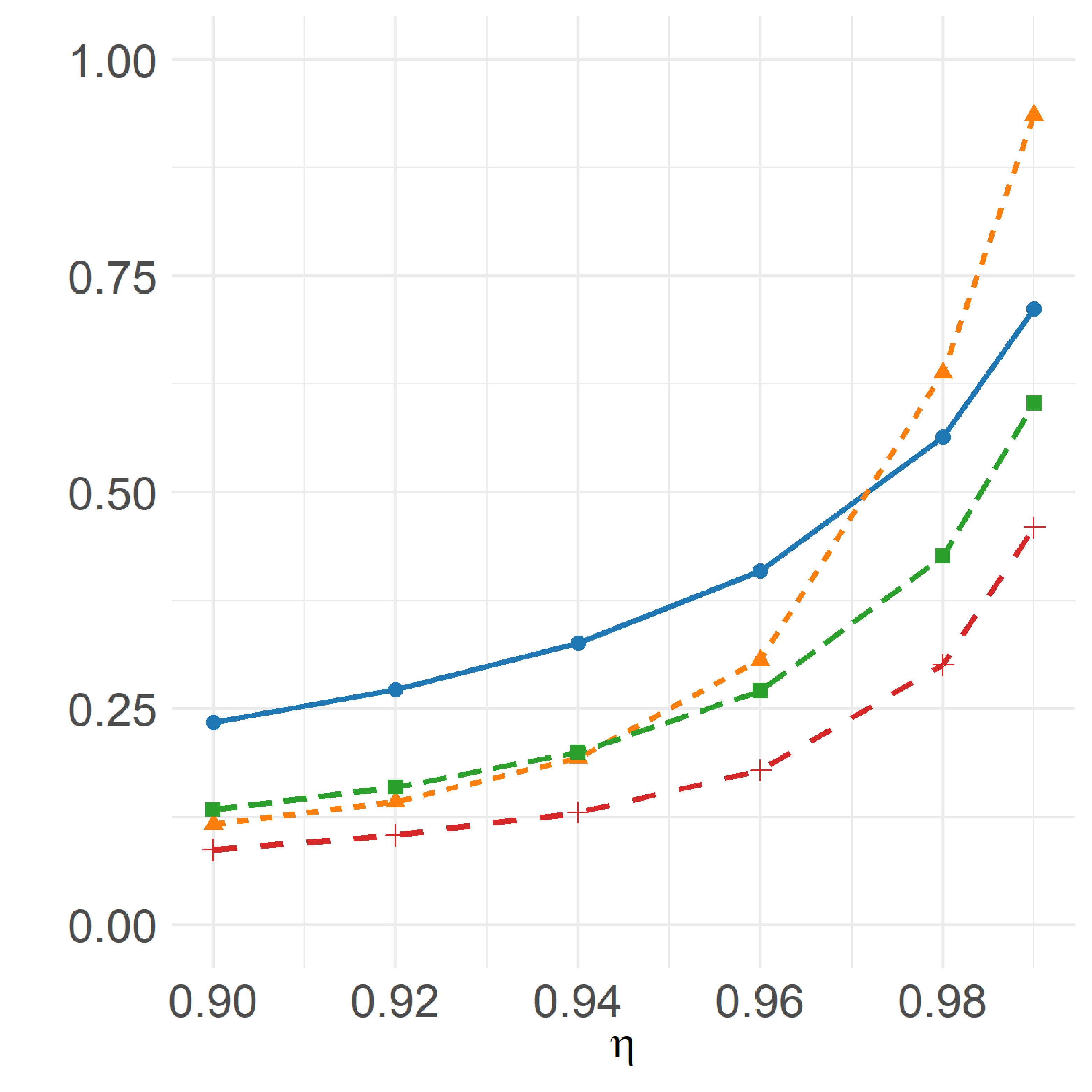}
		\end{minipage}
	}
	\caption{Hypothesis (V). Empirical local power of four tests with different values of $\eta$ and $p/n$. The solid, dashed, longdashed, and dotdashed curves are the empirical power functions of the central limit theorem test, directional test and two Skovgaard’s modifications \cite{skovgaard:2001}, respectively.  The extreme alternative setting is given in Section \ref{S4:simulation (III)-(VI)}. The six plots correspond to $p/n \in \{0.05, 0.1, 0.3, 0.5, 0.7, 0.9\}$, starting from top left and proceeding by row. }  
	\label{SMfig:power extreme case6}
\end{figure}

\begin{figure}[H]
	\centering
	\captionsetup{font=footnotesize}
	\subfigure{
		\begin{minipage}[b]{.3\linewidth}
			\centering
			\includegraphics[scale=0.0735]{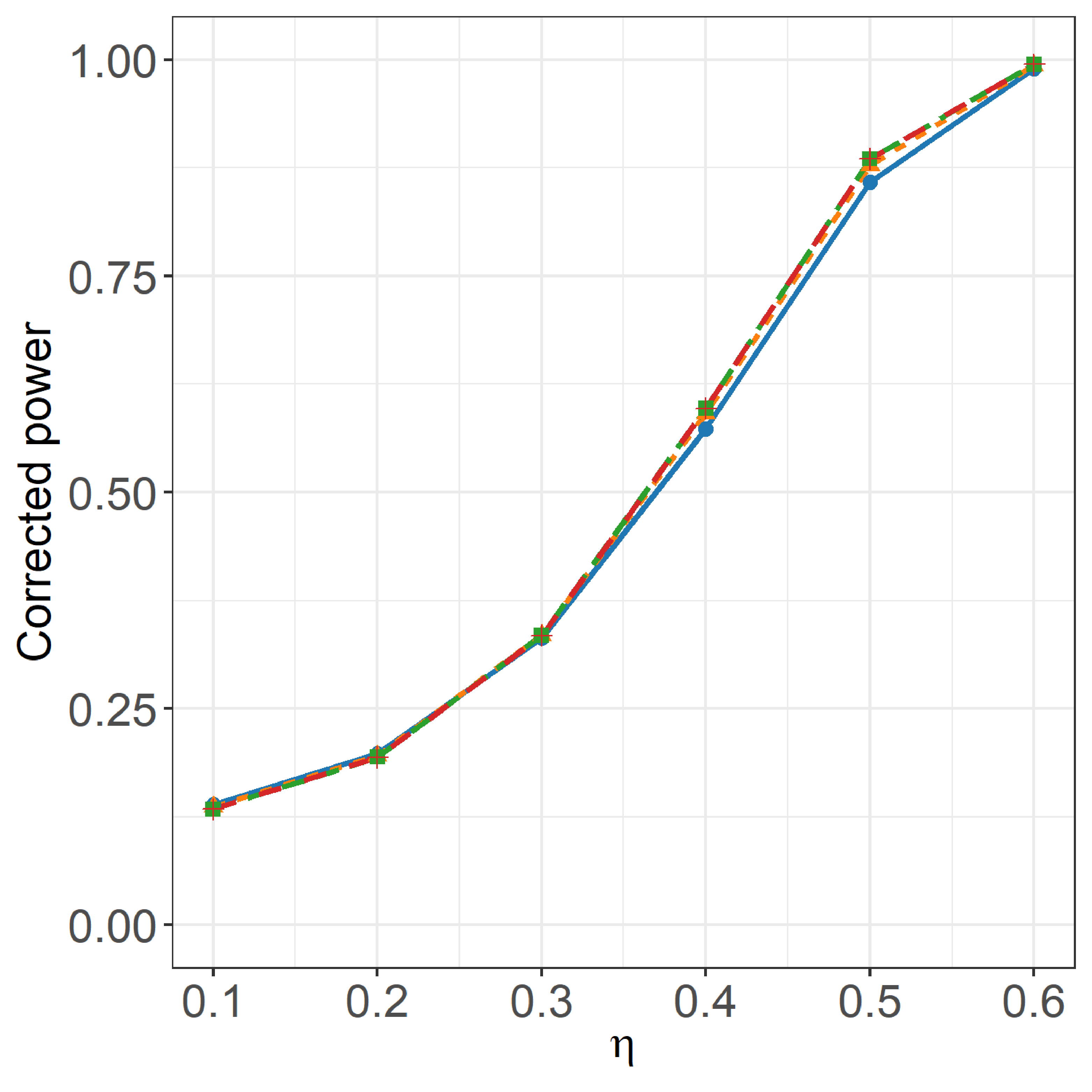}
		\end{minipage}
	}
	\subfigure{
		\begin{minipage}[b]{.3\linewidth}
			\centering
			\includegraphics[scale=0.0735]{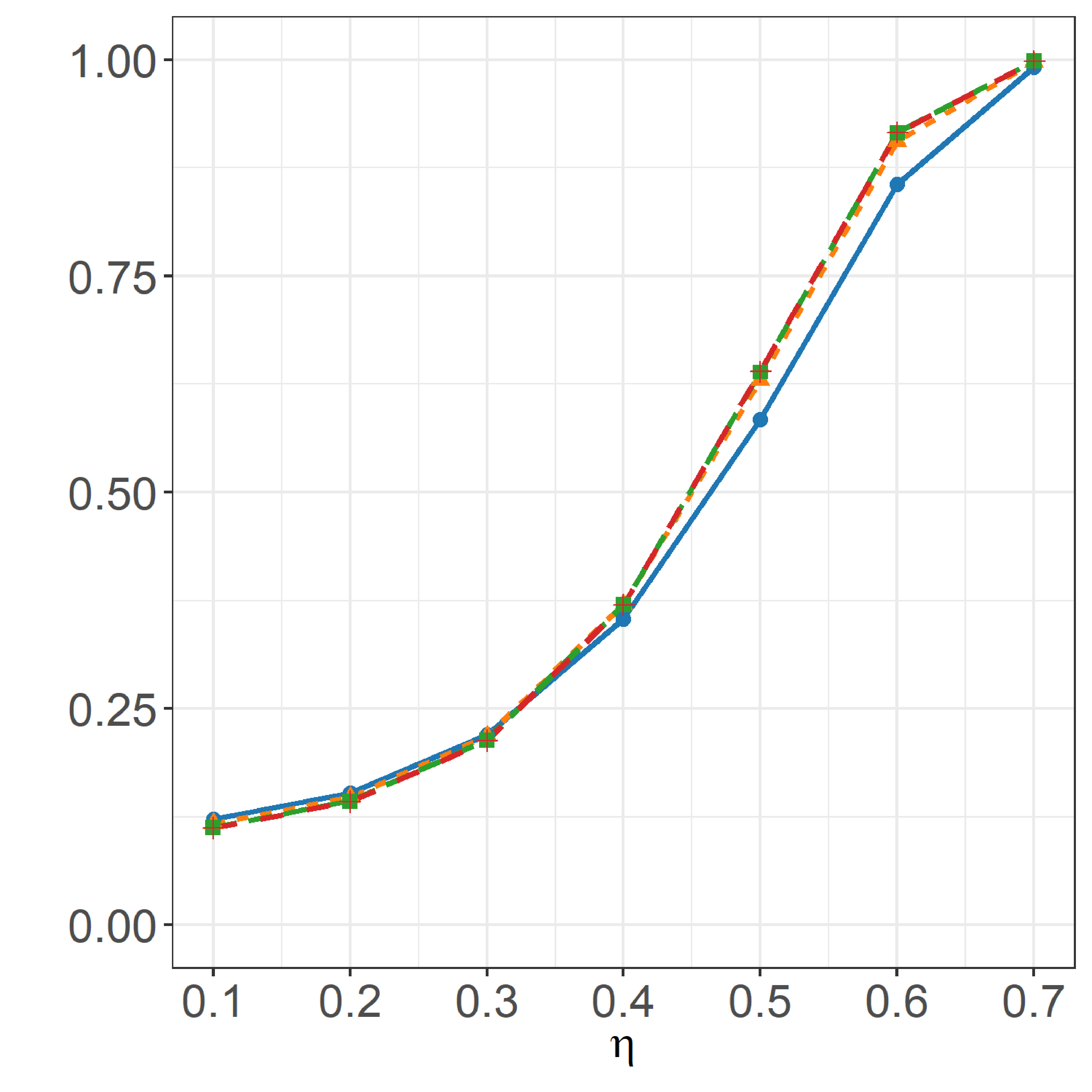}
		\end{minipage}
	}
	\subfigure{
		\begin{minipage}[b]{.3\linewidth}
			\centering
			\includegraphics[scale=0.0735]{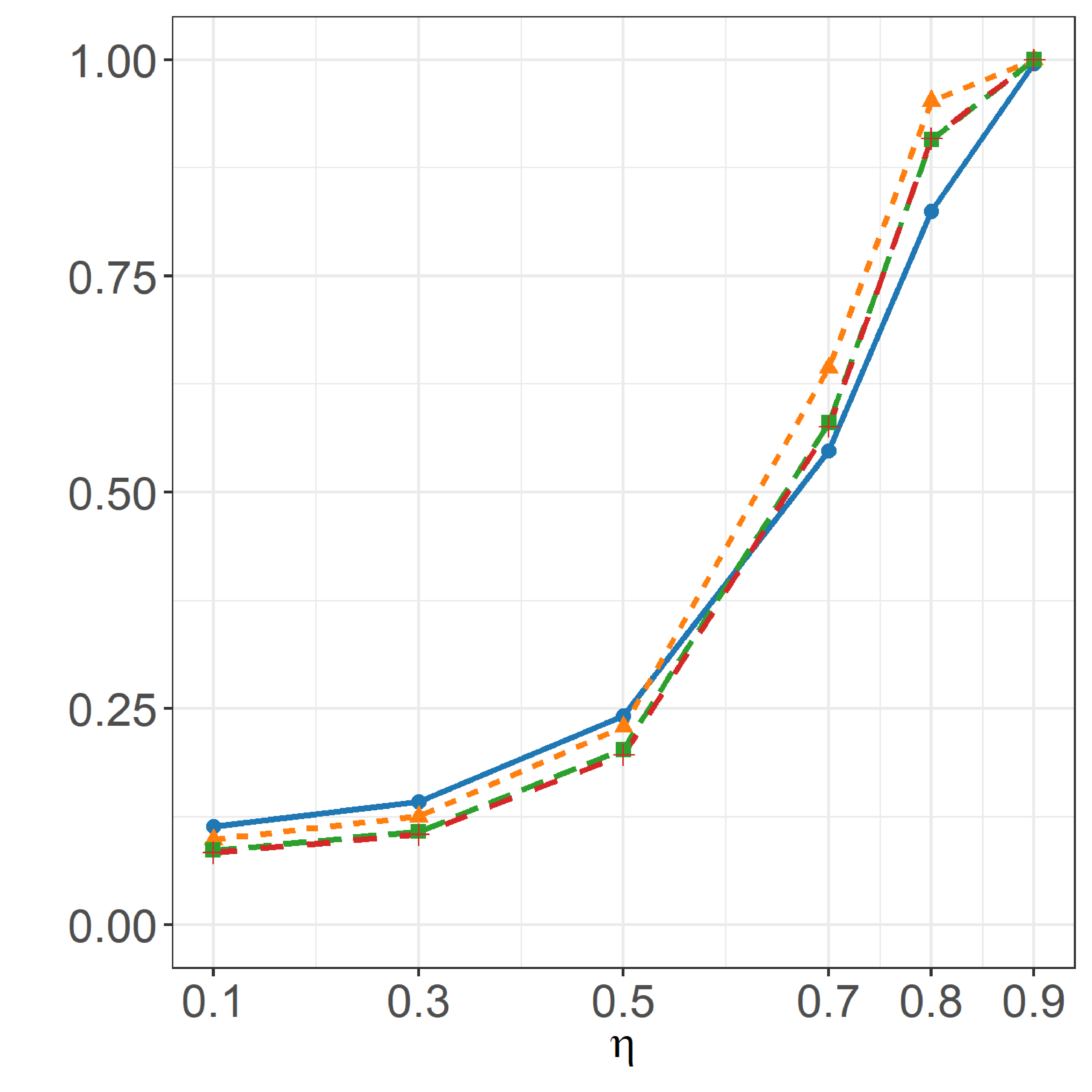}
		\end{minipage}
	}
	\subfigure{
		\begin{minipage}[b]{.3\linewidth}
			\centering
			\includegraphics[scale=0.0735]{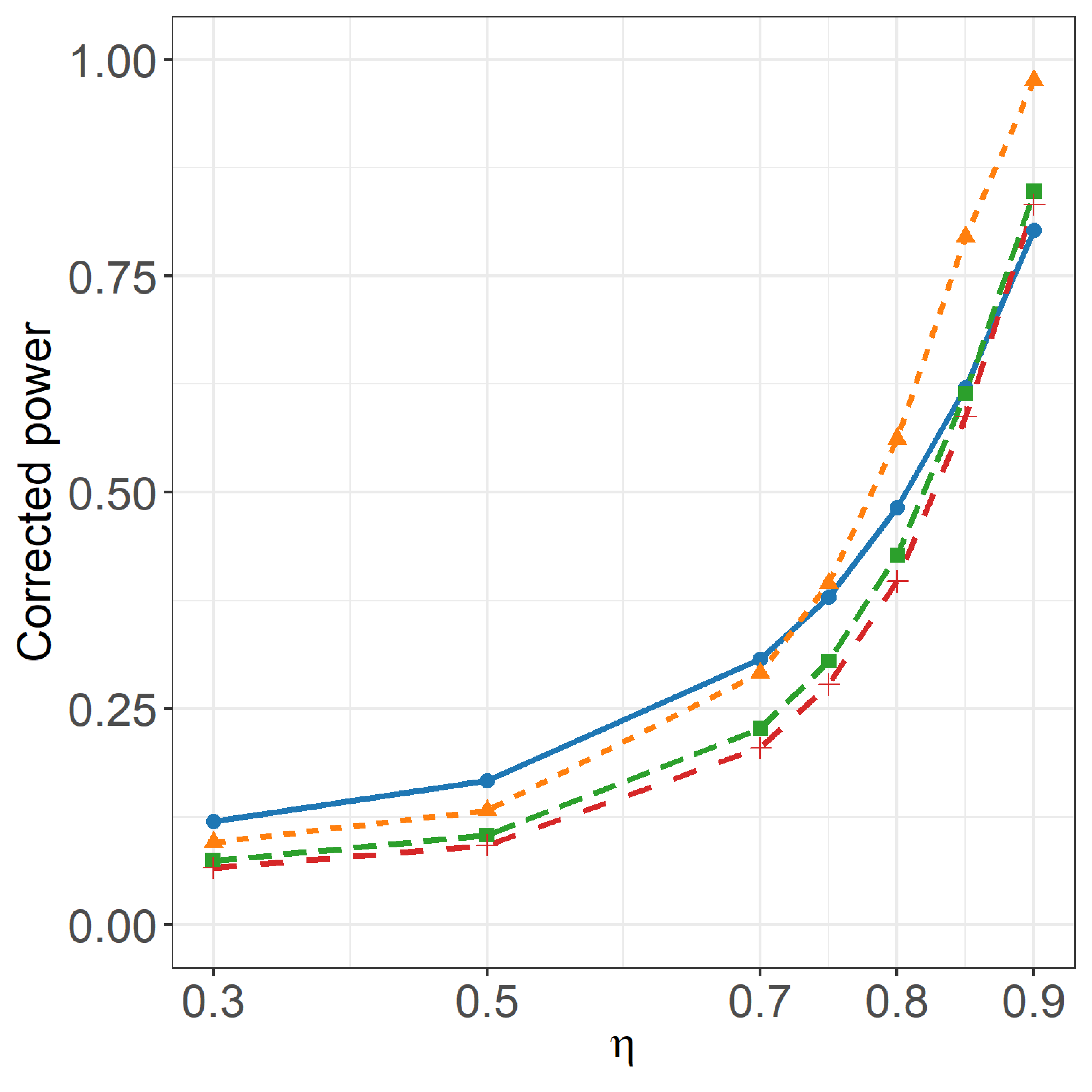}
		\end{minipage}
	}
	\subfigure{
		\begin{minipage}[b]{.3\linewidth}
			\centering
			\includegraphics[scale=0.0735]{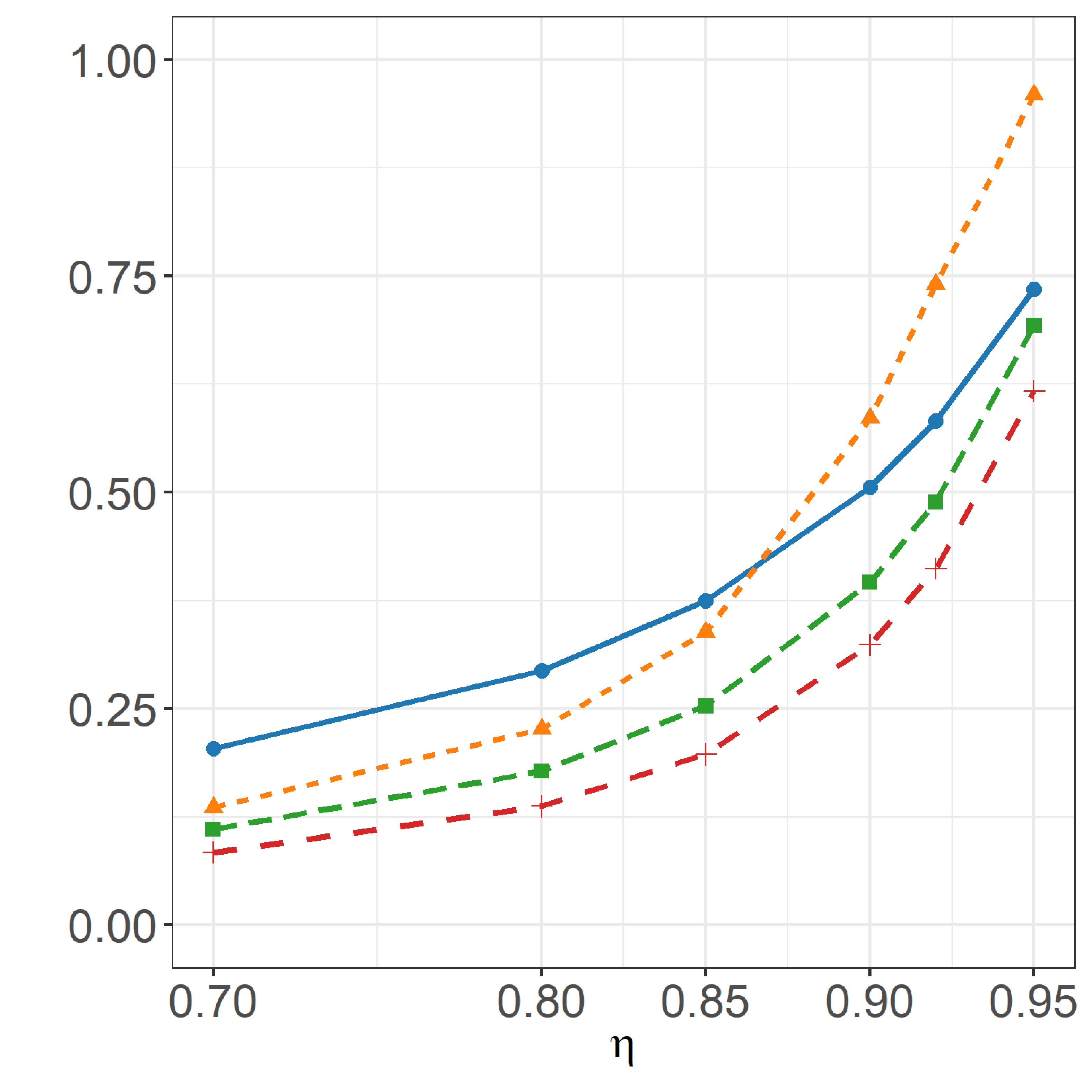}
		\end{minipage}
	}
	\subfigure{
		\begin{minipage}[b]{.3\linewidth}
			\centering
			\includegraphics[scale=0.0735]{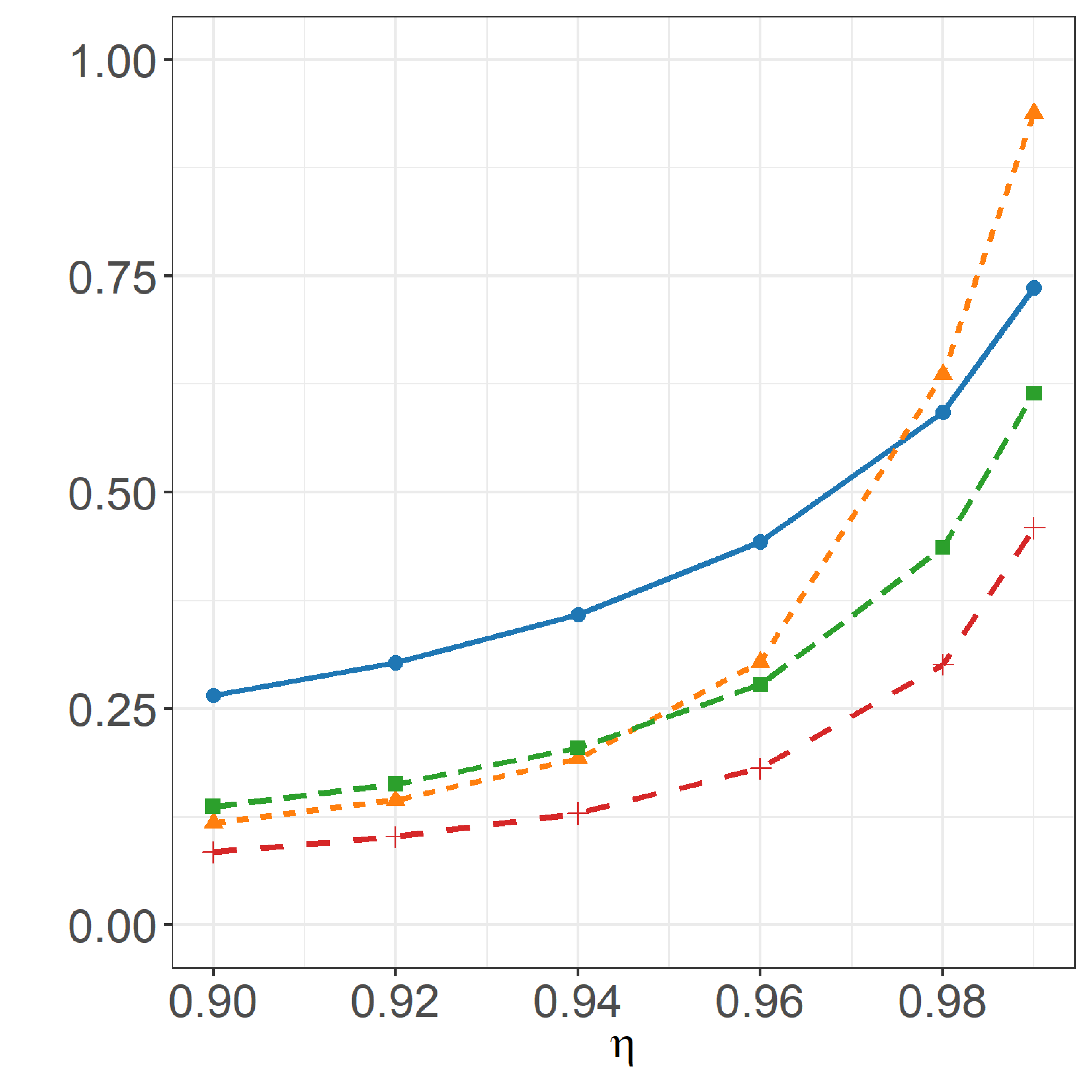}
		\end{minipage}
	}
	\caption{Hypothesis (VI). Empirical local power of four tests with different values of $\eta$ and $p/n$. The solid, dashed, longdashed, and dotdashed curves are the empirical power functions of the central limit theorem test, directional test and two Skovgaard’s modifications \cite{skovgaard:2001}, respectively.   The extreme alternative setting is given in Section \ref{S4:simulation (III)-(VI)}. The six plots correspond to $p/n \in \{0.05, 0.1, 0.3, 0.5, 0.7, 0.9\}$, starting from top left and proceeding by row.}  
	\label{SMfig:power extreme case5}
\end{figure}

We also examine the performance in terms of corrected power for hypotheses (III)--(VI).   Corrected power is obtained based on the corrected Type I error.  Three possible  choices for $\mu$ and $\Lambda^{-1}$ under the alternative hypothesis are considered. The central limit theorem test,  $W$ and $W_{BC}$  have the same corrected power since they use the same test statistic $W$, and result in different cutoff values for the corrected Type I error. 
The left column of Figure \ref{SMfigure:power} displays the corrected power in the alternative setting as  in \cite{jiang:2015}. The power of the directional test is comparable to the central limit theorem test when $p$ is moderate in hypothesis (IV).
The middle and right columns of Figure \ref{SMfigure:power} show the empirical power  for the local alternative hypothesis with various $\delta$ and fixed ratio $p/n=0.3$ and $0.7$. The power performance of the directional test is very close to that of the central limit theorem test, which is better than that of $W^*$ and $W^{**}$ with the ratio $p/n=0.3$. 
Figures \ref{SMfig:power extreme case1}-\ref{SMfig:power extreme case5} summarize the power examined in the alternative setting (3). The directional test has the best power performance, even in the setting when $p/n$ equals 0.9. 



\section{Numerical results based on  sample size $n=p+a$}\label{S5:n=p+2,based on 100,000 replications}

This section is devoted to examining the Monte Carlo simulations referred to the six hypotheses problems under two extreme settings where the dimension $p$ is very close to the sample size $n$. Specifically, we set $n=100$ and $p=95, 98$. The empirical  null distribution of the statistics and of the  $p$-value, and empirical Type I error at the nominal level $\alpha=0.05$ are evaluated based on $100,000$ replications.

\begin{figure}[H]
	\centering
	\captionsetup{font=footnotesize}
	\includegraphics[scale=0.08]{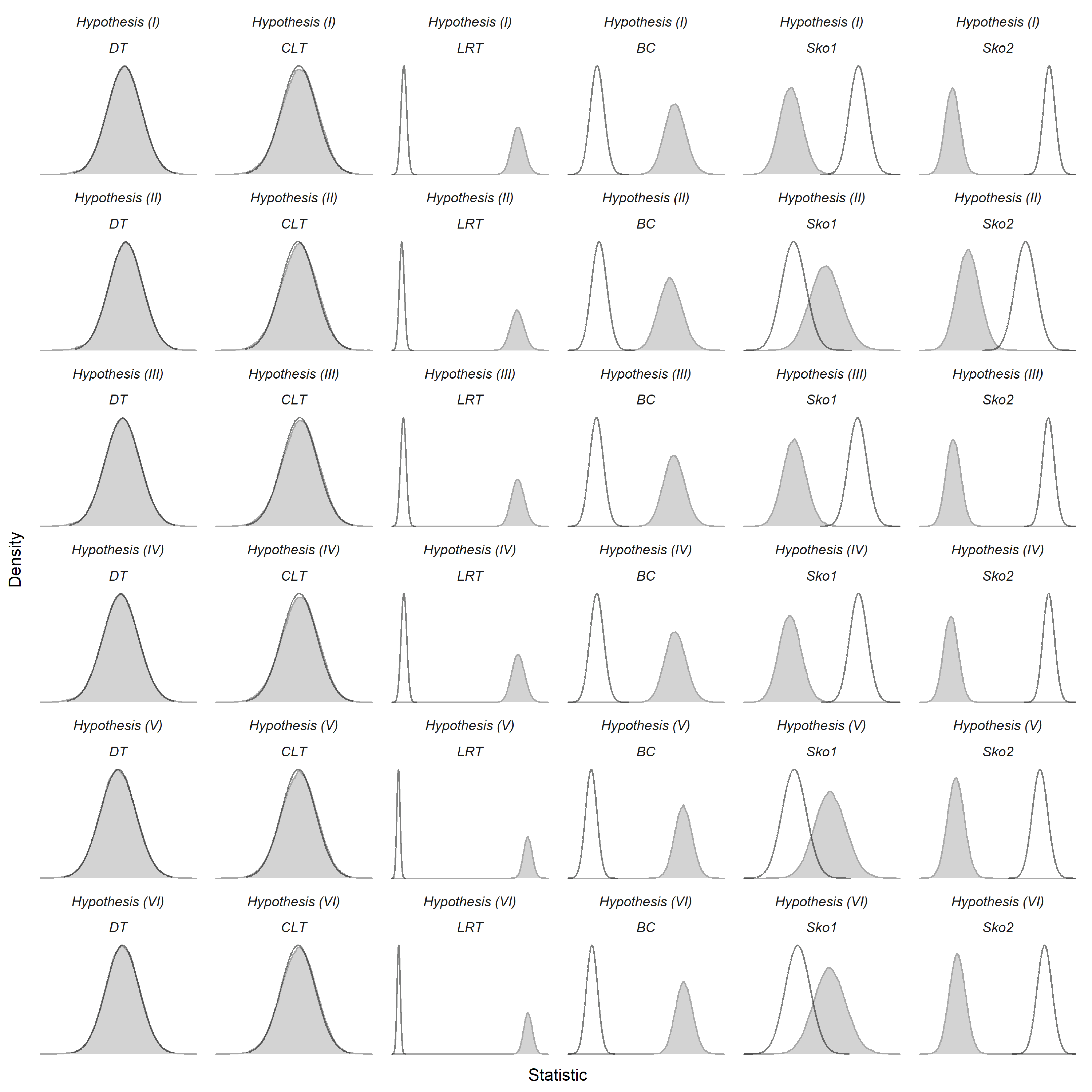}
	\caption{Comparison between   empirical (gray) and theoretical (black) null  distributions of the directional test (DT), the central limit theorem test (CLT), log-likelihood ratio test (LRT), Bartlett correction (BC),  and two Skovgaard's modifications \cite{skovgaard:2001} (Sko1 and Sko2, respectively),  with  $n=100$ and $p=95$.}
	\label{SMsimulation:density null p95}
\end{figure}

\begin{figure}[H]
	\centering
	\captionsetup{font=footnotesize}
	\includegraphics[scale=0.08]{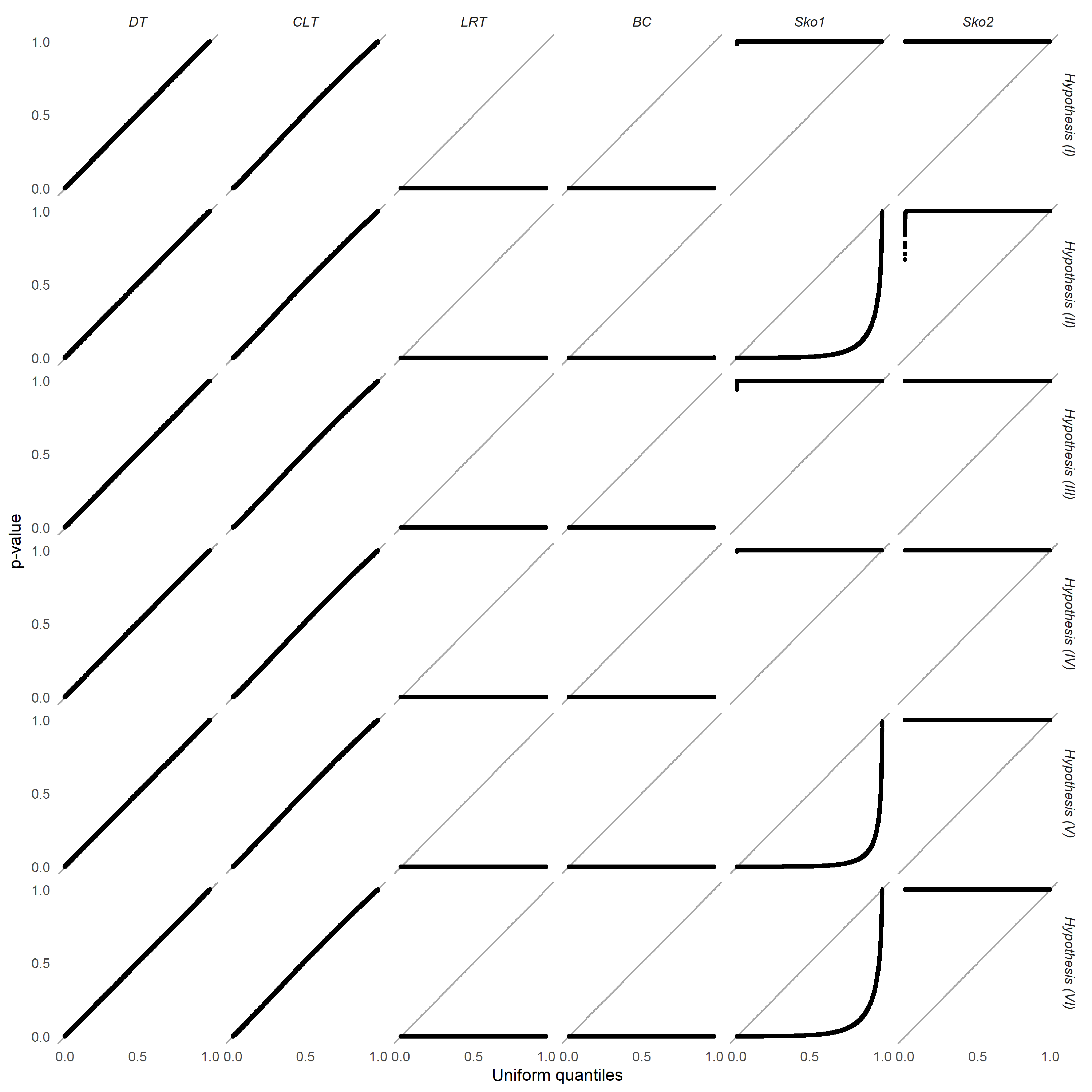}
	\caption{ Empirical null distribution of $p$-values for the directional test (DT), the central limit theorem test (CLT), log-likelihood ratio test (LRT), Bartlett correction (BC),  and two Skovgaard's modifications \cite{skovgaard:2001} (Sko1 and Sko2, respectively)  with $n=100$ and $p=95$, compared with the  $U(0,1)$ given by the gray diagonal}
	\label{SMsimulation:pvalue p95}
\end{figure}

\begin{figure}[H]
	\centering
	\captionsetup{font=footnotesize}
	\includegraphics[scale=0.08]{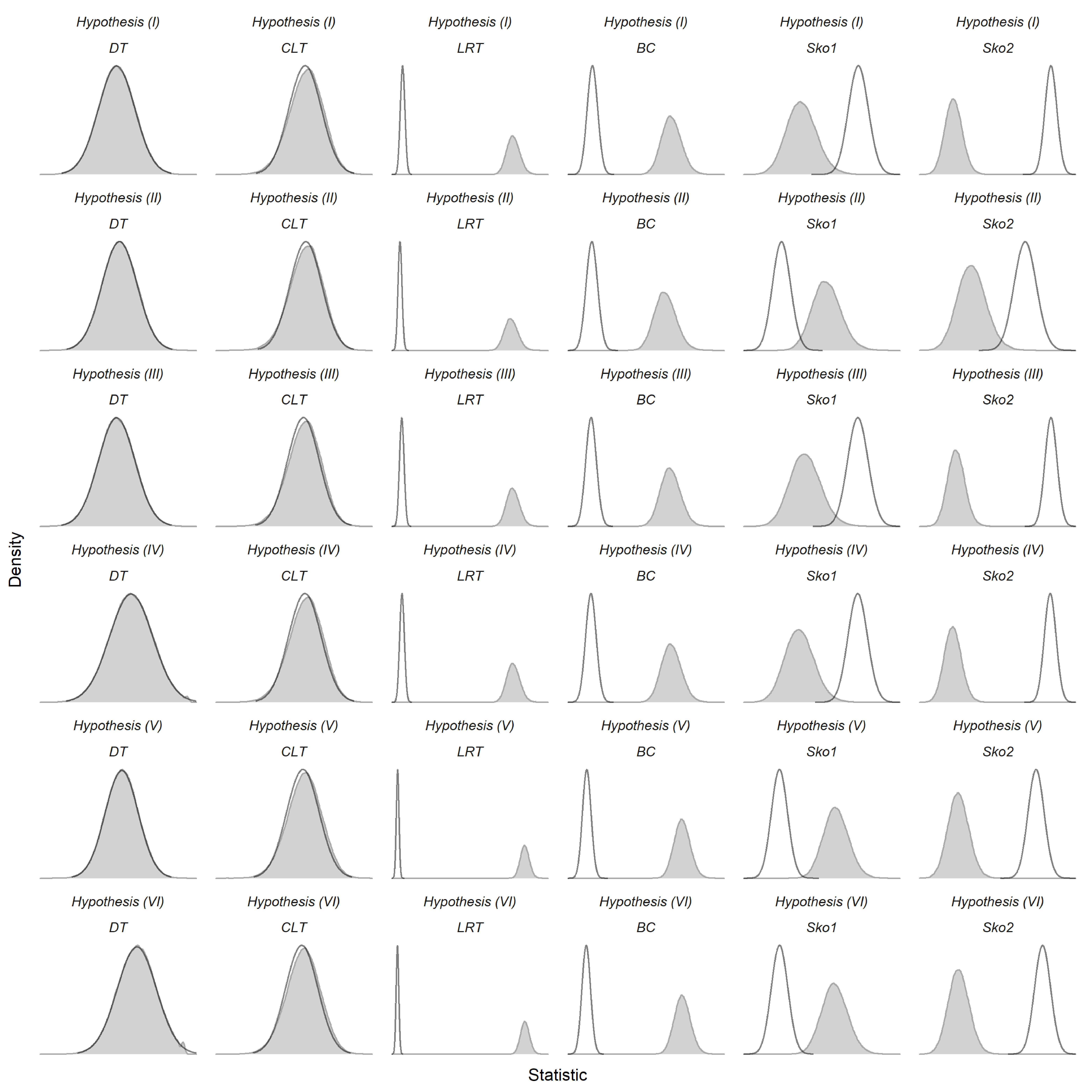}
	\caption{  Comparison between   empirical (gray) and theoretical (black) null  distributions of the directional test (DT), the central limit theorem test (CLT), log-likelihood ratio test (LRT), Bartlett correction (BC),  and two Skovgaard's modifications \cite{skovgaard:2001} (Sko1 and Sko2, respectively)  with  $n=100$ and $p=98$.}
	\label{SMsimulation:density null p98}
\end{figure}

\begin{figure}[H]
	\centering
	\captionsetup{font=footnotesize}
	\includegraphics[scale=0.08]{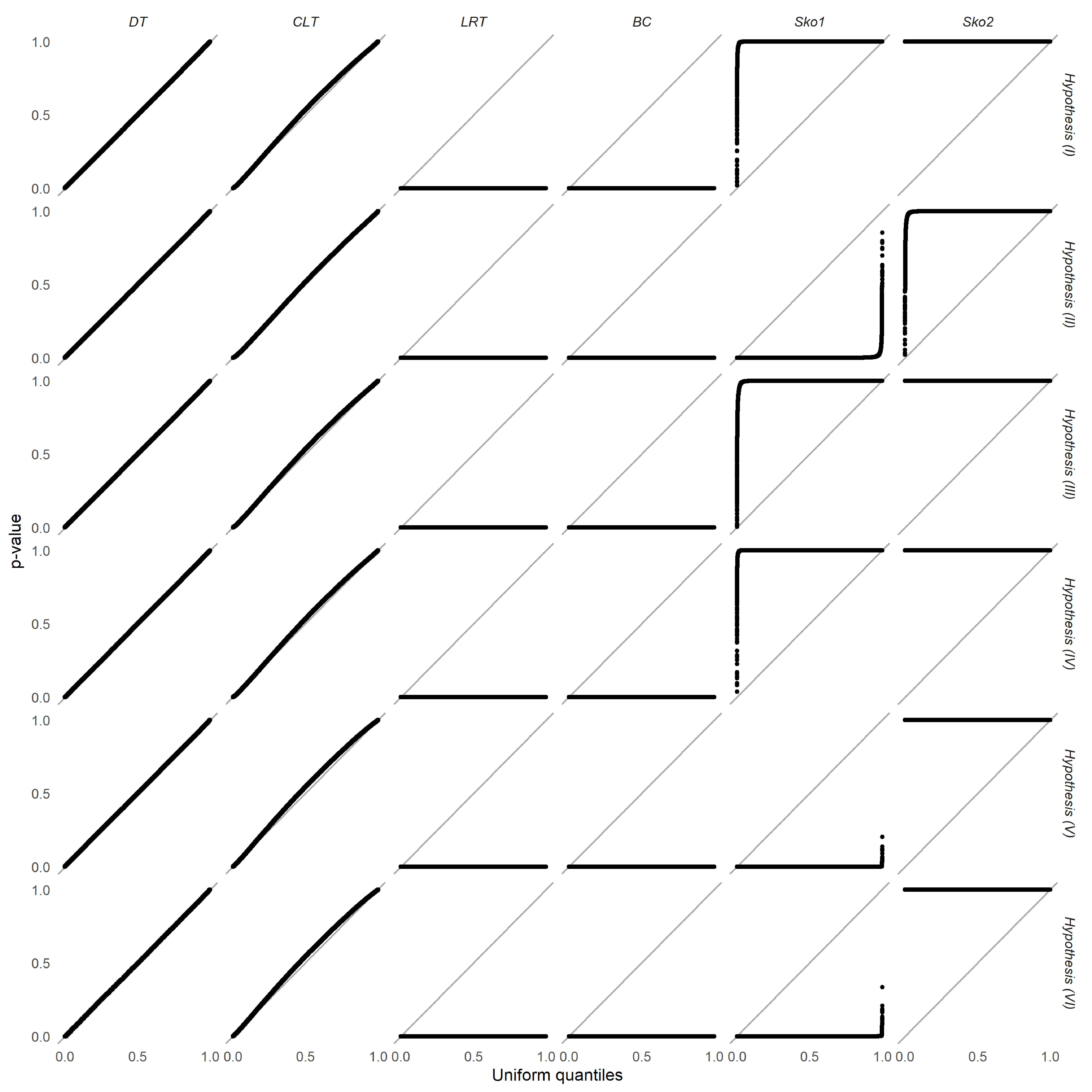}
	\caption{ 
		Empirical null  distribution of $p$-values for  the directional test (DT), the central limit theorem test (CLT), log-likelihood ratio test (LRT), Bartlett correction (BC),  and two Skovgaard's modifications \cite{skovgaard:2001} (Sko1 and Sko2, respectively) with $n=100$ and $p=98$, compared
		with the $U(0,1)$ given by the gray diagonal.
	}
	\label{SMsimulation:pvalue p98}
\end{figure}

Figures \ref{SMsimulation:density null p95}--\ref{SMsimulation:pvalue p98} show that  the directional test and central limit theorem test are superior to  the tests $W$, $W_{BC}$, $W^*$ and $W^{**}$.  Tables \ref{SMtable type I normal: p95} and \ref{SMtable type I normal: p98} report the estimated and corrected Type I error for $\alpha=0.05$.   The Type I errors of  the directional test and the central limit theorem test are very close to the nominal level, while the chi-square approximation for $W$, $W_{BC}$, $W^*$ and $W^{**}$ completely fails. In the most extreme scenario, when $p=98$, the directional test is overall more accurate than the central limit theorem test.

\section{Numerical results based on various sample sizes $n$}\label{S6:various sample sizes}

In this section, Monte Carlo simulations referred to the six hypotheses problems are considered with various sample sizes $n$  and fixed ratio $p/n=0.3$. We report the empirical null distribution of the statistics, of $p$-values, the 
empirical Type I errors at $\alpha=0.05$ and power based on $10,000$ replications. 
 The empirical distribution of the $p$-values under the null hypothesis is  instead compared to the  uniform distribution.

\begin{figure}[H]
	\centering
	\captionsetup{font=footnotesize}
	\includegraphics[scale=0.08]{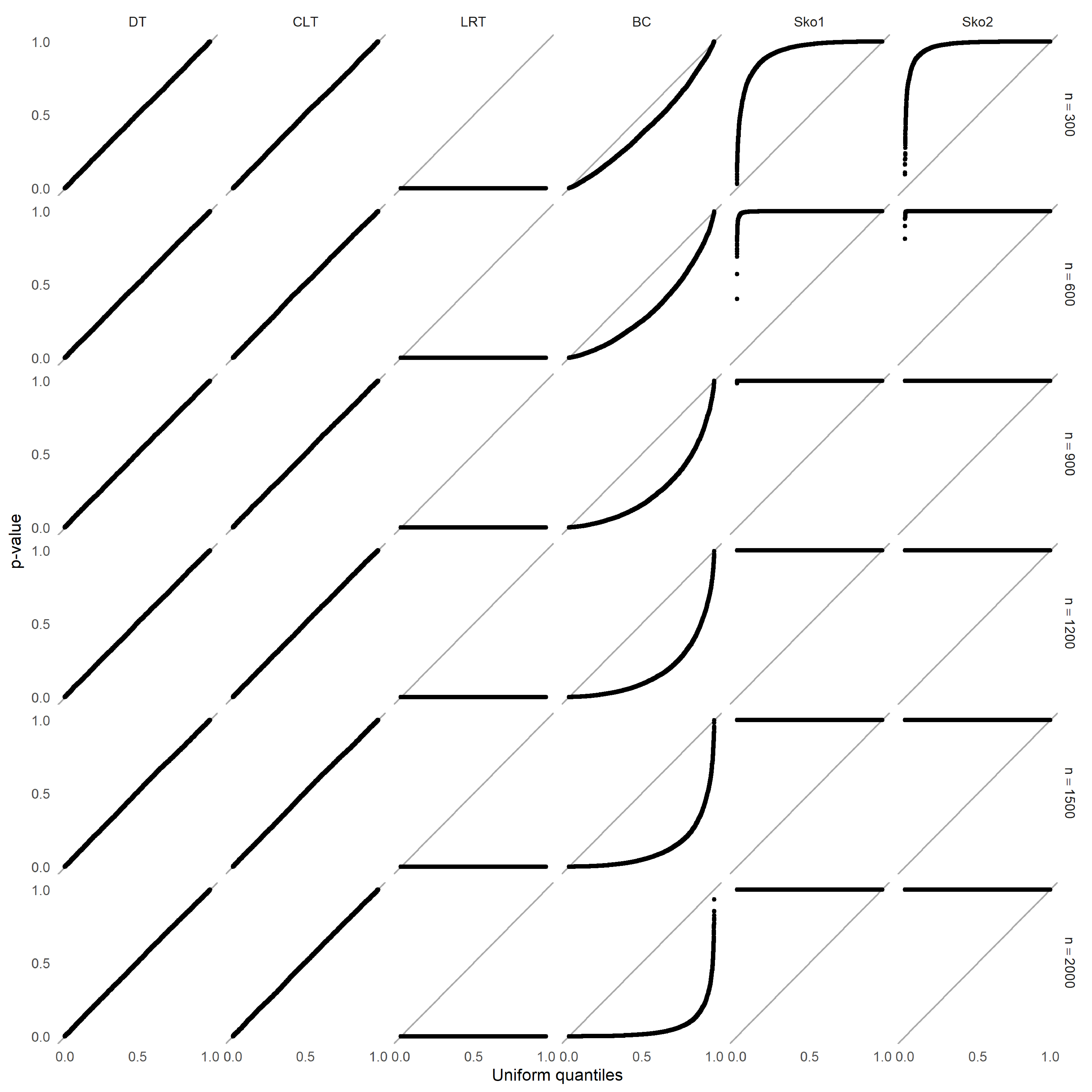}
	\caption{Hypothesis (I). Empirical null  distribution of $p$-values forthe directional test (DT), the central limit theorem test (CLT), log-likelihood ratio test (LRT), Bartlett correction (BC),  and two Skovgaard's modifications \cite{skovgaard:2001} (Sko1 and Sko2, respectively) in settings with various sample sizes $n_i$, $i \in \{1,\dots,k\}$ and fixed ratio $p/n_i=0.3$, compared with the $U(0,1)$ given by the gray diagonal.}
	\label{SMsimulation:pvalue null case4 n}
\end{figure}

\begin{figure}[H]
	\centering
	\captionsetup{font=footnotesize}
	\includegraphics[scale=0.08]{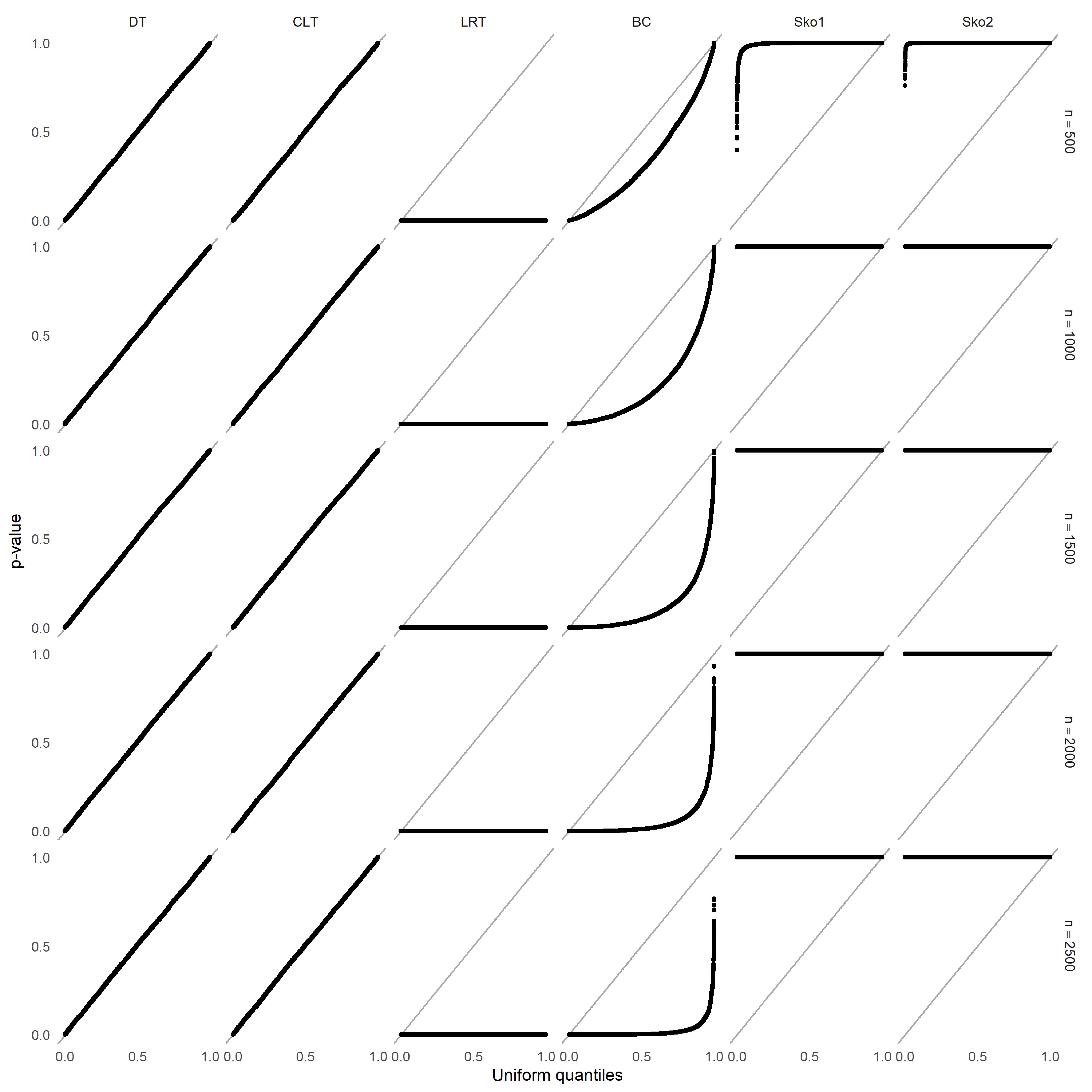}
	\caption{Hypothesis (II). Empirical null  distribution of $p$-values for the directional test (DT), the central limit theorem test (CLT), log-likelihood ratio test (LRT), Bartlett correction (BC),  and two Skovgaard's modifications \cite{skovgaard:2001} (Sko1 and Sko2, respectively) in settings with various sample sizes $n_i$, $i \in \{1,\dots, k\}$ and fixed ratio $p/n_i=0.3$, compared with the $U(0,1)$ given by the gray diagonal.}
	\label{SMsimulation:pvalue null case3 n}
\end{figure}

\begin{figure}[H]
	\centering
	\captionsetup{font=footnotesize}
	\includegraphics[scale=0.08]{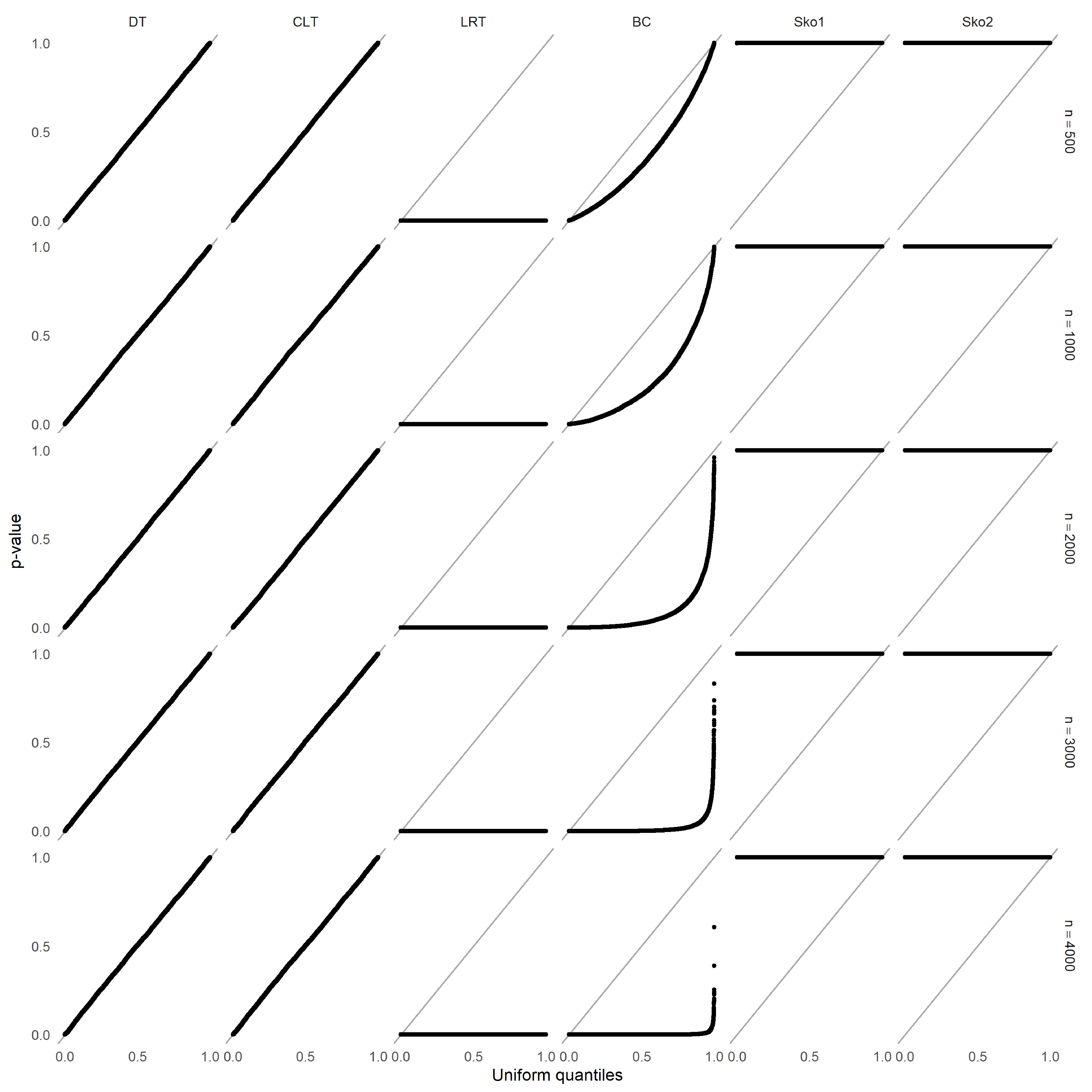}
	\caption{Hypothesis (III).   Empirical null  distribution of $p$-values for the directional test (DT), the central limit theorem test (CLT), log-likelihood ratio test (LRT), Bartlett correction (BC),  and two Skovgaard's modifications \cite{skovgaard:2001} (Sko1 and Sko2, respectively) in settings with various sample sizes $n$ and fixed ratio $p/n=0.3$, compared
		with the $U(0,1)$ given by the gray diagonal.}
	\label{SMsimulation:pvalue null case1 n}
\end{figure}

\begin{figure}[H]
	\centering
	\captionsetup{font=footnotesize}
	\includegraphics[scale=0.08]{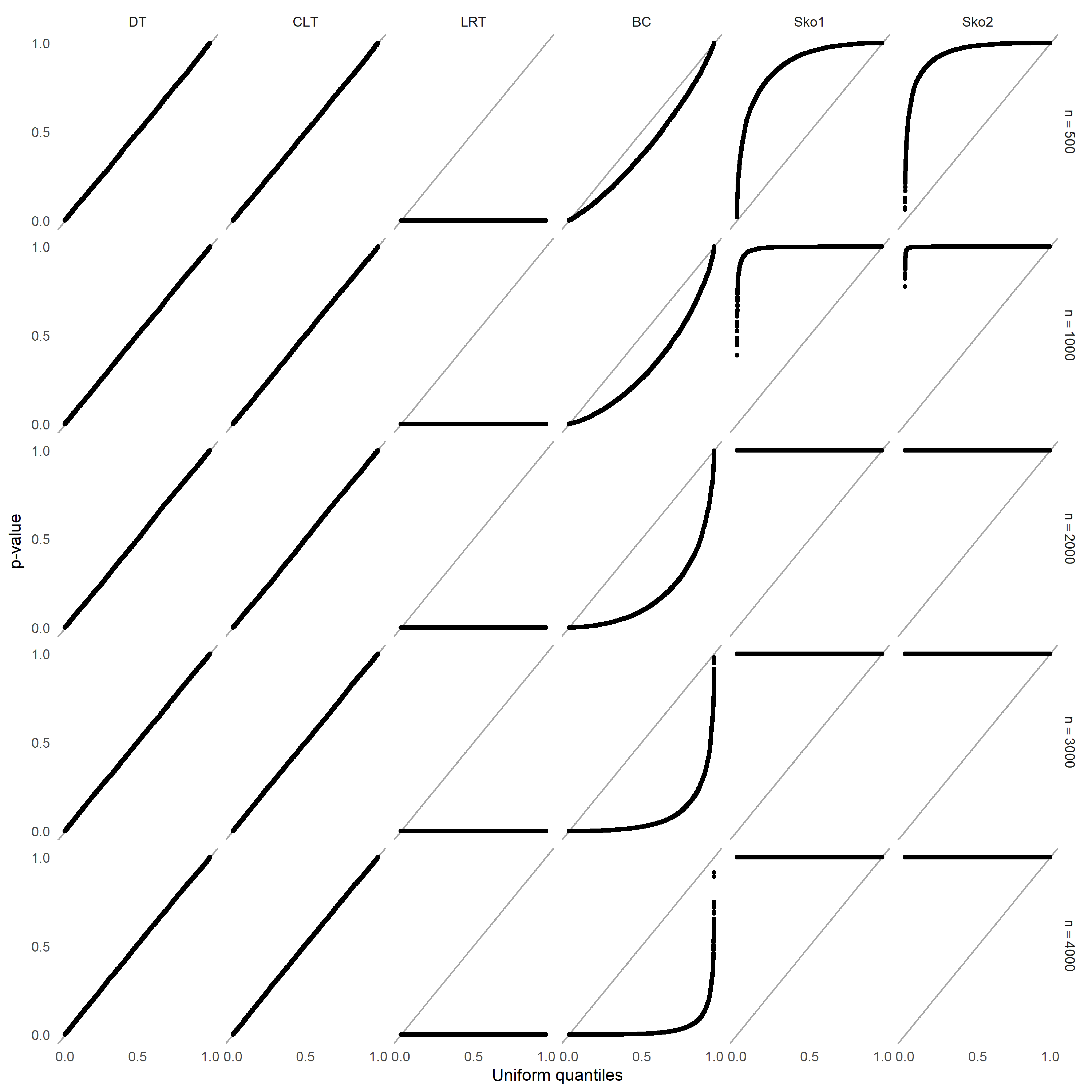}
	\caption{Hypothesis (IV). Empirical null  distribution of $p$-values for the directional test (DT), the central limit theorem test (CLT), log-likelihood ratio test (LRT), Bartlett correction (BC),  and two Skovgaard's modifications \cite{skovgaard:2001} (Sko1 and Sko2, respectively) in settings with various sample sizes $n$ and fixed ratio $p/n=0.3$, compared with the $U(0,1)$ given by the gray diagonal.}
	\label{SMsimulation:pvalue null case2 n}
\end{figure}

\begin{figure}[H]
	\centering
	\captionsetup{font=footnotesize}
	\includegraphics[scale=0.08]{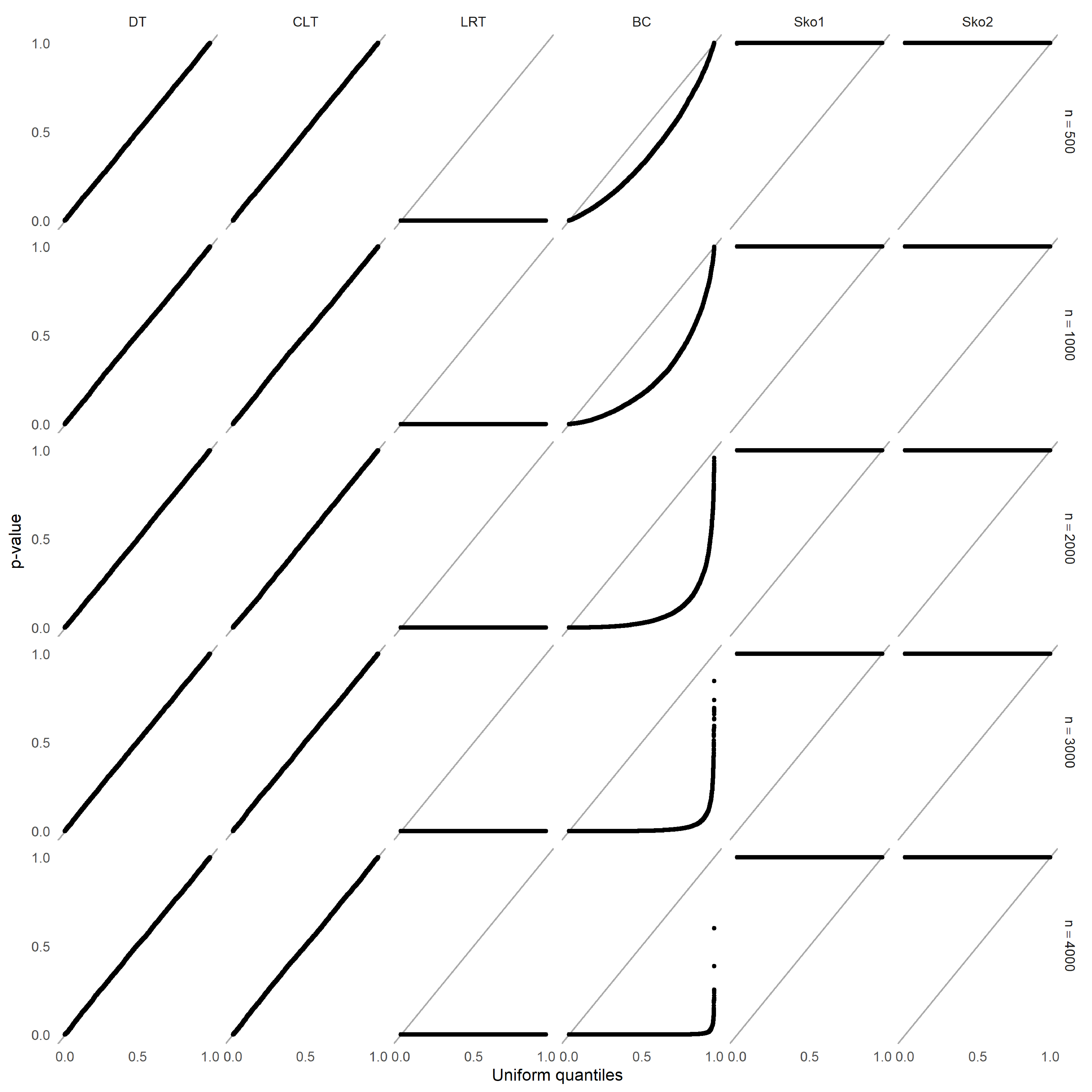}
	\caption{Hypothesis (V). Empirical null  distribution of $p$-values for the directional test (DT), the central limit theorem test (CLT), log-likelihood ratio test (LRT), Bartlett correction (BC),  and two Skovgaard's modifications \cite{skovgaard:2001} (Sko1 and Sko2, respectively) in settings with various sample sizes $n$ and fixed ratio $p/n=0.3$, compared with the $U(0,1)$ given by the gray diagonal.}
	\label{SMsimulation:pvalue null case6 n}
\end{figure}

\begin{figure}[H]
	\centering
	\captionsetup{font=footnotesize}
	\includegraphics[scale=0.08]{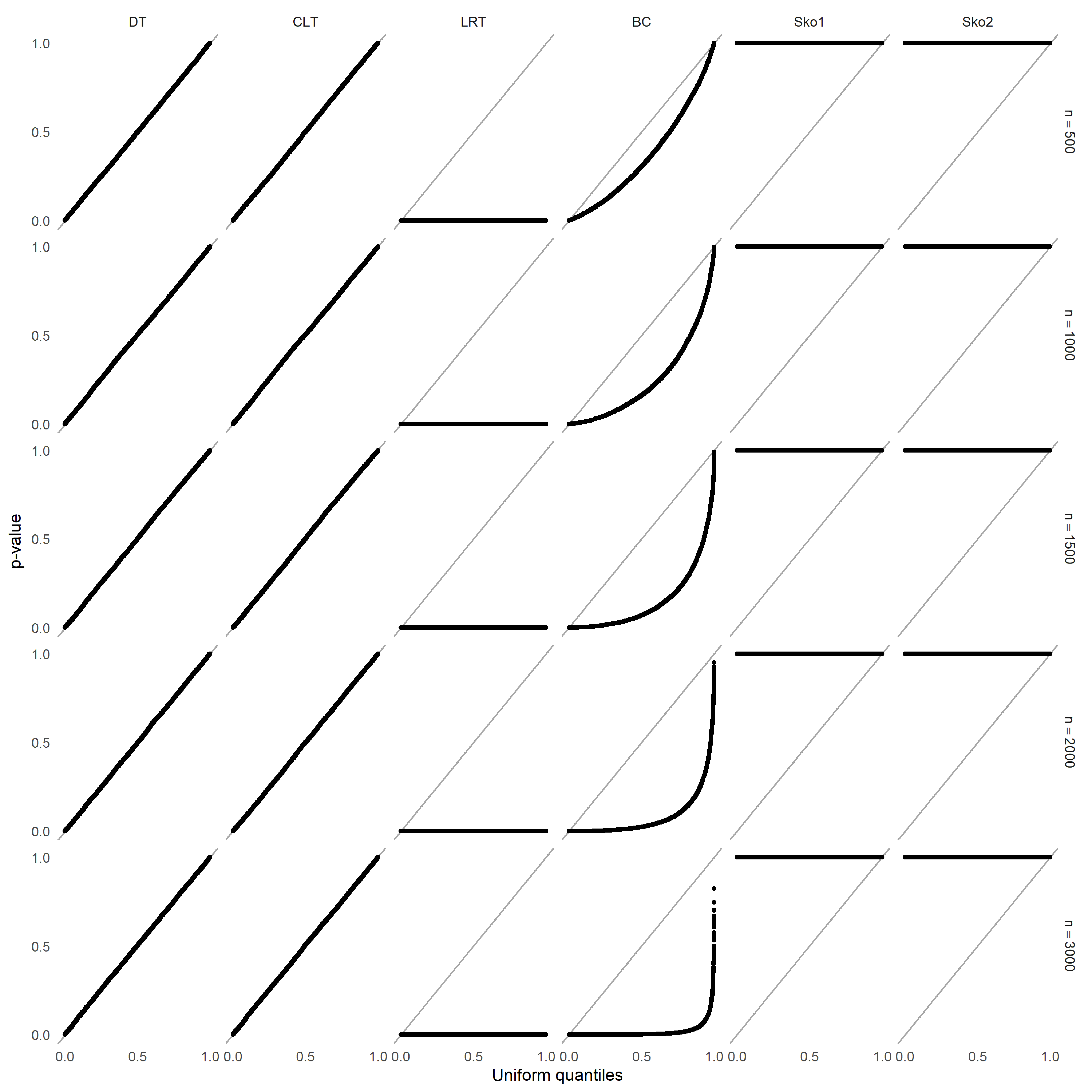}
	\caption{Hypothesis (VI). Empirical null  distribution of $p$-values for the directional test (DT), the central limit theorem test (CLT), log-likelihood ratio test (LRT), Bartlett correction (BC),  and two Skovgaard's modifications \cite{skovgaard:2001} (Sko1 and Sko2, respectively) in settings with various sample sizes $n$ and fixed ratio $p/n=0.3$, compared with the $U(0,1)$ given by the gray diagonal.}
	\label{SMsimulation:pvalue null case5 n}
\end{figure}

Figures \ref{SMsimulation:pvalue null case4 n}--\ref{SMsimulation:pvalue null case5 n} illustrate   the empirical  null distribution of the $p$-values.  The results show that the directional test and central limit theorem test maintain  extreme accuracy with various sample sizes. On the contrary, the chi-square approximations to $W$, $W_{BC}$, $W^*$ and $W^{**}$ are not reliable and get worse and worse as $n$ increases. Tables \ref{SMtable type I normal:case4 n}--\ref{SMtable type I normal:case6 n} report the estimated Type I error and corrected Type I error at the nominal level $\alpha = 0.05$. The directional test and central limit theorem test confirm an excellent performance for all  sample sizes. On the other hand, $W$, $W_{BC}$, $W^*$ and $W^{**}$ completely fail to control the Type I error when $p/n=0.3$.

We also investigate the performance in terms of corrected power for the six hypothesis problems.   Corrected power is obtained based on the corrected Type I error. The central limit theorem test, $W$ and $W_{BC}$  have the same corrected power since they use the same test statistic, and simply result in different cutoff values for the corrected Type I error.  

The alternative hypothesis settings for the six hypotheses problems are as follows, (I): $\Lambda_1^{-1} = \bI_p, \Lambda_2^{-1} = \Lambda_3^{-1} = 1.1\bI_p$; (II): $\mu_1 = 0_p, \mu_2 = \mu_3 = 0.05_p$ and  $\Lambda_1^{-1} = \bI_p, \Lambda_2^{-1} = \Lambda_3^{-1} = 1.05\bI_p$; (III): $\Lambda^{-1} = \diag(1.1, \dots, 1.1, 1,\dots, 1)$, where the number of diagonal entries equal to 1.1 is $\left\lceil {p/2} \right\rceil$; (IV): $\Lambda^{-1} = 0.005 1_p 1_p^T + (1-0.005) \bI_p$; (V):  $\Lambda^{-1} = \left(\lambda^{ij}\right)_{p\times p}$ where $\lambda^{ij} = 1$ for $i=j$, $\lambda^{ij} = 0.02$ for $0<|i-j|\le 3$, and $\lambda^{ij} = 0$ for $|i-j|>3$; (VI): $\mu = (0.02,\dots,0.02,0,\dots,0)$, where the number of diagonal entries equal to 0.02 is $\left\lceil {p/2} \right\rceil$, and the setup for $\Lambda^{-1}$ is as in hypothesis (V); .

\begin{figure}[H]
	\centering
	\captionsetup{font=footnotesize}
	\subfigure{
		\begin{minipage}[b]{.3\linewidth}
			\centering
			\includegraphics[scale=0.0735]{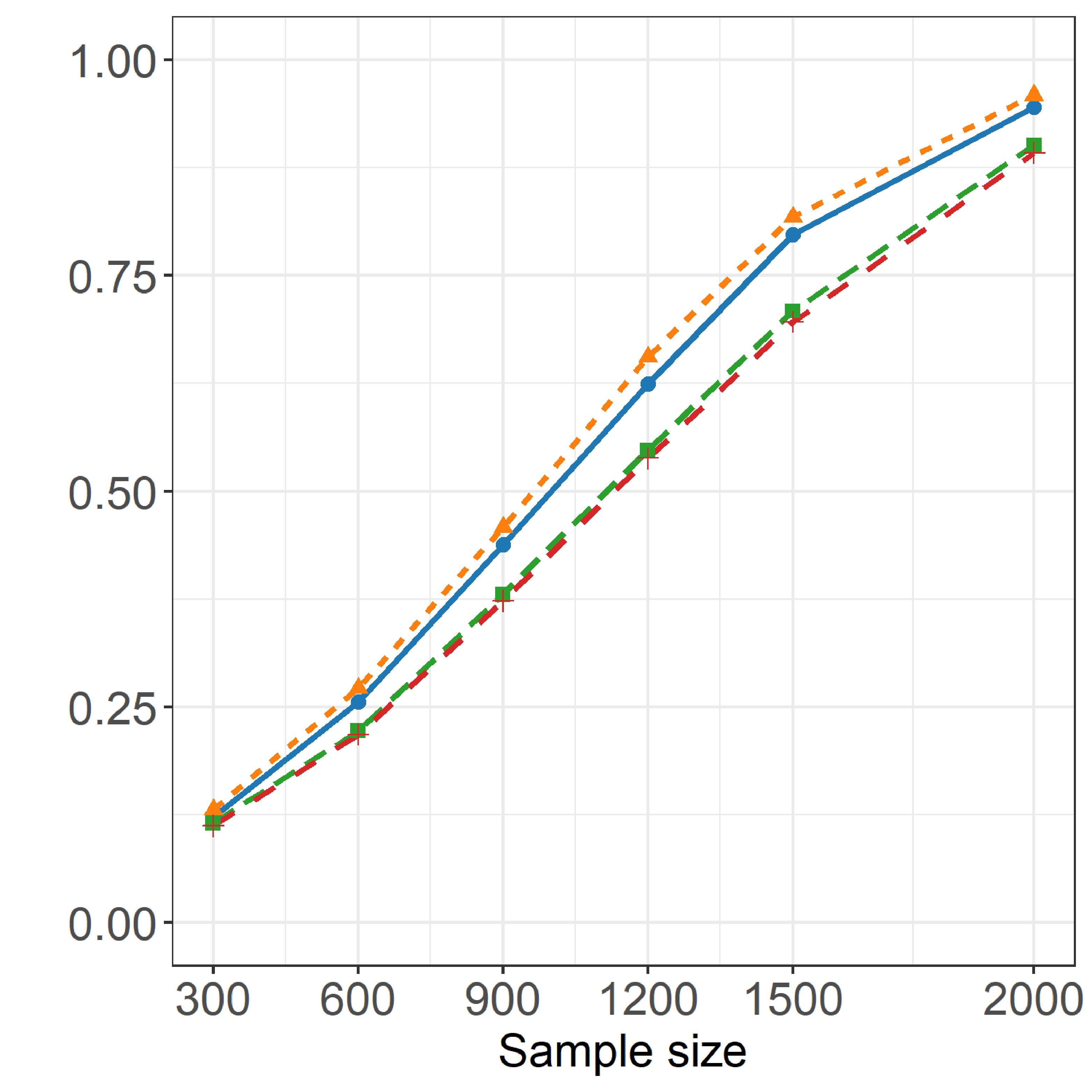}
		\end{minipage}
	}
	\subfigure{
		\begin{minipage}[b]{.3\linewidth}
			\centering
			\includegraphics[scale=0.0735]{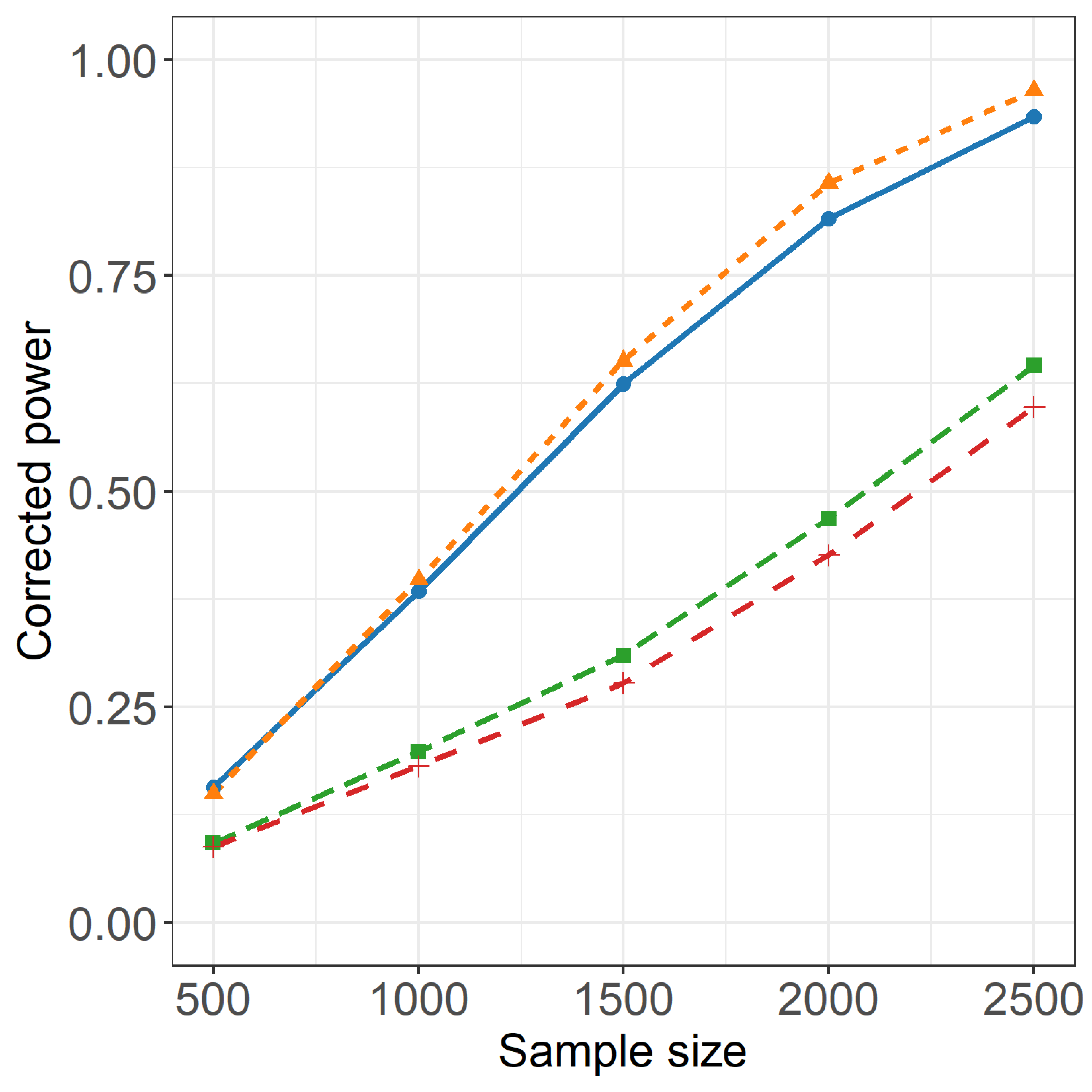}
		\end{minipage}
	}
	\subfigure{
		\begin{minipage}[b]{.3\linewidth}
			\centering
			\includegraphics[scale=0.0735]{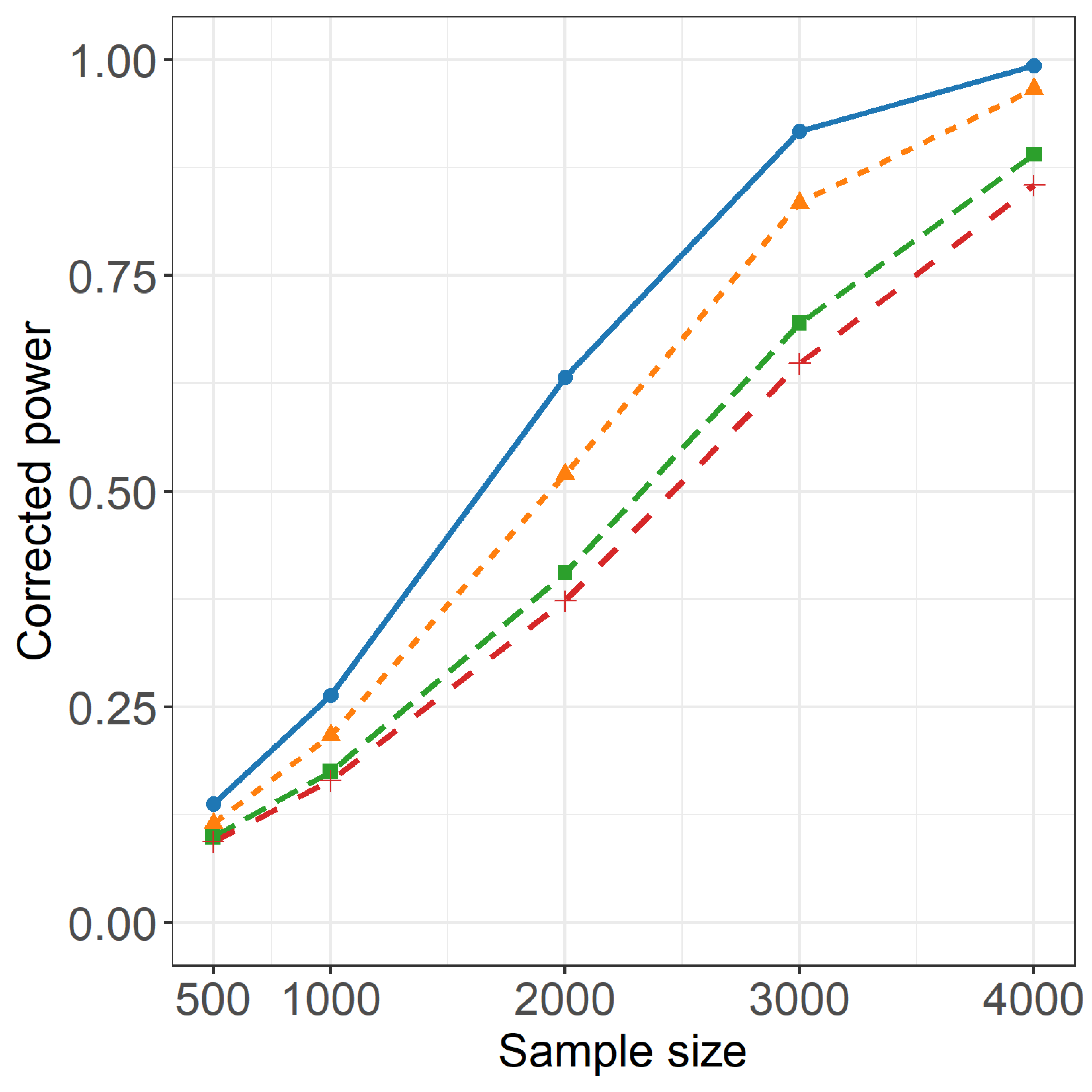}
		\end{minipage}
	}
	\subfigure{
		\begin{minipage}[b]{.3\linewidth}
			\centering
			\includegraphics[scale=0.0735]{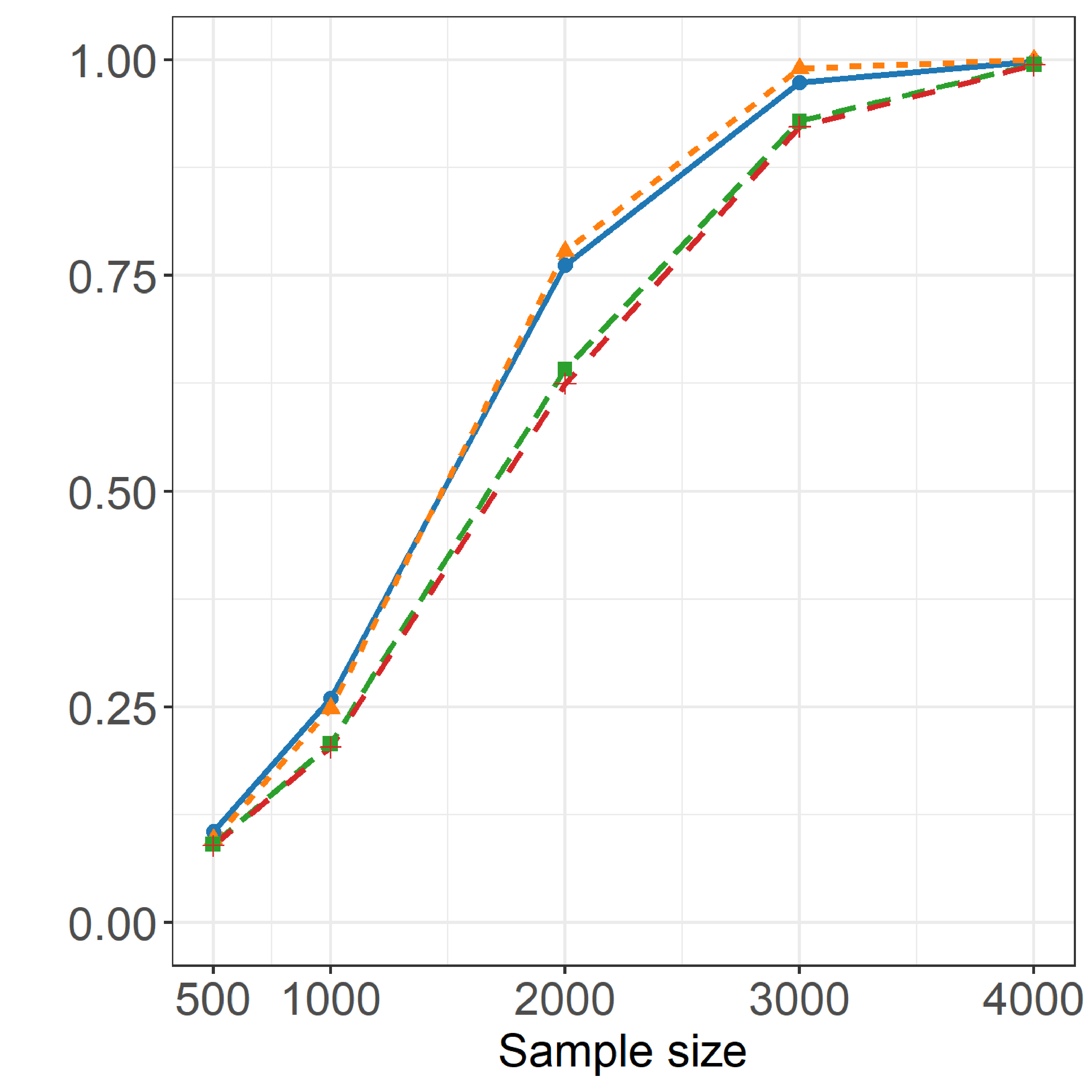}
		\end{minipage}
	}
	\subfigure{
		\begin{minipage}[b]{.3\linewidth}
			\centering
			\includegraphics[scale=0.0735]{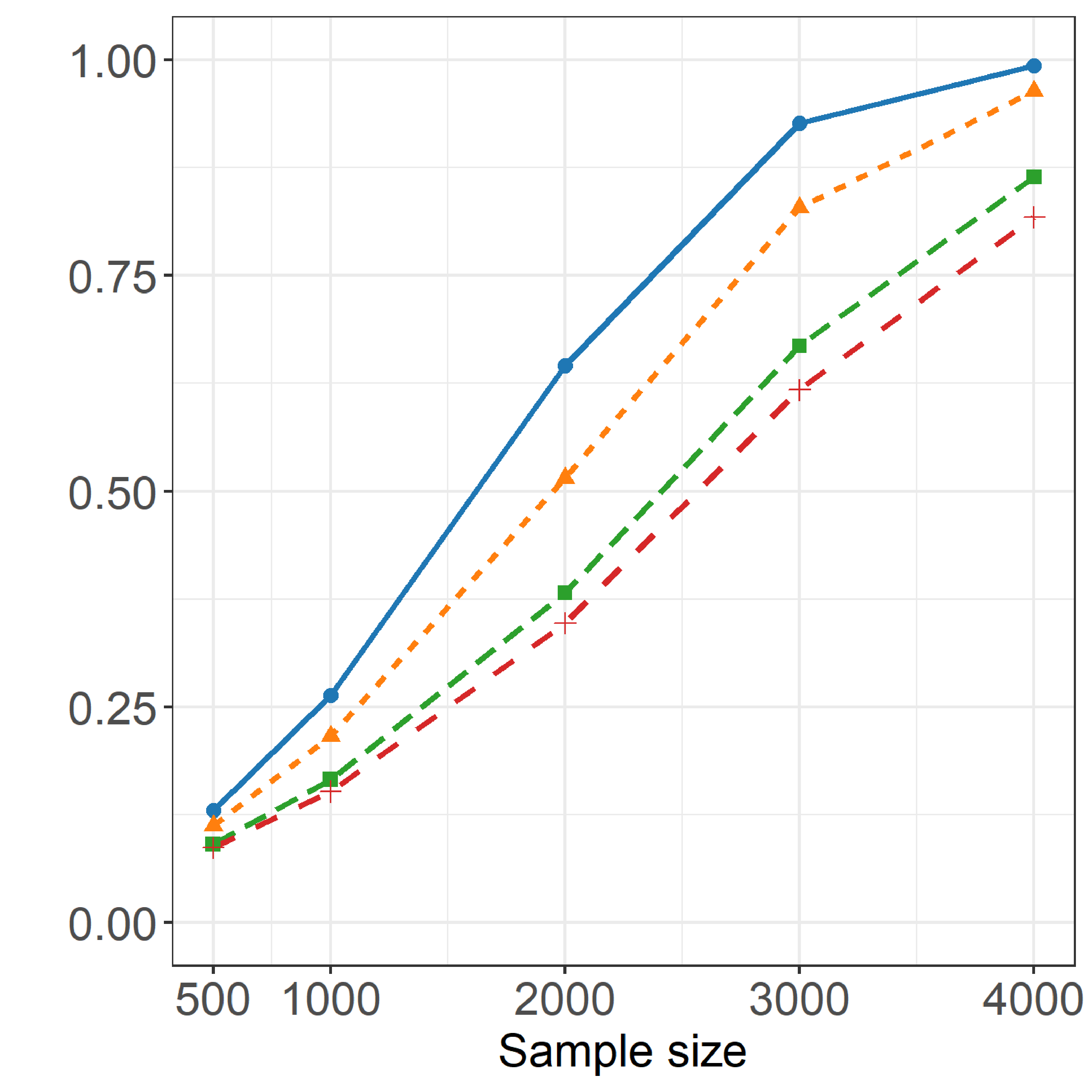}
		\end{minipage}
	}
	\subfigure{
		\begin{minipage}[b]{.3\linewidth}
			\centering
			\includegraphics[scale=0.0735]{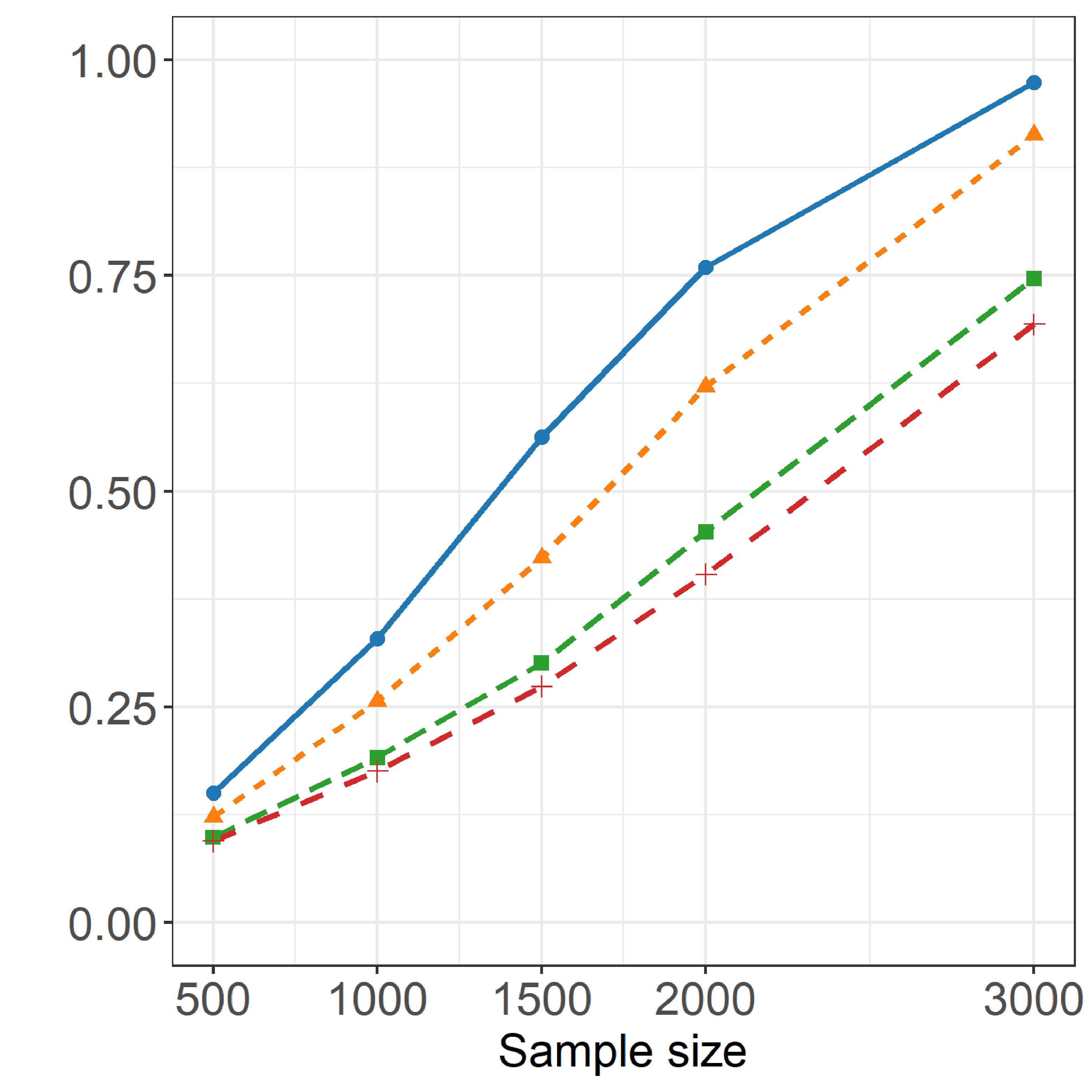}
		\end{minipage}
	}
	\caption{Hypotheses (I)-(VI). Empirical local power functions of four tests  with  various sample sizes $n$ and fixed ratio $p/n=0.3$. The solid, dashed, longdashed, and dotdashed curves are the empirical power functions of the central limit theorem test, directional test and two Skovgaard’s modifications \cite{skovgaard:2001}, respectively.  The six plots correspond to six hypotheses (I)-(VI), starting from top left and proceeding by row.}
	\label{SMFigures: local power vary sample size}
\end{figure}

Figure \ref{SMFigures: local power vary sample size} shows that the corrected power of the directional test is comparable to that of the central limit theorem test, having overall the best  performance across settings. The directional test is slightly more powerful than  the central limit theorem test in hypotheses (I) and (II). The statistics $W^*$ and $W^{**}$ have the lowest power.




\setlength{\tabcolsep}{3.4mm}{
	\begin{table}[H]
		\centering
		\caption{{Hypotheses (I)--(II). Corrected probability of Type I error for the directional test (DT), central limit theorem test (CLT), log-likelihood ratio test (LRT), Bartlett correction (BC) and two Skovgaard's modifications \cite{skovgaard:2001} (Sko1 and Sko2, respectively)  at nominal level $\alpha = 0.05$}}
		{\begin{tabular}{cccccccc}
				\hline
				Hypothesis	&  $p/n_i $ & DT & CLT  & LRT & BC & Sko1 & Sko2 \\
				\hline
				(V) &  0.05 & 0.050 & 0.026 & 0.040 & 0.050 & 0.052  & 0.052 \\ 
				&	0.1 & 0.051 & 0.037 & 0.022 & 0.051 & 0.061 & 0.062 \\ 
				&	0.3 & 0.049 & 0.043 & 0.000 & 0.036 & 0.173 & 0.229 \\ 
				&	0.5 & 0.050 & 0.046 & 0.000 & 0.008 & 0.533 & 0.884 \\ 
				&	0.7 & 0.050 & 0.047 & 0.000 & 0.000 & 0.798 & 1.000 \\ 
				&	0.9 & 0.051 & 0.046 & 0.000 & 0.000 & 0.037 & 1.000 \\ 		 
				(VI) & 0.05 & 0.051 & 0.039 & 0.037 & 0.051 & 0.055 & 0.055 \\ 
				&	0.1 & 0.052 & 0.045 & 0.018 & 0.051 & 0.067 & 0.068 \\ 
				&	0.3 & 0.049 & 0.045 & 0.000 & 0.036 & 0.216 & 0.291 \\ 
				&	0.5 & 0.050 & 0.047 & 0.000 & 0.007 & 0.638 & 0.941 \\ 
				&	0.7 & 0.050 & 0.048 & 0.000 & 0.000 & 0.882 & 1.000 \\ 
				&	0.9 & 0.051 & 0.047 & 0.000 & 0.000 & 0.078 & 1.000 \\
				\hline
			\end{tabular}
		}
		\label{SMsimulation: table type I normal:cases}
\end{table}}


\setlength{\tabcolsep}{5mm}{
	\begin{table}[H]
		\centering
		\caption{\textit{Hypothesis (III). Empirical probability of Type I error for the directional test (DT), central limit theorem test (CLT), log-likelihood ratio test (LRT), Bartlett correction (BC) and two Skovgaard's  modifications \cite{skovgaard:2001} (Sko1 and Sko2, respectively)  at nominal level $\alpha = 0.05$. Two empirical probability of Type I error are considered: estimated Type I error (top panel) and corrected Type I error (bottom panel)}}
		{\begin{tabular}{ccccccc}
				\hline
				$p/n$ &  DT  & CLT  & LRT  &BC &Sko1 &Sko2\\
				   \hline
				0.05 & 0.050 & 0.057 & 0.060 & 0.050 & 0.048 & 0.048 \\ 
				0.1 & 0.049 & 0.056 & 0.084 & 0.050 & 0.037 & 0.036 \\ 
				0.3 & 0.050 & 0.054 & 0.613 & 0.063 & 0.001 & 0.001 \\ 
				0.5 & 0.051 & 0.053 & 1.000 & 0.154 & 0.000 & 0.000 \\ 
				0.7 & 0.050 & 0.053 & 1.000 & 0.697 & 0.000 & 0.000 \\ 
				0.9 & 0.050 & 0.054 & 1.000 & 1.000 & 0.000 & 0.000 \\ 
				\multicolumn{7}{c}{Corrected Type I error} \\ 
				0.05 & 0.050 & 0.041 & 0.041 & 0.050 & 0.052 & 0.053 \\ 
				0.1 & 0.051 & 0.043 & 0.028 & 0.050 & 0.067 & 0.068 \\ 
				0.3 & 0.050 & 0.046 & 0.000 & 0.039 & 0.422 & 0.498 \\ 
				0.5 & 0.050 & 0.047 & 0.000 & 0.011 & 0.987 & 0.999 \\ 
				0.7 & 0.050 & 0.047 & 0.000 & 0.000 & 1.000 & 1.000 \\ 
				0.9 & 0.050 & 0.046 & 0.000 & 0.000 & 1.000 & 1.000 \\
				\hline
			\end{tabular}
		}
		\label{SM:table type I normal:case1}
\end{table}}

\setlength{\tabcolsep}{5mm}{
	\begin{table}[H]
		\centering
		\caption{\textit{Hypothesis (IV). Empirical probability of Type I error for the directional test (DT), central limit theorem test (CLT), log-likelihood ratio test (LRT), Bartlett correction (BC) and two Skovgaard's  modifications \cite{skovgaard:2001} (Sko1 and Sko2, respectively)  at nominal level $\alpha = 0.05$. Two empirical probability of Type I error are considered: estimated Type I error (top panel) and corrected Type I error (bottom panel)} }
		{\begin{tabular}{ccccccc}
				\hline
				$p/n$ &  DT  &CLT  & LRT  &BC & Sko1 &Sko2\\
				 			\hline
				0.05 & 0.050 & 0.068 & 0.061 & 0.050 & 0.049 & 0.049 \\ 
				0.1 & 0.051 & 0.062 & 0.087 & 0.051 & 0.048 & 0.048 \\ 
				0.3 & 0.051 & 0.056 & 0.657 & 0.059 & 0.026 & 0.020 \\ 
				0.5 & 0.051 & 0.053 & 1.000 & 0.112 & 0.009 & 0.001 \\ 
				0.7 & 0.050 & 0.054 & 1.000 & 0.481 & 0.006 & 0.000 \\ 
				0.9 & 0.050 & 0.055 & 1.000 & 1.000 & 0.214 & 0.000 \\ 
				\multicolumn{7}{c}{Corrected Type I error} \\ 
				0.05 & 0.050 & 0.030 & 0.041 & 0.050 & 0.051 & 0.051 \\ 
				0.1 & 0.049 & 0.037 & 0.026 & 0.049 & 0.052 & 0.052 \\ 
				0.3 & 0.049 & 0.044 & 0.000 & 0.042 & 0.088 & 0.110 \\ 
				0.5 & 0.049 & 0.046 & 0.000 & 0.019 & 0.181 & 0.397 \\ 
				0.7 & 0.049 & 0.046 & 0.000 & 0.000 & 0.235 & 0.930 \\ 
				0.9 & 0.050 & 0.045 & 0.000 & 0.000 & 0.004 & 1.000 \\
				\hline
			\end{tabular}
		}
		\label{SM:table type I normal:case2}
\end{table}}


\setlength{\tabcolsep}{5mm}{
	\begin{table}[H]
		\centering
		\captionsetup{font=footnotesize}
		\caption{\textit{Hypothesis (V). Empirical probability of Type I error for the directional test (DT), central limit theorem test (CLT), log-likelihood ratio test (LRT), Bartlett correction (BC) and two Skovgaard's modifications \cite{skovgaard:2001} (Sko1 and Sko2, respectively)  at nominal level $\alpha = 0.05$. Two empirical probability of Type I error are considered: estimated Type I error (top panel) and corrected Type I error (bottom panel)} }
		{\begin{tabular}{ccccccc}
				\hline
				$p/n$ &  DT  &CLT  & LRT  &BC & Sko1 &Sko2\\
				 			\hline
				0.05 & 0.050 & 0.052 & 0.061 & 0.050 & 0.049 & 0.049 \\ 
				0.1 & 0.050 & 0.053 & 0.087 & 0.050 & 0.041 & 0.041 \\ 
				0.3 & 0.050 & 0.052 & 0.629 & 0.061 & 0.002 & 0.002 \\ 
				0.5 & 0.050 & 0.053 & 1.000 & 0.147 & 0.000 & 0.000 \\ 
				0.7 & 0.049 & 0.053 & 1.000 & 0.677 & 0.000 & 0.000 \\ 
				0.9 & 0.050 & 0.055 & 1.000 & 1.000 & 0.000 & 0.000 \\ 
				\multicolumn{7}{c}{Corrected Type I error} \\ 
				0.05 & 0.050 & 0.047 & 0.041 & 0.050 & 0.051 & 0.051 \\ 
				0.1 & 0.050 & 0.047 & 0.027 & 0.050 & 0.060 & 0.061 \\ 
				0.3 & 0.050 & 0.048 & 0.000 & 0.040 & 0.344 & 0.412 \\ 
				0.5 & 0.050 & 0.047 & 0.000 & 0.012 & 0.969 & 0.998 \\ 
				0.7 & 0.050 & 0.047 & 0.000 & 0.000 & 1.000 & 1.000 \\ 
				0.9 & 0.050 & 0.046 & 0.000 & 0.000 & 1.000 & 1.000 \\
				\hline
			\end{tabular}
		}
		\label{SM:table type I normal:case6}
\end{table}}

\setlength{\tabcolsep}{5mm}{
	\begin{table}[H]
		\centering
		\captionsetup{font=footnotesize}
		\caption{\textit{Hypothesis (VI). Empirical probability of Type I error for the directional test (DT), central limit theorem test (CLT), log-likelihood ratio test (LRT), Bartlett correction (BC) and two SKovgaard's modifications \cite{skovgaard:2001} (Sko1 and Sko2, respectively)  at nominal level $\alpha = 0.05$. Two empirical probability of Type I error are considered: estimated Type I error (top panel) and corrected Type I error (bottom panel)}} 
		{\begin{tabular}{ccccccc}
				\hline
				$p/n$ &  DT  &CLT  & LRT  &BC & Sko1 &Sko2\\
				   \hline
				0.05 & 0.050 & 0.052 & 0.061 & 0.050 & 0.049 & 0.049 \\ 
				0.1 & 0.050 & 0.053 & 0.087 & 0.050 & 0.041 & 0.041 \\ 
				0.3 & 0.050 & 0.052 & 0.629 & 0.061 & 0.002 & 0.002 \\ 
				0.5 & 0.050 & 0.053 & 1.000 & 0.147 & 0.000 & 0.000 \\ 
				0.7 & 0.049 & 0.053 & 1.000 & 0.677 & 0.000 & 0.000 \\ 
				0.9 & 0.050 & 0.055 & 1.000 & 1.000 & 0.000 & 0.000 \\ 
				\multicolumn{7}{c}{Corrected Type I error} \\ 
				0.05 & 0.050 & 0.047 & 0.041 & 0.050 & 0.051 & 0.051 \\ 
				0.1 & 0.050 & 0.047 & 0.027 & 0.050 & 0.060 & 0.061 \\ 
				0.3 & 0.050 & 0.048 & 0.000 & 0.040 & 0.344 & 0.412 \\ 
				0.5 & 0.050 & 0.047 & 0.000 & 0.012 & 0.969 & 0.998 \\ 
				0.7 & 0.050 & 0.047 & 0.000 & 0.000 & 1.000 & 1.000 \\ 
				0.9 & 0.050 & 0.046 & 0.000 & 0.000 & 1.000 & 1.000 \\
				\hline
			\end{tabular}
		}
		\label{SM:table type I normal:case5}
\end{table}}


\setlength{\tabcolsep}{5mm}{
	\begin{table}[H]
		\centering
		\captionsetup{font=footnotesize}
		\caption{{Empirical probability of Type I error for the directional test (DT), central limit theorem test (CLT), log-likelihood ratio test (LRT), Bartlett correction (BC) and two Skovgaard's modifications \cite{skovgaard:2001} (Sko1 and Sko2, respectively) at nominal level $\alpha = 0.05$, with $n=100$ and $p=95$. Two empirical probability of Type I error are considered, estimated Type I error (top panel) and corrected Type I error (bottom panel)}}
		{\begin{tabular}{ccccccc}
				\hline
				Hypothesis  &  DT  &CLT  & LRT  &BC & Sko1 &Sko2\\
				   \hline
				(I) & 0.050 & 0.055 & 1.000 & 1.000 & 0.000 & 0.000 \\ 
				(II) & 0.050 & 0.056 & 1.000 & 1.000 & 0.764 & 0.000 \\  
				(III) & 0.050 & 0.056 & 1.000 & 1.000 & 0.000 & 0.000 \\
				(IV) & 0.050 & 0.057 & 1.000 & 1.000 & 0.000 & 0.000 \\ 
				(V) & 0.050 & 0.054 & 1.000 & 1.000 & 0.838 & 0.000 \\ 
				(VI) & 0.050 & 0.053 & 1.000 & 1.000 & 0.751 & 0.000 \\
				\multicolumn{7}{c}{Corrected Type I error} \\ 
				
				(I) & 0.050 & 0.045 & 0.000 & 0.000 & 1.000 & 1.000 \\ 
				(II) & 0.050 & 0.044 & 0.000 & 0.000 & 0.000 & 1.000 \\ 
				(III) & 0.050 & 0.044 & 0.000 & 0.000 & 1.000 & 1.000\\
				(IV) & 0.050 & 0.044 & 0.000 & 0.000 & 1.000 & 1.000 \\
				(V) & 0.050 & 0.046 & 0.000 & 0.000 & 0.000 & 1.000 \\  
				(VI) & 0.050 & 0.046 & 0.000 & 0.000 & 0.000 & 1.000 \\ 
				\hline
			\end{tabular}
		}
		\label{SMtable type I normal: p95}
\end{table}}

\setlength{\tabcolsep}{5mm}{
	\begin{table}[H]
		\centering
		\captionsetup{font=footnotesize}
		\caption{{Empirical probability of Type I error for the directional test (DT), central limit theorem test (CLT), log-likelihood ratio test (LRT), Bartlett correction (BC) and two Skovgaard's modifications \cite{skovgaard:2001} (Sko1 and Sko2, respectively) at nominal level $\alpha = 0.05$, with $n=100$ and $p=98$. Two empirical probability of Type I error are considered, estimated Type I error (top panel) and corrected Type I error (bottom panel)}}
		{\begin{tabular}{ccccccc}
				\hline
				Hypothesis &   DT  &CLT  & LRT  &BC & Sko1 & Sko2\\
				   \hline
				(I) & 0.050 & 0.057 & 1.000 & 1.000 & 0.000 & 0.000 \\ 
				(II) & 0.050 & 0.062 & 1.000 & 1.000 & 0.989 & 0.000 \\ 
				(III) & 0.050 & 0.056 & 1.000 & 1.000 & 0.000 & 0.000 \\ 
				(IV) & 0.053 & 0.057 & 1.000 & 1.000 & 0.000 & 0.000 \\  
				(V) & 0.050 & 0.051 & 1.000 & 1.000 & 1.000 & 0.000 \\ 
				(VI) & 0.058 & 0.051 & 1.000 & 1.000 & 1.000 & 0.000 \\ 
				\multicolumn{7}{c}{Corrected Type I error} \\ 
				(I) & 0.050 & 0.043 & 0.000 & 0.000 & 0.999 & 1.000 \\ 
				(II) & 0.050 & 0.038 & 0.000 & 0.000 & 0.000 & 0.995 \\ 
				(III) & 0.050 & 0.043 & 0.000 & 0.000 & 0.995 & 1.000 \\ 
				(IV) & 0.048 & 0.043 & 0.000 & 0.000 & 1.000 & 1.000 \\
				(V) & 0.050 & 0.049 & 0.000 & 0.000 & 0.000 & 1.000 \\ 
				(VI) & 0.048 & 0.049 & 0.000 & 0.000 & 0.000 & 1.000 \\
				\hline 
			\end{tabular}
		}
		\label{SMtable type I normal: p98}
\end{table}}



\setlength{\tabcolsep}{5mm}{
	\begin{table}[H]
		\centering
		\captionsetup{font=footnotesize}
		\caption{{Hypothesis (I). Empirical probability of Type I error for the directional test (DT), central limit theorem test (CLT), log-likelihood ratio test (LRT), Bartlett correction (BC) and two Skovgaard's modifications \cite{skovgaard:2001} (Sko1 and Sko2, respectively) at nominal level $\alpha = 0.05$ with $p/n=0.3$. Two empirical probability of Type I error are considered, estimated Type I error (top panel) and corrected Type I error (bottom panel)}}
		{\begin{tabular}{ccccccc}
				\hline
				$n_i$ &  DT  &CLT  & LRT  &BC & Sko1 & Sko2\\
				\hline
				300 & 0.050 & 0.052 & 1.000 & 0.097 & 0.000 & 0.000 \\ 
				600 & 0.049 & 0.051 & 1.000 & 0.169 & 0.000 & 0.000 \\ 
				900 & 0.047 & 0.048 & 1.000 & 0.265 & 0.000 & 0.000 \\ 
				1200 & 0.049 & 0.050 & 1.000 & 0.392 & 0.000 & 0.000 \\ 
				1500 & 0.051 & 0.049 & 1.000 & 0.516 & 0.000 & 0.000 \\ 
				2000 & 0.051 & 0.052 & 1.000 & 0.733 & 0.000 & 0.000 \\ 
				\multicolumn{7}{c}{Corrected Type I error} \\ 
				300 & 0.050 & 0.047 & 0.000 & 0.022 & 0.652 & 0.830 \\ 
				600 & 0.051 & 0.049 & 0.000 & 0.010 & 0.992 & 1.000 \\ 
				900 & 0.053 & 0.053 & 0.000 & 0.004 & 1.000 & 1.000 \\ 
				1200 & 0.052 & 0.051 & 0.000 & 0.001 & 1.000 & 1.000 \\ 
				1500 & 0.050 & 0.051 & 0.000 & 0.000 & 1.000 & 1.000 \\ 
				2000 & 0.049 & 0.048 & 0.000 & 0.000 & 1.000 & 1.000\\
				\hline  
		\end{tabular}}
		\label{SMtable type I normal:case4 n}
\end{table}}


\setlength{\tabcolsep}{5mm}{
	\begin{table}[H]
		\centering
		\captionsetup{font=footnotesize}
		\caption{{Hypothesis (II). Empirical probability of Type I error for the directional test (DT), central limit theorem test (CLT), log-likelihood ratio test (LRT), Bartlett correction (BC) and two Skovgaard's modifications \cite{skovgaard:2001} (Sko1 and Sko2, respectively) at nominal level $\alpha = 0.05$ with $p/n=0.3$. Two empirical probability of Type I error are considered, estimated Type I error (top panel) and corrected Type I error (bottom panel)}}
		{\begin{tabular}{ccccccc}
				\hline
				$n_i$ & DT  &CLT  & LRT  &BC & Sko1 & Sko2\\
				\hline
				500 & 0.053 & 0.054 & 1.000 & 0.147 & 0.000 & 0.000 \\ 
				1000 & 0.047 & 0.049 & 1.000 & 0.315 & 0.000 & 0.000 \\ 
				1500 & 0.050 & 0.049 & 1.000 & 0.521 & 0.000 & 0.000 \\ 
				2000 & 0.050 & 0.052 & 1.000 & 0.735 & 0.000 & 0.000 \\ 
				2500 & 0.046 & 0.047 & 1.000 & 0.888 & 0.000 & 0.000 \\  
				\multicolumn{7}{c}{Corrected Type I error} \\ 
				500 & 0.046 & 0.046 & 0.000 & 0.012 & 0.969 & 0.998 \\ 
				1000 & 0.053 & 0.051 & 0.000 & 0.003 & 1.000 & 1.000 \\ 
				1500 & 0.050 & 0.051 & 0.000 & 0.000 & 1.000 & 1.000 \\ 
				2000 & 0.050 & 0.048 & 0.000 & 0.000 & 1.000 & 1.000 \\ 
				2500 & 0.054 & 0.052 & 0.000 & 0.000 & 1.000 & 1.000 \\
				\hline
		\end{tabular}}
		\label{SMtable type I normal:case3 n}
\end{table}}

\setlength{\tabcolsep}{5mm}{
	\begin{table}[H]
		\centering
		\captionsetup{font=footnotesize}
		\caption{{Hypothesis (III).
				Empirical probability of Type I error for the directional test (DT), central limit theorem test (CLT), log-likelihood ratio test (LRT), Bartlett correction (BC) and two Skovgaard's modifications \cite{skovgaard:2001} (Sko1 and Sko2, respectively) at nominal level $\alpha = 0.05$ with $p/n=0.3$. Two empirical probability of Type I error are considered, estimated Type I error (top panel) and corrected Type I error (bottom panel)}}
		{\begin{tabular}{ccccccc}
				\hline
				$n$ & DT  &CLT  & LRT  &BC & Sko1 & Sko2\\
				\hline
				500 & 0.049 & 0.047 & 1.000 & 0.125 & 0.000 & 0.000 \\ 
				1000 & 0.049 & 0.050 & 1.000 & 0.250 & 0.000 & 0.000 \\ 
				2000 & 0.053 & 0.053 & 1.000 & 0.626 & 0.000 & 0.000 \\ 
				3000 & 0.049 & 0.051 & 1.000 & 0.903 & 0.000 & 0.000 \\ 
				4000 & 0.051 & 0.052 & 1.000 & 0.987 & 0.000 & 0.000 \\ 
				\multicolumn{7}{c}{Corrected Type I error} \\ 
				500 & 0.050 & 0.053 & 0.000 & 0.017 & 1.000 & 1.000 \\ 
				1000 & 0.051 & 0.050 & 0.000 & 0.004 & 1.000 & 1.000 \\ 
				2000 & 0.048 & 0.047 & 0.000 & 0.000 & 1.000 & 1.000 \\ 
				3000 & 0.051 & 0.049 & 0.000 & 0.000 & 1.000 & 1.000 \\ 
				4000 & 0.049 & 0.048 & 0.000 & 0.000 & 1.000 & 1.000 \\
				\hline
		\end{tabular}}
		\label{SMtable type I normal:case1 n}
\end{table}}


\setlength{\tabcolsep}{5mm}{
	\begin{table}[H]
		\centering
		\captionsetup{font=footnotesize}
		\caption{{Hypothesis (IV). Empirical probability of Type I error for the directional test (DT), central limit theorem test (CLT), log-likelihood ratio test (LRT), Bartlett correction (BC) and two Skovgaard's modifications \cite{skovgaard:2001} (Sko1 and Sko2, respectively) at nominal level $\alpha = 0.05$ with $p/n=0.3$. Two empirical probability of Type I error are considered, estimated Type I error (top panel) and corrected Type I error (bottom panel)}}
		{\begin{tabular}{ccccccc}
				\hline
				$n$ &   DT  &CLT  & LRT  &BC & Sko1 & Sko2\\
				\hline
				500 & 0.049 & 0.050 & 1.000 & 0.095 & 0.000 & 0.000 \\ 
				1000 & 0.049 & 0.050 & 1.000 & 0.163 & 0.000 & 0.000 \\ 
				2000 & 0.049 & 0.051 & 1.000 & 0.376 & 0.000 & 0.000 \\ 
				3000 & 0.050 & 0.049 & 1.000 & 0.631 & 0.000 & 0.000 \\ 
				4000 & 0.049 & 0.051 & 1.000 & 0.831 & 0.000 & 0.000 \\ 
				\multicolumn{7}{c}{Corrected Type I error} \\ 
				500 & 0.051 & 0.051 & 0.000 & 0.025 & 0.497 & 0.702 \\ 
				1000 & 0.052 & 0.050 & 0.000 & 0.010 & 0.951 & 0.997 \\ 
				2000 & 0.051 & 0.049 & 0.000 & 0.001 & 1.000 & 1.000 \\ 
				3000 & 0.050 & 0.051 & 0.000 & 0.000 & 1.000 & 1.000 \\ 
				4000 & 0.051 & 0.049 & 0.000 & 0.000 & 1.000 & 1.000\\
				\hline
			\end{tabular}
		}
		\label{SMtable type I normal:case2 n}
\end{table}}


\setlength{\tabcolsep}{5mm}{
	\begin{table}[H]
		\centering
		\captionsetup{font=footnotesize}
		\caption{{Hypothesis (V). Empirical probability of Type I error for the directional test (DT), central limit theorem test (CLT), log-likelihood ratio test (LRT), Bartlett correction (BC) and two Skovgaard's modifications \cite{skovgaard:2001} (Sko1 and Sko2, respectively) at nominal level $\alpha = 0.05$ with $p/n=0.3$. Two empirical probability of Type I error are considered, estimated Type I error (top panel) and corrected Type I error (bottom panel)} }
		{\begin{tabular}{ccccccc}
				\hline
				$n$ &   DT  &CLT  & LRT  &BC & Sko1 & Sko2\\
				\hline
				500 & 0.049 & 0.047 & 1.000 & 0.120 & 0.000 & 0.000 \\ 
				1000 & 0.048 & 0.050 & 1.000 & 0.247 & 0.000 & 0.000 \\ 
				2000 & 0.052 & 0.053 & 1.000 & 0.624 & 0.000 & 0.000 \\ 
				3000 & 0.049 & 0.051 & 1.000 & 0.901 & 0.000 & 0.000 \\ 
				4000 & 0.051 & 0.052 & 1.000 & 0.987 & 0.000 & 0.000 \\ 
				\multicolumn{7}{c}{Corrected Type I error} \\ 
				500 & 0.051 & 0.053 & 0.000 & 0.017 & 1.000 & 1.000 \\ 
				1000 & 0.052 & 0.049 & 0.000 & 0.004 & 1.000 & 1.000 \\ 
				2000 & 0.048 & 0.048 & 0.000 & 0.000 & 1.000 & 1.000 \\ 
				3000 & 0.051 & 0.049 & 0.000 & 0.000 & 1.000 & 1.000 \\ 
				4000 & 0.049 & 0.047 & 0.000 & 0.000 & 1.000 & 1.000\\
				\hline
			\end{tabular}
		}
		\label{SMtable type I normal:case6 n}
\end{table}}

\setlength{\tabcolsep}{5mm}{
	\begin{table}[H]
		\centering
		\captionsetup{font=footnotesize}
		\caption{\textit{Hypothesis (VI). Empirical probability of Type I error  for the directional test (DT), central limit theorem test (CLT), log-likelihood ratio test (LRT), Bartlett correction (BC) and two Skovgaard's modifications \cite{skovgaard:2001} (Sko1 and Sko2, respectively) at nominal level $\alpha = 0.05$ with $p/n=0.3$. Two empirical probability of Type I error are considered, estimated Type I error (top panel) and corrected Type I error (bottom panel)} }
		{\begin{tabular}{ccccccc}
				\hline
				$n$ &   DT  &CLT  & LRT  &BC & Sko1 & Sko2\\
				\hline
				500 & 0.049 & 0.049 & 1.000 & 0.127 & 0.000 & 0.000 \\ 
				1000 & 0.049 & 0.050 & 1.000 & 0.254 & 0.000 & 0.000 \\ 
				1500 & 0.054 & 0.052 & 1.000 & 0.436 & 0.000 & 0.000 \\ 
				2000 & 0.053 & 0.053 & 1.000 & 0.631 & 0.000 & 0.000 \\ 
				3000 & 0.050 & 0.052 & 1.000 & 0.904 & 0.000 & 0.000 \\ 
				\multicolumn{7}{c}{Corrected Type I error} \\ 
				500 & 0.051 & 0.051 & 0.000 & 0.017 & 1.000 & 1.000 \\ 
				1000 & 0.051 & 0.050 & 0.000 & 0.004 & 1.000 & 1.000 \\ 
				1500 & 0.046 & 0.049 & 0.000 & 0.001 & 1.000 & 1.000 \\ 
				2000 & 0.047 & 0.046 & 0.000 & 0.000 & 1.000 & 1.000 \\ 
				3000 & 0.050 & 0.048 & 0.000 & 0.000 & 1.000 & 1.000 \\ 
				\hline
			\end{tabular}
		}
		\label{SMtable type I normal:case5 n}
\end{table}}


\end{document}